\documentclass{memo-l}
\usepackage[mathscr]{eucal}
\usepackage{amsmath,amsfonts}

\parskip=\smallskipamount
\newtheorem{theorem}{Theorem}[chapter]

\newtheorem{proposition}[theorem]{Proposition}

\newtheorem{corollary}[theorem]{Corollary}
\newtheorem{lemma}[theorem]{Lemma}
\newtheorem{example}[theorem]{Example}

\newtheorem{remark}[theorem]{Remark}

\makeatletter \@addtoreset{equation}{chapter}
\makeatother

\numberwithin{section}{chapter} \numberwithin{equation}{chapter}

\newcommand{\BB}{{\mathbb B}}
\newcommand{\CC}{{\mathbb C}}
\newcommand{\NN}{{\mathbb N}}

\newcommand{\ZZ}{{\mathbb Z}}
\newcommand{\DD}{{\mathbb D}}
\newcommand{\RR}{{\mathbb R}}
\newcommand{\FF}{{\mathbb F}}
\newcommand{\TT}{{\mathbb T}}

\newcommand{\HH}{{\mathbb H}}
\newcommand{\KK}{{\mathbb K}}

\newcommand{\cA}{{\mathcal A}}

\newcommand{\cC}{{\mathcal C}}
\newcommand{\cD}{{\mathcal D}}
\newcommand{\cE}{{\mathcal E}}
\newcommand{\cF}{{\mathcal F}}
\newcommand{\cG}{{\mathcal G}}
\newcommand{\cH}{{\mathcal H}}
\newcommand{\cK}{{\mathcal K}}
\newcommand{\cL}{{\mathcal L}}
\newcommand{\cM}{{\mathcal M}}
\newcommand{\cN}{{\mathcal N}}
\newcommand{\cP}{{\mathcal P}}
\newcommand{\cR}{{\mathcal R}}
\newcommand{\cS}{{\mathcal S}}

\newcommand{\cU}{{\mathcal U}}
\newcommand{\cV}{{\mathcal V}}
\newcommand{\cY}{{\mathcal Y}}

\newcommand{\cW}{{\mathcal W}}

\newcommand{\rank}{\hbox{\rm{rank}}\,}

\newdimen\expt
\def\boxit#1{\setbox0\hbox{$\displaystyle{#1}$}
      \hbox{\lower.4\expt
 \hbox{\lower3\expt\hbox{\lower\dp0
      \hbox{\vbox{\hrule height.4\expt
 \hbox{\vrule width.4\expt\hskip3\expt
      \vbox{\vskip3\expt\box0\vskip2\expt}%
 \hskip3\expt\vrule width.4\expt}\hrule height.4\expt}}}}}}

\begin{document}
\frontmatter
 \pagestyle{myheadings}
\markboth{GELU POPESCU}{OPERATOR THEORY ON NONCOMMUTATIVE DOMAINS}
%\pagestyle{plain}

%\begin{flushright}
 % \it Date of this draft: \today
%\end{flushright}
%\bigskip

\title [Operator theory on noncommutative domains]
{Operator theory on noncommutative domains}
  \author{Gelu Popescu}

\address{Department of Mathematics, The University of Texas
at San Antonio \\ San Antonio, TX 78249, USA} \email{\tt
gelu.popescu@utsa.edu}

 \date{February 16, 2007}
\thanks{Research supported in part by an NSF grant}
\subjclass[2000]{Primary: 47A05; 47A56; 47A20;  46E40;  Secondary:
46L52; 46L07; 47A67; 47A63; 47A57; 47A60} \keywords{Multivariable
operator theory, Free holomorphic function,  Noncommutative domain,
Noncommutative variety, Fock space, Weighted shifts, Invariant
subspace, Hardy algebra, Cauchy transform, Poisson transform,  von
Neumann inequality, Bohr inequality,  Functional calculus, Wold
decomposition, Dilation, Characteristic function, Model theory,
Curvature, Commutant lifting, Interpolation}

\begin{abstract}
In this volume we study noncommutative domains $\cD_f\subset
B(\cH)^n$
 generated  by  positive regular free
holomorphic functions $f$ on $B(\cH)^n$, where $B(\cH)$ is the
algebra of all bounded linear operators on a Hilbert space $\cH$.

Each such a domain  has a  universal model $(W_1,\ldots, W_n)$ of
weighted shifts acting on the full Fock space with $n$ generators.
The  study of $\cD_f$ is close related to the study of the weighted
shifts $W_1,\ldots,W_n$, their joint invariant subspaces, and the
representations of the algebras they generate: the domain algebra
$\cA_n(\cD_f)$, the Hardy  algebra $F_n^\infty(\cD_f)$, and the
$C^*$-algebra $C^*(W_1,\ldots, W_n)$. A good part of this paper
deals with these issues. We also introduce the symmetric weighted
Fock space $F_s^2(\cD_f)$ and show that it can be identified with a
reproducing kernel Hilbert space. The algebra of all its
``analytic''  multipliers will play an important role in the
commutative case.

Free holomorphic functions, Cauchy transforms,  and Poisson
transforms on noncommutative domains $\cD_f$ are introduced and used
to provide an $F_n^\infty(\cD_f)$-functional calculus for completely
non-coisometric  elements of $\cD_f(\cH)$, and a free analytic
functional calculus for $n$-tuples of operators $(T_1,\ldots, T_n)$
with the joint spectral radius  $r_p(T_1,\ldots,T_n)<1$. Several
classical results from complex analysis have analogues in our
noncommutative setting of free holomorphic functions on $\cD_f$.

We associate with each $w^*$-closed two-sided ideal $J$ of the
algebra $F_n^\infty(\cD_f)$ a noncommutative variety
$\cV_{f,J}\subset \cD_f$. We develop a dilation theory and model
theory for $n$-tuples of operators $T:=(T_1,\ldots, T_n)$ in the
noncommutative domain $\cD_f$ (resp. noncommutative variety
$\cV_{f,J}$). We associate with each such an $n$-tuple of operators
a characteristic function $\Theta_{f,T}$ (resp. $\Theta_{f,T,J}$),
use it to provide a functional model,  and prove that it is a
complete unitary invariant for completely non-coisometric elements
of $\cD_f$ (resp. $\cV_{f,J}$). In particular,  we discuss the
commutative case when $T_iT_j=T_jT_i$, $i=1,\ldots,n$.

We introduce  two numerical invariants, the curvature and
$*$-curvature, defined on the noncommutative domain $\cD_p$, where
$p$ is positive regular noncommutative polynomial,  and present some
basic properties. We show that both curvatures can be express in
terms of the characteristic function $\Theta_{p,T}$.

  We present a
commutant lifting theorem for pure $n$-tuples of operators in
noncommutative domains $\cD_f$ (resp. varieties $\cV_{f,J}$) and
obtain   Nevanlinna-Pick and Schur-Carath\'eodory type interpolation
results.  We also obtain a corona theorem  for Hardy algebras
associated with $\cD_f$ (resp. $\cV_{f,J}$).

In the particular case when  $f=X_1+\cdots+X_n$, we recover several
results concerning the multivariable noncommutative (resp.
commutative) operator theory on the unit ball  $[B(\cH)^n]_1$.

\end{abstract}

\maketitle \setcounter{page}{4}

\tableofcontents

 %\clearpage
\chapter*{Introduction}
\pagenumbering{arabic}

The study of the operator unit ball
$$
[B(\cH)]_1:=\left\{T\in B(\cH):\ \|TT^*\|\leq 1\right\}
$$
has generated, in the last fifty years, the celebrated
Sz.-Nagy--Foia\c s theory of contractions on Hilbert spaces
\cite{SzF-book}. A key role in this theory is played by the model
operator associated with $[B(\cH)]_1$, i.e., the unilateral shift
$S$ acting on the Hardy space $H^2(\DD)$, where $\DD:=\{z\in\CC:\
|z|<1\}$. The study of the unit ball $[B(\cH)]_1$ is very close
related to the study of the unilateral shift, its invariant
subspaces, and the representations of the algebras it generates: the
disc algebra $A(\DD)$, the analytic Toeplitz algebra
$H^\infty(\DD)$, and the Toeplitz $C^*$-algebra $C^*(S)$. This
interplay between operator theory and harmonic analysis is
illustrated in the dilation and model theory of contractions, and
the profound implications in function theory, interpolation,
prediction theory, scattering theory, and linear system theory (see
\cite{SzF-book}, \cite{FF-book}, \cite{FFGK-book}, \cite{BaGR}).

In the noncommutative multivariable setting, the study of the
operator unit $n$-ball
$$
[B(\cH)^n]_1:=\left\{ (T_1,\ldots, T_n)\in B(\cH)^n:\
\|T_1T_1^*+\cdots +T_nT_n^*\|\leq 1\right\}
$$
has generated a {\it free analogue} of Sz.-Nagy--Foia\c s theory
(see \cite{F}, \cite{B}, \cite{Po-models}, \cite{Po-isometric},
\cite{Po-charact}, \cite{Po-multi}, \cite{Po-von}, \cite{Po-intert},
\cite{Po-funct}, \cite{Po-analytic}, \cite{Po-disc} and more
recently \cite{Po-poisson},  \cite{DKS}, \cite{BV}, \cite{BV3},
\cite{BBD}, \cite{Po-structure}, \cite{Po-curvature},
\cite{Po-nehari}, \cite{Po-entropy}, \cite{Po-unitary},
\cite{Po-charact-inv}, \cite{Po-varieties}, \cite{Po-varieties2},
\cite{BT2}).
In this case, the model  associated with $[B(\cH)^n]_1$ is the
$n$-tuple  $(S_1,\ldots, S_n)$ of  the left creation operators
acting on the full Fock space with $n$ generators $F^2(H_n)$. The
corresponding algebras are the noncommutative disc algebra $\cA_n$,
the noncommutative analytic Toeplitz algebra $F_n^\infty$, and the
Cuntz-Toeplitz algebra $C^*(S_1,\ldots, S_n)$. Introduced in
\cite{Po-von} in connection with a multivariable noncommutative von
Neumann inequality, the noncommutative analytic Toeplitz algebra
$F_n^\infty$
  has  been studied
    in several papers
\cite{Po-charact},  \cite{Po-multi},  \cite{Po-funct},
\cite{Po-analytic}, \cite{Po-disc}, \cite{Po-poisson},
 \cite{ArPo}, and recently in \cite{DP}, \cite{DP2},
  \cite{DP1},
  \cite{ArPo2}, \cite{Po-curvature},  \cite{DKP},  \cite{PPoS},
  \cite{Po-similarity}, and \cite{Po-holomorphic}.

Interpolation problems for the noncommutative analytic Toeplitz
algebra
 $F_n^\infty$ were first considered in \cite{Po-charact}, where we
 obtained a Sarason \cite{S} type interpolation result, and
in \cite{Po-analytic}, where   we
 obtained a
Schur-Carath\'eodory (\cite{Sc}, \cite{Ca}) type interpolation
theorem. In  1997,
 Arias and
 the author \cite{ArPo2} (see also \cite{Po-interpo})
  obtained a distance formula to
an arbitrary WOT-closed ideal in $F_n^\infty$ as well as a
Nevanlinna-Pick  (see \cite{N}) type interpolation theorem  for the
noncommutative analytic Toeplitz algebra
 $F_n^\infty$.
Using different methods, Davidson and Pitts  proved   these
  results  independently  in  \cite{DP}.
Let us mention that, recently, interpolation problems for
 $F_n^\infty$   and related interpolation problems on
  the unit  ball of
$\CC^n$  were  also considered in \cite{AMc2},
  \cite{BTV}, \cite{BB}, \cite{BB2},  \cite{Po-central},  \cite{Po-nehari},
    \cite{EP}, and \cite{Po-entropy}.

Applications of the noncommutative dilation  and model theory  to
prediction theory and entropy optimization were considered in
\cite{Po-structure} and \cite{Po-entropy}. It is well known that, in
the classical case, there is essentially a one-to-one correspondence
between any two of the following: operator model theory, scattering
theory, and unitary system theory. The same is true in the
noncommutative multivariable case. The model theory for row
contractions (elements of $[B(\cH)]_1$) (see \cite{Po-models},
\cite{Po-isometric}, and \cite{Po-charact}) lead Ball and Vinikov
\cite{BV} to their multivariable setting for Lax-Phillips scattering
and conservative linear systems. Along this line, we also mention
the paper  \cite{BGM} of Ball, Groenewald, and Malakorn.

In \cite{Po-holomorphic}, we developed a theory of holomorphic
functions in several noncommuting (free) variables and  provide a
framework for  the study of arbitrary
 $n$-tuples of operators. This theory enhances our program to develop a {\it free}
  analogue of
  Sz.-Nagy--Foia\c s theory \cite{SzF-book}, for row contractions.
  A free analytic functional calculus was introdued in \cite{Po-holomorphic}
  and studied in connection with
 Hausdorff derivations,  noncommutative Cauchy and Poisson transforms,
 and von Neumann inequalities.
 In a related area of research,
  we remark
 the work  of Helton, McCullough, Putinar, and Vinnikov, on symmetric noncommutative  polynomials (\cite{He}, \cite{He-C},
 \cite{He-C-P},
 \cite{He-C-P2}, \cite{He-C-V}), which contains a new type of
 engineering application.

 We should also remark that, in recent years, many results
 concerning the theory of row contractions were extended by Muhly
 and Solel (\cite{MuSo1}, \cite{MuSo5}, \cite{MuSo6}, \cite{MuSo7})
 to representations of tensor algebras over $C^*$-correspondences and Hardy algebras.

Motivated by the profound impact that  the   operator theory of the
unit ball $[B(\cH)^n]_1$, $n\geq 1$,  has had in several areas of
research, we  study,  in this paper, more general noncommutative
domains of the form
$$
\cD_f(\cH):=\left\{(X_1,\ldots, X_n)\in B(\cH)^n: \
\left\|\sum_{\alpha\in \FF_n^+, |\alpha|\geq 1} a_\alpha X_\alpha
X_\alpha^*\right\|\leq 1\right\},
$$
where  $f:=\sum_{\alpha\in \FF_n^+, |\alpha|\geq 1} a_\alpha
X_\alpha$, $a_\alpha\in \CC$,  is any {\it positive regular free
holomorphic function} on $B(\cH)^n$. Such an $f$ has positive
coefficients with $a_{g_i}\neq 0$, $i=1,\ldots,n$, and
$\limsup_{k\to\infty}
\left(\sum_{|\alpha|=k}|a_\alpha|^2\right)^{1/2k}< \infty. $ Here,
$\FF_n^+$ is the free semigroup with $n$  generators $g_1,\ldots,
g_n$ and $|\alpha|$ stands for the length of the word $\alpha\in
\FF_n^+$. For each $\alpha:=g_{i_1}\cdots g_{i_k}\in \FF_n^+$,
$X_\alpha$ denotes the product $X_{i_1}\cdots X_{i_k}$.

In Chapter 1, we associate with each noncommutative domain $\cD_f$
an essentially unique $n$-tuple $(W_1,\ldots, W_n)$ of weighted
shifts acting on the full Fock space  with $n$ generators
$F^2(H_n)$. This will play the role of {\it universal model} for the
elements of $\cD_f(\cH)$. Several of its properties will be
presented throughout the paper.

We introduce, in Section \ref{Noncommutative}, the {\it domain
algebra} $\cA_n(\cD_f)$, as the norm closure of all polynomials in
$W_1,\ldots, W_n$, and the identity, and characterize its completely
contractive (resp. completely bounded) representations. In
particular, we obtain a von Neumann  \cite{vN} type inequality,
i.e.,
$$\|q(T_1,\ldots, T_n)\|\leq \|q(W_1,\ldots, W_n)\|
$$
for any $(T_1,\ldots, T_n)\in \cD_f(\cH)$ and any polynomial $q$ in
$n$ noncommutative indeterminates. We also identify, in this
section, the set of all characters of $\cA_n(\cD_f)$ with the set
$$
\cD_f(\CC):=\left\{ (\lambda_1,\ldots, \lambda_n)\in \CC^n:\
\sum_{|\alpha|\geq 1} a_\alpha |\lambda_\alpha|^2\leq 1\right\}.
$$

In Section \ref{domain algebra}, we introduce the Hardy type algebra
$F_n^\infty(\cD_f)$ (resp. $R_n^\infty(\cD_f)$)  and provide some
basic properties including the fact that it is semisimple and
coincides with its double commutant. This algebra is the analogue of
$H^\infty(\DD)$ (see \cite{G} and \cite{H}), for the noncommutative
domain $\cD_f$, and it will play an important role in our
investigation.

In the next section,  we show that  the Hardy algebra
$F_n^\infty(\cD_f)$ is the $w^*$-(resp. WOT-, SOT-) closure of all
polynomials in $W_1,\ldots W_n$,  and the identity. Based on
noncommutative Poisson transforms associated with $\cD_f$, we obtain
an $F_n^\infty(\cD_f)$-functional calculus for {\it completely
non-coisometric} (c.n.c.) $n$-tuple of operators in the
noncommutative domain $\cD_f(\cH)$. More precisely, we prove that,
for each c.n.c. $n$-tuple $T:=(T_1,\ldots, T_n)\in \cD_f(\cH)$,
there exists a unital completely contractive homomorphism
$$
\Phi_{f,T}: F_n^\infty(\cD_f)\to B(\cH)
$$
which is $WOT$-continuous and $\Phi_{f,T}(W_\beta)=T_\beta$ for all
$\beta\in \FF_n^+$.

An important extension of this result to noncommutative varieties of
$\cD_f$ is obtained in  the next section. Let us explain a little
bit more.
 Each $w^*$-closed two-sided
ideal $J$ of the Hardy algebra $F_n^\infty(\cD_f)$ generates a
noncommutative variety $\cV_{f,J}\subset \cD_f$ which, represented
on a Hilbert space $\cH$,  is defined by
$$
\cV_{f,J}(\cH):=\left\{(X_1,\ldots, X_n)\in \cD_f(\cH):\
\varphi(X_1,\ldots, X_n)=0 \text{ for all } \varphi\in J\right\},
$$
where $\varphi(X_1,\ldots, X_n)$ is defined by an appropriate
functional calculus. The {\it universal model} associated with
$\cV_{f,J}$ is the $n$-tuple $(B_1,\ldots, B_n)$ of {\it constrained
weighted shifts} defined by $B_i:=P_{\cN_J} W_i |_{\cN_J}$,
$i=1,\ldots, n$, where  $(W_1,\ldots, W_n)$ is the universal model
associated with $\cD_f$ and $\cN_J:=F^2(H_n)\ominus \overline{J
F^2(H_n)}$. The Hardy algebra associated with $\cV_{f,J}$ is the
compression of $F_n^\infty(\cD_f)$ to the subspace $\cN_J$, and it
is denoted by $F_n^\infty(\cV_{f,J})$. It turns out that
$F_n^\infty(\cV_{f,J})$ is the $w^*$-closure of all polynomials in
$B_1,\ldots, B_n$ and the identity.

In Section \ref{Functional II}, we show that, for each
 c.n.c. $n$-tuple $T:=(T_1,\ldots, T_n)\in \cV_{f,J}(\cH)$, there
exists a unital completely contractive homomorphism
$$
\Psi_{f,T,J}: F_n^\infty(\cV_{f,J})\to B(\cH)
$$
which is $WOT$-continuous and $\Psi_{f,T, J}(B_\beta)=T_\beta$ for
all $\beta\in \FF_n^+$.

Section \ref{Symmetric}  is devoted to the weighted shifts
$(W_1,\ldots, W_n)$ associated with the noncommutative domain
$\cD_f$. We determine all the eigenvectors of $W_1^*,\ldots, W_n^*$,
and  prove that the right joint spectrum $\sigma_r(W_1,\ldots, W_n)$
coincides with $\cD_f(\CC)$. This will enable us to
 identify  the $w^*$-continuous
multiplicative functionals of the Hardy algebra $F_n^\infty(\cD_f)$.

We  introduce the symmetric weighted Fock space $F_s^2(\cD_f)$
associated with the noncommutative domain $\cD_f$ and prove that it
can be identified with $H^2(\cD_f^\circ(\CC))$, a Hilbert space  of
holomorphic functions defined on
$$
\cD_f^\circ(\CC):=\left\{ (\lambda_1,\ldots, \lambda_n)\in \CC^n: \
\sum_{|\alpha|\geq 1} a_\alpha |\lambda_\alpha|^2<1\right\}.
$$
Moreover, we show that $F_s^2(\cD_f)$ is the reproducing kernel
Hilbert space with reproducing kernel $K_f:\cD_f^\circ(\CC)\times
\cD_f^\circ(\CC)$ defined by
$$
K_f(\mu,\lambda):=\frac{1}{1-\sum_{|\alpha|\geq 1} a_\alpha
\mu_\alpha \overline{\lambda}_\alpha},\qquad \mu,\lambda\in
\cD_f^\circ(\CC).
$$
In the end of this section, we identify the algebra of all
multipliers of the Hilbert space $H^2(\cD_f^\circ(\CC))$ with
$F_n^\infty(\cV_{f,J_c})$, the $w^*$-closed algebra generated by the
operators $L_i:=P_{F_s^2(\cD_f)}W_i|_{F_s^2(\cD_f)}$,\ $i=1,\ldots,
n$, and the identity, and prove that it is reflexive.
 This algebra will play an important role in the
commutative case.

In Chapter 2, we introduce the algebra $Hol(\cD_f)$ of all free
holomorphic functions on the noncommutative domain $\cD_f$ and prove
a version of the maximum principle \cite{Co}. We identify the domain
algebra $\cA_n(\cD_f)$ and the Hardy algebra $F_n^\infty(\cD_f)$
with subalgebras of free holomorphic functions on $\cD_f$.

When $p$ is a positive regular noncommutative polynomial, we point
out  two Banach algebras of free holomorphic functions,
$H^\infty(\cD_p^\circ)$ and $A(\cD_p^\circ)$, which can be
identified with the Hardy algebra $F_n^\infty(\cD_p)$ and domain
algebra $\cA_n(\cD_p)$, respectively, and show that the elements of
these algebras can be seen as boundary functions for noncommutative
Poisson transforms on $\cD_p^\circ$. For example, we show that if
$u$ is a free holomorphic function on $\cD_p^\circ$, then there
exists $f\in F_n^\infty(\cD_p)$ with
$$u=P[f]\quad \text{ if and only if }\quad
\sup_{0\leq r<1}\|u(rW_1,\ldots, rW_n)\|<\infty. $$
 Moreover, in
this case, $u(rW_1,\ldots, rW_n)\to f$, as $r\to1$, in the
$w^*$-topology (or strong operator topology).

In Section \ref{Schwartz}, we obtain versions of Schwarz lemma
\cite{Co} and Bohr's inequality \cite{Bo} for the Hardy algebra
$F_n^\infty(\cD_f)$. Our version of Bohr's inequality states, in
particular, that if $\varphi=\sum_{\beta\in \FF_n^+}c_\beta W_\beta$
is in $F_n^\infty(\cD_f)$, then
$$
\sum_{\beta\in \FF_n^+} |c_\beta||\lambda_\beta|\leq \|\varphi\|
\quad \text{ for all }\quad (\lambda_1,\ldots, \lambda_n)\in
\cD_{f,1/3}(\CC), $$
 where
$$
\cD_{f,1/3}(\CC):=\left\{(\lambda_1,\ldots, \lambda_n)\in \CC^n:\
(3\lambda_1,\ldots, 3\lambda_n)\in \cD_f(\CC)\right\}.
$$
Notice that if $n=1$ and $f=X$ we obtain the classical Bohr's
inequality \cite{Bo}.

In Section \ref{Weierstrass}, we obtain Weierstrass and Montel type
theorems \cite{Co} for the algebra of free holomorphic functions on
$\cD_f$. This enables us to introduce a metric on $Hol(\cD_f)$ with
respect to which it becomes a complete metric space.

A noncommutative Cauchy transform  is introduced, in Section
\ref{Cauchy}, and used to provide a free analytic functional
calculus $\Phi_{p,T}: Hol(\cD_p)\to B(\cH)$ for $n$-tuples of
operators $T:=(T_1,\ldots, T_n)\in B(\cH)^n$ with joint spectral
radius $r_p(T_1,\ldots, T_n)<1$. We show that the  free analytic
functional calculus is continuous and unique. In the last part of
this section, we present multivariable commutative versions of the
free analytic functional calculus and Bohr's inequality.

Chapter 3 is devoted to dilation theory, model theory, and unitary
invariants on noncommutative domains. In Section \ref{Weighted}, we
obtain a Beurling \cite{Be} type characterization of the invariant
subspaces under the weighted shifts  $(W_1,\ldots, W_n)$ associated
with the noncommutative domain $\cD_f$. Similar results are deduced
for the model shifts $(B_1,\ldots, B_n)$ associated with the
noncommutative variety $\cV_{f,J}\subset \cD_f$. In Section
\ref{C*-algebras}, we develop a dilation theory for $n$-tuples of
operators in the noncommutative domain $\cD_f(\cH)$, or  in the
noncommutative variety $\cV_{f,J}(\cH)$, defined by
$$
\cV_{f,J}(\cH):=\left\{(X_1,\ldots, X_n)\in \cD_f(\cH):\
q(X_1,\ldots, X_n)=0, \text{ for all } q\in \cP_J\right\},
$$
where $J$ is a $w^*$-closed two-sided ideal of $F_n^\infty(\cD_f)$
generated by a set $\cP_J$  of noncommutative polynomials.

The  case when $(T_1,\ldots, T_n)\in \cD_f(\cH)$ and \
$T_iT_j=T_jT_i$, \ $i,j=1,\ldots n$, is obtained   when the  ideal
$J$ is generated by the  commutators $W_iW_j-W_jW_i$, \
$i,j=1,\ldots n$. In this commutative  case, and assuming that $f$
is a positive regular polynomial, we find some of the results of
S.~Pott \cite{Pot}. More particulary, if $f=X_1+\cdots +X_n$, we
recover some of the results obtained by Drury \cite{Dr},  Arveson
\cite{Arv}, and the author \cite{Po-poisson}.

In Section \ref{Characteristic}, we associate with each $n$-tuple of
$T:=(T_1,\ldots, T_n)\in \cD_f(\cH)$ a characteristic function
$\Theta_{f,T}$, which is a multi-analytic operator with respect to
the universal model $(W_1,\ldots, W_n)$, i.e.,
$$
\Theta_{f,T}\in R_n^\infty(\cD_f)\bar\otimes B(\cD_{C(T)^*},
\cD_{C(T)}),
$$
where $\cD_{C(T)^*}$ and $\cD_{C(T)}$ are some defect spaces
associated with  $T$.
 We prove
that the characteristic function is a {\it complete unitary
invariant} for completely non-coisometric elements of $\cD_f$, and
provide a model. Similar results are obtained for the {\it
constrained characteristic function} associated with noncommutative
varieties $\cV_{f,J}\subset \cD_f$.

In particular, we discuss the commutative case, when the $n$-tuple
of operators $T:=(T_1,\ldots, T_n)\in \cD_f(\cH)$  is completely
non-coisometric  and $T_iT_j=T_jT_i$ for all $i,j=1,\ldots n$. More
particularly, if $f$ is a positive regular polynomial and $T$ is
pure, we recover the results from \cite{BS}.

In the last section of this chapter, we introduce the curvature and
the $*$-curvature associated with $n$-tuples of operators in the
noncommutative domain $\cD_p$, where $p$ is a positive regular
noncommutative polynomial. We prove the existence of these numerical
invariants and present some basic properties. We also show that both
curvatures can be expressed only in terms of the characteristic
function $\Theta_{p,T}$.

The particular case when $p_e:=a_1 X_1+\cdots +a_n X_n$ is treated
in greater details. In this case, the curvature invariant $\text{\rm
curv}_{p_e}$ is defined on the noncommutative ellipsoid
$$
\cD_{p_e}(\cH):=\left\{ (X_1,\ldots, X_n)\in B(\cH)^n: \ a_1
X_1X_1^*+\cdots + a_nX_nX_n^*\leq I\right\}.
$$

We show that the range of the curvature $\text{\rm curv}_{p_e}$
coincides with $[0,\infty)$, and that $\text{\rm curv}_{p_e}$ can
detect the {\it pure} $n$-tuples of operators $(T_1,\ldots, T_n)\in
\cD_{p_e}$ which are unitarily equivalent with the model operator
$(W_1\otimes I_\cK,\ldots, W_n\otimes I_\cK)$, where $\cK$ is finite
dimensional. We remark that in the particular case when  $a_1=\cdots
=a_n=1$, we recover some of the results from \cite{Po-curvature},
\cite{Kr}, \cite{Po-similarity},  and \cite{Po-varieties}.

Chapter 4 deals with commutant lifting theorems  and  applications.
We provide in Section \ref{Interpolation}, a Sarason type \cite{S}
commutant lifting theorem  for pure $n$-tuples of operators in
noncommutative  domains $\cD_f$ or in  noncommutative varieties
$\cV_{f,J}$.  This  theorem also extends the corresponding result
obtained by Arias  \cite{Ar1}.

As consequences, we obtain Nevanlinna-Pick \cite{N} and
Schur-Carath\'eodory (\cite{Ca}, \cite{Sc}) type interpolation
results. In particular, we obtain the following Nevanlinna-Pick
interpolation result. If $\lambda_1,\ldots, \lambda_k$ are distinct
points in $\cD_f^\circ(\CC)$ and $A_1,\ldots, A_k$ are in  $B(\cK)$,
then there exists $\psi\in F_n^\infty(\cD_f)\bar\otimes B(\cK)$ such
that
$$\|\psi\|\leq 1\quad \text{ and }\quad \psi(\lambda_j)=A_j, \
j=1,\ldots, k, $$
 if and only if the operator matrix
$$
\left[(I_\cK-A_iA_j^*)K_f(\lambda_i,\lambda_j)\right]_{k\times k}
$$
is positive semidefinite.

In the last section of this paper  we  present, as an application of
our  commutant lifting theorem, a corona type result \cite{G} for a
class of Hardy algebras associated with the noncommutative domain
$\cD_f$ and the noncommutative variety $\cV_{f,J}$. In particular,
we show that if $\varphi_1,\ldots, \varphi_n\in F_n^\infty(\cD_f)$,
then there exist operators $g_1,\ldots, g_k\in F_n^\infty(\cD_f)$
such that
$$
\varphi_1 g_1+\cdots +\varphi_k g_k=I
$$
if and only if there exists $\delta>0$ such that
$$
\varphi_1 \varphi^*_1+\cdots +\varphi_k \varphi_k^*\geq \delta I.
$$
A commutative version of this result for the algebra
$F_n^\infty(\cV_{f,J_c})$ is obtained as well.

Now,  a few more remarks concerning the results of the present paper
are necessary. First, we mention that  many results of this paper
remain true if  the noncommutative domain $\cD_f(\cH)$   is replaced
by domains of type $\lambda I+\cD_f(\cH)$, where
$\lambda:=(\lambda_1,\ldots,\lambda_n)\in \CC^n$ and $\lambda
I:=(\lambda_1 I,\ldots,\lambda_nI)$. We  remark  that we can obtain,
as particular cases, all the results obtained by Arias and the
author in \cite{ArPo1}, where we considered interpolation problems
for Hardy algebras associated with a certain class of weighted
shifts on Fock spaces.

On the other hand, when $f=X_1+\cdots + X_n$, we recover several
results concerning the multivariable noncommutative (resp.
commutative) operator theory on the unit ball $[B(\cH)^n]_1$.
However,  there are many results for $[B(\cH)^n]_1$ which remain
open problems for our more general noncommutative domains generated
by free holomorphic functions on $B(\cH)^n$. For instance,  it
remains an open problem whether  the model theory for c.n.c
$n$-tuples of operators  in the noncommutative domain $\cD_f(\cH)$
can  be extended to the  class of c.n.u (completely nonunitary)
$n$-tuples, as in the classical Sz.-Nagy-Foias theory or the work of
Ball-Vinnikov in the particular case when $f=X_1+\cdots + X_n$.  A
similar question concerns the $F_n^\infty(\cD_f)$-functional
calculus. Other open problems will be mentioned
 throughout the paper.

It would be interesting to see if there are analogues of our results
for  noncommutative polydomains  $\cD_{f_1}(\cH)\times \cdots \times
\cD_{f_k}(\cH)$ and certain subvarieties determined by
noncommutative polynomials, extending in this way previous results
obtained (in the commutative case) by Curto and Vasilescu
\cite{CV2}, \cite{CV3}.

Finally, we add that several results of this paper are  extended, in
a forthcoming paper \cite{Po-Berezin},  to noncommutative domains
$$
{\bold D}^m_g(\cH):=\left\{ (X_1,\ldots, X_n)\in B(\cH)^n: \
\sum_{\alpha\in \FF_n^+} c_\alpha X_\alpha X_\alpha^*\geq 0\right\},
$$
where $g=\sum_{\alpha\in \FF_n^+} c_\alpha X_\alpha =(1-f)^m$, $m\in
\NN$, and $f$ is any positive regular free holomorphic function on
$B(\cH)^n$. In the commutative case,  we  recover several results
obtained  in \cite{Ag1}, \cite{Ag2}, \cite{M}, \cite{MV}, \cite{Va},
\cite{CV1},  and \cite{Pot}.

\bigskip
\chapter{Operator algebras associated with noncommutative domains}

      \section{The noncommutative domain $\cD_f$ and a universal
      model
} \label{Noncommutative}

In this section, we associate with each  positive regular free
holomorphic function $f$ on $B(\cH)^n$ a noncommutative domain
$\cD_f(\cH)\subset B(\cH)^n$ and an essentially unique $n$-tuple
$(W_1,\ldots, W_n)$ of weighted shifts, which will play the role of
the {\it universal model} for the elements  of $\cD_f(\cH)$.

      Let $H_n$ be an $n$-dimensional complex  Hilbert space with orthonormal
      basis
      $e_1$, $e_2$, $\dots,e_n$, where $n\in\{1,2,\dots\}$.        We consider
      the full Fock space  of $H_n$ defined by
      $$F^2(H_n):=\bigoplus_{k\geq 0} H_n^{\otimes k},$$
      where $H_n^{\otimes 0}:=\CC 1$ and $H_n^{\otimes k}$ is the (Hilbert)
      tensor product of $k$ copies of $H_n$.
      Define the left creation
      operators $S_i:F^2(H_n)\to F^2(H_n), \  i=1,\dots, n$,  by
      $$
       S_i\varphi:=e_i\otimes\varphi, \quad  \varphi\in F^2(H_n),
      $$
      and  the right creation operators
      $R_i:F^2(H_n)\to F^2(H_n)$, \  $i=1,\dots, n$,  by
      $
       R_i\varphi:=\varphi\otimes e_i$, \ $ \varphi\in F^2(H_n)$.

The    algebra   $F_n^\infty$
  and  its norm closed version,
  the noncommutative disc
 algebra  $\cA_n$,  were introduced by the author   \cite{Po-von} in connection
   with a multivariable noncommutative von Neumann inequality.
$F_n^\infty$  is the algebra of left multipliers of $F^2(H_n)$  and
can be identified with
 the
  weakly closed  (or $w^*$-closed) algebra generated by the left creation operators
   $S_1,\dots, S_n$  acting on   $F^2(H_n)$,
    and the identity.
     The noncommutative disc algebra $\cA_n$ is
    the  norm closed algebra generated by
   $S_1,\dots, S_n$,
    and the identity. For basic properties concerning
    the  noncommutative analytic Toeplitz algebra   $F_n^\infty$
  we refer to
\cite{Po-charact},  \cite{Po-multi},  \cite{Po-funct},
\cite{Po-analytic}, \cite{Po-disc}, \cite{Po-poisson},
 \cite{ArPo},
  \cite{DP1}, \cite{DP2},   \cite{DP}, \cite{DKP},
  and \cite{ArPo2}.

Let $\FF_n^+$ be the unital free semigroup on $n$ generators
$g_1,\ldots, g_n$ and the identity $g_0$.  The length of $\alpha\in
\FF_n^+$ is defined by $|\alpha|:=0$ if $\alpha=g_0$  and
$|\alpha|:=k$ if
 $\alpha=g_{i_1}\cdots g_{i_k}$, where $i_1,\ldots, i_k\in \{1,\ldots, n\}$.
If $(X_1,\ldots, X_n)\in B(\cH)^n$, where $B(\cH)$ is the algebra of
all bounded linear operators on the Hilbert space $\cH$,    we
denote $X_\alpha:= X_{i_1}\cdots X_{i_k}$  and $X_{g_0}:=I_\cH$.

Let $f(X_1,\ldots, X_n):= \sum_{\alpha\in \FF_n^+} a_\alpha X_\alpha$,
\ $a_\alpha\in \CC$,  be a free holomorphic function on $B(\cH)^n$  with
radius of convergence strictly positive. As shown in \cite{Po-holomorphic},
this condition is equivalent  to
\begin{equation} \label{lisu}
\limsup_{k\to\infty} \left( \sum_{|\alpha|=k} |a_\alpha|^2\right)^{1/2k}<\infty.
\end{equation}
Throughout this paper, we
  assume that  $a_\alpha\geq 0$ for any $\alpha\in \FF_n^+$, \ $a_{g_0}=0$,
 \ and  $a_{g_i}>0$, $i=1,\ldots, n$.
 A function $f$ satisfying  all these conditions on the coefficients is
 called a {\it positive regular free holomorphic function on}
 $B(\cH)^n$.
  Define  the noncommutative domain
$$
\cD_f=\cD_f(\cH):=\left\{ (X_1,\ldots, X_n)\in B(\cH)^n: \
\left\|\sum_{|\alpha|\geq 1} a_\alpha X_\alpha
X_\alpha^*\right\|\leq I_\cH\right\},
$$
where the convergence of the series is in the weak operator
topology. When the Hilbert space $\cH$ is understood, we use the
notation $\cD_f$. Due to the Schwarz type lemma for free holomorphic
functions on the open unit ball of $B(\cH)^n$ (see
\cite{Po-holomorphic}), there exists $r>0$ such that $\cD_f(\cH)$
contains the operatorial $n$-ball
$$
[B(\cH)^n]_r:=\left\{ (X_1,\ldots, X_n)\in B(\cH)^n: \ \|X_1
X_1^*+\cdots + X_n X_n^*\|^{1/2}<r\right\}
$$
  and  $f(rS_1,\ldots, rS_n)$ is a strict contraction in the noncommutative disc algebra
  $\cA_n$, i.e.,  $\|f(rS_1,\ldots, rS_n)\|<1$.
Therefore, the operator $I-f(rS_1,\ldots, rS_n)$ is invertible with
its inverse
 $(I-f(rS_1,\ldots, rS_n))^{-1}$  in $\cA_n\subset F_n^\infty$
 and has the ``Fourier representation''
$g(rS_1,\ldots, rS_n)= \sum_{\alpha\in \FF_n^+} b_\alpha
r^{|\alpha|} S_\alpha $ for some constants $b_\alpha\in \CC$. Hence,
 and  using the fact that $r^{|\alpha|}b_\alpha=P_\CC S_\alpha^*
g(rS_1,\ldots, rS_n)(1)$, we  deduce that
\begin{equation*}
\begin{split}
g(rS_1,\ldots, rS_n)&= I+ f(rS_1,\ldots, rS_n)+ f(rS_1,\ldots, rS_n)^2+\cdots \\
&=I+\sum_{k=1}^\infty \sum_{|\gamma_1|\geq 1,\ldots, |\gamma_k|\geq
1}
 a_{\gamma_1}\cdots a_{\gamma_k}
r^{|\gamma_1|+\cdots +|\gamma_k|}S_{\gamma_1\cdots \gamma_k}\\
&=I+\sum_{m=1}^\infty \sum_{|\alpha|=m}\left(\sum_{j=1}^{|\alpha|}
\sum_{{\gamma_1\cdots \gamma_j=\alpha }\atop {|\gamma_1|\geq 1,\ldots,
 |\gamma_j|\geq 1}} a_{\gamma_1}\cdots a_{\gamma_j} \right) r^{|\alpha|}
  S_\alpha.
\end{split}
\end{equation*}
Due to the uniqueness of the Fourier representation of the elements
in $F_n^\infty$,  we have
\begin{equation}\label{b_alpha}
b_{g_0}=1 \quad \text{ and }\quad b_\alpha= \sum_{j=1}^{|\alpha|}
\sum_{{\gamma_1\cdots \gamma_j=\alpha }\atop {|\gamma_1|\geq
1,\ldots, |\gamma_j|\geq 1}} a_{\gamma_1}\cdots a_{\gamma_j}   \quad
\text{ if } \ |\alpha|\geq 1.
\end{equation}
Since
 $$b_\alpha= a_{g_{i_1}}\cdots a_{g_{i_k}}+   \text{\it  positive
 terms},
 $$
  for any $\alpha=g_{i_1}\cdots g_{i_k}$
 and
 $a_{g_{i_1}}>0, \cdots, a_{g_{i_k}}>0$,   relation \eqref{b_alpha}
 implies  $b_\alpha>0$ for any  $\alpha\in \FF_n^+$.
Notice that any term of the sum
$$
b_\alpha b_\beta =\sum_{j=1}^{|\alpha|} \sum_{k=1}^{|\beta|}
\sum_{{\gamma_1\cdots \gamma_j=\alpha }\atop {|\gamma_1|\geq
1,\ldots, |\gamma_j|\geq 1}}\left( \sum_{{\sigma_1\cdots
\sigma_k=\beta }\atop {|\sigma_1|\geq 1,\ldots, |\sigma_k|\geq 1}}
a_{\gamma_1}\cdots a_{\gamma_j}
  a_{\sigma_1}\cdots a_{\sigma_k} \right)
$$
is also a term of  the sum
$$
b_{\alpha \beta}= \sum_{p=1}^{|\alpha|+|\beta|}
\sum_{{\epsilon_1\cdots \epsilon_p=\alpha\beta}
\atop{|\epsilon_1|\geq 1,\ldots,
 |\epsilon_p|\geq 1}} a_{\epsilon_1}\cdots a_{\epsilon_p}.
$$
 Since $a_\alpha\geq 0$, $|\alpha|\geq 1$, we deduce that

  \begin{equation}
\label{bbb}
  b_\alpha b_\beta \leq b_{\alpha\beta},\quad \text{ for any }\alpha,\beta\in
   \FF_n^+.
\end{equation}
On the other hand, since
$(I- f(rS_1,\ldots, rS_n)) g(rS_1,\ldots, rS_n)=I,
$
we deduce that
$$
I+\sum_{m=1}^\infty \sum_{|\gamma|=m} r^{|\gamma|}\left(b_\gamma-
\sum_{{\beta \alpha=\gamma}\atop{\alpha,\beta\in \FF_n^+,
|\beta|\geq 1}}a_\beta b_\alpha\right) S_\gamma = I
$$
and, therefore,
\begin{equation}\label{sum1}
b_\gamma-\sum_{\beta \alpha=\gamma, |\beta|\geq 1}a_\beta b_\alpha=0
\quad \text{ if } \ |\gamma|\geq 1.
\end{equation}
Similarly,  since  $g(rS_1,\ldots, rS_n)(I- f(rS_1,\ldots,
rS_n))=I$, we  deduce that
 \begin{equation}\label{sum2}
b_\gamma-\sum_{ \alpha\beta=\gamma, |\beta|\geq 1}a_\beta b_\alpha=0
 \quad \text{ if } \ |\gamma|\geq 1.
\end{equation}
%

%On the othe hand, if $p<|\gamma|$, then
%\begin{equation} \label{-ab}
%\sum_{\alpha \beta =\gamma, 1\leq |\alpha|\leq p} a_\alpha %b_\beta
%\leq 0\quad \text{ and }\quad  \sum_{ \beta \alpha %=\gamma, 1\leq
%|\alpha|\leq p} a_\alpha b_\beta \geq 0.
%\end{equation}

We define the diagonal operators $D_i:F^2(H_n)\to F^2(H_n)$,
$i=1,\ldots, n$, by setting
$$
D_ie_\alpha=\sqrt{\frac{b_\alpha}{b_{g_i \alpha}}} e_\alpha,\quad
 \alpha\in \FF_n^+.
$$
Due to relations   \eqref{b_alpha} and \eqref{bbb}, we have
$b_{g_i\alpha} \geq b_{g_i} b_\alpha= a_{g_i} b_\alpha$ and, consequently,

$$
\|D_i\|=\sup_{\alpha\in \FF_n^+} \sqrt{\frac{b_\alpha}{b_{g_i \alpha}}}=
\frac{1}{\sqrt{a_{g_i}}}, \quad i=1,\ldots,n.
$$
Now we define the {\it weighted left creation  operators}
$W_i:F^2(H_n)\to F^2(H_n)$, $i=1,\ldots, n$,  associated with the
 noncommutative domain $\cD_f$  by setting $W_i=S_iD_i$, where
 $S_1,\ldots, S_n$ are the left creation operators on the full
 Fock space $F^2(H_n)$.
Therefore, we have
$$W_i e_\alpha=\frac {\sqrt{b_\alpha}}{\sqrt{b_{g_i \alpha}}}
e_{g_i \alpha}, \quad \alpha\in \FF_n^+,
$$
where the coefficients  $b_\alpha$, $\alpha\in \FF_n^+$, are given
by relation \eqref{b_alpha}. Let's show that $(W_1,\ldots, W_n)\in
\cD_f(F^2(H_n))$, i.e.,
\begin{equation}\label{WWinD}
\sum_{|\beta|\geq 1} a_\beta W_\beta W_\beta^*\leq I.
\end{equation}
A simple calculation reveals that
\begin{equation}\label{WbWb}
W_\beta e_\gamma=
\frac {\sqrt{b_\gamma}}{\sqrt{b_{\beta \gamma}}}
e_{\beta \gamma} \quad
\text{ and }\quad
W_\beta^* e_\alpha =\begin{cases}
\frac {\sqrt{b_\gamma}}{\sqrt{b_{\alpha}}}e_\gamma& \text{ if }
\alpha=\beta\gamma \\
0& \text{ otherwise }
\end{cases}
\end{equation}
 for any $\alpha, \beta \in \FF_n^+$.
Notice also that $\|W_\beta\|=\frac{1}{\sqrt{b_\beta}}$, $\beta\in \FF_n^+$.
 Indeed, due to inequality \eqref{bbb}, if
$f=\sum_{\alpha\in \FF_n^+} c_\alpha e_\alpha\in F^2(H_n)$,
then
$$
\|W_\beta f\|=\left( \sum_{\alpha\in \FF_n^+} |c_\alpha|^2
\frac{b_\alpha}{b_{\beta\alpha}}\right)^{1/2}\leq \frac{1}{\sqrt{b_\beta}}
\|f\|.
$$
Since $W_\beta 1= \frac{1}{\sqrt{b_\beta}}$, the assertion follows.
 Now, using  relation \eqref{WbWb}, we deduce that
$$
W_\beta W_\beta^* e_\alpha =\begin{cases}
\frac {{b_\gamma}}{{b_{\alpha}}} e_\alpha & \text{ if } \alpha=\beta\gamma\\
0& \text{ otherwise. }
\end{cases}
$$
Notice that
$\left(I-\sum\limits_{1\leq |\beta|\leq N} a_\beta W_\beta W_\beta^*\right)
 e_\alpha= C_{N, \alpha} e_\alpha,
$
 where $C_{N, \alpha}=1$ if $\alpha=g_0$, and
$$
C_{N, \alpha}=1- \sum_{\beta\gamma=\alpha, 1\leq |\beta|\leq N}
\frac{a_\beta b_\gamma}{b_\alpha}
\qquad \text{if } \quad |\alpha|\geq 1.
$$
Due to relation \eqref{sum1}, if $1\leq|\alpha|\leq N$, we have
$C_{N, \alpha}=0$.
Using  the same relation and the fact that
\begin{equation*} b_\alpha-\sum_{\beta \gamma =\alpha, 1\leq |\beta|\leq N}
 a_\beta b_\gamma\leq b_\alpha,\quad |\alpha| \geq 1,
\end{equation*}
we deduce that  $0\leq C_{N, \alpha}\leq 1$,  whenever $|\alpha|>N$.
  On the  other hand,
notice that if  $1\leq N_1\leq N_2\leq |\alpha|$, then
\begin{equation}\label{bga}
b_\alpha-\sum_{\beta \gamma =\alpha, 1\leq |\beta|\leq N_2} a_\beta b_\gamma\leq
b_\alpha-\sum_{\beta \gamma =\alpha, 1\leq |\beta|\leq N_1} a_\beta b_\gamma.
\end{equation}
Consequently, $\left\{I-\sum\limits_{1\leq |\beta|\leq N} a_\beta W_\beta
 W_\beta^*\right\}_{N=1}^\infty$ is a decreasing
sequence of positive diagonal operators which  converges in the
strong operator topology to $P_\CC$, the orthogonal projection of
$F^2(H_n)$ onto $1\otimes \CC$. Therefore, we have
\begin{equation}\label{proj}
I-\sum_{ |\beta|\geq 1} a_\beta W_\beta W_\beta^*=P_\CC.
\end{equation}
This also shows that $(W_1,\ldots, W_n)\in \cD_f(F^2(H_n))$. Now,
since
\begin{equation}\label{projW*}
P_\CC W_\beta^* e_\alpha =\begin{cases}
\frac {1}{\sqrt{b_{\beta}}}   & \text{ if } \alpha=\beta\\
0& \text{ otherwise, }
\end{cases}
\end{equation}
we have $\sum_{\beta\in \FF_n^+} b_\beta W_\beta P_\CC W_\beta^*
e_\alpha= e_\alpha$.  Therefore,
\begin{equation}
\label{WWWW}
\sum_{\beta\in \FF_n^+}b_\beta W_\beta\left(I-\sum_{ |\alpha|\geq 1}
a_\alpha W_\alpha W_\alpha^*\right)   W_\beta^* =I,
\end{equation}
where the convergence is in the strong operator topology.

We can also define the {\it weighted right creation operators}
$\Lambda_i:F^2(H_n)\to F^2(H_n)$ by setting
$\Lambda_i:= R_i G_i$, $i=1,\ldots, n$,  where $R_1,\ldots, R_n$ are
 the right creation operators on the full Fock space $F^2(H_n)$ and
 each  diagonal operator $G_i$  is defined by
$$
G_ie_\alpha=\sqrt{\frac{b_\alpha}{b_{ \alpha g_i}}} e_\alpha,\quad
 \alpha\in \FF_n^+,
$$
where the coefficients $b_\alpha$, $\alpha\in \FF_n^+$, are given
by relation \eqref{b_alpha}.
In this case, we have
\begin{equation}\label{WbWb-r}
\Lambda_\beta e_\gamma=
\frac {\sqrt{b_\gamma}}{\sqrt{b_{ \gamma \tilde\beta}}}
e_{ \gamma \tilde \beta} \quad
\text{ and }\quad
\Lambda_\beta^* e_\alpha =\begin{cases}
\frac {\sqrt{b_\gamma}}{\sqrt{b_{\alpha}}}e_\gamma& \text{ if }
\alpha=\gamma \tilde \beta \\
0& \text{ otherwise }
\end{cases}
\end{equation}
 for any $\alpha, \beta \in \FF_n^+$, where $\tilde \beta$ denotes
 the reverse of $\beta=g_{i_1}\cdots g_{i_k}$, i.e.,
 $\tilde \beta=g_{i_k}\cdots g_{i_1}$.
Using relations \eqref{b_alpha}, \eqref{bbb}, and \eqref{sum2},
we deduce that
$$\left(I-\sum\limits_{1\leq |\beta|\leq N}
a_{\tilde\beta} \Lambda_\beta \Lambda_\beta^*\right) e_\alpha=
C_{N, \alpha} e_\alpha,
$$
 where $C_{N, \alpha}=1$ if $\alpha=g_0$, and
$$
C_{N, \alpha}=1- \sum_{\gamma\tilde\beta=\alpha,
 1\leq |\tilde\beta|\leq N} \frac{a_{\tilde\beta} b_\gamma}{b_\alpha}
\qquad \text{if } \quad |\alpha|\geq 1.
$$
As in the case of weighted left creation operators, one can show
that \begin{equation} \label{proj2} I-\sum\limits_{|\beta|\geq 1}
a_{\tilde\beta} \Lambda_\beta
 \Lambda_\beta^*=P_\CC
 \end{equation}
and,
since
$$
P_\CC \Lambda_\beta^* e_\alpha =\begin{cases}
\frac {1}{\sqrt{b_{\alpha}}}   & \text{ if } \alpha=\tilde\beta\\
0& \text{ otherwise, }
\end{cases}
$$
we have
\begin{equation*}
\sum_{\beta\in \FF_n^+}b_{\tilde\beta} \Lambda_\beta
\left(I-\sum_{ |\alpha|\geq 1} a_{\tilde\alpha} \Lambda_\alpha
\Lambda_\alpha^*\right)   \Lambda_\beta^* =I,
\end{equation*}
where the convergence is in the strong operator topology.

Notice that if  $f(X_1,\ldots, X_n):=\sum_{|\alpha|\geq 1} a_\alpha
X_\alpha$ is a positive regular free holomorphic function on
$B(\cH)^n$, then so is $\tilde{f}(X_1,\ldots,
X_n):=\sum_{|\alpha|\geq 1} a_{\tilde \alpha} X_\alpha$, where
$\tilde \alpha$ denotes the reverse of $\alpha$. We also denote by
$(W_1^{(f)},\ldots, W_n^{(f)})$ the weighted left creation operators
$(W_1,\ldots, W_n)$ associated with $\cD_f$. The notation
$(\Lambda_1^{(f)},\ldots, \Lambda_n^{(f)})$ is now  clear.

\begin{proposition}\label{tilde-f}
Let $(W_1^{(f)},\ldots, W_n^{(f)})$ (resp. $(\Lambda_1^{(f)},\ldots,
\Lambda_n^{(f)})$) be the weighted left (resp. right) creation
operators associated with the noncommutative domain  $\cD_f$. Then
the following statements hold:
\begin{enumerate}
\item[(i)] $(W_1^{(f)},\ldots, W_n^{(f)})\in \cD_f(F^2(H_n));$
\item[(ii)] $(\Lambda_1^{(f)},\ldots,
\Lambda_n^{(f)})\in \cD_{\tilde f}(F^2(H_n));$
\item[(iii)] $U^* \Lambda_i^{(f)} U=W_i^{(\tilde f)}$,
$i=1,\ldots,n$, where $U\in B(F^2(H_n))$ is the unitary operator
defined by $U e_\alpha:= e_{\tilde \alpha}$, $\alpha\in \FF_n^+.$
\end{enumerate}
\end{proposition}
\begin{proof}
Items (i) and (ii) follow from relation \eqref{WWinD} and
\eqref{proj2}, respectively. Using relation \eqref{WbWb} when $f$ is
replaced by $\tilde f$, we obtain
$$
W_i^{(\tilde f)} e_\gamma=\frac{\sqrt{b_{\tilde
\gamma}}}{\sqrt{b_{\tilde\gamma g_i}}} e_{g_i\gamma}.
$$
On the other hand, using relation \eqref{WbWb-r}, we deduce that
$$
U^*\Lambda_i^{(f)} U e_\gamma =U^*\left(\frac{\sqrt{b_{\tilde
\gamma}}}{\sqrt{b_{\tilde\gamma g_i}}} e_{\tilde\gamma g_i}\right)=
\frac{\sqrt{b_{\tilde \gamma}}}{\sqrt{b_{\tilde\gamma g_i}}}
e_{g_i\gamma}.
$$
Therefore,  $U^* \Lambda_i^{(f)} U=W_i^{(\tilde f)}$,
$i=1,\ldots,n$. The proof is complete.
\end{proof}

Now, we ca prove that the weighted left creation operators
associated with noncommutative domains are essentially unique.

\begin{theorem}
\label{unitar-equiv} Let $(W_1,\ldots, W_n)$ and $(W_1',\ldots,
W_n')$ be the weighted left creation operators  associated with the
noncommutative domains $\cD_f$ and $\cD_{f'}$, respectively. Then
there exists a unitary operator $U\in B(F^2(H_n))$ such that
$$
UW_i=W_i'U,\quad i=1,\ldots, n,
$$
if and only if $f=f'$.
\end{theorem}
\begin{proof}
One implication is obvious. Let $[u_{\alpha,
\beta}]_{\alpha,\beta\in \FF_n^+}$ be the matrix representation of
the unitary operator $U$, i.e. $u_{\alpha, \beta}:=\left< Ue_\beta,
e_\alpha\right>$, $\alpha,\beta\in \FF_n^+$. Using the definition of
the weighted left creation operators, we deduce that
$$
\left<UW_ie_\alpha,
e_\beta\right>=\frac{\sqrt{b_\alpha}}{\sqrt{b_{g_i \alpha}}}\left<U
e_{g_i\alpha}, e_\beta\right>=\frac{\sqrt{b_\alpha}}{\sqrt{b_{g_i
\alpha}}} u_{\beta, g_i\alpha}
$$
and
\begin{equation*}
\begin{split}
\left<W_i'U e_\alpha, e_\beta\right>&=\left<U e_\alpha, {W_i'}^*
e_\beta\right>\\
&=\begin{cases}\frac{\sqrt{b_\gamma'}}{\sqrt{b_{\beta}'}},&\quad
\text{ if } \beta=g_i\gamma\\
0,&\quad \text{ otherwise}
\end{cases}
\end{split}
\end{equation*}
for any $\alpha,\beta\in \FF_n^+$ and $i=1,\ldots, n$. Therefore, if
$\beta=g_i \gamma$ for some $\gamma\in \FF_n^+$ and $i=1,\ldots, n$,
then the equality $UW_i=W_i'U$ implies
\begin{equation}\label{rad-rad}
\frac{\sqrt{b_{g_i\gamma}'}}{\sqrt{b_{g_i \alpha}}} u_{g_i\gamma,
g_i\alpha}= \frac{\sqrt{b_{\gamma}'}}{\sqrt{b_{\alpha}}} u_{\gamma,
\alpha}.
\end{equation}
If $\beta$ is a word which does not start with the generator $g_i$,
then the same equality implies $\frac{\sqrt{b_\alpha}}{\sqrt{b_{g_i
\alpha}}} u_{\beta, g_i\alpha}=0,$ which implies
\begin{equation}
\label{u=0} u_{\beta, g_i\alpha}=0,\quad \text{ if $\beta$ does not
start with $g_i$}.
\end{equation}
In particular, $u_{g_0, \omega}=0$ for any $\omega\in \FF_n^+$ with
$|\omega|\geq 1$. Since $U$ is a unitary operator, we must have
$u_{g_0,g_0}\neq 0$. Similarly, using relations  \eqref{rad-rad} and
\eqref{u=0}, we deduce that $u_{\beta,\gamma}=0$ if $\beta\neq
\gamma$ and $u_{\beta,\beta}\neq 0$ for any $\beta\in \FF_n^+$. Let
$U^*=[c_{\alpha,\beta}]_{\alpha,\beta\in \FF_n^+}$ be the matrix
representation of $U^*$, where
$\overline{c}_{\beta,\alpha}=u_{\alpha,\beta}$. Since $U^*
W_i'=W_iU^*$ for $i=1,\ldots,n$, similar arguments as those used to
obtain relations \eqref{rad-rad} and \eqref{u=0}, imply
\begin{equation}
\label{rad-rad'} \frac{\sqrt{b_{g_i\gamma}}}{\sqrt{b_{g_i \alpha}'}}
c_{g_i\gamma, g_i\alpha}=
\frac{\sqrt{b_{\gamma}}}{\sqrt{b_{\alpha}'}} c_{\gamma, \alpha}.
\end{equation}
Using relations \eqref{rad-rad} and \eqref{rad-rad'}, we obtain
\begin{equation*}
\begin{split}
u_{g_i\gamma, g_i\alpha}&=\frac{\sqrt{b_{g_i
\alpha}}}{\sqrt{b_{g_i\gamma}'}}
\frac{\sqrt{b_{\gamma}'}}{\sqrt{b_{\alpha}}} u_{\gamma, \alpha}\\
&=\overline{c}_{g_i\alpha, g_i\gamma} =\frac{\sqrt{b_{g_i
\gamma}'}}{\sqrt{b_{g_i\alpha}}}
\frac{\sqrt{b_{\alpha}}}{\sqrt{b_{\gamma}'}}\overline{c}_{\alpha,
\gamma}\\
&=\frac{\sqrt{b_{g_i \gamma}'}}{\sqrt{b_{g_i\alpha}}}
\frac{\sqrt{b_{\alpha}}}{\sqrt{b_{\gamma}'}} u_{\gamma,\alpha}
\end{split}
\end{equation*}
for any $\alpha, \gamma\in \FF_n^+$ and $i=1,\ldots,n$. Hence, we
deduce that
$$
\frac{\sqrt{b_{g_i \alpha}}}{\sqrt{b_{g_i\gamma}'}}
\frac{\sqrt{b_{\gamma}'}}{\sqrt{b_{\alpha}}} = \frac{\sqrt{b_{g_i
\gamma}'}}{\sqrt{b_{g_i\alpha}}}
\frac{\sqrt{b_{\alpha}}}{\sqrt{b_{\gamma}'}}
$$
if $ u_{\gamma,\alpha}\neq 0$, and, consequently,
$\frac{b_{g_i\alpha}}{b_\alpha}=\frac{b_{g_i\gamma}'}{b_\gamma'}$.
Since we already know that $u_{\alpha,\alpha}\neq 0$, we deduce that
$$
\frac{b_{g_i\alpha}}{b_\alpha}=\frac{b_{g_i\alpha}'}{b_\alpha'}
$$
for any $\alpha\in \FF_n^+$ and $i=1,\ldots,n$. Now, using the fact
that $b_{g_0}=b_{g_0}'=1$ ( see \eqref{b_alpha}), one  can easily
show that $b_\omega=b_\omega'$ for any $\omega\in \FF_n^+$. Based on
the results
  preceding the theorem, we have
\begin{equation*}
\begin{split}
(I-f(rS_1,\ldots, rS_n))^{-1}&=g(rS_1,\ldots, rS_n)=\sum_{\omega\in
\FF_n^+} b_\alpha r^{|\omega|} S_\omega\\
&=\sum_{\omega\in \FF_n^+} b_\alpha' r^{|\omega|}
S_\omega=g'(rS_1,\ldots, rS_n)\\
&=(I-f'(rS_1,\ldots, rS_n))^{-1}.
\end{split}
\end{equation*}
Hence, we have $f(rS_1,\ldots, rS_n)=f'(rS_1,\ldots, rS_n)$ which,
due to the uniqueness of the Fourier representation of the  elements
of the noncommutative analytic algebra $F_n^\infty$, implies $f=f'$.
This completes the proof.
\end{proof}

Now, we can prove the following similarity result.

\begin{theorem}
\label{similar-equiv} Let $(W_1,\ldots, W_n)$ and $(W_1',\ldots,
W_n')$ be the weighted left creation operators  associated with the
noncommutative domains $\cD_f$ and $\cD_{f'}$, respectively. Then
there exists an invertible operator $X\in B(F^2(H_n))$ such that
$$
W_i=X^{-1}W_i'X,\quad i=1,\ldots, n,
$$
if and only if there exist positive constants $C_1$ and $C_2$ such
that

\begin{equation}\label{C1C2}
 0<C_1\leq \frac{b_\alpha}{b_\alpha'}\leq C_2
\quad \text{ for all } \ \alpha\in \FF_n^+,
\end{equation}
where the coefficients $b_\alpha$ and $b_\alpha'$ are given by
relation \eqref{b_alpha}. Moreover the operator $X$ can be chosen
such that
$$
\max\left\{\|X\|, \|X^{-1}\| \right\}\leq \max\left\{ \frac{1}{C_1},
C_2\right\}.
$$
\end{theorem}
\begin{proof}
Assume that $X\in B(F^2(H_n))$ is a positive invertible operator
such that $XW_i=W_i' X$, $i=1,\ldots,n$, and let
$X=[x_{\alpha,\beta}]_{\alpha,\beta\in \FF_n^+}$ be its matrix
representation. For each $\alpha\in \FF_n^+$ and $i=1,\ldots,n$, we
denote
\begin{equation}
\label{weig} w(g_i,
\alpha):=\frac{\sqrt{b_\alpha}}{\sqrt{b_{g_i\alpha}}},
\end{equation}
where the coefficients $b_\alpha$ are given by equation
\eqref{b_alpha}. Using the definition of the weighted left creation
operators, we have $\left<XW_ie_\alpha, e_\beta\right>= w(g_i,
\alpha) x_{\beta, g_i\alpha}$ and
\begin{equation*}
\begin{split}
\left<W_i'X e_\alpha, e_\beta\right>&=\left<X e_\alpha, {W_i'}^*
e_\beta\right>\\
&=\begin{cases} w'(g_i,\gamma)x_{\gamma,\alpha},&\quad
\text{ if } \beta=g_i\gamma\\
0,&\quad \text{ otherwise}.
\end{cases}
\end{split}
\end{equation*}
Since $XW_i=W_i' X$, we deduce that $x_{\beta, g_i\alpha}=0$ if
$\beta$ does not start with the generator $g_i$, and
\begin{equation}
\label{xgg}
x_{g_i\gamma,g_i\alpha}=\frac{w'(g_i,\gamma)}{w(g_i,\alpha)}
x_{\gamma,\alpha}\qquad \text{ if } \ \beta=g_i\gamma.
\end{equation}
Hence, $x_{g_0,\alpha}=0$ if $|\alpha|\geq 1$. Since $X$ is
invertible, we must have $x_{g_0,g_0}\neq 0$. Using relations
\eqref{weig} and \eqref{xgg}, we deduce that, for any
$\omega=g_{i_1}g_{i_2}\cdots g_{i_k}\in \FF_n^+$,
\begin{equation*}
\begin{split}
x_{\omega,\omega}&= \frac{w'(g_{i_1},g_{i_2}\cdots g_{i_k})}
{w(g_{i_1},g_{i_2}\cdots g_{i_k})} \frac{w'(g_{i_2},g_{i_3}\cdots
g_{i_k})} {w(g_{i_2},g_{i_3}\cdots g_{i_k})}\cdots \frac{w'(
g_{i_k},
g_0)} {w( g_{i_k}, g_0)}x_{g_0,g_0}\\
&=\frac{\frac{\sqrt{b_{g_0}}}{ \sqrt{b_{\omega}}}}
{\frac{\sqrt{b_{g_0}'}}{\sqrt{b_{\omega}'}}}x_{g_0,g_0}=
\frac{\sqrt{b_\omega'}}{\sqrt{b_\omega}}x_{g_0,g_0}.
\end{split}
\end{equation*}
Hence, we obtain
$$
\frac{\sqrt{b_\omega'}}{\sqrt{b_\omega}}
=\frac{x_{\omega,\omega}}{x_{g_0,g_0}}=\frac{\left<Xe_\omega,
e_\omega\right>}{x_{g_0,g_0}} \leq \frac{\|X\|}{x_{g_0,g_0}}.
$$
Consequently, there is a constant $C_1>0$ such that $C_1\leq
\frac{\sqrt{b_\omega}}{\sqrt{b_\omega'}}$ for any $\omega\in
\FF_n^+$. Similarly, if we work with $X^{-1}$ instead of $X$ and
with the equality $X^{-1} W_i'=W_iX^{-1}$, $i=1,\ldots, n$, we find
a constant $C_2>0$ such that
$\frac{\sqrt{b_\omega}}{\sqrt{b_\omega'}}\leq C_2$ for any
$\omega\in \FF_n^+$.

Conversely, assume that condition \eqref{C1C2} holds. Define the
diagonal operator $D\in B(F^2(H_n))$ by setting $De_\alpha= d_\alpha
e_\alpha$, $\alpha\in \FF_n^+$. As in the first part of the proof,
the relation $DW_iW_i'D$, $i=1,\ldots,n$, holds if and only if
$d_\alpha=\frac{\sqrt{b_\alpha'}}{\sqrt{b_\alpha}} d_{g_0}$,
$\alpha\in \FF_n^+$. We can choose $d_{g_0}=1$. It is clear that
condition \eqref{C1C2} implies that the operator $D$ is invertible.
Moreover, we have
$$
\|D\|=\sup_{\alpha\in
\FF_n^+}\frac{\sqrt{b_\alpha'}}{\sqrt{b_\alpha}}\leq \frac{1}{C_1}
\quad  \text{ and }\quad \|D^{-1}\|=\sup_{\alpha\in
\FF_n^+}\frac{\sqrt{b_\alpha}}{\sqrt{b_\alpha'}}\leq C_2, $$
 which
proves the last part of the theorem.
\end{proof}

We remark that  results similar to Theorem \ref{unitar-equiv} and
Theorem \ref{similar-equiv} can be obtained for the  weighted right
creation operators $(\Lambda_1,\ldots, \Lambda_n)$. We should
mention that, in the particular case when $n=1$, we find some
results from Shields' survey paper \cite{Sh}.

%\bigskip

      \section{The domain algebra $\cA_n(\cD_f)$ and its
      representations
  }
\label{domain algebra}

We introduce the  domain  algebra $\cA_n(\cD_f)$ associated with the
noncommutative domain $\cD_f$ to be the norm closure of all
polynomials in the weighted left creation operators  $W_1,\ldots,
W_n$ and the identity. Using noncommutative Poisson transforms
associated with $\cD_f$, we characterize the completely contractive
(resp. completely bounded) representations of the domain algebra
$\cA_n(\cD_f)$. In particular, we identify the set of all characters
of $\cA_n(\cD_f)$ with the set $\cD_f(\CC)\subset \CC^n$.

 Using the weighted right creation operators
 associated with $\cD_f$, one can   define    the corresponding
  domain algebra ${\cR}_n(\cD_f)$ and obtain similar results.

Let $f(X_1,\ldots, X_n):=\sum_{|\alpha|\geq 1} a_\alpha X_\alpha$ be
a positive regular free holomorphic function on $B(\cH)^n$,
 and let $T:=(T_1,\ldots, T_n)$ be an $n$-tuple of operators in the
noncommutative domain $ \cD_f(\cH)$, i.e.,
$\sum\limits_{|\alpha|\geq 1} a_\alpha T_\alpha T_\alpha^*\leq
I_\cH$.
 Define the positive linear mapping
$$\Phi_{f,T}:B(\cH)\to
B(\cH) \quad \text{  by  } \quad
\Phi_{f,T}(X)=\sum\limits_{|\alpha|\geq 1} a_\alpha T_\alpha
XT_\alpha^*,
$$
where the convergence is in the weak operator topology.
Since $\Phi_{f,T}(I)\leq I$ and $\Phi_{f,T}(\cdot)$ is a positive linear map,
 it is easy to see that $\{\Phi_{f,T}^m(I)\}_{m=1}^\infty$ is a decreasing
  sequence of positive operators and, consequently,
$Q_{f,T}:=\text{\rm SOT-}\lim\limits_{m\to\infty} \Phi_{f,T}^m(I)$ exists.
We say that $T$  is of class $C_{\cdot 0}$ (or {\it pure}) if
~$\text{\rm SOT-}\lim \limits_{m\to\infty} \Phi_{f,T}^m(I)=0$.

We remark that, for any $(T_1,\ldots, T_n)\in \cD_f(\cH)$ and $0\leq
r<1$, the $n$-tuple $(rT_1,\ldots, rT_n)\in \cD_f(\cH)$ is of class
$C_{\cdot 0}$. Indeed, it is enough to see  that
$$\Phi_{f,rT}^m(I)\leq r^m \Phi_{f,T}^m(I)\leq r^m I$$ for any  $m\in
\NN$. Notice also that if $\|\Phi_{f,T}(I)\|<1$, then $T$ is of
class $C_{\cdot 0}$. This is due to the fact that
$\|\Phi_{f,T}^m(I)\|\leq \|\Phi_{f,T}(I)\|^m$.

Following \cite{Po-poisson} (see also \cite{ArPo1}), we define the
Poisson kernel associated with the $n$-tuple $T:=(T_1,\ldots,
T_n)\in \cD_f(\cH)$ to be  the operator $K_{f,T}:\cH\to
F^2(H_n)\otimes \overline{\Delta_{f,T}(\cH)}$ defined by
\begin{equation}
\label{Po-ker}
 K_{f,T}h=\sum_{\alpha\in \FF_n^+} \sqrt{b_\alpha}
e_\alpha\otimes \Delta_{f,T} T_\alpha^* h,\quad h\in \cH,
\end{equation}
where
$\Delta_{f,T}:=\left(I- \Phi_{f,T}(I) \right)^{1/2}$.
Due to relation \eqref{b_alpha}, we  deduce that
\begin{equation*}
\begin{split}
&\left< \sum_{\beta\in \FF_n^+} b_\beta T_\beta \Delta_{f,T}^2T_\beta^*h,
h\right>
 =
\left<\Delta_{f,T}^2 h,h\right>+
\sum_{m=1}^\infty \sum_{|\beta|=m}\left<
b_\beta T_\beta \Delta_{f,T}^2T_\beta^*h, h\right>\\
& =\left<\Delta_{f,T}^2 h,h\right>+ \sum_{m=1}^\infty
\sum_{|\beta|=m}\left<\left(
\sum_{j=1}^{|\beta|}\sum_{{\gamma_1\cdots
\gamma_j=\beta}\atop{|\gamma_1|\geq 1,\ldots, |\gamma_j|\geq 1}}
a_{\gamma_1}\cdots a_{\gamma_j}\right) T_{\gamma_1\cdots \gamma_j}
\Delta_{f,T}^2 T_{\gamma_1\cdots \gamma_j}^*h,h\right>\\
& =\left<\Delta_{f,T}^2 h,h\right>+
\sum_{k=1}^\infty \left<\Phi_{f,T}^k(I-\Phi_{f,T}(I))h,h\right>\\
& =\|h\|^2 -\lim_{m\to\infty} \left<\Phi_{f,T}^m(I)h,h\right>
\end{split}
\end{equation*}
for any $h\in \cH$.
Therefore, we have
\begin{equation}\label{I-Q}
\sum_{\beta\in \FF_n^+} b_\beta T_\beta
\Delta_{f,T}^2T_\beta^*=I_\cH-Q_{f,T}.
\end{equation}
Due to the above calculations, we have
$$
\|K_{f,T}h\|= \sum_{\beta\in \FF_n^+} b_\beta \left<T_\beta
\Delta_{f,T}^2T_\beta^*h, h\right>=\|h\|^2-\|Q^{1/2}_{f,T}h\|^2
$$
for any $h\in \cH$. Therefore $K_{f,T}$ is a contraction and
\begin{equation}
\label{K*K} K_{f,T}^* K_{f,T}=I_\cH-Q_{f,T}.
\end{equation}
 On the other hand, a simple calculation reveals that
\begin{equation}
\label{ker-inter}
K_{f,T} T_i^*=(W_i^*\otimes I_\cH)K_{f,T},\quad i=1,\ldots, n,
\end{equation}
where $(W_1,\ldots, W_n)$ is the $n$-tuple of weighted left creation operators
  associated with the noncommutative domain $\cD_f$.
Moreover, $K_{f,T}$ is an isometry if and only if $T=(T_1,\ldots,
T_n)\in \cD_f(\cH)$  is of class $C_{\cdot 0}$. In this
   case,  we have
$$
K_{f,T}^*[W_\alpha W_\beta^*\otimes I_\cH)K_{f,T}=T_\alpha T_\beta^*,
\quad \alpha,\beta\in \FF_n^+.
$$
 Due to Proposition \ref{tilde-f},
 the $n$-tuple  $W:=(W_1,\ldots, W_n)$ is in the noncommutative domain
  $\cD_f(F^2(H_n))$ and, due to relation \eqref{WWWW}, we have $K_{f,W}^* K_{f,W}=I_\cH$,
   which implies that $W$ is of class $C_{\cdot 0}$.

  We identify $M_m(B(\cH))$, the set of
$m\times m$ matrices with entries in $B(\cH)$, with
$B( \cH^{(m)})$, where $\cH^{(m)}$ is the direct sum of $m$ copies
of $\cH$.
Thus we have a natural $C^*$-norm on
$M_m(B(\cH))$. If $X$ is an operator space, i.e., a closed subspace
of $B(\cH)$, we consider $M_m(X)$ as a subspace of $M_m(B(\cH))$
with the induced norm.
Let $X, Y$ be operator spaces and $u:X\to Y$ be a linear map. Define
the map
$u_m:M_m(X)\to M_m(Y)$ by
$$
u_m ([x_{ij}]_{m\times m}):=[u(x_{ij})]_{m\times m}.
$$
We say that $u$ is completely bounded   if
$$
\|u\|_{cb}:=\sup_{m\ge1}\|u_m\|<\infty.
$$
If $\|u\|_{cb}\leq1$
(resp. $u_m$ is an isometry for any $m\geq1$) then $u$ is completely
contractive (resp. isometric),
 and if $u_m$ is positive for all $m$, then $u$ is called
 completely positive. For more information  on completely bounded maps
  and  the classical von Neumann inequality \cite{vN}, we refer
 to \cite{Pa-book} and \cite{Pi}.

The next result  extends Theorem 3.7 and Theorem 3.8  from
\cite{Po-poisson}.

\begin{theorem}\label{Poisson-C*}
Let $T:=(T_1,\ldots, T_n)$ be an $n$-tuple of operators in the
 noncommutative domain $ \cD_f(\cH)$. Then there is
    a unital completely contractive linear map
$$
\Psi_{f,T}: \overline{\text{\rm  span}} \{ W_\alpha W_\beta^*;\
\alpha,\beta\in \FF_n^+\}\to B(\cH)
$$
such that
 $$
\Psi_{f,T}(W_\alpha W_\beta^*)=T_\alpha T_\beta^*, \quad
\alpha,\beta\in \FF_n^+.
$$
Moreover, the Poisson transform $\Psi_{f,T}$ satisfies the equality
\begin{equation}\label{Po-tran}
\Psi_{f,T}(g)=\lim_{r\to 1} K_{f,rT}^* (g\otimes I_\cH)K_{f,rT},
\end{equation}
where the limit exists in the norm topology of $B(\cH)$.
\end{theorem}

\begin{proof}
 The $n$-tuple $rT:=(rT_1,\ldots, rT_n)\in
\cD_f(\cH)$ is of class $C_{\cdot 0}$ for each $0<r<1$.
Consequently, $K_{f,rT}$ is an isometry and  we can use relation
\eqref{ker-inter} to deduce that
\begin{equation}\label{K*WWK}
K_{f,rT}^*[W_\alpha W_\beta^*\otimes I_\cH)K_{f,rT}=
r^{|\alpha|+|\beta|}T_\alpha T_\beta^*, \quad \alpha,\beta\in
\FF_n^+.
\end{equation}
Hence, we have
$$
 \left\|\sum_{\alpha,\beta\in \Lambda} c_{\alpha, \beta}r^{|\alpha|+|\beta|} T_\alpha
 T_\beta^*\right\|\leq \left\|\sum_{\alpha,\beta\in\Lambda}
  c_{\alpha, \beta} W_\alpha
 W_\beta^*\right\|
 $$
 for any
   finite subset $\Lambda$  of $\FF_n^+$ and
   $c_{\alpha,\beta}\in \CC$, $\alpha, \beta\in
  \Lambda$.
 Taking $r\to 1$,  we deduce that
 $$
 \left\|\sum_{\alpha,\beta\in \Lambda} c_{\alpha, \beta} T_\alpha
 T_\beta^*\right\|\leq \left\|\sum_{\alpha,\beta\in\Lambda}
  c_{\alpha, \beta} W_\alpha
 W_\beta^*\right\|.
 $$
If $g\in \overline{\text{\rm  span}} \{ W_\alpha W_\beta^*;\
\alpha,\beta\in \FF_n^+\}$, we define
\begin{equation}\label{Poi}
\Psi_{f,T}(g):=\lim_{k\to\infty} q_k(T_i,T_i^*),
\end{equation}
where the limit is in the norm topology and
$\{q_k(W_i,W_i^*)\}_{k=1}^\infty$ is a sequence of ``polynomials''
of the form $\sum_{\alpha,\beta\in\Lambda} c_{\alpha, \beta}
W_\alpha W_\beta^*$ such that $\|q_k(W_i,W_i^*)-g\|\to 0$, as
$k\to\infty$. Due to  the inequality above, the operator
$\Psi_{f,T}(g)$ is well-defined and
\begin{equation*}
\|\Psi_{f,T}(g)\|\leq \|g\|.
\end{equation*}
Now, the   first part of the theorem follows from  an operator
matrix version of relation \eqref{K*WWK} and the considerations
above.

 To prove the second part of the theorem, let $g\in
\overline{\text{\rm span}} \{ W_\alpha W_\beta^*;\ \alpha,\beta\in
\FF_n^+\}$ and let $\{q_k(W_i,W_i^*)\}_{k=1}^\infty$ be a sequence
of polynomials in $\overline{\text{\rm  span}} \{ W_\alpha
W_\beta^*;\ \alpha,\beta\in \FF_n^+\}$ such that
$\|q_k(W_i,W_i^*)-g\|\to 0$, as $k\to\infty$. Due to relation
\eqref{Poi}, we have $\Psi_{f,rT}(g):=\lim_{k\to\infty}
q_k(rT_i,rT_i^*)$ in norm. On the other hand, relation \eqref{K*WWK}
implies
$$
 q_k(rT_i,rT_i^*)=K_{f,rT}^*(q_k(W_i,W_i^*)\otimes I_\cH)K_{f,rT}
$$
for any $k\in \NN$ and $r\in (0,1)$. Consequently, we have
\begin{equation}
\label{I-tensor} \Psi_{f,rT}(g)=K_{f,rT}^* (g\otimes I_\cH)K_{f,rT}.
\end{equation}

 Now, let $\epsilon>0$ and $m\in \NN$ such that
$\|q_m(W_i,W_i^*)-g\|<\frac{\epsilon}{3}. $ Due to the first part of
the theorem, we have
$$
\|\Psi_{f,T}(g)-q_m(T_i,T_i^*)\|\leq
\|g-q_m(W_i,W_i^*)\|<\frac{\epsilon}{3}
$$
and
$$
\|\Psi_{f,rT}(g)-q_m(rT_i,rT_i^*)\|\leq
\|g-q_m(W_i,W_i^*)\|<\frac{\epsilon}{3}.
$$
On the other hand, since $q_m(W_i,W_i^*)$  has a finite number of
terms, there exists $\delta\in (0,1)$ such that
$$
\|q_m(rT_i,rT_i^*)-q_m(T_i,T_i^*)\|<\frac{\epsilon}{3}
$$
for any $r\in (\delta,1)$.
Now, using these inequalities and  relation \eqref{I-tensor}, we
deduce that
\begin{equation*}
\begin{split}
\left\|\Psi_{f,T}(g)-K_{f,rT}^* (g\otimes I_\cH)K_{f,rT}\right\| &=
\left\|
\Psi_{f,T}(g)- \Psi_{f,rT}(g)\right\|\\
&\leq \left\| \Psi_{f,T}(g)-q_m(T_i,
T_i^*)\right\|\\
&\quad+\left\| q_m(T_i, T_i^*)- q_m(rT_i,
rT_i^*)\right\|\\
&\quad+ \left\| q_m(rT_i, rT_i^*)- \Psi_{f,
rT}(g)\right\|\\
&<\frac{\epsilon}{3}+ \frac{\epsilon}{3} +
\frac{\epsilon}{3}=\epsilon
\end{split}
\end{equation*}
for any $r\in (\delta, 1)$. This proves \eqref{Po-tran} and
completes the proof.
\end{proof}

\begin{corollary}\label{cc-rep}
Let $(A_1,\dots,A_n)$ be in $B(\cH)^n$. Then $(A_1,\dots,A_n)$ is in
the noncommutative domain $\cD_f(\cH)$ if and only if  the map
$$
\Phi:\cA_n\to B(\cH),\quad\Phi(p(W_1,\ldots, W_n))=p(A_1,\dots,A_n)
$$
is completely contractive.
\end{corollary}

An $n$-tuple  of operators $(A_1,\ldots,A_n)\in B(\cH)^n$ will be
called completely polynomially bounded with respect to the
noncommutative disc algebra $\cA_n(\cD_f)$
  if there is a constant $C>0$ such that, for all $m\in \NN$
and all
$m\times  m$ matrices $[p_{ij}(W_1,\ldots, W_n)]$ with entries in
$\cP_n(W_1,\ldots, W_n)$, the set of all polynomials in $W_1,\ldots,
W_n$ and the identity, we have
$$\|[p_{ij}(A_1,\dots,A_n)]\|_{M_m(B(\cH))}\leq C\|[p_{ij}(W_1,\ldots,
W_n)]\|_{M_m(\cA_n(\cD_f))}.
$$

In \cite{Pa1}, V.~Paulsen proved that an operator $A\in B(\cH)$ is
similar to a contraction if and only if   it is completely
polynomially bounded. In what follows we will extend this result to
our setting.
\begin{theorem} Let $(A_1,\dots,A_n)$ be in $B(\cH)^n$. Then
there is an $n$-tuple of operators $(T_1,\dots,T_n)\in \cD_f(\cH)$
and an invertible operator $X$ such that
$$
A_i=X^{-1}T_i X,\text{\quad for any\ }i=1,\dots,n.
$$
 if and only if
the $n$-tuple $(A_1,\dots,A_n)$ is  completely polynomially bounded
with respect to the noncommutative  domain algebra $\cA_n(\cD_f)$.
\end{theorem}
\begin{proof}
Using Theorem \ref{Poisson-C*} and Paulsen's similarity result
\cite{Pa1}, the result follows.
\end{proof}

\begin{corollary} A representation $\Phi:\cA_n(\cD_f)\to B(\cH)$ is
completely   bounded if and only if  it is given by
$\Phi(W_i)=XT_iX^{-1}, i=1,\dots,n$ where $(T_1,\dots,T_n)\in
\cD_f(\cH)$  and $X$ is an invertible operator.
\end{corollary}

Now, we identify the characters of the noncommutative
  domain
algebra $\cA_n(\cD_f)$. Let $\lambda=(\lambda_1,\dots,\lambda_n)$ be
in $\cD_f(\CC)$, where
$$\cD_f(\CC):=\left\{(\lambda_1,\dots,\lambda_n)\in \CC^n:\
\sum_{|\alpha|\geq 1} a_\alpha |\lambda_\alpha|^2\leq 1\right\},
$$
  and define the
evaluation functional
$$
\Phi_\lambda:\cP_n(W_1,\ldots, W_n)\to
\CC,\quad\Phi_\lambda(p(W_1,\ldots,
W_n))=p(\lambda_1,\dots,\lambda_n).
$$
According to  Theorem \ref{Poisson-C*},   we have
$$
|p(\lambda_1,\dots,\lambda_n)|=\|p(\lambda_1I_{\CC},\dots,\lambda_n
I_{\CC})\|\leq
\|p(W_1,\dots,W_n)\|.
$$
Hence, $\Phi_\lambda$ has a unique extension to the domain algebra
$\cA_n(\cD_f)$. Therefore
$\Phi_\lambda$ is a character of $\cA_n(\cD_f)$.
\begin{theorem}
Let $M_{\cA_n(\cD_f)}$ be the set of all characters of
$\cA_n(\cD_f)$. Then the map
$$\Psi: \cD_f(\CC)\to M_{\cA_n(\cD_f)},\quad
\Psi(\lambda)=\Phi_\lambda $$
 is a homeomorphism of $\cD_f(\CC)$
onto
$M_{\cA_n(\cD_f)}$.
\end{theorem}

\begin{proof} First, notice  that $\Psi$ is one-to-one. Indeed, let
$\lambda:=(\lambda_1,\ldots,\lambda_n)$ and
$\mu:=(\mu_1,\ldots,\mu_n)$ be in $ \cD_f(\CC)$ and assume that
$ \Psi(\lambda)=\Psi(\mu)$.  Then, we have
$$
\lambda_i=\Phi_\lambda(W_i)=\Phi_\mu(W_i)=\mu_i,\text{\ for  \
}i=1,\ldots,n,
$$
which implies $\lambda=\mu$.
Now, assume that $\Phi:\cA_n(\cD_f)\to\CC$ is a character. Setting
$\lambda_i:=\Phi(W_i)$,\ $i=1,\ldots,n$, we  deduce that
$$
\Phi(p(W_1,\ldots,W_n))=p(\lambda_1,\ldots,\lambda_n)
$$
for any polynomial $p(W_1,\ldots,W_n)$ in $\cA_n(\cD_f)$.
Since $\Phi$ is a character it follows that it is completely
contractive. Applying Corollary \ref{cc-rep} in the particular case
when
$A_i:=\lambda_i I_{\CC}$,\ $i=1,\ldots,n$,  it follows that
$(\lambda_1 I_{\Bbb{C}},\ldots,\lambda_n I_{\CC})\in \cD_f(\CC)$.
Moreover, since
$$
\Phi(p(W_1,\ldots,S_n))=p(\lambda_1,\ldots,\lambda_n)=
\Phi_\lambda(p(W_1,\ldots,W_n))
$$
for any polynomial $p(W_1,\ldots,W_n)$ in $\cA_n(\cD_f)$,
   we must have $\Phi=\Phi_\lambda$.

Suppose now that
$\lambda^\alpha:=(\lambda_1^\alpha,\dots,\lambda_n^\alpha);\
\alpha\in J$, is
a net in $\cD_f(\CC)$ such that
$\lim_{\alpha\in
J}\lambda^\alpha=\lambda:=(\lambda_1,\dots,\lambda_n)$.
It is clear  that
$$
\lim_{\alpha\in J}\Phi_{\lambda^\alpha}(p(W_1,\ldots,W_n))=
\lim_{\alpha\in J} p(\lambda_1^\alpha,\ldots,\lambda_n^\alpha)=
p(\lambda_1,\ldots,\lambda_n)=\Phi_\lambda(p(W_1,\ldots,W_n))
$$
for every  polynomial  $p(W_1,\ldots,W_n)$. Since the set of all
polynomials $\cP_n(W_1,\ldots, W_n)$  is  dense
in $\cA_n(\cD_f)$  and $\sup_{\alpha\in J}
\|\Phi_{\lambda^\alpha}\|\leq1$,
 it follows that $\Psi$
is continuous. On the other hand, since  both $ \cD_f(\CC)$ and
$M_{\cA_n}$ are compact Hausdorff
spaces and $\Psi$ is also one-to-one and onto,  the result follows.
The proof is complete.
\end{proof}

Due to Proposition \ref{tilde-f}, it is clear that
$\cR_n(\cD_f)=U^*\cA_n(\cD_{\tilde f})U$. Consequently, all the
results of this section can be written for the algebra
$\cR_n(\cD_f)$.

\pagebreak

\section{The Hardy algebra $F_n^\infty(\cD_f)$}
\label{Hardy}

 In this section, we introduce the Hardy algebra
$F_n^\infty(\cD_f)$ (resp. $R_n^\infty(\cD_f)$) associated with the
noncommutative domain $\cD_f$  and prove some basic properties. In
particular, we show that $F_n^\infty(\cD_f)$ is semisimple and
coincides with its double commutant.

Let $\varphi(W_1,\ldots, W_n)=\sum\limits_{\beta\in \FF_n^+} c_\beta
W_\beta $ be a formal sum with the property  that $\sum_{\beta\in
\FF_n^+} |c_\beta|^2 \frac{1}{b_\beta}<\infty$, where the
coefficients $b_\beta$, $\beta\in \FF_n^+$, are given by relation
 \eqref{b_alpha}.
Using relations \eqref{bbb} and \eqref{WbWb}, one can see that
$\sum\limits_{\beta\in \FF_n^+} c_\beta W_\beta (p)\in F^2(H_n)$ for
any $p\in \cP$, where $\cP$ is the set of all polynomial in
$F^2(H_n)$. Indeed, for each $\gamma\in \FF_n^+$,  we have
$\sum\limits_{\beta\in \FF_n^+} c_\beta W_\beta(e_\gamma)=
\sum\limits_{\beta\in \FF_n^+}
c_\beta\sqrt{\frac{b_\gamma}{b_{\beta\gamma}}} e_{\beta \gamma}$
and, due to relation \eqref{bbb}, we deduce that
$$
\sum_{\beta\in \FF_n^+}|c_\beta|^2 \frac{b_\gamma}{b_{\beta\gamma}}
\leq \sum_{\beta\in \FF_n^+} |c_\beta|^2 \frac{1}{b_\beta}<\infty.
$$
If
$$
\sup_{p\in\cP, \|p\|\leq 1} \left\|\sum\limits_{\beta\in \FF_n^+}
c_\beta W_\beta (p)\right\|<\infty,
$$
then there is a unique bounded operator acting on $F^2(H_n)$, which
we denote by $\varphi(W_1,\ldots, W_n)$, such that
$$
\varphi(W_1,\ldots, W_n)p=\sum\limits_{\beta\in \FF_n^+} c_\beta
W_\beta (p)\quad \text{ for any } \ p\in \cP.
$$
The set of all operators $\varphi(W_1,\ldots, W_n)\in B(F^2(H_n))$
satisfying the above mentioned properties is denoted by
$F_n^\infty(\cD_f)$. When $f=X_1+\cdots +X_n$, it coincides with the
noncommutative analytic Toeplitz algebra $F_n^\infty$, which was
introduced in \cite{Po-von} in connection with a noncommutative
multivariable von Neumann inequality. As in this particular case,
one can prove that $F_n^\infty(\cD_f)$ is a Banach algebra, which we
call Hardy algebra associated with the noncommutative domain
$\cD_f$. In Section \ref{Functional}, we will show that
$F_n^\infty(\cD_f)$ is the $WOT$-closure (resp. $SOT$-closure,
$w^*$-closure) of all polynomial in $W_1,\ldots, W_n$ and the
identity.

 In a similar manner,  using the weighted right creation
operators
 associated with $\cD_f$, one can   define    the corresponding
     the Hardy algebra $R_n^\infty(\cD_f)$.
More precisely, if $g(\Lambda_1,\ldots,
\Lambda_n)=\sum\limits_{\beta\in \FF_n^+} c_{\tilde\beta
}\Lambda_\beta $ is a formal sum with the property  that
$\sum_{\beta\in \FF_n^+} |c_\beta|^2 \frac{1}{b_\beta}<\infty$,
where the
 coefficients $b_\alpha$, $\alpha\in \FF_n^+$, are given by relation
 \eqref{b_alpha}, and such that
$$
\sup_{p\in\cP, \|p\|\leq 1} \left\|\sum\limits_{\beta\in \FF_n^+}
 c_{\tilde\beta} \Lambda_\beta (p)\right\|<\infty,
$$
then there is a unique bounded operator on $F^2(H_n)$, which we
 denote by $g(\Lambda_1,\ldots, \Lambda_n)$, such that
$$
g(\Lambda_1,\ldots, \Lambda_n)p=\sum\limits_{\beta\in \FF_n^+}
c_{\tilde\beta} \Lambda_\beta (p)\quad \text{ for any } \ p\in \cP.
$$
The set of all operators $g(\Lambda_1,\ldots, \Lambda_n)\in
B(F^2(H_n))$
 satisfying the above mentioned properties is denoted by $R_n^\infty(\cD_f)$.

\begin{proposition}\label{tilde-f2}
Let ~$(W_1^{(f)},\ldots, W_n^{(f)}) ($resp.
$(\Lambda_1^{(f)},\ldots, \Lambda_n^{(f)}))$ be the weighted left
(resp. right) creation operators associated with the noncommutative
domain  $\cD_f$. Then the following statements hold:
\begin{enumerate}
\item[(i)] $F_n^\infty(\cD_f)'=U^*(F_n^\infty(\cD_{\tilde
f}))U=R_n^\infty(\cD_f)$, where $'$ stands for the commutant and
$U\in B(F^2(H_n))$ is the unitary operator defined by ~$U
e_\alpha=e_{\tilde\alpha}$, $\alpha\in \FF_n^+$;
\item[(ii)] $F_n^\infty(\cD_f)''=F_n^\infty(\cD_f)$ and
$R_n^\infty(\cD_f)''=R_n^\infty(\cD_f)$.
\end{enumerate}
\end{proposition}
\begin{proof}
Due to Proposition \ref{tilde-f}, we have
$U^*(F_n^\infty(\cD_{\tilde f}))U=R_n^\infty(\cD_f)$. Since
$W_i^{(f)} \Lambda_j^{(f)}=\Lambda_j^{(f)} W_i^{(f)}$ for any
$i,j=1,\ldots,n$, it is clear that $R_n^\infty(\cD_f)\subseteq
F_n^\infty(\cD_f)'$. To prove the reverse inclusion, let $A\in
F_n^\infty(\cD_f)'$. Since $A(1)\in F^2(H_n)$, we have
$A(1)=\sum_{\beta\in \FF_n^+} c_{\tilde
\beta}\frac{1}{\sqrt{b_{\tilde \beta}}} e_{\tilde \beta}$ for some
coefficients $\{c_\beta\}_{\FF_n^+}$ with $\sum_{\beta\in\FF_n^+}
|c_\beta|^2 \frac{1}{b_\beta}<\infty$. On the other hand, since $
AW_i^{(f)}=W_i^{(f)}A$ for $i=1,\ldots,n$, relations \eqref{WbWb}
and \eqref{WbWb-r} imply
\begin{equation*}
\begin{split}
Ae_\alpha &=\sqrt{b_\alpha}AW_\alpha(1)=\sqrt{b_\alpha}W_\alpha
A(1)\\
&=\sum_{\beta\in \FF_n^+} c_{\tilde \beta}
\frac{\sqrt{b_\alpha}}{\sqrt{b_{\alpha \tilde \beta}}} e_{\alpha
\tilde \beta}\\
&=\sum_{\beta\in \FF_n^+} c_{\tilde \beta} \Lambda_\beta(e_\alpha).
\end{split}
\end{equation*}
Therefore, $A(q)=\sum_{\beta\in \FF_n} c_{\tilde \beta}
\Lambda_\beta(q)$ for any polynomial $q$ in in the full Fock space
$F^2(H_n)$. Since $A$ is a bounded operator, $g(\Lambda_1,\ldots,
\Lambda_n):=\sum_{\beta\in \FF_n} c_{\tilde \beta} \Lambda_\beta$ is
in $R_n^\infty(\cD_f)$ and $A=g(\Lambda_1,\ldots, \Lambda_n)$.
Therefore, $R_n^\infty(\cD_f)= F_n^\infty(\cD_f)'$. The item (ii)
follows easily applying part (i). This completes the proof.
\end{proof}

As in the particular case when $f=X_1+\cdots +X_n$ (see
 \cite{Po-analytic}, \cite{Po-central}), one can
 use Proposition \ref{tilde-f2}
  to show that
if $M:F^2(H_n)\otimes \cH\to F^2(H_n)\otimes \cK$ is an operator
such that
 $M(W_i\otimes I_\cH)=(W_i\otimes I_\cK)M$ for any $i=1,\ldots,n$, then
 $M\in R_n^\infty(\cD_f)\bar\otimes B(\cH,\cK)$ and has a unique formal
 Fourier expansion $\sum_{\alpha\in \FF_n^+} \Theta_{(\alpha)}\otimes
 \Lambda_\alpha$ for some  coefficients $\Theta_{(\alpha)}\in B(\cH,\cK)$,
  with the property that  $M$ acts  like its Fourier representation on
   vector-valued polynomials in $F^2(H_n)\otimes \cH$. The operator
   $M$ is called multi-analytic with respect to $W_1,\ldots, W_n$.
   Similar results can be obtained for multi-analytic operators
   with respect to $\Lambda_1,\ldots, \Lambda_n$

In what follows, we employ some ideas from \cite{DP1} to extend some
of the results obtained by  Davidson and Pitts,  to our more general
setting.

\begin{theorem}\label{prop-F-infty} The following statements hold:
\begin{enumerate}
\item[(i)]The Hardy algebra $F_n^\infty(\cD_f)$ is inverse closed.
\item[(ii)] The only normal elements in $F_n^\infty(\cD_f)$ are the
scalars.
\item[(iii)] Every element  $A\in F_n^\infty(\cD_f)$
has
 its spectrum
$\sigma (A)\neq \{0\}$ and it is  injective.
\item[(iv)] The algebra $F_n^\infty(\cD_f)$  contains no non-trivial
idempotents and no non-zero quasinilpotent elements.
\item[(v)] The algebra $F_n^\infty(\cD_f)$ is semisimple.
\item[(vi)] If  $A\in F_n^\infty(\cD_f)$, $n\geq 2$, then
$\sigma (A)=\sigma_e (A)$.
\end{enumerate}
\end{theorem}

\begin{proof}
The first item   follows from Proposition \ref{tilde-f2} part (ii),
using the fact that  $F_n^\infty(\cD_f)''=F_n^\infty(\cD_f)$. To
prove (ii), let $\varphi(W_1,\ldots, W_n):=\sum_{\alpha\in \FF_n^+}
c_\alpha W_\alpha$ be in $F_n^\infty(\cD_f)$. Using relation
\eqref{WbWb}, we have
$$\left<\varphi(W_1,\ldots, W_n)1,1\right>=c_0 \ \text{
and } \ \varphi(W_1,\ldots, W_n)^* 1=\bar{c}_0 1. $$ By normality,
$\varphi(W_1,\ldots, W_n)1=c_0 1$. Using now relation
\eqref{WbWb-r}, we have $\Lambda_\beta (1)=\frac{1}{\sqrt{b_{\tilde
\beta}}} e_{\tilde \beta}$ for $\alpha\in \FF_n^+$. Since
$\varphi(W_1,\ldots, W_n)\Lambda_i=\Lambda_i \varphi(W_1,\ldots,
W_n)$, $i=1,\ldots,n$, we deduce that
\begin{equation*}
\begin{split}
\varphi(W_1,\ldots, W_n) e_\beta&= \sqrt{b_\beta}
\varphi(W_1,\ldots, W_n)\Lambda_{\tilde \beta}(1)\\
&=\sqrt{b_\alpha} \Lambda_{\tilde \beta} \varphi(W_1,\ldots,
W_n)(1)\\
&=c_0\sqrt{b_\beta}\Lambda_{\tilde \beta}(1)= c_0 e_\beta.
\end{split}
\end{equation*}
Therefore, $\varphi(W_1,\ldots, W_n)=c_0 I$. Now we prove (iii).
Assume that  $\varphi(W_1,\ldots, W_n)$ is a nonzero element in
$F_n^\infty(\cD_f)$. As above, we have
$$\varphi(W_1,\ldots, W_n)
e_\beta=\sqrt{b_\alpha} \Lambda_{\tilde \beta} \varphi(W_1,\ldots,
W_n)(1)$$ for any $\beta\in \FF_n^+$. This implies
$\varphi(W_1,\ldots, W_n) e_\beta =\sum_{\alpha\in \FF_n^+} c_\alpha
\frac{1}{\sqrt{b_\alpha}} e_\alpha\neq 0$. Choose an element
$\gamma\in \FF_n^+$ of minimal length $m:=|\gamma|$ such that
$c_\gamma\neq 0$. Using relation \eqref{WbWb}, we can prove by
induction that, for each $k\in \NN$,
$$
\varphi(W_1,\ldots, W_n)^k(1)=c_\gamma^k \frac{1}{\sqrt{\gamma^k}}
e_{\gamma^k}+\sum_{\beta\neq \gamma^k, |\beta|\geq |\gamma|^k}
d_{\beta,k} e_\beta
$$
for some constants $d_{\beta,k}$. Hence,  we deduce that
\begin{equation*}
\begin{split}
\left\|\varphi(W_1,\ldots, W_n)\right\|^{1/k}&\geq
\left|\left<\varphi(W_1,\ldots, W_n)^k1,
e_{\gamma^k}\right>\right|^{1/k}\\
&=|c_\gamma|\left(\frac{1}{\sqrt{b_{\gamma^k}}}\right)^{1/k}.
\end{split}
\end{equation*}
Due to \cite{Po-holomorphic} and the results of Section
\ref{Noncommutative}, since $g(rS_1,\ldots, rS_n)=\sum_{\alpha\in
\FF_n^+} b_\alpha r^{|\alpha|} S_\alpha$ is convergent for some
$r>0$, the radius of convergence $R$ of the power series
$\sum_{\alpha\in \FF_n^+} b_\alpha X_\alpha$ is strictly positive
and $\frac{1}{R}=\limsup\limits_{k\to\infty}\left(\sum_{|\alpha|=
k}|b_\alpha|^2\right)^{1/2k}$. Therefore, we have
$$
\limsup\limits_{k\to\infty}\left( |b_{\gamma^k}|^2\right)^{1/2k}
\leq\limsup\limits_{k\to\infty}\left[\left(\sum_{|\alpha| =
mk}|b_\alpha|^2\right)^{1/(2km)}\right]^m\leq \frac{1}{R^m}.
$$
Hence, there exists $N\in \NN$ such that $|b_{\gamma^k}|^{1/k}\leq
\frac{1}{R^m}$ for any $k\geq N$. Combining this with the above
inequalities, we obtain $\|\varphi(W_1,\ldots, W_n)\|^{1/k}\geq
|c_\gamma|R^{m/2}$ for $k\geq N$. Consequently, the spectral radius
of the operator $\varphi(W_1,\ldots, W_n)$ is greater than  or equal
to $|c_\gamma|R^{m/2}>0$, which shows that
$\sigma(\varphi(W_1,\ldots, W_n))\neq \{0\}$.

To prove injectivity of $\varphi(W_1,\ldots, W_n)$, let
$x:=\sum_{\alpha\in \FF_n^+} d_\alpha e_\alpha\in F^2(H_n)$ be
nonzero and let $\sigma\in \FF_n^+$ be a word of minimal length such
that $d_\sigma\neq 0$. It is clear that
\begin{equation*}
\begin{split}
\left<\varphi(W_1,\ldots, W_n)x,e_{\gamma\sigma}\right>
&=\left<c_\gamma d_\sigma W_\gamma(e_\sigma),
e_{\gamma\sigma}\right>\\
&=c_\gamma d_\sigma\left<
\frac{\sqrt{b_\gamma}}{\sqrt{b_{\gamma\sigma}}} e_{\gamma\sigma},
e_{\gamma\sigma}\right>\\
&=c_\gamma
d_\sigma\frac{\sqrt{b_\gamma}}{\sqrt{b_{\gamma\sigma}}}\neq 0.
\end{split}
\end{equation*}
Therefore $\varphi(W_1,\ldots, W_n)$  is injective. Now, it is clear
that (iii)$\implies$(iv) and (iv)$\implies$(v). It remains to prove
(vi). The implication $\sigma_e (A)\subseteq \sigma (A)$ is obvious.
To prove $\sigma (A)\subset\sigma_e (A)$, it is enough to show that
$0\in \sigma (A)$ implies $0\in\sigma_e (A)$. Since $A$ is
injective, we may assume that $0\in \sigma (A)$ and $\text{\rm
range}\,A$ is  a closed and proper subspace of $F^2(H_n)$. Notice
that due to the fact that the subspace $\sum_
{|\alpha|=k}\nolimits^\oplus \Lambda_\alpha(F^2(H_n))$ is invariant
under $A:=\varphi(W_1,\ldots, W_n)$, we have
\begin{equation}
\label{int1} A\left(\sum_ {|\alpha|=k}\nolimits^\oplus
\Lambda_\alpha(F^2(H_n))\right)\bigcap \text{\rm span} \left\{
e_\beta:\ |\beta|\leq k-1\right\}=\emptyset.
\end{equation}
Now, we prove  that
\begin{equation}
\label{int2} A\left(\sum_ {|\alpha|=k}\nolimits^\oplus
\Lambda_\alpha(F^2(H_n))\right)\bigcap \text{\rm span} \left\{
e_\beta:\ |\beta|= k\right\}=\{0\}.
\end{equation}
Let $x:=\sum\limits_{|\alpha|=k}a_\alpha e_\alpha
=A\left(\sum\limits_{|\alpha|=k} \Lambda_\alpha y_{(\alpha)}\right)$
for some vectors $y_{(\alpha)}\in F^2(H_n)$. If $\beta\in \FF_n^+$
and $|\beta|=k$, then, due to relation \eqref{WbWb-r}, we have
\begin{equation*}
\begin{split}
a_\beta \frac{1}{\sqrt{b_\beta}}&=\Lambda_{\tilde \beta}^*
x=\Lambda_{\tilde\beta} A\left(\sum\limits_{|\alpha|=k}
\Lambda_\alpha y_{(\alpha)}\right)\\
&=\sum\limits_{|\alpha|=k}\Lambda_{\tilde \beta}^* \Lambda_\alpha
Ay_{(\alpha)}=\Lambda_{\tilde \beta}^* \Lambda_{\tilde \beta}
Ay_{(\alpha)}\\
&=D_{\tilde \beta}Ay_{(\alpha)}
\end{split}
\end{equation*}
where $D_{\tilde \beta}:=\Lambda_{\tilde \beta}^* \Lambda_{\tilde
\beta}$ is the diagonal operator on $F^2(H_n)$ satisfying the
relation $D_{\tilde
\beta}(e_\alpha)=\frac{b_\alpha}{b_{\alpha\beta}}e_\alpha$, \
$\alpha\in \FF_n^+$. Notice that the constant $1$ is not in the
range of $A$ because $1$ is cyclic for $\{\Lambda_\alpha:\ \alpha\in
\FF_n^+\}$ and $A\Lambda_i=\Lambda_i A$, $i=1,\ldots,n$, which would
contradict that the range of $A$ is proper. Consequently, we must
have $Ay_{(\tilde\beta)}=\sum_{\alpha\in \FF_n^+} d_\alpha
e_\alpha$, where $d_\alpha=0$ for any $\alpha\in \FF_n^+$, or at
least one $\gamma\in \FF_n^+$, $|\gamma|\geq 1$ is such that
$d_\gamma\neq 0$. In the first case, the above calculations show
that $a_\beta=0$. In the second case, we have
$$a_\beta\frac{1}{\sqrt{b_\beta}}=D_{\tilde \beta}\left(\sum_{\alpha\in \FF_n^+} d_\alpha
e_\alpha\right) =\sum_{\alpha\in \FF_n^+} d_\alpha
\frac{b_\alpha}{b_{\alpha \beta}} e_\beta.
$$
Since $ d_\alpha \frac{b_\alpha}{b_{\alpha \beta}}\neq 0$, we got a
contradiction. Therefore, we have $a_\beta=0$ for any $\beta\in
\FF_n^+$ with $|\beta|=k$. This implies $x=0$, which proves relation
\eqref{int2}. Using relations \eqref{int1} and \eqref{int2}, one can
see that
$$
\dim \left[\text{\rm range}\, A\bigcap \text{\rm span}\,\{e_\alpha:\
|\alpha|\leq k\}\right]\leq 1+n+n^2+\cdots +n^{k-1}.
$$
Hence the range of $A$  has codimension at least $n^k$ for any $k\in
\NN$. This  shows that  the range of $A$ has infinite codimension
and, consequently, $0\in \sigma_e(A)$. The  proof is complete.
\end{proof}

\bigskip

\section{Functional calculus for
 $n$-tuples of operators  in $\cD_f$} \label{Functional}

In this section,  we obtain an $F_n^\infty(\cD_f)$- functional
calculus for completely noncoisometric  $n$-tuples of operators in
the noncommutative domain $\cD_f$.

 Let $T:=(T_1,\ldots, T_n)$ be an $n$-tuple of
operators in the noncommutative domain $ \cD_f(\cH)$, i.e.,
$\sum\limits_{|\alpha|\geq 1} a_\alpha T_\alpha T_\alpha^*\leq
I_\cH$, and define the positive linear mapping
$$\Phi_{f,T}:B(\cH)\to
B(\cH)\quad \text{ by }\quad \Phi_{f,T}(X)=\sum\limits_{|\alpha|\geq
1} a_\alpha T_\alpha XT_\alpha^*,
$$
where the convergence is in the weak operator topology. The
$n$-tuple $T \in \cD_f(\cH)$ is called {\it completely
non-coisometric }(c.n.c) with respect to $\cD_f(\cH)$ if there is no
vector $h\in \cH$, $h\neq 0$, such that
$$
\left< \Phi_{f,T}^m(I)h,h\right>=\|h\|^2 \quad \text{ for any } \
m=1,2,\ldots.
$$
According to Section \ref{domain algebra}, we have
$$
\|K_{f,T} h\|^2=\|h\|^2-\|Q_{f,T}^{1/2}h\|^2,\quad h\in \cH,
$$
where $Q_{f,T}:=\text{\rm SOT-}\lim_{m\to\infty}\Phi_{f,T}^m(I)$ and
$K_{f,T}$ is the Poisson kernel associated with $T$ and $\cD_f$. As
a consequence, one can easily prove the following.

\begin{proposition}\label{cnc}
An $n$-tuple of operators $T:=(T_1,\ldots, T_n)\in \cD_f(\cH)$ is
c.n.c. if and only if the noncommutative Poisson kernel ~$K_{f,T}$
is one-to-one.
\end{proposition}

We recall from Section \ref{domain algebra} the following property
of the Poisson kernel:
\begin{equation}
\label{KTWK} K_{f,T} T_i^*=(W_i^*\otimes I_\cH) K_{f,T}, \quad
i=1,\ldots, n,
\end{equation}
wherer $(W_1,\ldots, W_n)$ are the weighted left creation operators
associated  with $\cD_f$.
 Now, we can prove the main result of this
section.

\begin{theorem}
\label{funct-calc} Let $T:=(T_1,\ldots, T_n)\in B(\cH)^n$ be a
completely non-coisometric $($c.n.c.$)$
 $n$-tuple of operators in the noncommutative domain $\cD_f(\cH)$. Then
\begin{equation}
\label{sot}
g(T_1,\ldots, T_n):=\text{\rm SOT-}\lim_{r\to 1} g_r(T_1,\ldots, T_n)
\end{equation}
exists in the strong operator topology of $B(\cH)$ for any
$g=g(W_1,\ldots, W_n)\in F_n^\infty(\cD_f)$, and the mapping
$$\Phi:F_n^\infty(\cD_f)\to B(\cH)\quad \text{defined by} \quad
 \Phi(g):=g(T_1,\ldots, T_n)$$
has the following properties:
\begin{enumerate}
\item[(i)] $\Phi$ is    WOT-continuous (resp.
SOT-continuous)  on bounded sets;
 \item[(ii)]
$\Phi$ is a unital completely contractive homomorphism;
\item[(iii)]
$\Phi$ coincides with the Poisson transform
$\Psi_{f,T}:F_n^\infty(\cD_f)\to B(\cH)$ defined by
$$
\Psi_{f,T}(g):=\text{\rm SOT-}\lim_{r\to 1} K_{f,rT}^*(g\otimes
I_\cH)K_{f,rT}.
$$
 \end{enumerate}
In particular, if $T$ is of class $C_{\cdot 0}$, then
$$
\Phi(g):= g(T_1,\ldots, T_n)= K_{f,T}^*(g\otimes I_\cH)K_{f,T}
$$
for any $g\in F_n^\infty(\cD_f)$, and $\Phi$ is $w^*$-continuous.
Moreover,
$$
g(T_1,\ldots, T_n)^*=\text{\rm SOT-}\lim_{r\to 1} g_r(T_1,\ldots,
T_n)^*.
$$
 \end{theorem}

\begin{proof} First, we shall prove that the limit in \eqref{sot}
exists for the $n$-tuple of weighted left creation operators
$W:=(W_1,\ldots, W_n)\in \cD_f(F^2(H_n))$ and, moreover, that
\begin{equation}\label{c0sot}
g(W_1,\ldots,W_n)=\text{\rm SOT-}\lim_{t\to 1} g_t(W_1,\ldots, W_n)
\end{equation}
for any   ~$g(W_1,\ldots, W_n):=\sum\limits_{\beta\in \FF_n^+}
 c_\beta W_\beta\in F_n^\infty(\cD_f)$.
 According to Section \ref{Noncommutative}, the operators
$\{W_\beta\}_{|\beta|=k}$ have orthogonal ranges and
$\|W_\beta\|=\frac{1}{\sqrt{b_\beta}}$, $\beta\in \FF_n^+$.
Consequently, $\left\|\sum\limits_{|\beta|=k} b_\beta W_\beta
W_\beta^*\right\|\leq 1$ for any $k=0,1,\ldots$. Hence, and using
the fact that $\sum\limits_{\beta\in \FF_n^+}
|c_\beta|^2\frac{1}{b_\beta}<\infty$,
 we deduce that, for $0<t<1$,

\begin{equation*}
 \begin{split}
 \sum_{k=0}^\infty t^k \left\|\sum_{|\beta|=k} c_\beta W_\beta\right\|
 &\leq
 \sum_{k=0}^\infty t^k \left(\sum_{|\beta|=k} |c_\beta|^2\frac{1}
 {b_\beta}\right)^{1/2}\left\|\sum_{|\beta|=k} b_\beta
  W_\beta W_\beta^*\right\|^{1/2}  \\
 &\leq
  \frac{1}{1-t} \left(\sum_{\beta\in \FF_n^+} |c_\beta|^2\frac{1}{b_\beta}
   \right)^{1/2} < \infty,
 \end{split}
 \end{equation*}
which proves that
 $g_t(W_1,\ldots, W_n)$ is in the  noncommutative domain  algebra
$\cA_n(\cD_f)$ and
 \begin{equation}
\label{limm}
 \lim_{m\to \infty}\sum_{|\alpha|\leq m}t^{|\alpha|} c_\alpha W_\alpha
  =g_t(W_1,\ldots, W_n),
 \end{equation}
 where the convergence is in the operator norm.

 Fix now  $\gamma, \sigma, \epsilon\in \FF_n^+$ and denote
$p(W_1,\ldots,W_n):=\sum\limits_{\beta\in \FF_n^+, |\beta|\leq
|\gamma|} c_\beta W_\beta$. Since $W_\beta^* e_\gamma =0$ for any
$\beta\in \FF_n^+$ with $|\beta|>|\gamma|$,  we have
$$
g_r(W_1,\ldots, W_n)^* e_\alpha =p_r(W_1,\ldots, W_n)^* e_\alpha$$
for any $\alpha\in \FF_n^+$ with $|\alpha|\leq |\gamma|$ and any
$r\in [0,1]$. On the other hand, using relation \eqref{KTWK} when
$T=rW$, we obtain
$$
K_{f,rW}p_r(W_1,\ldots, W_n)^*=[p(W_1,\ldots, W_n)^*\otimes
I_{F^2(H_n)}]K_{f,rW}
$$
for any $r\in[0,1)$. Using all these facts,  careful calculations
reveal that
\begin{equation*}
\begin{split}
\left<K_{f,rW}g_r(W_1,\right.&\left.\ldots, W_n)^*e_\gamma,
e_\sigma\otimes e_\epsilon\right> \\
&=\left<K_{f,rW}p_r(W_1,\ldots, W_n)^*e_\gamma,
e_\sigma\otimes e_\epsilon\right>\\
&=\left<[(p(W_1,\ldots, W_n)^*\otimes
I_{F^2(H_n)})]K_{f,rW}e_\gamma,
e_\sigma\otimes e_\epsilon\right>\\
&=\sum_{\beta\in \FF_n^+} r^{|\beta|}
\sqrt{b_\beta}\left<p(W_1,\ldots, W_n)^* e_\beta,e_\sigma\right>
\left< W_\beta^*e_\gamma,  \Delta_{f,rW} e_\epsilon\right>\\
&=\sum_{\beta\in \FF_n^+} r^{|\beta|}
\sqrt{b_\beta}\left<g(W_1,\ldots, W_n)^* e_\beta,e_\sigma\right>
\left< W_\beta^*e_\gamma,  \Delta_{f,rW} e_\epsilon\right>\\
&= \left<[g(W_1,\ldots, W_n)^*\otimes I_{F^2(H_n)}] K_{f,rW}
e_\gamma, e_\sigma\otimes e_\epsilon\right>
\end{split}
\end{equation*}
for any $r\in [0,1)$ and $\gamma, \sigma,\epsilon\in \FF_n^+$.
Hence, since $g_r(W_1,\ldots, W_n)$ and  $g(W_1,\ldots, W_n)$  are
bounded operators, we deduce that
$$
K_{f,rW}g_r(W_1,\ldots, W_n)^*=[g(W_1,\ldots, W_n)^*\otimes
I_{F^2(H_n)}]K_{f,rW}.
$$
Consequently, since  the $n$-tuple $rW:=(rW_1,\ldots, rW_n)\in
\cD_f(F^2(H_n))$ is of class $C_{\cdot 0}$, the Poisson kernel
$K_{f,rW}$ is an isometry and therefore,
\begin{equation}
\label{gr} \|g_r(W_1,\ldots, W_n)\|\leq \|g(W_1,\ldots, W_n)\|\quad
\text{ for any } r\in [0,1).
\end{equation}
 Hence, and due to the fact that
$g(W_1,\ldots,W_n)p= \lim\limits_{r\to 1}g_r(W_1,\ldots, W_n)p$ for
any polynomial $p$ in $F^2(H_n)$, we deduce  relation \eqref{c0sot}.

Let $\{r_k\}_{k=1}^\infty$ be an increasing sequence  of positive
numbers which converges to $1$. Since $g_{r_k}(W_1,\ldots,
W_n)=\sum\limits_{m=1} \sum\limits_{|\alpha|=m} r_k^{|\alpha|}
c_\alpha W_\alpha$, where the convergence is in the uniform topology
(see \eqref{limm}), we can find a polynomial $p_k(W_1,\ldots, W_n)$
such that
$$
\|g_{r_k}(W_1,\ldots, W_n)-p_k(W_1,\ldots, W_n)\|\leq
\frac{1}{k}\quad \text{ for any } \ k=1,2,\ldots.
$$
 Using relation \eqref{c0sot}, it is easy to see that
$$g(W_1,\ldots,W_n)=\text{\rm SOT-}\lim\limits_{k\to
\infty}p_k(W_1,\ldots, W_n).
$$
Since  the WOT and $w^*$ topology coincide on bounded sets, the
principle of uniform boundedness implies
$$g(W_1,\ldots,W_n)=\text{\rm WOT-}\lim\limits_{k\to
\infty}p_k(W_1,\ldots, W_n)=w^*\text{-}\lim\limits_{k\to
\infty}p_k(W_1,\ldots, W_n).
$$
 On the other hand,  since $F_n^\infty (\cD_f)=R_n^\infty (\cD_f)'$,
 it is clear that $F_n^\infty(\cD_f)$ is WOT-closed.
 Now, one can easily
   see that $F_n^\infty(\cD_f)$ is the  SOT-(resp. WOT-, $w^*$-) closure of
     all polynomials in
$W_1,\ldots, W_n$ and the identity.

The next step is to  prove that, if $T:=(T_1,\ldots, T_n)\in
\cD_f(\cH)$, then
 \begin{equation}\label{Kfrt}
 K_{f,rT}^* (g(W_1,\ldots, W_n)\otimes I_\cH)=g_r(T_1,\ldots, T_n)
  K_{f,rT}^*
 \end{equation}
 for any $g(W_1,\ldots, W_n) \in F_n^\infty(\cD_f)$ and $0<r<1$.
To this end, notice that the  relation \eqref{KTWK} implies
 \begin{equation}\label{eq-ker2}
 K_{f,rT}^*[p(W_1,\ldots, W_n)\otimes I_\cH]=p(rT_1,\ldots, rT_n)
  K_{f,rT}^*
 \end{equation}
for any polynomial $p(W_1,\ldots, W_n)$ and $r\in [0,1)$.
  Since $rT:=(rT_1,\ldots, rT_n)\in \cD_f(\cH)$, relation \eqref{limm} and
 Theorem \ref{Poisson-C*}   imply
  $$
 \lim_{m\to \infty} \sum_{|\alpha|\leq m}
 t^{|\alpha|} r^{|\alpha|} c_\alpha T_\alpha
   =g_t(rT_1,\ldots, rT_n) \quad \text{ for any } r,t\in [0,1),
 $$
 where the convergence is in the operator norm topology.
  Using  relation \eqref{eq-ker2}, when
 $p(W_1,\ldots, W_n):=
 \sum\limits_{k=0}^m \sum\limits_{|\alpha|=k}t^{|\alpha|} c_\alpha
 W_\alpha$,
 and taking the limit as
 $m\to \infty$, we get
 \begin{equation}
 \label{eq-ker3}
 K_{f,rT}^* [g_t(W_1,\ldots, W_n)\otimes I_\cH]=g_t(rT_1,\ldots, rT_n)
 K_{f,rT}^*.
 \end{equation}
 Now, let us prove  that
\begin{equation}\label{lim-t}
 \lim_{t\to 1} g_t(rT_1,\ldots, rT_n)=g_r(T_1,\ldots, T_n),
 \end{equation}
 where the convergence is in the operator norm topology.
  Notice that, if  $\epsilon>0$,  there is   $m_0\in \NN$
 such that $\left(\sum\limits_{k=m_0}^\infty r^k\right) \|g(W_1,\ldots,
  W_n)(1)\|<\frac {\epsilon} {2}$.
  Since $(T_1,\ldots, T_n)\in \cD_f(\cH)$,  Theorem \ref{Poisson-C*} implies
  $$\left\| \sum_{|\beta|=k} b_\alpha T_\beta T_\beta^*\right\|\leq
  \left\| \sum_{|\beta|=k} b_\alpha W_\beta W_\beta^*\right\|\leq 1.
  $$
  Now, we can deduce that
\begin{equation*}
 \begin{split}
 \sum_{k=m_0}^\infty r^k \left\|\sum_{|\beta|=k} c_\beta T_\beta\right\|
 &\leq
 \sum_{k=m_0}^\infty r^k \left(\sum_{|\beta|=k} |c_\alpha|^2\frac{1}
 {b_\beta}\right)^{1/2}\left\|\sum_{|\beta|=k} b_\beta T_\beta
 T_\beta^*\right\|^{1/2}  \\
 &\leq
 \left(\sum_{k=m_0}^\infty r^k\right) \|g(W_1,\ldots, W_n)(1)\|
  <\frac {\epsilon} {2}.
 \end{split}
 \end{equation*}
 Consequently, we have
 \begin{equation*}
 \begin{split}
 \left\|\sum_{k=0}^\infty \sum_{|\alpha|=k}
  t^{|\alpha|} r^{|\alpha|} c_\alpha T_\alpha \right.&-\left.
  \sum_{k=0}^\infty \sum_{|\alpha|=k}
    r^{|\alpha|} c_\alpha T_\alpha\right\|\\
    &\leq
    \left\| \sum_{k=1}^{m_0-1} \sum_{|\beta|
    =k}r^k(t^k-1)c_\beta T_\beta\right\|+ \epsilon\\
    &\leq \sum_{k=1}^{m_0-1} r^k(t^k-1)\|g(W_1,\ldots, W_n)(1)\|
    +\epsilon.
 \end{split}
 \end{equation*}
 Now, it clear that there exists $0<\delta<1$ such that
 $$\sum\limits_{k=1}^{m_0-1} r^k(t^k-1)\|g(W_1,\ldots, W_n)(1)\|
 <\epsilon\quad \text{ for any } \ t\in (\delta, 1),
 $$
  which proves \eqref{lim-t}.
Since the map $Y\mapsto Y\otimes I_\cH$ is SOT-continuous on bounded
sets, relations \eqref{c0sot} and \eqref{gr} imply that
$$
\text{\rm SOT-}\lim_{r\to 1}g_r(W_1,\ldots, W_n )\otimes
I_\cH=g(W_1,\ldots, W_n )\otimes I_\cH.
$$
Using relation
   \eqref{lim-t}  and  taking  the limit, as $t\to 1$, in
 \eqref{eq-ker3}, we   deduce   \eqref{Kfrt}.

Since  the operator $g_r(W_1,\ldots, W_n)$ is in the  domain algebra
$\cA_n(\cD_f)$ and the sequence $\sum_{k=0}^m \sum_{|\alpha|=k}
c_\alpha r^{|\alpha|} W_\alpha $ is convergent in norm to
$g_r(W_1,\ldots, W_n)$, as $m\to\infty$, Theorem \ref{Poisson-C*}
and relation \eqref{gr} imply that $\sum\limits_{k=0}^m
\sum\limits_{|\alpha|=k} c_\alpha r^{|\alpha|} T_\alpha $ is
convergent in norm to $g_r(T_1,\ldots, T_n)$, as $m\to\infty$, and
\begin{equation}
\label{vnvn} \|g_r(T_1,\ldots, T_n)\|\leq \|g_r(W_1,\ldots, W_n)\|
\leq \|g(W_1,\ldots, W_n)\|
\end{equation}
for any $r\in [0,1)$. Since $T_i K_{f,T}^*=K_{f,T}^*(W_i\otimes
I_\cH)$,\ $i=1,\ldots, n$, we deduce that

\begin{equation}
\label{frf} g_r(T_1,\ldots, T_n) K_{f,T}^*
=K_{f,T}^*(g_r(W_1,\ldots, W_n )\otimes I_\cH) \qquad \text{ for any
} \ r\in [0,1).
\end{equation}

Taking $r\to 1$ in relation \eqref{frf}, we deduce that the map $A:
\text{\rm range} \,K_{f,T}^*\to \cH$ given by $Ay:=\lim\limits_{r\to
1}g_r(T_1,\ldots,T_n)y$ is well-defined, linear, and
$$
\|AK_{f,T}^*\varphi\|\leq \limsup_{r\to 1} \|g_r(W_1,\ldots,
W_n)\|\|K_{f, T}^*\varphi\|\leq \|g(W_1,\ldots, W_n)\|\|K_{f,
T}^*\varphi\|
$$
   for any $\varphi\in
 F^2(H_n)\otimes \cH$.

Now, assume that $T=(T_1,\ldots, T_n)\in \cD_f(\cH)$ is c.n.c..
According to Proposition \ref{cnc}, the Poisson kernel $K_{f,T}$ is
one-to-one and therefore the range of $K_{f,T}^*$ is dense in $\cH$.
Consequently, the map $A$ has a unique extension to a bounded linear
operator on $\cH$, denoted also by $A$, with  $\|A\|\leq
\|g(W_1,\ldots, W_n)\|$.
Let us show that
$$
\lim_{r\to 1} g_r(T_1,\ldots, T_n)h =Ah\quad \text{ for any }\ h\in
\cH.
$$
To this end,  let $\{y_k\}_{k=1}^\infty$  be a sequence of vectors
in the range of $K_{f,T}^*$, which converges to $y$.
Due to relation \eqref{vnvn}, we have
\begin{equation*}
\begin{split}
\|Ah-g_r(T_1,\ldots, T_n)h\|&\leq \|Ah-A
y_k\|+\|Ay_k-g_r(T_1,\ldots, T_n)y_k\|\\
&\qquad \qquad  +\|g_r(T_1,\ldots,
T_n)y_k-g_r(T_1,\ldots, T_n)h\|\\
&\leq 2\|g(W_1,\ldots, W_n)\| \|h-y_k\|+\|Ay_k-g_r(T_1,\ldots,
T_n)y_k\|.
\end{split}
\end{equation*}
Now, since $\lim\limits_{r\to 1} g_r(T_1,\ldots, T_n)y_k =Ay_k$, our
assertion follows. Denoting  the operator $A$ by $g(T_1,\ldots,
T_n)$, we complete the proof of  relation \eqref{sot}.

Due to relation \eqref{Kfrt}, we have
\begin{equation}
\label{anot} g_r(T_1,\ldots, T_n)=K_{f,rT}^* [g(W_1,\ldots,
W_n)\otimes I_\cH]K_{f,rT},
\end{equation}
which together with \eqref{sot} imply part (iii) of the theorem.
Now, using \eqref{sot}, the fact that $\|g_r(W_1,\ldots, W_n)\|\leq
\|g(W_1,\ldots, W_n)\|$, and taking $r\to 1$ in relation
\eqref{frf}, we obtain
\begin{equation}
\label{gKKg} g(T_1,\ldots, T_n)K_{f,T}^*=K_{f,T}^*[g(W_1,\ldots,
W_n)\otimes I_\cH].
\end{equation}

To prove part (i), let $\{g_i(W_1,\ldots, W_n)\}$   be a bounded net
in $F_n^\infty(\cD_f)$  which is WOT (resp. SOT) convergent to an
element $g(W_1,\ldots, W_n)\in F_n^\infty(\cD_f)$. Then, due to
standard facts in functional analysis,  the net $g_i(W_1,\ldots,
W_n)\otimes I_\cH$ is WOT (resp. SOT) convergent to $g(W_1,\ldots,
W_n)\otimes I_\cH$. On the other hand, as shown above,
$\|g_i(T_1,\ldots, T_n)\|\leq \|g_i(W_1,\ldots, W_n)\|$ and,
consequently, the net $\{g_i(T_1,\ldots, T_n)\}$ is bounded. Now,
one can use  relation \eqref{gKKg}, the fact that the range of
$K_{f,T}^*$ is dense in $\cH$, and a standard approximation argument
to deduce that $g_i(T_1,\ldots, T_n)$ is WOT (resp. SOT) convergent
to $g(T_1,\ldots, T_n)$.

 To prove part (ii), notice that since $K_{f,rT}$ is an isometry,
  relation \eqref{anot} implies
  $$
   \left\|\left[(f_{ij})_r(T_1,\ldots, T_n)\right]_{k\times k}\right\|
   \leq
\left\|\left[f_{ij}(W_1,\ldots, W_n)\right]_{k\times k}\right\|
$$
for any operator-valued matrix $\left[f_{ij}(W_1,\ldots,
W_n)\right]_{k\times k}$ in $M_k(F_n^\infty(\cD_f))$ and $r\in
[0,1)$. Due to the first part of this theorem,
 $\text{\rm SOT-}\lim\limits_{r\to 1}(f_{ij})_r(T_1,\ldots, T_n)
 =f_{ij}(T_1,\ldots, T_n)$
 for each $i,j=1,\ldots, k$. Now, one can easily deduce that
$$
   \left\|\left[f_{ij}(T_1,\ldots, T_n)\right]_{k\times k}\right\|
   \leq
\left\|\left[f_{ij}(W_1,\ldots, W_n)\right]_{k\times k}\right\|.
$$
 Since $\Phi$ is  a homomorphism  on polynomials in $F_n^\infty(\cD_f)$,
      the WOT-continuity
 of $\Phi$ and
the WOT-density of the polynomials  in $F_n^\infty(\cD_f)$ shows
that $\Phi$  is also a homomorphism on $F_n^\infty(\cD_f)$.

Now, assume that  $T$ is of class $C_{\cdot 0}$.  Then the Poisson
kernel  $K_{f,T}$ is an isometry and the map
$$
\Psi_{f,T}(g(W_1,\ldots, W_n)):=K_{f,T}^* (g(W_1,\ldots, W_n)\otimes
I_\cH) K_{f,T},
$$
where $g(W_1,\ldots, W_n)\in F_n^\infty(\cD_f)$,
 coincides with  $\Phi$ on  polynomials. Since the polynomials are
 sequentially
 WOT-dense in $F_n^\infty(\cD_f)$, and $\Psi_{f,T}$ and $\Phi$ are
 WOT- continuous on bounded sets, one can use the principle of
 uniform boundedness to deduce that $\Psi_{f,T}=\Phi$. The $w^*$
 continuity of $\Phi$ is now obvious.

Finally, to prove the last part of the theorem, notice that if\
  $ g(W_1,\ldots, W_n)=
\sum\limits_{\alpha\in \FF_n^+}c_\alpha W_\alpha \in
F_n^\infty(\cD_f)$, then using relation \eqref{c0sot}, we deduce
that
    for every $\beta\in \FF_n^+$  and
$y\in F^2(H_n) $, we have

\begin{equation*}
\begin{split}
\left< g(W_1,\ldots, W_n)^*e_\beta ,y\right> &=
\left<e_\beta ,g(W_1,\ldots, W_n)y\right>\\
&=\left<e_\beta , \sum_{\alpha\in \FF_n^+, |\alpha|\leq
|\beta|} c_\alpha W_\alpha y\right>\\
&= \left< \sum_{\alpha\in \FF_n^+, |\alpha|\leq
|\beta|}\overline{c}_\alpha W_\alpha^*  e_\beta , y\right>.
\end{split}
\end{equation*}
 Therefore,
$$
g(W_1,\ldots, W_n)^*e_\beta  =\sum_{\alpha\in \FF_n^+, |\alpha|\leq
|\beta|} \overline{c}_\alpha W_\alpha^*  e_\beta
$$
and
$$
g(rW_1,\ldots, rW_n)^*e_\beta  = \left(\sum_{\alpha\in \FF_n^+,
|\alpha|\leq |\beta|}r^{|\alpha|}\overline{c}_\alpha W_\alpha^*
\right)e_\beta .
$$
Using the last two equalities, we obtain
\begin{equation*}
\label{lim2} \lim_{r\to 1}g(rW_1,\ldots, rW_n)^*e_\beta
=g(W_1,\ldots, W_n)^*e_\beta
\end{equation*}
for any $\beta\in \FF_n^+$. On the other hand,  using relation
\eqref{gr}
  and   the fact that the closed span of all vectors
$e_\alpha $,  $\beta\in \FF_n^+$,  coincides with $F^2(H_n)$, we
deduce (using standard arguments) that
\begin{equation}\label{so-e}
\text{\rm SOT-}\lim_{r\to 1}g(rW_1,\ldots, rW_n)^*=g(W_1,\ldots,
W_n)^*.
\end{equation}
If $T$ is of class $C_{\cdot 0}$, then
$$
g(T_1,\ldots, T_n):=K_{f,T}^* (g(W_1,\ldots, W_n)\otimes I_\cH)
K_{f,T}
$$
and
$$
g(rT_1,\ldots, rT_n):=K_{f,T} (g(rW_1,\ldots, rW_n)\otimes I_\cH)
K_{f,T}.
$$
Using again that $\|g(rT_1,\ldots, rT_n)\|\leq \|g(T_1,\ldots,
T_n)\|$ and relation \eqref{so-e}, we complete the   proof.
\end{proof}

\begin{corollary}\label{w*-density}
The Hardy algebra $F_n^\infty(\cD_f)$ is the $w^*$- (resp. WOT-,
SOT-) closure of all polynomials   in $W_1,\ldots, W_n$ and the
identity.
\end{corollary}

\bigskip

\section{The noncommutative variety  $\cV_{f,J}$ and a functional calculus}
\label{Functional II}

 In this section, we obtain a functional calculus for $n$-tuples of
 operators in noncommutative varieties $\cV_{f,J}$ generated by
 $w^*$-closed two sided ideals of the Hardy algebra
 $F_n^\infty(\cD_f)$.

Let $J$ be a $w^*$-closed  two-sided ideal of $F^\infty_n(\cD_f)$.
 Denote by $F_n^\infty(\cV_{f,J})$   the $w^*$-closed algebra
generated
 by the operators
 $B_i:=P_{\cN_J} W_i |\cN_J$, \ $i=1,\ldots, n$, and the identity, where
 $$
 \cN_J:= F^2(H_n)\ominus \cM_J\quad \text{ and }\quad
 \cM_J:=\overline{ J F^2(H_n)}.
 $$
 Notice that
 $\cM_J:=\overline{\{\varphi (1):\ \varphi \in J \}}$
  and
 $
 \cN_J=\bigcap\limits_{\varphi\in J}\ker \varphi^*.
 $

The following result is an extension   of  Theorem 4.2 from
\cite{ArPo1} to our more general setting. Since the proof follows
the same lines, we shall omit  it.

\begin{theorem}
\label{F/J} Let $J$ be a $w^*$-closed  two-sided ideal of the Hardy
algebra $F^\infty_n(\cD_f)$. Then  the map
$$\Gamma:F_n^\infty(\cD_f)/J\to
B(\cN_J) \quad \text{ defined by } \quad \Gamma(\varphi+J)=P_{\cN_J}
\varphi|_{\cN_J}
$$
is a $w^*$-continuous, completely isometric representation.
\end{theorem}

Since the set of all polynomials in $W_1,\ldots, W_n$  and the
identity is $w^*$-dense in $F_n^\infty(\cD_f)$, Theorem \ref{F/J}
implies that $P_{\cN_J}F_n^\infty(\cD_f)|_{\cN_J}$ is  a
$w^*$-closed subalgebra  of $B(\cN_J)$ and, moreover,
$F_n^\infty(\cV_{f,J})=P_{\cN_J}F_n^\infty(\cD_f)|_{\cN_J}$.
 Now, one can easily deduce that
  \begin{equation}\label{w-formulas}
 F_n^\infty(\cV_{f,J})  =\{g(B_1,\ldots, B_n):
 \ g(W_1,\ldots, W_n)\in F_n^\infty(\cD_f)\},
 \end{equation}
where $g(B_1,\ldots, B_n)$ is defined by the
 $F_n^\infty(\cD_f)$-functional calculus of Section \ref{Functional}.

Let $\cH$ be a separable complex Hilbert space,
 $J$ be a $w^*$-closed  two-sided ideal of $F^\infty_n(\cD_f)$, and let
$T:=(T_1,\ldots, T_n)\in \cD_f(\cH)$   be a c.n.c. $n$-tuple of
operators. We say that $T$ is   in the
 noncommutative variety $\cV_{f,J}(\cH)$ if
  $$
\varphi(T_1,\ldots, T_n)=0\quad \text{for any } \ \varphi\in J.
$$
The {\it constrained Poisson kernel} associated with $T\in
\cV_{f,J}(\cH)$ is the operator
$$K_{f,T,J}: \cH\to \cN_J\otimes
\overline{\Delta_{f,T}(\cH)}\quad \text{ defined by }\quad
K_{f,T,J}:=(P_{\cN_J}\otimes I_\cH) K_{f,T}, $$
 where $K_{f,T}$ is
the Poisson kernel given by \eqref{Po-ker}.

The next result provides a functional calculus for c.n.c. $n$-tuples
of operators in  noncommutative varieties  $\cV_{f,J}(\cH)$.

\begin{theorem}
\label{funct-calc2} Let $T:=(T_1,\ldots, T_n)\in B(\cH)^n$ be a
c.n.c. $n$-tuple of operators in the
 noncommutative variety
  $\cV_{f,J}(\cH)$,
 where $J$ is  a $w^*$-closed two-sided
 ideal of $F_n^\infty(\cD_f)$.
Then the map
$$\Psi_{f,T,J}:F_n^\infty(\cV_{f,J})\to B(\cH), \quad
\Psi_{f,T,J}(g(B_1,\ldots, B_n)):=g(T_1,\ldots, T_n), $$
 where
$g(T_1,\ldots, T_n)$ is defined by the
$F_n^\infty(\cD_f)$-functional calculus, has the following
properties:
\begin{enumerate}
\item[(i)]
 $\Psi_{f,T,J} $ is WOT-continuous (resp. SOT-continuous) on bounded sets;
 \item[(ii)]
 $\Psi_{f,T,J}$ is   a unital completely contractive
homomorphism.
\end{enumerate}
Moreover,  the Poisson transform $\Psi_{f,T,J} $ is uniquely
determined by the relation
$$
\Psi_{f,T,J} (g(B_1,\ldots, B_n))
K_{f,T,J}^*=K_{f,T,J}^*[g(B_1,\ldots, B_n)\otimes I_\cH],
$$
 where $K_{f,T,J}$ is the constrained  Poisson kernel associated
with $T\in \cV_{f,J}(\cH)$.

 In particular, if $T$ is of class $C_{\cdot 0}$, then
  the
Poisson transform
$$\Lambda_{f,T,J}:B(\cN_J)\to B(\cH), \quad
\Lambda_{f,T,J}(X):=K_{f,T,J}^*(X\otimes I)K_{f,T,J}$$
 is
$w^*$-continuous  and has  the properties:
\begin{enumerate}
\item[(iii)]
 $
\Lambda_{f,T,J}(g(B_1,\ldots, B_n))=g(T_1,\ldots, T_n) $ ~for any
$g(W_1,\ldots, W_n)\in F_n^\infty(\cD_f)$;
\item[(iv)]
$\Lambda_{f,T,J}(B_\alpha B_\beta^*)=T_\alpha T_\beta^*$ ~for any
$\alpha,\beta\in \FF_n^+$.
\end{enumerate}
 \end{theorem}

\begin{proof}
Note first the the map $\Psi_{f,T,J}$ is well-defined. Indeed, let
$g(W_1,\ldots, W_n)$ and $ \psi(W_1,\ldots, W_n)$ be in
$F_n^\infty(\cD_f)$ and assume that $g(B_1,\ldots,
B_n)=\psi(B_1,\ldots, B_n)$. Using  the
$F_n^\infty(\cD_f)$-functional calculus, we obtain
$$P_{\cN_J}(g(W_1,\ldots, W_n)-\psi(W_1,\ldots, W_n))|_{\cN_J}=0.
$$
Due to Theorem \ref{F/J}, we have $g(W_1,\ldots,
W_n)-\psi(W_1,\ldots, W_n)\in J$. On the other hand, since
$(T_1,\ldots, T_n)$ is a $J$-constrained $n$-tuple of operators,
i.e.,  $\varphi(T_1,\ldots, T_n)=0$ for $\varphi\in J$,  we deduce
that $g(T_1,\ldots, T_n)=\psi(T_1,\ldots, T_n)$, which proves our
assertion.

The next step is to show that $g(T_1,\ldots, T_n)$ satisfies and is
uniquely determined by the equation
$$
g(T_1,\ldots, T_n)K_{f,T,J}^*=K_{f,T,J}^*[g(B_1,\ldots, B_n)\otimes
I_\cH]
$$
for any $g(B_1,\ldots, B_n)\in F_n^\infty(\cV_{f,J})$,  where
$K_{f,T,J}^*$ is the constrained  Poisson kernel associated with
$T\in \cV_{f,J}(\cH)$.

According to Theorem \ref{funct-calc} (see relation \eqref{gKKg}),
we have
 \begin{equation}
 \label{ker-p}
 \left<(\varphi(W_1,\ldots, W_n)^*\otimes I_\cH)K_{f,T}h,1\otimes k\right>=
 \left<K_{f,T}\varphi(T_1,\ldots, T_n)^*h,1\otimes k\right>
 \end{equation}
 for any $\varphi(W_1,\ldots, W_n)\in F_n^\infty(\cD_f)$ and $h,k\in \cH$.
 Note that if $\varphi(W_1,\ldots, W_n)\in J$, then
 $\varphi(T_1,\ldots, T_n)=0$, and relation \eqref{ker-p} implies
 $$\left<K_{f,T}h, \varphi(W_1,\ldots, W_n)(1)\otimes k\right>=0$$ for
 any $h,k\in \cH$.
 Taking into account the definition of $\cM_J$, we deduce that
 $K_{f,T}(\cH)\subseteq \cN_J\otimes \overline{\Delta_{f,T}(\cH)}$.
  This implies
 \begin{equation}
 \label{KT}
 K_{f,T}=\left(P_{ \cN_J}\otimes I_\cH\right) K_{f,T}=K_{f,T,J}.
 \end{equation}
 Since $J$ is a left ideal of $F_n^\infty(\cD_f)$, the subspace
  $\cN_J$ is
  invariant
 under each operator $W_1^*,\ldots, W_n^*$ and, therefore,
  $B_\alpha=P_{ \cN_J}W_\alpha|\cN_J$,
 $\alpha\in \FF_n^+$.
 On the other hand, since $(B_1,\ldots, B_n)\in \cD_f(\cN_J)$ is of class
  $C_{\cdot 0}$, we can use
 the $F_n^\infty(\cD_f)$-functional calculus of Theorem \ref{funct-calc}
 to deduce that
 \begin{equation}
 \label{fnn}
 g(B_1,\ldots, B_n)=P_{ \cN_J}g(W_1,\ldots, W_n)|\cN_J
 \end{equation}
 for any $g(W_1,\ldots, W_n)\in F_n^\infty(\cD_f)$.
 Taking into account relations \eqref{gKKg},
  \eqref{KT}, and \eqref{fnn},
 we obtain
 \begin{equation*}
 \begin{split}
 K_{f,T,J}g(T_1,\ldots, T_n)^*&=\left(P_{ \cN_J}\otimes I_\cH\right)
 K_{f,T} g(T_1,\ldots, T_n)^*\\
 &=\left(P_{ \cN_J}\otimes I_\cH\right)
 [g(W_1,\ldots, W_n)^*\otimes I_\cH]\left(P_{ \cN_J}\otimes
  I_\cH\right)K_{f,T}\\
 &=
 \left[\left(P_{ \cN_J}g(W_1,\ldots, W_n)|
 \cN_J\right)^*\otimes I_\cH\right] K_{f,T,J}\\
&=
 \left[g(B_1,\ldots, B_n)^*\otimes I_\cH\right]K_{f,T,J}.
 \end{split}
 \end{equation*}
Therefore, we have
\begin{equation}\label{KFFK}
 g(T_1,\ldots, T_n)K_{f,T,J}^*=K_{f,T,J}^*\left[g(B_1,\ldots, B_n)
 \otimes I_\cH\right]
\end{equation}
 for any $g(B_1,\ldots, B_n)\in F_n^\infty(\cV_{f,J})$.

On the other hand, since $T:=(T_1,\ldots, T_n)\in \cD_f(\cH)$ is
c.n.c. the Poisson kernel $K_{f,T}:\cH \to F^2(H_n)\otimes
\overline{\Delta_{f,T}(\cH)}$ is one-to-one (see Proposition
\ref{cnc}). Since $K_{f,T}(\cH)\subset \cN_J\otimes
\overline{\Delta_{f,T}(\cH)}$ and $K_{f,T,J}=\left(P_{ \cN_J}\otimes
I_\cH\right) K_{f,T}$, we deduce that the constrained Poisson kernel
$K_{f,T,J}$ is also one-to-one. Consequently, the range of
$K_{f,T,J}^*$ is dense in $\cH$ and the equation \eqref{KFFK}
uniquely determines the operator $g(T_1,\ldots,T_n)$.

If $\varphi\in J$ and $g(W_1,\ldots, W_n)\in F_n^\infty(\cD_f)$,
then Theorem \ref{funct-calc} implies
$$
\|g(T_1,\ldots, T_n)\|= \|(g+\varphi)(T_1,\ldots, T_n)\|\leq
\|g(W_1,\ldots, W_n)+\varphi(W_1,\ldots, W_n)\|.
$$
Hence, and using now Theorem \ref{F/J}, we obtain
\begin{equation*}
\begin{split}
 \|g(T_1,\ldots,
T_n)\|&\leq\text{\rm dist}(g(W_1,\ldots, W_n), J)
=\|P_{\cN_J}g(W_1,\ldots, W_n)|_{\cN_J}\|\\
&=\|g(B_1,\ldots, B_n)\|.
\end{split}
\end{equation*}
Similarly, we can show that $\Psi_{f,T,J}$ is a completely
contractive map.

 To prove part (i), let $\{g_i(B_1,\ldots, B_n)\}$
be a bounded net in $F_n^\infty(\cV_{f,J})$  which is WOT (resp.
SOT) convergent to an element $g(B_1,\ldots, B_n)$ in the algebra
$F_n^\infty(\cV_{f,J})$. Then, due to standard facts in functional
analysis, the net $g_i(B_1,\ldots, B_n)\otimes I_\cH$ is WOT (resp.
SOT) convergent to $g(B_1,\ldots, B_n)\otimes I_\cH$.

On the other hand, as shown above, $\|g_i(T_1,\ldots, T_n)\|\leq
\|g_i(B_1,\ldots, B_n)\|$ and, consequently, the net
$\{g_i(T_1,\ldots, T_n)\}$ is bounded. Now, one can use  relation
\eqref{KFFK}, the fact that the range of $K_{f,T,J}^*$ is dense in
$\cH$, and a standard approximation argument to deduce that
$g_i(T_1,\ldots, T_n)$ is WOT (resp. SOT) convergent to
$g(T_1,\ldots, T_n)$.

 Le us  prove part (ii). First, note that we have
 already proved that $\Psi_{f,T,J}$ is a completely contractive map.
 Due to relation \eqref{KFFK},  the map $\Psi_{f,T,J}$ is  a homomorphism
  on polynomials in $B_1,\ldots, B_n$ and the identity.
      Using the WOT-continuity
 of $\Psi_{f,T,J}$ and
the WOT-density of the polynomials  in $F_n^\infty(\cV_{f,J})$, one
can deduce that $\Psi_{f,T,J}$ is also a homomorphism on
$F_n^\infty(\cV_{f,J})$.

Now, notice that if $T$ is of class $C_{\cdot 0}$, then the Poisson
kernel  $K_{f,T,J}$ is an isometry and, due to \eqref{KFFK}, we have
$$
\Psi_{f,T,J}(g(B_1,\ldots, B_n)):=K_{f,T,J}^* (g(B_1,\ldots,
B_n)\otimes I_\cH) K_{f,T,J}
$$
for any $g(B_1,\ldots, B_n)\in F_n^\infty(\cV_{f,J})$.
 The rest of the proof is now clear.
\end{proof}

\begin{theorem}\label{poiss-poly}
Let $J$ be a $w^*$-closed two-sided
 ideal of $F_n^\infty(\cD_f)$ generated by some homogenous
 polynomials $q_1,\ldots q_d$, and
let $T:=(T_1,\ldots, T_n)$ be an  arbitrary $n$-tuple of operators
in the
 noncommutative variety $\cV_{f,J}(\cH)$, i.e., $q_j(T_1,\ldots,
 T_n)=0$ for  $j=1,\ldots, n$.
Then the constrained Poisson transform
$$\Psi_{f,T,J}:\overline{\text{\rm span}}\{B_\alpha B_\beta^*:\
\alpha,\beta\in \FF_n^+\}\to B(\cH)
$$
defined by
$$
\Psi_{f,T,J}(X):=\lim_{r\to 1} K_{f,rT,J}^*(X\otimes
I_\cH)K_{f,rT,J},
$$
where the limit exists in the operator norm, has the following
properties:
\begin{enumerate}
\item[(i)] $\Psi_{f,T,J}(B_\alpha B_\beta^*)=T_\alpha T_\beta^*$ for
any $\alpha,\beta\in \FF_n^+$.
\item[(ii)] $\Psi_{f,T,J}$ is a unital completely contractive linear
map.
\end{enumerate}
If, in addition,  $T$ is a c.n.c. $n$-tuple of operators in
$\cD_f(\cH)$, then
$$
g(T_1,\ldots, T_n)= \text{\rm SOT-}\lim_{r\to 1}
K_{f,rT,J}^*(g(B_1,\ldots, B_n)\otimes I_\cH)K_{f,rT,J}
$$
for  any \ $g(B_1,\ldots, B_n)\in F_n^\infty(\cV_{f,J})$, where
$g(T_1,\ldots, T_n)$ is defined by the
$F_n^\infty(\cD_f)$-functional calculus.
\end{theorem}

\begin{proof}
As  in the proof of  Lemma 2.11 from \cite{Po-unitary}, we can use
the fact that $J$ is generated by homogenous polynomials to show
that $K_{f,rT}(\cH)\subset \cN_J \otimes \cH$ and
$$
K_{f,rT,J}g_r(T_1,\ldots, T_n)^*=(g(B_1,\ldots, B_n)^*\otimes
I_\cH)K_{f,rT,J}
$$
for any $0<r<1$ and $g(B_1,\ldots, B_n)\in \cW(B_1,\ldots, B_n)$.
Since $rT$ is of class $C_{\cdot 0}$, the constrained Poisson kernel
$K_{f,rT,J}$ is an isometry and, consequently,
\begin{equation}
\label{KfTJ} g_r(T_1,\ldots, T_n)= K_{f,rT,J}^*(g(B_1,\ldots,
B_n)\otimes I_\cH)K_{f,rT,J}.
\end{equation}
Using this relation,  the proof of (i) and (ii) is similar to that
of Theorem \ref{Poisson-C*}, we should omit it. If
 $T$ is a c.n.c. $n$-tuple of operators in
$\cD_f(\cH)$, then Theorem \ref{funct-calc} shows that
$g(T_1,\ldots, T_n)=\text{SOT-}\lim\limits_{r\to 1} g_r(T_1,\ldots,
T_n)$. Using again relation \ref{KfTJ}, we deduce the last part of
the theorem.
\end{proof}

\bigskip

\section{Weighted shifts, symmetric weighted Fock spaces, and  multipliers}
 \label{Symmetric}

In this section we find the joint right spectrum of the weighted
left creation operators $(W_1,\ldots, W_n)$ associated with the
noncommutative domain $\cD_f$ and identify the $w^*$-continuous
multiplicative  linear functional of the Hardy algebra
$F_n^\infty(\cD_f)$. We introduce the symmetric weighted Fock space
$F_s^2(\cD_f)$ and identify the algebra $F_n^\infty(\cV_{f,J_c})$ of
all its multipliers, which turns out to be   reflexive.

 Let $f=\sum_{|\alpha|\geq 1} a_\alpha X_\alpha$ be a
positive regular free holomorphic function on $B(\cH)^n$ and define
$$\cD_f^\circ(\CC):=\left\{ \lambda=(\lambda_1,\ldots, \lambda_n)\in
\CC^n:\ \sum_{|\alpha|\geq 1} a_\alpha |\lambda_\alpha|^2<1\right\},
$$
where  $\lambda_\alpha:=\lambda_{i_1}\cdots \lambda_{i_m}$ if
$\alpha=g_{i_1}\cdots g_{i_m}\in \FF_n^+$, and $\lambda_{g_0}$=1.

\begin{theorem}\label{eigenvectors}
Let $(W_1,\ldots, W_n)$ (resp. $(\Lambda_1,\ldots, \Lambda_n)$) be
the weighted left (resp. right)  creation operators associated with
the noncommutative domain $\cD_f$. The eigenvectors for
$W_1^*,\ldots, W_n^*$ (resp. $\Lambda_1^*,\ldots, \Lambda_n^*$) are
precisely the vectors
$$
z_\lambda:= \sum_{\beta\in \FF_n^+} \sqrt{b_\beta}
\overline{\lambda}_\beta e_\beta \in F^2(H_n)\quad \text{ for }\
\lambda=(\lambda_1,\ldots, \lambda_n)\in \cD_f^\circ(\CC),
$$
where the coefficients $b_\beta$, $\beta\in \FF_n^+$, are defined by
relation \eqref{b_alpha}. They satisfy the equations
$$
W_i^*z_\lambda=\overline{\lambda}_i z_\lambda, \quad
\Lambda_i^*z_\lambda=\overline{\lambda}_i z_\lambda \qquad \text{
for } \ i=1,\ldots,n.
$$
Each vector $z_\lambda$ is cyclic for $R_n^\infty(\cD_f)$ and
$$
z_\lambda=\left( I-\sum_{|\alpha|\geq 1} a_{\tilde \alpha}
\overline{\lambda}_\alpha \Lambda_\alpha\right)^{-1}(1),
$$
where $\tilde \alpha$ denotes the reverse of $\alpha$.

If $\lambda:=(\lambda_1,\ldots, \lambda_n)\in \cD_f^\circ(\CC)$ and
$\varphi(W_1,\ldots, W_n):=\sum_{\beta\in \FF_n^+} c_\beta W_\beta$
is in the Hardy algebra $F_n^\infty(\cD_f)$, then $\sum_{\beta\in
\FF_n^+} |c_\beta||\lambda_\beta|<\infty$ and the map
$$\Phi_\lambda:F_n^\infty(\cD_f)\to \CC, \qquad
 \Phi_\lambda(\varphi(W_1,\ldots, W_n)):=\varphi(\lambda),
$$
is $w^*$-continuous and multiplicative. Moreover,
$\varphi(W_1,\ldots,
W_n)^*z_\lambda=\overline{\varphi(\lambda)}z_\lambda$ and
$$
\varphi(\lambda)=\left<\varphi(W_1,\ldots, W_n)1,
z_\lambda\right>=\left< \varphi(W_1,\ldots, W_n)u_\lambda,
u_\lambda\right>,
$$
where $u_\lambda:=\frac{z_\lambda}{\|z_\lambda\|}$.
\end{theorem}

\begin{proof}
First, notice that if $\lambda=(\lambda_1,\ldots, \lambda_n)\in
\cD_f^\circ(\CC)$, then $\lambda$ is  of class $C_{\cdot 0}$ with
respect to  $\cD_f(\CC)$. Using relation \eqref{I-Q} in our
particular case, we get
$$
\left(1-\sum_{|\alpha|\geq 1} a_\alpha
|\lambda_\alpha|^2\right)\left(\sum_{\beta\in \FF_n^+} b_\beta
|\lambda_\beta|^2\right)= 1.
$$
Consequently, the vector $z_\lambda$ is in $F^2(H_n)$ and
$$\|z_\lambda\|=\frac{1}{\sqrt{1-\sum_{|\alpha|\geq 1} a_\alpha
|\lambda_\alpha|^2}}.
$$
Due to  relation \eqref{WbWb}, we have
$$
W_i^* e_\alpha =\begin{cases} \frac
{\sqrt{b_\gamma}}{\sqrt{b_{\alpha}}}e_\gamma& \text{ if }
\alpha=g_i\gamma \\
0& \text{ otherwise. }
\end{cases}
$$
A simple computation shows that $W_i^*
z_\lambda=\overline{\lambda}_{i} z_\lambda$ for $i=1,\ldots, n$.
Similarly, one can use relation \eqref{WbWb-r} to prove that
$\Lambda_i^*z_\lambda=\overline{\lambda}_i z_\lambda$ for
$i=1,\ldots, n$.

Conversely, let $z=\sum_{\beta\in \FF_n^+} c_\beta e_\beta \in
F^2(H_n)$ and assume that $W_i^*z=\overline{\lambda}_i z$,
$i=1,\ldots,n$, for some $n$-tuple $(\lambda_1,\ldots, \lambda_n)\in
\CC^n$.  Using the definition of the weighted left creation
operators $W_1,\ldots, W_n$, we deduce that
\begin{equation*}
\begin{split}
c_\alpha&=\left< z, e_\alpha\right>=\left< z, \sqrt{b_\alpha}
W_\alpha(1)\right>\\
&=\sqrt{b_\alpha}\left< W_\alpha^* z,1\right>=
\sqrt{b_\alpha}\overline{\lambda}_\alpha\left< z,1\right>\\
&=c_0\sqrt{b_\alpha}\overline{\lambda}_\alpha
\end{split}
\end{equation*}
 for any $\alpha\in \FF_n^+$, whence
  $z= a_0\sum_{\beta\in \FF_n^+} \sqrt{b_\beta}
\overline{\lambda}_\beta e_\beta$.
 Since $z\in F^2(H_n)$, we must have $\sum_{\beta\in \FF_n^+} b_\beta
|\lambda_\beta|^2<\infty$. On the other hand, relation
\eqref{b_alpha} implies
$$
\sum_{j=0}^m \left(\sum_{|\alpha|\geq 1} a_\alpha
|\lambda_\alpha|^2\right)^j\leq \sum_{\beta\in \FF_n^+} b_\beta
|\lambda_\beta|^2<\infty
$$
for any $m\in \NN$. Letting $m\to\infty$ in the relation above, we
must have $\sum_{|\alpha|\geq 1} a_\alpha |\lambda_\alpha|^2<1$,
whence $(\lambda_1,\ldots, \lambda_n)\in \cD_f^\circ(\CC)$. A
similar result can be proved for the weighted right creation
operators $\Lambda_1,\ldots, \Lambda_n$ if one uses relation
\eqref{WbWb-r}.

Since $\sum_{|\alpha|\geq 1} a_{\tilde \alpha}
  \Lambda_\alpha \Lambda_\alpha^*$ is SOT-convergent (see Section
  \ref{Noncommutative}) and

  \begin{equation*}
  \begin{split}
  \left\|\sum_{|\alpha|\geq 1} a_{\tilde \alpha}
  \overline{\lambda}_{\tilde\alpha} \Lambda_\alpha\right\|&\leq
\left\|\sum_{|\alpha|\geq 1} a_{\tilde \alpha}
  \Lambda_\alpha \Lambda_\alpha^*\right\|
  \left(\sum_{|\alpha|\geq 1} a_\alpha
  |\lambda_\alpha|^2\right)^{1/2}\\
  &\leq\left(\sum_{|\alpha|\geq 1} a_\alpha
  |\lambda_\alpha|^2\right)^{1/2}<1,
\end{split}
  \end{equation*}
the operator $\left( I-\sum_{|\alpha|\geq 1} a_{\tilde \alpha}
\overline{\lambda}_{\tilde\alpha} \Lambda_\alpha\right)^{-1}$ is
well-defined. Due to the results of Section \ref{Noncommutative}, we
have
$$
\left( I-\sum_{|\alpha|\geq 1} a_{\tilde \alpha}
\overline{\lambda}_{\tilde\alpha} \Lambda_\alpha\right)^{-1}=
\sum_{\beta\in \FF_n^+}  b_{\tilde
    \beta}\overline{\lambda}_{\tilde \beta} \Lambda_\beta.
    $$
    Hence, and using relation \eqref{WbWb-r}, we obtain
    \begin{equation*}
    z_\lambda=\sum_{\beta\in
\FF_n^+}  b_{\tilde
    \beta}\overline{\lambda}_{\tilde \beta} \Lambda_\beta (1)=
    \left( I-\sum_{|\alpha|\geq 1} a_{\tilde \alpha}
\overline{\lambda}_\alpha \Lambda_\alpha\right)^{-1}(1),
    \end{equation*}
which implies that $z_\lambda$ is a cyclic  vector for
$R_n^\infty(\cD_f)$.

Now, let us prove the last part of the theorem. Since
$\varphi(W_1,\ldots, W_n)=\sum_{\beta\in \FF_n^+} c_\beta W_\beta$
is in the Hardy algebra $F_n^\infty(\cD_f)$, we have $\sum_{\beta\in
\FF_n^+} |c_\beta|^2 \frac{1}{b_\beta}<\infty$ (see Section
\ref{Hardy}). As shown above, if $\lambda=(\lambda_1,\ldots,
\lambda_n)\in \cD_f^\circ(\CC)$, then $\sum_{\beta\in \FF_n^+}
b_\beta |\lambda_\beta|^2<\infty$. Applying Cauchy's inequality, we
have
$$
\sum_{\beta\in \FF_n^+} |c_\beta||\lambda_\beta|\leq
\left(\sum_{\beta\in \FF_n^+} |c_\beta|^2
\frac{1}{b_\beta}\right)^{1/2} \left(\sum_{\beta\in \FF_n^+} b_\beta
|\lambda_\beta|^2\right)^{1/2}<\infty.
$$
Note  also that
\begin{equation*}
\begin{split}
\left<\varphi(W_1,\ldots, W_n)1, z_\lambda\right>&= \left<
\sum_{\beta\in \FF_n^+}c_\beta \frac{1}{\sqrt{b_\beta}} e_\beta,
\sum_{\beta\in \FF_n^+} \sqrt{b_\beta} \overline{\lambda}_\beta
e_\beta \right>\\
&=\sum_{\beta\in\FF_n^+} c_\beta \lambda_\beta
=\varphi(\lambda_1,\ldots, \lambda_n).
\end{split}
\end{equation*}
Now, for each $\beta\in \FF_n^+$, we have
\begin{equation*}
\begin{split}
\left<\varphi(W_1,\ldots, W_n)^* z_\lambda,
\frac{1}{\sqrt{b_\alpha}}e_\beta\right>&= \left< z_\lambda,
\varphi(W_1,\ldots, W_n) W_\beta (1)\right>\\
&=\overline {\lambda_\beta \varphi(\lambda)} =\left<\overline {
\varphi(\lambda)} z_\lambda, \frac{1}{\sqrt{b_\alpha}}
e_\beta\right>.
\end{split}
\end{equation*}
Hence, we deduce that
\begin{equation}
\label{fi*z} \varphi(W_1,\ldots, W_n)^* z_\lambda=\overline {
\varphi(\lambda)} z_\lambda.
\end{equation}
One can easily see that
\begin{equation*}
\begin{split}
\left< \varphi(W_1,\ldots, W_n)u_\lambda, u_\lambda\right>&=
\frac{1}{\|z_\lambda\|^2}\left< z_\lambda,\varphi(W_1,\ldots, W_n)^*
z_\lambda\right>\\
&=\frac{1}{\|z_\lambda\|^2}\left< z_\lambda,\overline {
\varphi(\lambda)} z_\lambda\right>= \varphi(\lambda).
\end{split}
\end{equation*}
The fact that  the map $\Phi_\lambda$ is multiplicative and
$w^*$-continuous is now obvious.  This completes the proof.
\end{proof}

 We recall   that the joint  right spectrum
  $\sigma_r(T_1,\ldots, T_n)$ of an   $n$-tuple
 $(T_1,\ldots, T_n)$ of operators
   in $B(\cH)$ is the set of all $n$-tuples
    $(\lambda_1,\ldots, \lambda_n)$  of complex numbers such that the
     right ideal of $B(\cH)$  generated by the operators
     $\lambda_1I-T_1,\ldots, \lambda_nI-T_n$ does
      not contain the identity operator.
We recall \cite{Po-unitary} that
  $(\lambda_1,\ldots, \lambda_n)\notin \sigma_r(T_1,\ldots, T_n)$
  if and only if  there exists $\delta>0$ such that
  $\sum\limits_{i=1}^n (\lambda_iI-T_i)
  (\overline{\lambda}_iI-T_i^*)\geq \delta I$.

\begin{proposition}\label{right-spec} If ~$(W_1,\ldots, W_n)$~ are the
weighted left creation operators associated with the noncommutative
domain $\cD_f$, then the right joint spectrum $\sigma_r(W_1,\ldots,
W_n)$ coincide with $\cD_f(\CC)$.
\end{proposition}

\begin{proof}
Let $(\lambda_1,\ldots, \lambda_n)\in \sigma_r(W_1,\ldots, W_n)$.
 Since the left ideal of $B(F^2(H_n))$ generated by the operators
 $W_1^*-\overline{\lambda}_1 I, \ldots,
 W_n^*-\overline{\lambda}_n I$ does not contain the identity,
 there is a pure
  state
 $\varphi$ on $B(F^2(H_n))$ such that
 $\varphi(X(W_i^*-\overline{\lambda}_i I))=0$
  for
 any $X\in B(F^2(H_n))$
 and $i=1,\ldots, n$. In particular, we have
 $\varphi(W_i)= \lambda_i=\overline{
 \varphi(W_i^*)}$ and
 $$
 \varphi(W_\alpha W_\alpha^*)=\overline{\lambda}_\alpha
  \varphi(W_\alpha)=|\lambda_\alpha|^2,\qquad \alpha\in \FF_n^+.
 $$
 Hence, using relation \eqref{WWinD},  we infer that
 \begin{equation*}
 \begin{split}
   \sum_{1\leq |\alpha|\leq m} a_\alpha|\lambda_\alpha|^2  &=
   \varphi\left(\sum_{1\leq |\alpha|\leq m} a_\alpha W_\alpha
   W_\alpha^*\right)
   \leq \left\|\sum_{1\leq |\alpha|\leq m} a_\alpha W_\alpha
   W_\alpha^*\right\|\\
   &\leq \left\|\sum_{ |\alpha|\geq 1} a_\alpha W_\alpha
   W_\alpha^*\right\|=1
\end{split}
 \end{equation*}
 for any $m=1,2,\ldots.$ Therefore $\sum_{|\alpha|\geq 1} a_\alpha
 |\lambda_\alpha|^2\leq 1$, which proves  that
 $(\lambda_1,\ldots, \lambda_n)\in \cD_f(\CC)$.
Now, let  $\mu:=(\mu_1,\ldots, \mu_n)\in \cD_f(\CC)$  and assume
that
 there is $\delta>0$ such that
 $$\sum_{i=1}^n \|(W_i-\mu_iI)^*h\|^2\geq \delta \|h\|^2 \quad
 \text{ for all } \ h\in F^2(H_n).
 $$
 Taking $h=z_\lambda:=\sum_{\beta\in \FF_n^+}
 \sqrt{b_\beta}\overline{\lambda}_\beta e_\beta$,
 $\lambda\in \cD_f^\circ(\CC)$, in the
 inequality above, and using  Theorem \ref{eigenvectors}, we obtain
 $$
 \sum_{i=1}^n |\lambda_i-\mu_i|^2\geq \delta\quad \text{ for all }\
 \lambda\in \cD_f^\circ(\CC),
 $$
 which is a contradiction. Due to the remarks preceding the theorem,
 we deduce that $\mu:=(\mu_1,\ldots, \mu_n)\in \sigma_r(W_1,\ldots,
 W_n)$, which completes the proof.
\end{proof}

Now, we define the symmetric weighted Fock space associated with the
noncommutative domain $\cD_f$.

We need a few definitions. For each $\lambda=(\lambda_1,\ldots,
\lambda_n)$ and each $n$-tuple ${\bf k}:=(k_1,\ldots, k_n)\in
\NN_0^n$, where $\NN_0:=\{0,1,\ldots \}$, let $\lambda^{\bf
k}:=\lambda_1^{k_1}\cdots \lambda_n^{k_n}$. For each ${\bf k}\in
\NN_0^n$, we denote
$$
\Lambda_{\bf k}:=\{\alpha\in \FF_n^+: \ \lambda_\alpha =\lambda^{\bf
k} \text{ for all } \lambda\in \CC^n\}.
$$
For each ${\bf k}\in \NN_0^n$, define the vector
$$
w^{\bf k}:=\frac{1}{\gamma_{\bf k}} \sum_{\alpha \in \Lambda_{\bf
k}} \sqrt{b_\alpha} e_\alpha\in F^2(H_n), \quad  \text{ where } \
\gamma_{\bf k}:=\sum_{\alpha\in \Lambda_{\bf k}} b_\alpha
$$
 and the
coefficients $b_\alpha$, $\alpha\in \FF_n^+$, are defined by
relation \eqref{b_alpha}. Note that the set  $\{w^{\bf k}:\ {\bf
k}\in \NN_0^n\}$ consists  of orthogonal vectors in $F^2(H_n)$ and
$\|w^{\bf k}\|=\frac{1}{\sqrt{\gamma_{\bf k}}}$. We denote by
$F_s^2(\cD_f)$ the closed span of these vectors, and call it the
symmetric weighted  Fock space associated with the noncommutative
domain $\cD_f$.

\begin{theorem}\label{symm-Fock} Let $f=\sum_{|\alpha|\geq 1} a_\alpha X_\alpha$ be a
positive regular free holomorphic function on $B(\cH)^n$ and  let
$J_c$ be the $w^*$-closed two-sided ideal of the Hardy algebra
$F_n^\infty(\cD_f)$ generated by the commutators
$$W_iW_j-W_j W_i,\qquad i,j=1,\ldots, n.
$$
 Then the following statements hold:
 \begin{enumerate}
 \item[(i)]
 $
 F_s^2(\cD_f)=\overline{\text{\rm span}}\{z_\lambda: \ \lambda\in
\cD^\circ_f(\CC)\}=\cN_{J_c}:=F^2(H_n)\ominus \overline{J_c(1)}. $
\item[(ii)] The symmetric weighted Fock space $F_s^2(\cD_f)$ can be
identified with the Hilbert space $H^2(\cD_f^\circ(\CC))$ of all
functions $\varphi:\cD^\circ_f(\CC)\to \CC$ which admit a power
series representation $\varphi(\lambda)=\sum_{{\bf k}\in \NN_0^n}
c_{\bf k} \lambda^{\bf k}$ with
$$
\|\varphi\|_2=\sum_{{\bf k}\in \NN_0^n}|c_{\bf
k}|^2\frac{1}{\gamma_{\bf k}}<\infty.
$$
More precisely, every  element  $\varphi=\sum_{{\bf k}\in \NN_0^n}
c_{\bf k} w^{\bf k}$ in $F_s^2(\cD_f)$  has a functional
representation on $\cD_f^\circ(\CC)$ given by
$$
\varphi(\lambda):=\left<\varphi, z_\lambda\right>=\sum_{{\bf k}\in
\NN_0^n} c_{\bf k} \lambda^{\bf k}, \quad \lambda=(\lambda_1,\ldots,
\lambda_n)\in \cD^\circ_f(\CC),
$$
and
$$
|\varphi(\lambda)|\leq
\frac{\|\varphi\|_2}{\sqrt{1-\sum_{|\alpha|\geq 1} a_\alpha
|\lambda_\alpha|^2}},\quad \lambda=(\lambda_1,\ldots, \lambda_n)\in
\cD_f^\circ(\CC).
$$
\item[(iii)]
 The mapping $K_f:\cD_f^\circ(\CC)\times
\cD_f^\circ(\CC)\to \CC$ defined by
$$
K_f(\mu,\lambda):= \left<z_\lambda,
z_\mu\right>=\frac{1}{1-\sum_{|\alpha|\geq 1} a_\alpha \mu_\alpha
\overline{\lambda}_\alpha}\quad \text{ for all }\ \lambda,\mu\in
\cD_f^\circ(\CC)
$$
is positive definite.
\end{enumerate}
\end{theorem}

\begin{proof}
First, notice that $z_\lambda=\sum_{{\bf k}\in \NN_0^n}
\overline{\lambda}^{\bf k} \gamma_{\bf k} w^{\bf k}$, $\lambda \in
\cD^\circ_f(\CC)$, and therefore
 $$\overline{\text{\rm span}}\{z_\lambda: \ \lambda\in
\cD^\circ_f(\CC)\}\subseteq F_s^2(\cD_f).
$$
Now, we prove  that $w^{\bf k}\in \cN_{J_c}:=F^2(H_n)\ominus
\overline{J_c(1)}$. Indeed, relation \ref{WbWb} implies
\begin{equation}
\label{WWWWW}
W_\gamma(W_jW_i-W_iW_j)W_\beta(1)=\frac{1}{\sqrt{b_{\gamma
g_jg_i\beta}}} e_{\gamma g_jg_i\beta}- \frac{1}{\sqrt{b_{\gamma
g_ig_j\beta}}} e_{\gamma g_ig_j\beta}.
\end{equation}
Hence, we deduce that
$$
\left<\sum_{\alpha \in \Lambda_{\bf k}} \sqrt{b_\alpha} e_\alpha,
W_\gamma(W_jW_i-W_iW_j)W_\beta(1)\right>=0
$$
for any ${\bf k}\in \NN_0^n$. This shows that $w^{\bf k}\in
\cN_{J_c}$ and proves our assertion. Consequently,
$F_s^2(\cD_f)\subseteq \cN_{J_c}$. To complete the proof of part
(i), it is enough to show that
$$
\overline{\text{\rm span}}\{z_\lambda: \ \lambda\in
\cD^\circ_f(\CC)\}=\cN_{J_c}.
$$
Assume that there is a vector $ x:=\sum_{\beta\in \FF_n^+} c_\beta
e_\beta\in \cN_{J_c}$ and $x\perp z_\lambda$ for all $\lambda\in
\cD_f^\circ(\CC)$. Then
\begin{equation*}
\begin{split}
\left<\sum_{\beta\in \FF_n^+} c_\beta e_\beta, z_\lambda\right>&=
\sum_{\beta\in \FF_n^+} c_\beta\sqrt{b_\beta} \lambda_\beta\\
&=\sum_{{\bf k}\in \NN_0^n}\left(\sum_{\beta\in\Lambda_{\bf k}}
c_\beta \sqrt{b_\beta}\right)\lambda^{\bf k}=0
\end{split}
\end{equation*}
for any $\lambda\in \cD^\circ_f(\CC)$. Since $\cD^\circ_f(\CC)$
contains  an open ball in $\CC^n$, we deduce that
\begin{equation}\label{sigma=0}
\sum_{\beta\in \Lambda_{\bf k}} c_\beta \sqrt{b_\beta}=0 \quad
\text{ for all } \ {\bf k}\in \NN_0^n.
\end{equation}
Fix $\beta_0\in \Lambda_{\bf k}$ and let $\beta\in \Lambda_{\bf k}$
be such that $\beta$ is obtained from $\beta_0$ by transposing just
two generators. So we can assume that $\beta_0=\gamma g_j g_i\omega$
and $\beta=\gamma g_i g_j\omega$ for some $\gamma,\omega\in \FF_n^+$
and $i\neq j$, $i,j=1,\ldots,n$. Since $x\in
\cN_{J_c}=F^2(H_n)\ominus\overline{J_c(1)}$, we must have
$$
\left<x,W_\gamma(W_jW_i-W_iW_j)W_\omega(1)\right>=0,
$$
which implies
$\frac{c_{\beta_0}}{\sqrt{b_{\beta_0}}}=\frac{c_\beta}{\sqrt{b_{\beta}}}$.
Since any element $\gamma\in \Lambda_{\bf k}$ can be obtained from
$\beta_0$  by successive  transpositions, repeating the above
argument, we deduce that
$$
\frac{c_{\beta_0}}{\sqrt{b_{\beta_0}}}=\frac{c_\gamma}{\sqrt{b_{\gamma}}}
\quad \text{ for all }\ \gamma\in \Lambda_{\bf k}.
$$
Setting $t:=\frac{c_{\beta_0}}{\sqrt{b_{\beta_0}}}$, we have
$c_\gamma=t\sqrt{b_\gamma}$, $\gamma\in \Lambda_{\bf k}$, and
relation  \eqref{sigma=0} implies  $t=0$ (remember that
$b_\beta>0$). Therefore, $c_\gamma=0$ for any $\gamma\in
\Lambda_{\bf k}$ and ${\bf k}\in \NN_0^n$, so $x=0$. Consequently,
we have $\overline{\text{\rm span}}\{z_\lambda: \ \lambda\in
\cD^\circ_f(\CC)\}=\cN_{J_c}. $

Now, let us prove part (ii) of the theorem. Note that
\begin{equation*}
\begin{split}
\left<w^{\bf k},z_\lambda\right>&=\frac{1}{\gamma_{\bf k}} \left<
\sum_{\beta\in\Lambda_{\bf k}}\sqrt{b_\beta} e_\beta,
z_\lambda\right>\\
&=\frac{1}{\gamma_{\bf k}}\sum_{\beta\in\Lambda_{\bf k}} b_\beta
\lambda_\beta =\lambda^{\bf k}
\end{split}
\end{equation*}
for any $\lambda\in \cD_f^\circ (\CC)$ and ${\bf k}\in \NN_0^n$.
Hence, every  element  $\varphi=\sum_{{\bf k}\in \NN_0^n} c_{\bf k}
w^{\bf k}$ in $F_s^2(\cD_f)$  has a functional representation on
$\cD_f^\circ(\CC)$ given by
$$
\varphi(\lambda):=\left<\varphi, z_\lambda\right>=\sum_{{\bf k}\in
\NN_0^n} c_{\bf k} \lambda^{\bf k}, \quad \lambda=(\lambda_1,\ldots,
\lambda_n)\in \cD^\circ_f(\CC),
$$
and
$$|\varphi(\lambda)|\leq \|\varphi\|_2 \|z_\lambda\|=
\frac{\|\varphi\|_2}{\sqrt{1-\sum_{|\alpha|\geq 1} a_\alpha
|\lambda_\alpha|^2}}.
$$
We recall that $\|z_\lambda\|$ was calculated in the proof of
Theorem \ref{eigenvectors}. The identification of $F_s^2(\cD_f)$
with  $H^2(\cD_f^\circ(\CC))$  is now clear.

As in the proof of Theorem \ref{eigenvectors}, we deduce that
$$
\left( I-\sum_{|\alpha|\geq 1} a_{\tilde \alpha}
\overline{\lambda}_{\tilde\alpha} \Lambda_\alpha\right)^{-1}=
\sum_{\beta\in \FF_n^+}  b_{\tilde
    \beta}\overline{\lambda}_{\tilde \beta} \Lambda_\beta
    $$
    if $
 (\lambda_1,\ldots, \lambda_n)\in
\cD_f^\circ(\CC)$.
 Now, the $R_n^\infty(\cD_f)$-functional calculus  applied to the
 $n$-tuple
$(\mu_1,\ldots, \mu_n)\in\cD_f^\circ(\CC)\cap \cD_{\tilde{
f}}^\circ(\CC)$  and the equation above imply
$$
\sum_{\beta\in \FF_n^+} b_\beta \mu_\beta\overline{\lambda}_\beta =
\left( I-\sum_{|\alpha|\geq 1} a_{ \alpha}
 \mu_\alpha\overline{\lambda}_{\alpha}\right)^{-1}.
$$
Since
$$K_f(\mu,\lambda):= \left<z_\lambda, z_\mu\right>=
\sum_{\beta\in \FF_n^+} b_\beta \mu_\beta\overline{\lambda}_\beta,
$$
the result in part (iii) follows. The proof is complete.
\end{proof}

\begin{proposition} Let
$J_c$ be the $w^*$-closed two-sided ideal of the Hardy algebra
$F_n^\infty(\cD_f)$ generated by the commutators
$$W_iW_j-W_j W_i,\qquad i,j=1,\ldots, n.
$$
 Then the following statements hold:
 \begin{enumerate}
 \item[(i)] $J_c$ coincides with the $w^*$-closure  of the
 commutator ideal  of $F_n^\infty(\cD_f)$.
\item[(ii)]
$\overline{J_c(1)}$ coincides with
 $$\overline{\text{\rm
span}}\left\{\frac{1}{\sqrt{b_{\gamma g_jg_i\beta}}} e_{\gamma
g_jg_i\beta}- \frac{1}{\sqrt{b_{\gamma g_ig_j\beta}}} e_{\gamma
g_ig_j\beta}:\  \gamma, \beta\in \FF_n^+, i,j=1,\ldots,n\right\}.
$$
\end{enumerate}
\end{proposition}
\begin{proof}
One inclusion is obvious. Since $W_iW_j-W_jW_i\in J_c$ and every
permutation of $k$ objects is a product of transpositions, it is
clear that $W_\alpha W_\beta-W_\beta W_\alpha\in J_c$ for any
$\alpha,\beta\in \FF_n^+$. Consequently, $W_\gamma(W_\alpha
W_\beta-W_\beta W_\alpha)W_\omega\in J_c$ for any $\alpha,\beta,
\gamma,\omega\in \FF_n^+$. Since the polynomials in $W_1,\ldots,
W_n$ are $w^*$ dense in $F_n^\infty(\cD_f)$ (see Corollary
\ref{w*-density}), item (i) follows. Part (ii) is obvious if one
uses relation \eqref{WWWWW}.
\end{proof}

\begin{theorem}\label{w*-funct} A map $\varphi: F_n^\infty(\cD_f)\to
\CC$ is a $w^*$-continuous multiplicative linear functional  if and
only if there exists $\lambda\in \cD_f^\circ (\CC)$  such that
$$
\varphi(A)=\varphi_\lambda(A):=\left<Au_\lambda,
u_\lambda\right>,\quad A\in F_n^\infty(\cD_f),
$$
where $u_\lambda:=\frac{z_\lambda}{\|z_\lambda\|}$.
\end{theorem}

\begin{proof}
Note that $J:=\ker \varphi$ is a $w^*$-closed two-sided ideal of
$F_n^\infty(\cD_f)$  of  codimension one.  The subspace
$\cM_{J}:=\overline{J(1)}$  has codimension one and $\cM_J+\CC
1=F^2(H_n)$. Indeed, assume that there exists a vector $y\in
(\cM_J+\CC 1)^\perp$, $\|y\|=1$. Then we can choose a polynomial
$p:=p(W_1,\ldots, W_n)\in F_n^\infty(\cD_f)$ such that
$\|p(1)-y\|<1$. Since $p-\varphi(p)I\in \ker \varphi=J$, we have
$p(1)-\varphi(p)\in \cM_J$. Since $y\perp \cM_J+\CC 1$, we deduce
that
\begin{equation*}
\begin{split}
1&=\|y\|=\left<y-\varphi(p),y\right>\leq
|\left<h-p(1),y\right>|+|\left<p(1)-\varphi(p),y\right>|\\
&=|\left<h-p(1),y\right>|<1,
\end{split}
\end{equation*}
which is a contradiction. Therefore, $M_J$ is an invariant subspace
under $W_1,\ldots, W_n$  and has codimension one. Due to Theorem
\ref{eigenvectors}, there exists $\lambda\in \cD_f^\circ(\CC)$ such
that $\cM_J=\{z_\lambda\}^\perp$.

On the other hand, due to the same theorem the map  $\varphi:
F_n^\infty(\cD_f)\to \CC$  defined by
$\varphi_\lambda(A):=\left<Au_\lambda, u_\lambda\right>$, \ $ A\in
F_n^\infty(\cD_f)$, is a $w^*$-continuous multiplicative linear
functional, and $\varphi_\lambda(A)=\left<A (1), z_\lambda\right>$.
Hence, $J:=\ker \varphi\subseteq \ker \varphi_\lambda\subset
F_n^\infty(\cD_f)$.  Since both $\ker \varphi$ and $\ker
\varphi_\lambda$ are $w^*$-closed two-sided maximal ideals  of
$F_n^\infty(\cD_f)$ of codimension one, we must have $\ker
\varphi=\ker \varphi_\lambda$. Therefore, $\varphi=\varphi_\lambda$
and the proof is complete.
\end{proof}

\begin{proposition}\label{W(L)}
Let \, $F_n^\infty(\cV_{f,J_c})$ be the $w^*$-closed algebra
generated by operators
$$L_i:=P_{F_s^2(\cD_f)} W_i|_{F_s^2(\cD_f)}, \qquad i=1,\ldots, n
$$
and the identity, where $W_1,\ldots, W_n$ are the weighted left
creation operators associated with $\cD_f$. Then the following
statements hold:
\begin{enumerate}
\item[(i)] $F_n^\infty(\cV_{f,J_c})
=P_{F_s^2(\cD_f)} F_n^\infty(\cD_f)|_{F_s^2(\cD_f)}$;
\item[(ii)]  $F_n^\infty(\cV_{f,J_c})$ can be identified
 with the algebra of  all multipliers of the  Hilbert space $H^2(\cD_f^\circ(\CC))$.
\end{enumerate}
\end{proposition}

\begin{proof}
The first part of the proposition is a particular case of a result
mentioned in the beginning of  Section \ref{Functional II}. To prove
the last part, let $\varphi(W_1,\ldots, W_n)\in F_n^\infty(\cD_f)$
and $\varphi(L_1,\ldots, L_n)=P_{F_s^2(\cD_f)} \varphi(W_1,\ldots,
W_n)|_{F_s^2(\cD_f)}$. Since $z_\lambda\in F_s^2(\cD_f)$ for
$\lambda\in \cD_f^\circ(\CC)$, and  $\varphi(W_1,\ldots,
W_n)^*z_\lambda= \overline{\varphi(\lambda)} z_\lambda$, we have
\begin{equation*}
\begin{split}
[\varphi(L_1,\ldots, L_n) \psi](\lambda)&=\left<\varphi(L_1,\ldots,
L_n)\psi, z_\lambda\right>\\
&=\left<\varphi(W_1,\ldots,
W_n)\psi, z_\lambda\right>\\
&=\left< \psi,\varphi(W_1,\ldots, W_n)^*z_\lambda\right>\\
&=\left< \psi,\overline{\varphi(\lambda)} z_\lambda\right>
=\varphi(\lambda)\psi(\lambda)
\end{split}
\end{equation*}
for any $\psi\in F_s^2(\cD_f)$ and $\lambda\in \cD_f^\circ(\CC)$.
Therefore, the operators in $F_n^\infty(\cV_{f,J_c})$ are
``analytic'' multipliers of $F_s^2(\cD_f)$. Moreover,
$$
\|\varphi(L_1,\ldots, L_n)\|=\sup\{\|\varphi f\|_2:\ f\in
F_s^2(\cD_f), \,  \|f\|\leq 1\}.
$$
Conversely, suppose that $g=\sum_{{\bf k}\in \NN_0} c_{\bf k} w^{\bf
k}$ is a bounded multiplier, i.e., $M_g\in B(F_s^2(\cD_f))$. As in
\cite{DP1} (see Lemma 1.1), using Cesaro means,  one can find a
sequence of polynomials $p_m=\sum  c_{\bf k}^{(m)} w^{\bf k}$ such
that $M_{p_m}$ converges to $M_g$ in the strong operator topology
and, consequently, in  the $w^*$-topology. The latter is due to the
fact that a sequence of operators converges $w^*$ if and only if it
converges WOT. Since $M_{p_m}$ is a polynomial in $L_1,\ldots, L_n$,
we conclude that $M_g\in F_n^\infty(\cV_{f,J_c})$. In particular
$L_i$ is is the multiplier $M_{\lambda_i}$  by the coordinate
function. The proof is complete.
\end{proof}

  If $~A\in B(\cH)~$ then the set of all
invariant subspaces of $~A~$ is denoted by $\text{\rm Lat~}A$. For
any $~\cU\subset B(\cH)~$ we define
$$
\text{\rm Lat~}\cU=\bigcap_{A\in\cU}\text{\rm Lat~}A.
$$
 If $~\cS~$ is any collection of subspaces of $~\cH$,
then we define $\text{\rm Alg~}\cS$ by setting $$\text{\rm
Alg~}\cS:=\{A\in B(\cH):\ \cS\subset\text{\rm Lat~}A\}.$$ We recall
that the algebra $~\cU\subset B(\cH)~$ is reflexive if
$\cU=\text{\rm Alg Lat~}\cU.$

\begin{theorem}\label{reflexivity}
The algebra $F_n^\infty(\cV_{f,J_c})$ is reflexive.
\end{theorem}
\begin{proof} Let $X\in B(F_s^2(\cD_f))$ be an operator   that leaves
invariant all the invariant subspaces under each operator
$L_1,\ldots, L_n$. Due to Theorem \ref{eigenvectors} and Proposition
\ref{W(L)}, we have $L_i^*z_\lambda=\overline{\lambda}_i z_\lambda$
for any $\lambda\in \cD_f^\circ(\CC)$ and $i=1,\ldots,n$. Since
$X^*$ leaves invariant all the invariant subspaces under
$L_1^*,\ldots, L_n^*$, the vector $z_\lambda$ must be an eigenvector
for $X^*$. Consequently, there is a function
$\varphi:\cD_f^\circ(\CC)\to \CC$ such that
$X^*z_\lambda=\overline{\varphi(\lambda)} z_\lambda$ for any
$\lambda\in \cD_f^\circ(\CC)$. Notice that, if $f\in F_s^2(\cD_f)$,
then, due to Theorem \ref{eigenvectors}, $Xf$ has the functional
representation
$$
(Xf)(\lambda)=\left<Xf,z_\lambda\right>=\left<f,X^*z_\lambda\right>=
\varphi(\lambda)f(\lambda)\quad \text{ for all }\ \lambda\in
\cD_f^\circ(\CC).
$$
In particular, if $f=1$, then  the the functional representation  of
$X(1)$  coincide with $\varphi$. Consequently, $\varphi$ admits a
power series representation  on $\cD_f^\circ(\CC)$ and can be
identified with $X(1)\in F_s^2(\cD_f)$. Moreover, the equality above
shows that $\varphi f\in H^2(\cD_f^\circ(\CC))$ for any $f\in
F_s^2(\cD_f)$. Applying Proposition \ref{W(L)}, we deduce that
$X=M_\varphi\in F_n^\infty(\cV_{f,J_c})$. This completes the proof.
\end{proof}

We remark that, in  the particular case when $f=X_1+\cdots +X_n$, we
recover some of the results obtained by Arias and the author,
Davidson and Pitts, and Arveson (see  \cite{ArPo}, \cite{Po-disc},
\cite{DP1}, and \cite{Arv}), where some of the ideas originated.

\chapter{Free holomorphic functions on noncommutative domains}

\section{Free holomorphic functions and Poisson transforms}
\label{Free}

We introduce the algebra $Hol(\cD_f)$  of all free holomorphic
functions on the noncommutative domain $\cD_f$ and prove a maximum
principle. We identify the domain algebra $\cA_n(\cD_f)$ and the
Hardy algebra $F_n^\infty(\cD_f)$ with subalgebras of free
holomorphic functions on $\cD_f$. When $p$ is a positive regular
noncommutative polynomial, we point out two Banach algebras of free
holomorphic functions $H^\infty(\cD_p^\circ)$ and $A(\cD_p^\circ)$,
which can be identified with the Hardy algebra $F_n^\infty(\cD_p)$
and domain algebra $\cA_n(\cD_p)$, respectively. The elements of
these algebras can be seen as boundary functions for noncommutative
Poisson transforms on $\cD_p^\circ$.

A formal power series $g(Z_1,\ldots, Z_n):=\sum\limits_{\alpha\in
\FF_n^+} c_\alpha Z_\alpha$ is called  free holomorphic function on
the noncommutative domain $\cD_f$ if the series
$$
g(rW_1,\ldots, rW_n):=\sum_{k=0}^\infty \sum_{|\alpha|=k} c_\alpha
r^{|\alpha|} W_\alpha <\infty
$$
is convergent   in the operator norm topology for any  $ r\in
[0,1)$, where $(W_1,\ldots, W_n)$ is the universal model associated
with the noncommutative domain $\cD_f$. We remark that if
$f=X_1+\cdots +X_n$, then our definition is equivalent to the one
for free holomorphic function (see \cite{Po-holomorphic}) on the
unit ball $[B(\cH)^n]_1$

We denote by $Hol(\cD_f)$ the set of all free holomorphic functions
on  $\cD_f$.  Since $g(rW_1,\ldots, rW_n)$ in  the noncommutative
domain algebra $\cA_n(\cD_f)$, one can easily see that $Hol(\cD_f)$
is an algebra.

\begin{theorem}\label{free-ho}
A formal power series $g(Z_1,\ldots, Z_n):=\sum\limits_{\alpha\in
\FF_n^+} c_\alpha Z_\alpha$ is a  free holomorphic function on
$\cD_f$ if and only if
\begin{equation}\label{radius}
\limsup_{k\to \infty} \left(\sum_{|\beta|=k}
|c_\beta|^2\frac{1}{b_\beta}\right)^{1/2k}\leq 1,
\end{equation}
where the coefficients $b_\beta$, $\alpha\in \FF_n^+$,  are given by
relation \eqref{b_alpha}.

In particular, if $\cH$ is infinite dimensional  and
$f=p:=\sum\limits_{1\leq |\alpha|\leq m} a_\alpha X_\alpha$ is a
positive regular free polynomial,  then the above condition is
equivalent to the convergence in the operator norm
 of the series
$$
g(X_1,\ldots, X_n):=\sum_{k=0}^\infty \sum_{|\alpha|=k} c_\alpha
X_\alpha\qquad  \text{for all } \ (X_1,\ldots, X_n)\in
\cD^\circ_p(\cH),
$$
  where
$$
\cD_p^\circ(\cH):=\left\{(X_1,\ldots, X_n)\in B(\cH)^n:\
\left\|\sum_{1\leq |\alpha|\leq m} a_\alpha X_\alpha
X_\alpha^*\right\|<1\right\}.
$$
\end{theorem}
\begin{proof}
Assume that  relation \eqref{radius} holds and let $0<r<\rho<1$.
Then, we can find $m_0\in \NN$ such that
$$
\sum_{|\beta|=k}
|c_\beta|^2\frac{1}{b_\beta}<\left(\frac{1}{\rho}\right)^{2k}\quad
\text{ for any } \ k\geq m_0.
$$
Since $\left\|\sum\limits_{|\beta|=k} b_\beta W_\beta
W_\beta^*\right\|\leq 1$ (see the proof of Theorem
\ref{funct-calc}), we deduce that
\begin{equation*}
\begin{split}
\left\|\sum_{|\beta|=k} c_\beta r^{|\beta|}W_\beta \right\|&\leq
r^k\left(\sum_{|\beta|=k}
|c_\beta|^2\frac{1}{b_\beta}\right)^{1/2}\left\|\sum_{|\beta|=k}
b_\beta W_\beta W_\beta^*\right\|^{1/2}\\
&\leq \left(\frac{r}{\rho}\right)^{k}
\end{split}
\end{equation*}
for any $k=0,1,\ldots$. Hence,
$$
\sum_{k=0}^\infty \left\|\sum_{|\beta|=k} c_\beta r^{|\beta|}W_\beta
\right\|\leq \sum_{k=0}^\infty
\left(\frac{r}{\rho}\right)^{k}<\infty,
$$
which shows that $g(Z_1,\ldots, Z_n)$ is a free holomorphic function
on $\cD_f$. To prove the converse, assume that
$$
\limsup_{k\to \infty} \left(\sum_{|\beta|=k}
|c_\beta|^2\frac{1}{b_\beta}\right)^{1/2k}>\gamma>1.
$$
Then there is $k$ as large as needed such that
$\sum\limits_{|\beta|=k} |c_\beta|^2\frac{1}{b_\beta}>\gamma^{2k}$.
Fix $r$ such that $\frac{1}{\gamma}<r<1$ and notice that, using
relation \eqref{WbWb}, we obtain
\begin{equation*}
\begin{split}
\left\|\sum_{|\beta|=k} c_\beta r^{|\beta|}W_\beta \right\|&\geq
\left\|\sum_{|\beta|=k} c_\beta r^{|\beta|}W_\beta (1)\right\|\\
&=r^k\left(\sum_{|\beta|=k}
|c_\beta|^2\frac{1}{b_\beta}\right)^{1/2} >(r\gamma)^k.
\end{split}
\end{equation*}
Since $r\gamma>1$, the series $\sum\limits_{k=0}^\infty
 \sum\limits_{|\beta|=k} c_\beta r^{|\beta|}W_\beta$ is divergent and
 therefore $g(Z_1,\ldots, Z_n)$ is not a free holomorphic function.

Now, we prove the last part of the theorem. Assume that relation
\eqref{radius} holds and let $(X_1\ldots, X_n)\in \cD_p^\circ(\cH)$.
Since $p$ is a polynomial,  $\cD_p^\circ(\cH)$ is an open set in
$B(\cH)^n$ and consequently there exists $r\in (0,1)$ such that
$$
(X_1\ldots, X_n)\in\cD_{p,r}(\cH):=\left\{(Y_1\ldots, Y_n)\in
B(\cH)^n: \ \left\|\sum_{1\leq |\alpha\leq m}
a_\alpha\frac{1}{r^{|\alpha|}} Y_\alpha Y_\alpha^*\right\|\leq
1\right\}.
$$
Therefore, $(\frac{1}{r} X_1,\ldots,\frac{1}{r}X_n)\in \cD_p(\cH)$
and, due to Theorem \ref{Poisson-C*}, we have
$$
\left\|\sum_{|\alpha|=k}c_\alpha X_\alpha\right\|\leq
\left\|\sum_{|\alpha|=k}c_\alpha r^{|\alpha|}W_\alpha\right\|.
$$
Hence, the series
$\sum\limits_{k=0}^\infty\left\|\sum\limits_{|\alpha|=k}c_\alpha
X_\alpha\right\|$ is convergent. The converse follows taking into
account that $(rW_1,\ldots, rW_n)\in \cD_p^\circ(F^2(H_n))$ for
$0<r<1$, and, since $\cH$ is infinite dimensional,
 the series
$\sum\limits_{k=0}^\infty \sum\limits_{|\alpha|=k} c_\alpha
r^{|\alpha|} W_\alpha $
 converges   in the operator norm topology for any $r\in [0,1)$.
 The proof is complete.
\end{proof}

For $0<r<1$, we define $\cD_{f,r}=\cD_{f,r}(\cH)\subset B(\cH)^n$ by
setting
$$
\cD_{f,r}(\cH):=\left\{ (X_1,\ldots, X_n)\in B(\cH)^n:\
\left(\frac{1}{r} X_1,\ldots,\frac{1}{r}X_n\right)\in
\cD_f(\cH)\right\}.
$$
Notice that $(X_1,\ldots, X_n)\in \cD_f$ if and only if
$(rX_1,\ldots, rX_n)\in \cD_{f,r}$.

\begin{proposition}\label{pro}
Let $g(Z_1,\ldots, Z_n):=\sum\limits_{\alpha\in \FF_n^+} c_\alpha
Z_\alpha$ be  a  free holomorphic function on  $\cD_f$.
\begin{enumerate}
\item[(i)] If $0<r_1<r_2<1$, then $\cD_{f,r_1}\subset
\cD_{f,r_2}\subset\cD_{f}$ and
$$
\|g(r_1W_1,\ldots, r_1 W_n)\|\leq \|g(r_2W_1,\ldots, r_2 W_n)\|.
$$
\item[(ii)] If $0<r<1$, then the map \ $G:\cD_{f,r}(\cH)\to B(\cH)$
defined by
$$G(X_1,\ldots, X_n):=\sum_{k=0}^\infty \sum_{|\alpha|=k}
c_\alpha X_\alpha,\qquad (X_1,\ldots, X_n) \in \cD_{f,r}(\cH),
$$
is continuous and \ $\|G(X_1,\ldots, X_n)\|\leq \|G(rW_1,\ldots,
rW_n)\|$.  Moreover, the series defining $G$  converges uniformly on
$\cD_{f,r}$ in the operator norm topology.
\end{enumerate}
\end{proposition}

\begin{proof}
Since $\varphi(W_1,\ldots, W_n):=\sum\limits_{k=0}^\infty
\sum\limits_{|\alpha|=k} c_\alpha  r_2^{|\alpha|} W_\alpha$ is in
the domain algebra $ \cA_n(\cD_f)$, Theorem \ref{Poisson-C*} implies
$$\|\varphi(rW_1,\ldots, rW_n)\|\leq
\|\varphi(W_1,\ldots, W_n)\|\quad \text{ for any } r\in [0,1).
$$
Taking $r:=\frac{r_1}{r_2}$, we complete the proof of part (i).

To prove (ii), notice that $g_r(W_1,\ldots, W_n)\in \cA_n(\cD_f)$
and $\left(\frac{1}{r} X_1,\ldots,\frac{1}{r}X_n\right)\in
\cD_f(\cH)$. Using  again Theorem \ref{Poisson-C*}, we obtain
$$
\|G(X_1,\ldots, X_n)\|=\left\|g_r\left(\frac{1}{r}
X_1,\ldots,\frac{1}{r}X_n\right)\right\|\leq \|G(rW_1,\ldots,
rW_n)\|
$$
and
$$
\sum_{k=0}^\infty \left\|\sum_{|\alpha|=k}c_\alpha
X_\alpha\right\|\leq \sum_{k=0}^\infty
\left\|\sum_{|\alpha|=k}c_\alpha r^{|\alpha|}W_\alpha\right\|.
$$
Now, one can easily complete the  proof.
\end{proof}

Now we can identify the Hardy algebra $F_n^\infty(\cD_f)$ with a
subalgebra  of free holomorphic function on $\cD_f$.

\begin{theorem}\label{f-infty}
Let $G(Z_1,\ldots, Z_n):=\sum\limits_{\alpha\in \FF_n^+} c_\alpha
Z_\alpha$ be a free holomorphic function on the   operatorial domain
$\cD_f$. Then the following statements are equivalent:
\begin{enumerate}
\item[(i)] $g(W_1,\ldots, W_n):=\sum\limits_{\alpha\in \FF_n^+} c_\alpha W_\alpha$
 is in the noncommutatice Hardy algebra $F_n^\infty(\cD_f)$;
\item[(ii)]$\sup\limits_{0\leq r<1}\|G(rW_1,\ldots, rW_n)\|<\infty$.
\end{enumerate}
In this case,  we have
\begin{equation}\label{many eq}
\|g(W_1,\ldots, W_n)\|=\sup_{0\leq r<1}\|G(rW_1,\ldots, rW_n)\|=
\lim_{r\to 1}\|G(rW_1,\ldots, rW_n)\|.
\end{equation}
\end{theorem}

\begin{proof} Assume (i) holds.
 Since $g(W_1,\ldots, W_n)\in F_n^\infty(\cD_f)$,  Theorem \ref{funct-calc}
 implies
\begin{equation}\label{f-inf}
\|G(rW_1,\ldots, rW_n)\|=\|g(rW_1,\ldots, rW_n)\| \leq
\|g(W_1,\ldots, W_n)\|
\end{equation}
for any $r\in [0,1)$. Therefore, (i)$\implies$(ii). To prove that
(ii)$\implies$(i), assume that  (ii) holds. Consequently,  due to
relation \eqref{WbWb}, we have
\begin{equation*}
\begin{split}
\sum_{\alpha\in \FF_n^+} r^{2|\alpha|} |a_\alpha|^2
\frac{1}{b_\alpha}&=
\left\|\sum_{\alpha\in \FF_n^+} r^{|\alpha|} a_\alpha W_\alpha(1)\right\|\\
&\leq \sup\limits_{0\leq r<1}\|G(rW_1,\ldots, rW_n)\|<\infty
\end{split}
\end{equation*}
for any $0\leq r<1$. Hence, $\sum\limits_{\alpha\in \FF_n^+}
|a_\alpha|^2\frac{1}{b_\alpha}<\infty$, which implies that $
g(W_1,\ldots, W_n)p$ is in $F^2(H_n)$ for any polynomial $p\in
F^2(H_n)$. Now assume that
  $g(W_1,\ldots, W_n)\notin F_n^\infty(\cD_f)$.
 Due to the definition of $F_n^\infty(\cD_f)$, given an arbitrary positive
  number $M$, there exists a polynomial $q\in F^2(H_n)$ with $\|q\|=1$ such that
$$
\|g(W_1,\ldots, W_n)q\|>M.
$$
Since $\|g_r(W_1,\ldots, W_n)(1)-g(W_1,\ldots, W_n)(1)\|\to 0$ as
$r\to 1$, we have
$$\|g(W_1,\ldots, W_n)q-g_r(W_1,\ldots, W_n)q\| \to 0,
\quad \text{ as }\ r\to 1.
$$
Therefore, there is $r_0\in (0,1)$ such that $ \|g_{r_0}(W_1,\ldots,
W_n)  q\|> M. $ Hence,
$$
\|g_{r_0}(W_1,\ldots, W_n)\| >M.
$$
Since $M>0$ is arbitrary, we deduce that
$$
\sup_{0\leq r<1}\|g(rW_1,\ldots, rW_n)\|=\infty,
$$
which is a contradiction. Consequently, (i)$\Longleftrightarrow $
(ii).

We  prove  now the last part of the theorem. If $g(W_1,\ldots,
W_n)\in F_n^\infty(\cD_f)$ and $\epsilon>0$, then there exists a
polynomial $q\in F^2(H_n)$ with $\|q\|=1$ such that
$$
\|g(W_1,\ldots, W_n)q\|>\|g(W_1,\ldots, W_n)\|-\epsilon.
$$
Due to  Theorem \ref{funct-calc}, there is $r_0\in (0,1)$ such that
$$\|g_{r_0}(W_1,\ldots, W_n)q\|>\|g(W_1,\ldots, W_n)\|-\epsilon.
$$
Using now relation \eqref{f-inf}, we deduce that
$$
\sup_{0\leq r<1}\|g(rW_1,\ldots, rW_n)\|=\|g(W_1,\ldots, W_n)\|.
$$
 Due to Proposition \ref{pro},  the   function
 $[0,1)\ni r\to \|g(rW_1,\ldots, rW_n)\|\in \RR^+$
is increasing. Hence, and using  the above  equality, we complete
the proof.
\end{proof}

We can identify the domain algebra $\cA_n(\cD_f)$ with a subalgebra
of $Hol(\cD_f)$.

\begin{theorem}\label{A-infty} Let $G:=\sum\limits_{\alpha\in \FF_n^+}
c_\alpha Z_\alpha$ be a free holomorphic function on the
noncommutative domain $\cD_f$. Then the following statements are
equivalent:
\begin{enumerate}
\item[(i)] $g(W_1,\ldots,W_n):=\sum\limits_{\alpha\in \FF_n^+} c_\alpha W_\alpha$
is in the noncommutative domain algebra $\cA_n(\cD_f)$;
\item[(ii)]$\{G(rW_1,\ldots, rW_n)\}_{0\leq r<1}$ is convergent in
the operator norm as $r\to 1$.
\end{enumerate}
\end{theorem}

\begin{proof} The implication (i)$\implies$ (ii) is due to Theorem
\ref{Poisson-C*}. Conversely, assume item (ii) holds.   Since
$G(rW_1,\ldots, rW_n)\in \cA_n(\cD_f)$, which is an operator
algebra, there exists $\varphi(W_1,\ldots, W_n)\in \cA_n(\cD_f)$
such that
$$
\varphi(W_1,\ldots, W_n)=\lim\limits_{r\to
\infty} G(rW_1,\ldots, rW_n)
$$
 in the operator norm topology. On the other
hand, due to Theorem \ref{f-infty}, we deduce that $g(W_1,\ldots,
W_n):=\sum\limits_{\alpha\in \FF_n} a_\alpha W_\alpha$ is in $
F_n^\infty(\cD_f)$. Since $g(rW_1,\ldots, rW_n)=G(rW_1,\ldots,
rW_n)$  for  $0\leq r<1$, and
$$g(W_1,\ldots, W_n)=\text{\rm
SOT-}\lim_{r\to 1} g(rW_1,\ldots, rW_n), $$ we conclude that
 $\varphi(W_1,\ldots, W_n)=g(W_1,\ldots,
W_n)$,
 which completes the proof.
  \end{proof}

 Define
 $$
 \partial \cD_f:=\left\{ (X_1,\ldots, X_n)\in \cD_f:\
 \left\|\sum_{|\alpha|\geq 1} a_\alpha X_\alpha
 X_\alpha^*\right\|=1\right\}.
 $$
Notice that $(rX_1,\ldots, rX_n)\in \partial \cD_{f,r}$ if and only
if $(X_1,\ldots, X_n)\in \partial\cD_f$.

 Here is our version of the maximum principle \cite{Co} for free holomorphic functions
 on noncommutative domains.

\begin{theorem}\label{max-mod2}
Let $G$ be a free  holomorphic function on $\cD_f(\cH)$, where $\cH$
is an infinite dimensional Hilbert space, and let $r\in [0,1)$. Then
\begin{equation*}
\begin{split}
\max\{\|G(X_1,\ldots, X_n)\|&:\ (X_1,\ldots, X_n)\in {\cD}_{f,r}\}\\
&=
\max\{\|G(X_1,\ldots, X_n)\|:\ (X_1,\ldots, X_n)\in \partial \cD_{f,r}\}\\
&=\|G(rW_1,\ldots, rW_n)\|.
\end{split}
\end{equation*}
\end{theorem}

\begin{proof}
 Due to
Proposition  \ref{pro},  we have
\begin{equation}
\label{ff} \|G(X_1,\ldots, X_n)\|\leq \|G(rW_1,\ldots, rW_n)\| \quad
\text{ for any } \ (X_1,\ldots, X_n)\in  {\cD}_{f,r}.
\end{equation}
 Since $\cH$ is infinite dimensional, there exists a subspace
 $\cK\subset \cH$ and a unitary operator $U:F^2(H_n)\to \cK$.
Define the operators
$$
V_i:=\left(\begin{matrix}
UW_iU^*&0\\
0&0
\end{matrix}\right), \quad i=1,\ldots,n,
$$
with respect to the orthogonal decomposition $\cH=\cK\oplus
\cK^{\perp}$, where $W_1,\ldots, W_n$  are the weighted left
creation operators associated with $\cD_f$. Notice that
$(rV_1,\ldots, rV_n)\in
\partial{\cD}_{f,r}$ and
$$
G(rV_1,\ldots, rV_n)= \left(\begin{matrix}
UG(rW_1,\ldots, rW_n)U^*&0\\
0&0
\end{matrix}\right).
$$
 Consequently,
$$
\|G(rV_1,\ldots, rV_n)\|= \|G(rW_1,\ldots, rW_n)\|.
$$
Hence, and using  inequality \eqref{ff}, we can  complete the proof.
 \end{proof}

Let us consider the particular case when $f$ is equal to a positive
regular  polynomial $p:=\sum\limits_{1\leq |\alpha|\leq m} a_\alpha
X_\alpha$ and $\cH$ is an infinite dimensional Hilbert space.
According to Theorem \ref{free-ho}, the algebra $Hol(\cD_p)$ of all
free holomorphic function on $\cD_p$ can be identified with the set
all functions $G:\cD_p^\circ(\cH)\to B(\cH)$ of the form
$$G(X_1,\ldots,X_n)=\sum_{k=0}^\infty\sum_{|\alpha|=k} c_\alpha
X_\alpha,\quad (X_1,\ldots, X_n)\in \cD_p^\circ(\cH),
$$
where the convergence is in the operator norm topology.
Let $H^\infty(\cD_p^\circ)$  denote the set of  all elements $G$ in
$Hol(\cD_p)$     such that
$$\|G\|_\infty:= \sup \|G(X_1,\ldots, X_n)\|<\infty,
$$
where the supremum is taken over all $n$-tuples $(X_1,\ldots,
X_n)\in \cD_p^\circ(\cH)$.
 We  denote by  $A(\cD_p^\circ)$  be the set of all  elements $G$
  in $Hol(\cD_p)$   such that the
mapping
$$\cD_p^\circ(\cH)\ni (X_1,\ldots, X_n)\mapsto G(X_1,\ldots, X_n)\in B(\cH)$$
 has a continuous extension to  ${\cD}_p$.
Using the results of  this section,  one can easily prove  that
$H^\infty(\cD_p^\circ)$ and $A(\cD_p^\circ)$ are Banach algebras
under pointwise multiplication and the norm $\|\cdot \|_\infty$,
which can be identified with the Hardy algebra $F_n^\infty(\cD_p)$
and the noncommutative domain  algebra $\cA_n(\cD_p)$, respectively.

Given an operator $A\in B(F^2(H_n))$, the noncommutative Poisson
transform  associated with the noncommutative domain $\cD_p$
generates a function
$$
P[A]:\cD_p^\circ(\cH)\to B(\cH)
$$
by setting
$$
P[A](X_1,\ldots, X_n):=\Psi_{f,X}(A)\quad  \text{ for any }\
X:=(X_1,\ldots, X_n)\in \cD_p^\circ(\cH),
$$
where $\Psi_{f,X}(A):=K_{f,X}^* (A\otimes I_\cH) K_{f,X}$ and
$K_{f,T}$ is the Poisson kernel associated with $f$ and $X\in
\cD_p^\circ(\cH)$ (see Section \ref{domain algebra}).

 Using some of the results  from Theorem \ref{Poisson-C*}, Theorem \ref{funct-calc}, Theorem
\ref{f-infty}, and  Theorem \ref{A-infty}, we can provide classes of
operators $A\in B(F^2(H_n))$ such that the mapping $P[A]$ is a free
holomorphic function  on $\cD_p^\circ(\cH)$. In this case, the
operator $A$ can be seen as the boundary function of the Poisson
transform $P[A]$. The following theorem can be easily deduced from
the above-mentioned results.
\begin{theorem}\label{behave} Let $\cH$ be an infinite dimensional
  Hilbert space and  $u$ be
 a free holomorphic function on $\cD_p^\circ(\cH)$, where $p$ is
 a positive regular noncommutative polynomial.
 Then the following statements hold.
\begin{enumerate}
\item[(i)]
There exists $f\in F_n^\infty(\cD_p)$ with $u=P[f]$ if and only if
$$\sup\limits_{0\leq r<1}\|u(rW_1,\ldots, rW_n)\|<\infty.
$$
 In this case, $u(rW_1,\ldots, rW_n)\to f$,  as $ r\to 1$,
  in the $w^*$-topology (or strong operator topology).

\item[(ii)]
There is $f\in \cA_n(\cD_p)$ with $u=P[f]$ if and only if
$\{u(rW_1,\ldots, rW_n)\}_{0\leq r<1}$ is convergent in norm, as
$r\to 1$. In this case, $u(rW_1,\ldots, rW_n)\to f$ in the operator
norm, as $r\to 1$.
\end{enumerate}
\end{theorem}

\bigskip

\section{Schwarz lemma and Bohr's inequality for
$F_n^\infty(\cD_f)$ } \label{Schwartz}

In this section we obtain versions of Schwarz lemma \cite{Co} and
Bohr's inequality \cite{Bo} for the  Hardy algebra
$F_n^\infty(\cD_f)$.

First, we  prove an analogue of Schwarz lemma from complex analysis,
in our multivariable  noncommutative  setting.

\begin{theorem}\label{Schwartz1}
Let $G=G(W_1,\ldots, W_n)=\sum\limits_{\alpha\in \FF_n^+} c_\alpha
W_\alpha$ be  in the noncommutative Hardy algebra
$F_n^\infty(\cD_f)$ with the properties $G(0)=0$, \    $ \|G\|\leq
1, $ and such that
\begin{equation}
\label{M-sup}
 M:=\sup\limits_{{\alpha\in \FF_n^+}\atop{i=1,\ldots,
n}}\left\{\sqrt{ \frac{b_{g_i\alpha}}{b_\alpha}} \right\}<\infty,
\end{equation}
where  the coefficients $b_\alpha$, $\alpha\in \FF_n^+$,  are given
by \eqref{b_alpha}. Then, for any c.n.c. $n$-tuple of operators
$(X_1,\ldots, X_n)\in \cD_f(\cH)$,
 $$\|G(X_1,\ldots, X_n)\| \leq
 {M}\|[X_1,\ldots,
X_n]\|.
$$
\end{theorem}

\begin{proof}
For each $ i=1,\ldots, n$, define the formal power series
$$
\Phi_i(Z_1,\ldots, Z_n):=\sum_{\alpha\in \FF_n^+} a_{g_i\alpha}
Z_\alpha.
$$
Since
$$
\sum_{|\alpha|=k} |a_{ g_i \alpha}|^2\leq \sum_{|\gamma|=k+1} |a_{
\gamma}|^2,\quad i=1,\ldots,n,
$$
we deduce that
$$
\limsup_{k\to\infty} \left( \sum_{|\alpha|=k} |a_{ g_i
\alpha}|^2\right)^{1/2k} \leq \limsup_{k\to\infty} \left(
\sum_{|\gamma|=k+1} |a_{\gamma}|^2 \right)^{\frac{1}{2(k+1)}}.
$$
Consequently, due to Theorem \ref{free-ho},  $\Phi_{i}$ represents a
free holomorphic function on  $\cD_f$. Since $G(0)=0$, we have
$$
G(Z_1,\ldots, Z_n)=\sum_{i=1}^n Z_i\left(\sum_{\alpha\in \FF_n^+}
a_{g_i \alpha} Z_\alpha\right)= \sum_{i=1}^n Z_i\Phi_{i}(Z_1,\ldots,
Z_n)
$$
and
\begin{equation}\label{F-Phi}
G(W_1,\ldots, W_n) =\sum_{i=1}^n     W_i \Phi_{i}(W_1,\ldots, W_n).
\end{equation}

According to Section \ref{Noncommutative}, we have $W_i=S_iD_i$,
$i=1,\ldots,n$, where $S_1,\ldots, S_n$ are the left creation
operators on $F^2(H_n)$ and $D_1,\ldots, D_n$ are the diagonal
operators defined by
$$ D_ie_\alpha:=\sqrt{\frac{b_\alpha}{b_{g_i
\alpha}}} e_\alpha,\quad \alpha\in \FF_n^+.
$$
The condition \eqref{M-sup} shows that the positive diagonal
operator $D_i$ is invertible and $\max\limits_{i=1,\ldots,n}
\|D_i^{-1}\|=M$. Since $W_1,\ldots, W_n$ have orthogonal ranges,
     we deduce that
\begin{equation}\label{M2}
\begin{split}
G(W_1,\ldots, W_n)^*G(W_1,\ldots, W_n) \geq \frac{1}{M^2}
\sum_{i=1}^n \Phi_{i}(W_1,\ldots, W_n)^* \Phi_{i}(W_1,\ldots, W_n).
\end{split}
\end{equation}

Now, due to  Theorem \ref{funct-calc},  we deduce that
\begin{equation}\label{M-P}
\left\|\left[\begin{matrix}
\Phi_{1}(X_1,\ldots, X_n)\\
\vdots\\
\Phi_{n}(X_1,\ldots, X_n)
\end{matrix}
\right]\right\|\leq \left\|\left[\begin{matrix}
\Phi_{1}(W_1,\ldots, W_n)\\
\vdots\\
\Phi_{n}(W_1,\ldots, W_n)
\end{matrix}
\right]\right\|
\end{equation}
for any c.n.c $n$-tuple of operators $(X_1,\ldots, X_n)\in
\cD_f(\cH)$. Consequently, using relations \eqref{F-Phi},
\eqref{M2}, and \eqref{M-P}, we obtain
\begin{equation*}
\begin{split}
\|G(X_1,\ldots, X_n)\| &= \left\|\sum_{i=1}^n  X_i
\Phi_{i}(X_1,\ldots, X_n)\right\| \\
&\leq \left\|[ X_1,\ldots, X_n]\right\| \left\|\left[\begin{matrix}
\Phi_{1}(X_1,\ldots, X_n)\\
\vdots\\
\Phi_{n}(X_1,\ldots, X_n)
\end{matrix}
\right]\right\| \\
&\leq \left\|\sum_{i=1} X_i X_i^*\right\|^{1/2} \left\| \sum_{i=1}^n
\Phi_{i}(W_1,\ldots, W_n)^*\Phi_{i}(W_1,\ldots, W_n)
\right\|^{1/2}\\
&\leq {M}\left\|\sum_{i=1}^n X_i X_i^*\right\|^{1/2}
\left\|G(W_1,\ldots, W_n)^*G(W_1,\ldots, W_n)\right\|\\
&\leq {M}\left\|\sum_{i=1} X_i X_i^*\right\|^{1/2}.
 \end{split}
\end{equation*}
 Therefore,
$$
\|G(X_1,\ldots, X_n)\|\leq  {M}\|[X_1,\ldots, X_n]\|\quad
 \text{ for any }  \ (X_1,\ldots, X_n)\in \cD_f(\cH),
$$
which  completes the proof.
 \end{proof}

Notice that, in particular, if the domain $\cD_p$ is generated  by
the polynomial ~$p=X_1+\cdots +X_n$, then $M=1$.

\bigskip

Let $f(z):=\sum\limits_{k=0}^\infty a_kz^k$ be an  analytic function
on the open unit disc $\DD:=\{z\in\CC: \ |z|<1\}$ such that
$\|f\|_\infty\leq 1$. Bohr's inequality \cite{Bo} asserts  that
$$\sum_{k=0}^\infty r^k |a_k|\leq 1\quad \text{  for }\
 0\leq r\leq \frac{1} {3}.
$$
The fact that $\frac{1} {3}$ is the best possible constant was
obtained independently by M. Riesz, Schur, and Wiener.  In recent
years, multivariable analogues of Bohr's inequality were considered
by several authors (see \cite{Ai}, \cite{BK},
 \cite{DT}, \cite{PPoS}, and \cite{Po-Bohr}).

In what follows we extend Bohr's inequality to the Hardy algebra
$F_n^\infty(\cD_f)$. First, let us recall the  following result
from{ \cite{PPoS}.

\begin{lemma}\label{mat-ine} Let $\cH$ and $\cK$ be Hilbert spaces, let A be
a bounded operator from $\cH$ to $\cK$, and let $z$ be a complex
number. Then
$$
\left[ \begin{matrix} z I_\cH & 0 \\ A & z I_\cK
\end{matrix}\right]
$$
 is a contraction if and only if $\|A\| \leq 1 -
|z|^2$.
\end{lemma}

We can use this result to prove a Wiener type inequality for the
 Hardy  algebra $F_n^\infty(\cD_f)$.

\begin{theorem}\label{Wiener} If
 ~$\varphi(W_1,\dots, W_n)=\sum_{\beta\in \FF_n^+} c_\beta W_\beta$
 is in  the noncommutative Hardy algebra $F_n^\infty(\cD_f)$ and
  $\|\varphi(W_1,\dots, W_n)\|\leq 1$, then
$$
\left(\sum_{|\beta|=k} |c_\beta|^2\frac{1}{b_\beta}\right)^{1/2}
\leq 1-|c_0|^2 \quad
$$
for any $k=1,2,\dots$.
 \end{theorem}
\begin{proof} Let $\cM$ be the subspace spanned by the vectors
 $1$ and $\{e_\alpha\}_{|\alpha|=k}$.
Due to the definition of the weighted left creation operators
$W_1,\ldots, W_n$ ( see \eqref{WbWb}), we have
$$ \left<P_\cM \varphi(W_1,\dots, W_n)
e_\alpha, 1\right>=  \begin{cases} 0 &\quad \text{ if }
|\alpha|=k\\
c_0&\quad \text{ if }  \alpha=g_0 \end{cases} $$
  and
  $$
\left<P_\cM \varphi(W_1,\dots, W_n)1, e_\alpha\right>=
\frac{c_\alpha}{\sqrt{b_\alpha}} \quad  \text{ if } |\alpha|=k.
 $$
 Therefore
 the operator
  $P_\cM \varphi(W_1,\dots, W_n)|_\cM$ is a contraction and
  its matrix with respect to the orthonormal basis
$\{1, \{e_\alpha\}_{|\alpha|=k}\}$ is
$$
\left[\begin{matrix} c_0&0\\
            A& c_0I_{n^k}
\end{matrix}\right],
$$
where $A$ is the column matrix with entries $
\frac{c_\alpha}{\sqrt{b_\alpha}}$, $|\alpha|=k$. Now, applying Lemma
\ref{mat-ine},
 we  complete the proof.
\end{proof}
We recall that, for each $r\in (0,1)$,
$$
\cD_{f,r}(\cH):=\left\{ (X_1,\ldots, X_n)\in B(\cH)^n:\
\left(\frac{1}{r} X_1,\ldots,\frac{1}{r}X_n\right)\in
\cD_f(\cH)\right\}.
$$
Now, we can prove the following.

\begin{theorem}\label{Bohr1} Let
~$\varphi(W_1,\dots, W_n)=\sum_{\beta\in
 \FF_n^+} c_\beta W_\beta$ be  in the noncommutative Hardy algebra
  $F_n^\infty(\cD_f)$  and such that
 $\|\varphi(W_1,\dots,
  W_n)\|\leq 1$. If  $(T_1,\dots, T_n)\in \cD_{f,r}(\cH)$, where
  $0\leq r<1$, then
$$
\left\|\sum_{\beta\in \FF_n^+} |c_\beta| T_{\beta}\right\| \leq
|c_0|+(1-|c_0|^2){r\over 1-r}.
$$
 \end{theorem}
\begin{proof}
 Due to the results of Section \ref{Free}, since
 $\varphi(W_1,\ldots, W_n)\in F_n^\infty(\cD_f)$, we have
 $\sum_{k=0}^\infty\left\|\sum_{|\beta|=k} |c_\beta|r^{|\beta|}
 W_\beta\right\|<\infty$.
Therefore, the series $\sum_{\beta\in \FF_n^+} |c_\beta|
r^{|\beta|}W_{\beta}$ converges in the operator norm for $0\leq r<
1$ and it is in the
  algebra $\cA_n(\cD_f)$. Since
$\left(\frac{1}{r}T_1,\dots, \frac{1}{r}T_n\right)\in \cD_f(\cH)$,
Proposition \ref{pro}, part(ii),   implies
$$
\left\|\sum_{\beta\in \FF_n^+} |c_\beta| T_{\beta}\right\|\leq
\left\|\sum_{\beta\in \FF_n^+} |c_\beta|r^{|\beta|}
W_{\beta}\right\|.
$$
On the other hand, according to Section \ref{Noncommutative}, the
operators $\{W_\beta\}_{|\beta|=k}$ have orthogonal ranges and
$\|W_\beta\|=\frac{1}{\sqrt{b_\beta}}$, $\beta\in \FF_n^+$.
Consequently, $\left\|\sum\limits_{|\beta|=k} b_\beta W_\beta
W_\beta^*\right\|\leq 1$ for any $k=0,1,\ldots$. Hence, and using
the fact that $\sum\limits_{\beta\in \FF_n^+}
|c_\beta|^2\frac{1}{b_\beta}<\infty$, which is due to the fact that
$\varphi(W_1,\ldots, W_n)\in F_n^\infty(\cD_f)$,
 we deduce that
\begin{equation*}
\begin{split}
\left\|\sum_{\beta\in \FF_n^+} |c_\beta|r^{|\beta|}
W_{\beta}\right\|&\leq \sum_{k=0}^\infty r^k\left\|\sum_{|\beta|=k}
|c_\beta| W_\beta\right\|  \\
&\leq\sum_{k=0}^\infty r^k \left(\sum_{|\beta|=k} |c_\beta|^2
\frac{1}{{b_\beta}} \right)^{1/2}\left\|\sum\limits_{|\beta|=k}
b_\beta W_\beta
W_\beta^*\right\|^{1/2}\\
&=\sum_{k=0}^\infty r^k \left(\sum_{|\beta|=k} |c_\beta|^2
\frac{1}{{b_\beta}} \right)^{1/2}.
\end{split}
\end{equation*}
Now, since  $\|f(W_1,\dots, W_n)\|\leq1$,  we can use  the Wiener
type inequality of Theorem \ref{Wiener} and combine it with  the
inequalities above to  complete the proof.
\end{proof}
Our version of Bohr's inequality for the Hardy algebra
$F_n^\infty(\cD_f)$ is the following.

\begin{theorem}\label{Bohr2}
If $\varphi(W_1,\dots, W_n)=\sum_{\beta\in \FF_n^+} c_\beta W_\beta$
is in $F_n^\infty(\cD_f)$,   then
$$
\left\|\sum_{\beta\in \FF_n^+} |c_\beta| T_{\beta}\right\|\leq
\|\varphi(W_1,\dots, W_n)\|
$$
for any $(T_1,\dots, T_n)\in \cD_{f,1/3}(\cH)$.
 In particular, if $\lambda:=(\lambda_1,\dots,
\lambda_n)\in \cD_{f,1/3}(\CC)$, then
$$
\sum_{\beta\in \FF_n^+} |c_\beta| |\lambda_{\beta}|\leq
\|\varphi(W_1,\dots, W_n)\|
$$
for any $\varphi(W_1,\dots, W_n)\in F_n^\infty(\cD_f)$.

Moreover the inequalities above  are  strict unless
$\varphi(W_1,\dots, W_n)$ is a multiple of the identity.
\end{theorem}
\begin{proof} Let  $s\in (0,1/3)$
and assume $\|\varphi(W_1,\dots, W_n)\|=1$. Using  the inequality
from Theorem \ref{Bohr1}, we have

$$\left\|\sum_{\beta\in \FF_n^+} |c_\beta| T_{\beta}\right\|\leq
|c_0|+(1-|c_0|^2) {s\over 1-s}\\
\leq  |c_0|+{1-|c_0|^2\over 2}\leq 1,
$$
for any $(T_1,\dots, T_n)\in \cD_{f,1/3}(\cH)$ and  $s\in (0,1/3)$.
To prove the last part of the theorem, assume that we have equality
in the inequality above and $\|\varphi(W_1,\dots, W_n)\|=1$. Then,
$$|c_0|+{1-|c_0|^2\over 2}=1.
$$
Hence, we deduce that $|c_0|=1$. Since
$$\|\varphi(W_1,\dots, W_n)(1)\|^2=\sum_{\beta\in \FF_n^+}
|c_\beta|^2 \frac{1}{b_\beta}\le1,
$$
we must have $c_\beta=0$ for any $\beta\in \FF_n^+$ with
$|\beta|\geq 1$. Therefore $\varphi(W_1,\dots, W_n)$ is a multiple
of the identity.

%To see that the constant 1/3 is the best possible, it is enough to consider
%functions of one variable only and observe that for such functions
%$\|f(S_1)\| = \|f(z)\|_{\infty}$ where
%the latter supremum is over the unit %disc.
\end{proof}

Whether $1/3$ is the best possible constant in Theorem \ref{Bohr2}
remains an open problem. The answer is affirmative  in the
particular case when $f=X_1+\cdots + X_n$ (see \cite{PPoS}).  We
remark that a commutative version of Bohr's inequality for
holomorphic functions on the Reinhardt domain $\cD_p^\circ(\CC)$
will be presented at the end of Section \ref{Cauchy}.

\bigskip

\section{Weierstrass and Montel theorems for the algebra $Hol(\cD_f)$}
 \label{Weierstrass}

In this section, we obtain Weierstrass and Montel type theorems for
the algebra   of free holomorphic functions on the  noncommutative
domain $\cD_f$. This enables us to introduce a metric on
$Hol(\cD_f)$ with respect to which it becomes a complete metric
space.

The first result of this section
 is a multivariable operatorial version of Weierstrass theorem \cite{Co}.

\begin{theorem}\label{Weierstrass1}
 Let $\{G_m\}_{m=1}^\infty\subset Hol(\cD_f)$ be a
 sequence of free holomorphic functions   such that, for each
 $r\in[0,1)$, the sequence $\{G_m(rW_1,\ldots, rW_n)\}_{m=1}^\infty$
 is convergent in the operator norm topology.
  Then there is a free holomorphic function   $G\in Hol(\cD_f)$
   such that
$G_m(rW_1,\ldots, rW_n)$ converges to $G(rW_1,\ldots, rW_n)$  for
any  $r\in [0,1)$.
\end{theorem}

\begin{proof} Let $G_m:=\sum\limits_{k=0}^\infty
\sum\limits_{|\alpha|=k} a_{\alpha}^{(m)} Z_\alpha$ and fix $r\in
(0,1)$.  Since
$$
G_m(rW_1,\ldots,
rW_n)=\sum\limits_{k=0}^\infty\sum\limits_{|\alpha|=k} r^{|\alpha|}
a_{\alpha}^{(m)} W_\alpha
$$
is in the noncommutative disc algebra $\cA_n(\cD_f)$ and
  the sequence of operators $\{G_m(rW_1,\ldots,
rW_n)\}_{m=1}^\infty$ is convergent in the operator norm of
$B(F^2(H_n))$,   there exists $G\in \cA_n(\cD_f)$ such that
\begin{equation}\label{Fm-to}
G_m(rW_1,\ldots, rW_n)\to G(W_1,\ldots, W_n), \quad \text{ as } \
m\to\infty.
\end{equation}
 Let $\sum\limits_{\alpha\in \FF_n^+}
 d_\alpha(r) W_\alpha$,
  $d_\alpha(r)\in \CC$, $\alpha\in \FF_n^+$, be the Fourier
  representation
  of $G$. Since $W_1,\ldots, W_n$
   have orthogonal ranges and using \eqref{WbWb}, we have
\begin{equation}
\label{d_alpha}
 \frac{d_\alpha(r)}{b_\alpha}=\left< W_\alpha^*
G(W_1,\ldots, W_n)(1),1\right>, \quad \alpha\in \FF_n^+.
\end{equation}
If $\lambda_{(\beta)}\in \CC$ for    $\beta\in \FF_n^+$ with
$|\beta|=k$,  we have
\begin{equation*}
\begin{split}
\Bigl|\Bigl<\sum_{|\beta|=k} &\lambda_{(\beta)}
W_\beta^*(G_m(rW_1,\ldots, rW_n)  -G(W_1,\ldots,
W_n))1,1\Bigr>\Bigr| \\
&\leq \|G_m(rW_1,\ldots, rW_n)-G(W_1,\ldots, W_n)\|
\left\|\sum_{|\beta|=k} \lambda_{(\beta)} W_\beta^*\right\|.
\end{split}
\end{equation*}
 Hence, using \eqref{d_alpha} and the fact that $W_1,\ldots, W_n$
 have orthogonal ranges and $\|W_\beta\|=\frac{1}{\sqrt{b_\beta}}$,
 $\beta\in \FF_n^+$, we have
\begin{equation*}
\begin{split}
&\left|\sum_{|\beta|=k}\frac{1}{b_\beta}
(r^ka_\beta^{(m)}-d_\beta(r))\lambda_{(\beta)}\right|\\
& \qquad\leq \|G_m(rW_1,\ldots, rW_n)-G(W_1,\ldots, W_n)\|
\left(\sum_{|\beta|=k}
\frac{1}{b_\beta}|\lambda_{(\beta)}|^2\right)^{1/2}.
\end{split}
\end{equation*}
for any $\lambda_{(\beta)}\in \CC$ with    $|\beta|=k$.
Consequently,
$$
\left(\sum_{|\beta|=k}\frac{1}{b_\beta}|r^ka_\beta^{(m)}-d_\beta(r)|^2\right)^{1/2}\leq
\|G_m(rW_1,\ldots, rW_n)-G(W_1,\ldots, W_n)\|
$$
for any $k=0,1,\ldots$. Now, since $\|G_m(rW_1,\ldots,
rW_n)-G(W_1,\ldots, W_n)\| \to 0$, as $m\to\infty$, the inequality
above  implies  $r^ka_\beta^{(m)}\to d_\beta(r)$, as $m\to\infty$,
for any $|\beta|=k$ and $k=0,1,\ldots$. Therefore,
$a_\beta:=\lim\limits_{m\to\infty}a_\beta^{(m)}$ exists and
$d_\beta(r)=r^k a_\beta$ for any $\beta\in \FF_n^+$ with $|\beta|=k$
and $k=0,1,\ldots$. Define the  formal power series $F:=
\sum_{\alpha\in \FF_n^+} a_\alpha Z_\alpha$ and let us prove  that
$F$ is a free holomorphic function on  the noncommutative domain
$\cD_f$.  The  calculations above show that, for each $r\in [0,1)$,

\begin{equation*}
\begin{split}
r^k &\left|\left(
\sum_{|\beta|=k}\frac{1}{b_\beta}|a_\beta^{(m)}|^2\right)^{1/2}-
\left(\sum_{|\beta|=k}\frac{1}{b_\beta}|a_\beta|^2\right)^{1/2}\right|\\
&\qquad \qquad\qquad\leq \|G_m(rW_1,\ldots, rW_n)-G(W_1,\ldots,
W_n)\|.
\end{split}
\end{equation*}
This shows that
\begin{equation}
\label{conv-coef} \sum_{|\beta|=k}\frac{1}{b_\beta}|a_\beta^{(m)}|^2
\to \sum_{|\beta|=k}\frac{1}{b_\beta}|a_\beta|^2,\quad \text{ as } \
m\to\infty,
\end{equation}
uniformly with respect to $k=0,1,\ldots$. Assume now  that
$\gamma>1$ and
$$
\limsup_{k\to\infty}
\left(\sum_{|\beta|=k}\frac{1}{b_\beta}|a_\beta|^2\right)^{1/2k}>\gamma.
$$
Then,  there is $k\in \NN$ as large as needed such that
\begin{equation}
\label{sup-ga}
\left(\sum_{|\beta|=k}\frac{1}{b_\beta}|a_\beta|^2\right)^{1/2}>\gamma^k.
\end{equation}
Let $\lambda$  be such that $1<\lambda< \gamma$ and let $\epsilon>0$
be with the property that $\epsilon<\gamma-\lambda$. It is clear
that $\epsilon<\gamma^k-\lambda^k$ for any $k=1,2,\ldots$. Now, the
convergence in  \eqref{conv-coef} implies that there exists
$K_\epsilon\in \NN$ such that
$$
\left|\left(\sum_{|\beta|=k}\frac{1}{b_\beta}|a_\beta^{(m)}|^2\right)^{1/2}
-
\left(\sum_{|\beta|=k}\frac{1}{b_\beta}|a_\beta|^2\right)^{1/2}\right|<
\epsilon
$$
for any  $m>K_\epsilon$ and any $k=0,1,\ldots$. Hence, and  due to
  \eqref{sup-ga}, we have
$$
\left(\sum_{|\beta|=k}\frac{1}{b_\beta}|a_\beta^{(m)}|^2\right)^{1/2}\geq
\gamma^k-\epsilon>\lambda^k
$$
for any  $m>K_\epsilon$ and some $k$ as large  needed.  Therefore,
$$
\limsup_{k\to \infty}
\left(\sum_{|\beta|=k}\frac{1}{b_\beta}|a_\beta^{(m)}|^2\right)^{1/2k}\geq
\lambda>1
$$
for $m\geq K_\epsilon$.   This shows that   $G_m$ is not a free
holomorphic function  on $\cD_f$ (see Theorem \ref{free-ho}), which
is a contradiction. Therefore,
$$
\limsup_{k\to\infty}\left(\sum_{|\beta|=k}\frac{1}{b_\beta}
|a_\beta|^2\right)^{1/2k}\leq 1
$$ and,  due to  Theorem \ref{free-ho}, $F$ is a
free holomorphic function on  $\cD_f$. Consequently, $F(rW_1,\ldots,
rW_n)=\sum\limits_{k=0}^\infty \sum\limits_{|\alpha|=k} r^{|\alpha|}
a_\alpha W_\alpha$ is convergent in the operator norm topology.
Since $G$ and $F(rW_1,\ldots, rG_n)$ have the same Fourier
coefficients, we must have $G=F(rW_1,\ldots, rW_n)$. Due to relation
\eqref{Fm-to}, for each $r\in [0,1)$,  we have
$$
\|G_m(rW_1,\ldots, rW_n)-F(rW_1,\ldots, rW_n)\|\to 0, \quad \text{
as }\ m\to \infty.
$$
  The proof is
complete.
 \end{proof}

We say that a set $\cG\subset Hol(\cD_f)$ is normal if each sequence
$\{G_m\}_{m=1}^\infty$  in $\cG$ has a subsequence
$\{G_{m_k}\}_{k=1}^\infty$ which converges to an element $G\in \cG$,
i.e.,   for any $r\in [0,1)$,
$$
\|G_{m_k}(rW_1,\ldots, rW_n)-G(rW_1,\ldots, rW_n)\|\to 0, \quad
\text{ as }\ k\to \infty.
$$
 The set $\cG$ is called  locally bounded if, for any $r\in[0,1)$,
 there exists $M>0$ such that
$\|f(rW_1,\ldots, rW_n)\|\leq M$ for any $f\in \cG$.

An important consequence of Theorem \ref{Weierstrass} is   the
following noncommutative version of Montel theorem (see \cite{Co}).
Since the proof is similar to that from \cite{Po-holomorphic}, we
shall omit it.

\begin{theorem}\label{Montel}
Let $\cG\subset Hol(\cD_f)$  be a family of free holomorphic
functions. Then the  following statements are equivalent:
\begin{enumerate}
\item[(i)]  $\cF$ is a  normal set.
\item[(ii)] $\cF$ is  locally bounded.
\end{enumerate}
\end{theorem}

  If $\varphi,\psi\in  Hol(\cD_f)$ and
$0<r<1$, we define
$$
d_r(\varphi,\psi):=\|\varphi(rW_1,\ldots, rW_n)-\psi(rW_1,\ldots,
rW_n)\|.
$$
Due to the maximum principle of Theorem \ref{max-mod2}, if  $\cH$ is
an infinite dimensional Hilbert space,  then
$$
d_r(\varphi,\psi)=\sup_{(X_1,\ldots, X_n)\in {\cD}_{f,r}(\cH)}
\|\varphi(X_1,\ldots, X_n)-\psi(X_1,\ldots, X_n)\|.
$$
 Let $0<r_m<1$ be such that $\{r_m\}_{m=1}^\infty$ is an increasing sequence
  convergent  to $1$.
For any $\varphi,\psi\in Hol(\cD_f)$, we define
$$
\rho (\varphi,\psi):=\sum_{m=1}^\infty \left(\frac{1}{2}\right)^m
\frac{d_{r_m}(\varphi,\psi)}{1+d_{r_m}(\varphi,\psi)}.
$$

  Using standards arguments,
    one can    show  that $\rho$ is a metric
on $Hol(\cD_f)$.

\begin{theorem}\label{complete-metric}
$\left(Hol(\cD_f), \rho\right)$  is a complete metric space.
 \end{theorem}

\begin{proof}

First, notice  that  if $\epsilon>0$, then  there exists $\delta>0$
and $m\in \NN$ such that, for any   $\varphi,\psi\in Hol(\cD_f)$,
 \ $d_{r_m}(\varphi, \psi)<\delta\implies \rho(\varphi,\psi)<\epsilon$.
Conversely, if $\delta>0$ and $m\in \NN$ are fixed, then there is
$\epsilon>0$ such that, for  any $\varphi,\psi\in Hol(\cD_f)$,  \ $
\rho(\varphi,\psi)<\epsilon \implies d_{r_m}(\varphi, \psi) <\delta.
$

 Now, let
$\{g_k\}_{k=1}^\infty\subset Hol(\cD_f)$ be a Cauchy sequence in the
metric $\rho$. An immediate consequence of the above observation  is
that
  $\{g_k(r_mW_1,\ldots, r_mW_n)\}_{k=1}^\infty $  is a Cauchy sequence
   in $B(F^2(H_n))$,
   for any
  $m=1,2,\ldots$. Consequently, for each $m=1,2\ldots$,
the sequence $\{g_k(r_mW_1,\ldots, r_mW_n)\}_{k=1}^\infty$ is
convergent in the operator norm of $B(F^2(H_n))$. According to
Theorem \ref{Weierstrass}, there is a free holomorphic function
$g\in Hol(\cD_f)$
   such that
$g_k(rW_1,\ldots, rW_n)$ converges to $g(rW_1,\ldots, rW_n)$  for
any $r\in [0,1)$. Using again the observation made at the beginning
of this proof, we deduce that $\rho(g_k, g)\to 0$, as $k\to\infty$,
which completes the proof.
 \end{proof}

We remark that, Theorem \ref{Montel} implies the following
compactness criterion for subsets of $Hol(\cD_f)$.

\begin{corollary}
A subset ~$\cG$ of  $(Hol(\cD_f), \rho)$ is compact if and only if
it is closed and locally bounded.
\end{corollary}

\bigskip

\section{Cauchy transforms and analytic functional
calculus for $n$-tuples of operators} \label{Cauchy}

In this section,  we define a   Cauchy transform $\cC_{f,T}$
 associated with any positive
regular free holomorphic function  $f$ on $B(\cH)^n$ and any
$n$-tuple of operators $T:=(T_1,\ldots, T_n)\in B(\cH)^n$
 with joint spectral radius
$r_f(T_1,\ldots, T_n)<1$.

This is used to provide   a  free analytic functional calculus for
such  $n$-tuples of operators.
       In particular,  if $p$ is a positive regular  noncommutative
       polynomial, we  prove  that
$$
g(T_1,\ldots, T_n)=\cC_{p,T}(g(W_1,\ldots, W_n)),\quad g(W_1,\ldots,
W_n)\in F_n^\infty (\cD_p),
$$
where  $g(T_1,\ldots, T_n)$ is  defined by the free analytic
functional calculus associated with the domain $\cD_p$. In the last
part of this section we obtain an analytic functional calculus and
Bohr type inequalities  in the  commutative case.

Let $f=\sum_{\alpha\in \FF_n^+} a_\alpha X_\alpha$, $\alpha\in \CC$,
be a positive regular free holomorphic on $B(\cH)^n$. For any
$n$-tuple of operators $T:=(T_1,\ldots, T_n)\in B(\cH)^n$
 such  that
 $ \left\|\sum\limits_{|\alpha|\geq 1}
a_{\alpha} T_{\alpha}T_{\alpha}^*\right\|<\infty$, we define the
joint spectral radius of $T$  with respect to the noncommutative
domain $\cD_f$ to be
$$
r_f(T_1,\ldots, T_n):=\lim_{m\to\infty}\|\Phi_{f,T}^m(I)\|^{1/2m},
$$
where the positive linear map $\Phi_{f,T}:B(\cH)\to B(\cH)$ is given
by
$$
\Phi_{f,T}(X):=\sum_{|\alpha|\geq 1} a_\alpha T_\alpha
XT_\alpha^*,\quad X\in B(\cK),
$$
and  the convergence is in the week operator topology.
In the particular case when $f:=X_1+\cdots +X_n$, we obtain the
usual definition of the joint operator radius for $n$-tuples of
operators.

Since $\sum\limits_{|\alpha|\geq 1} a_{\tilde \alpha} \Lambda_\alpha
\Lambda_\alpha^*$ is SOT convergent (see Section \ref
{Noncommutative}), one can easily see that the  series
$\sum\limits_{|\alpha|\geq 1} a_{\tilde\alpha} \Lambda_\alpha
\otimes T_{\tilde\alpha}^*$ is SOT-convergent  in $B(F^2(H_n)\otimes
\cH)$ and
$$
\left\|\sum_{|\alpha|\geq 1} a_{\tilde\alpha} \Lambda_\alpha \otimes
T_{\tilde\alpha}^*\right\|\leq \left\|\sum_{|\alpha|\geq 1}
a_{\tilde\alpha} \Lambda_\alpha
\Lambda_\alpha^*\right\|\left\|\sum_{|\alpha|\geq 1}
a_{\tilde\alpha} T_{\tilde\alpha}T_{\tilde\alpha}^*\right\|.
$$
We call the operator
$$ \sum_{|\alpha|\geq 1} a_{\tilde\alpha} \Lambda_\alpha \otimes
T_{\tilde\alpha}^*$$
 the {\it  reconstruction operator} associated with the $n$-tuple
  $T:=(T_1,\ldots, T_n)$ and the noncommutative domain $\cD_f$.
 Notice also that
\begin{equation}
\label{RRR}
 \left\|\left(\sum_{|\alpha|\geq 1} a_{\tilde\alpha}
\Lambda_\alpha \otimes T_{\tilde\alpha}^*\right)^m\right\|\leq
\left\|\Phi^m_{\tilde f,\Lambda}(I)\right\|^{1/2}
\left\|\Phi_{f,T}^m(I)\right\|^{1/2},\quad m\in \NN,
\end{equation}
where $\tilde f:=\sum\limits_{|\alpha|\geq 1} a_{\tilde
\alpha}X_\alpha$ and $\Phi_{\tilde
f,\Lambda}(Y):=\sum\limits_{|\alpha|\geq 1}
a_{\tilde\alpha}\Lambda_\alpha Y \Lambda_\alpha^*$. Hence, we deduce
that
\begin{equation}\label{RARARA}
r\left(\sum_{|\alpha|\geq 1} a_{\tilde\alpha} \Lambda_\alpha \otimes
T_{\tilde\alpha}^*\right)\leq r_{\tilde f}(\Lambda_1,\ldots,
\Lambda_n)r_f(T_1,\ldots, T_n),
\end{equation}
where $r(A)$ denotes the usual spectral radius of an operator $A$.
Due to the results of Section \ref{Noncommutative}, we have
$\left\|\Phi_{\tilde f,\Lambda}(I)\right\|\leq 1$, which implies
$r_{\tilde f}(\Lambda_1,\ldots, \Lambda_n)\leq 1$. Consequently, we
have
\begin{equation}
\label{r<r} r\left(\sum_{|\alpha|\geq 1} a_{\tilde\alpha}
\Lambda_\alpha \otimes T_{\tilde\alpha}^*\right)\leq r_f(T_1,\ldots,
T_n).
\end{equation}

Assume now that $T:=(T_1,\ldots, T_n)\in B(\cH)^n$ is an $n$-tuple
of operators with  $r_f(T_1,\ldots, T_n)<1$.
 We introduce the {\it Cauchy kernel} associated  with   $T $  to be
  the operator
\begin{equation}
\label{Cauc-inv} C_{f,T}(\Lambda_1,\ldots, \Lambda_n):=\left(
I-\sum_{|\alpha|\geq 1} a_{\tilde\alpha} \Lambda_\alpha \otimes
T_{\tilde\alpha}^* \right)^{-1},
\end{equation}
which is well-defined due to   relation \eqref{r<r}.
Notice that
$$\tilde f(Y_1,\ldots, Y_n)=\sum_{|\alpha|\geq 1} a_{\tilde\alpha} \Lambda_\alpha \otimes
T_{\tilde\alpha}^*
$$
where $\tilde f(Y_1,\ldots, Y_n):=\sum\limits_{|\alpha|\geq 1}
a_{\tilde\alpha} Y_\alpha$ and $Y_i:=\Lambda_i\otimes T_i^*$,
$i=1,\ldots, n$. Due to inequality \eqref{r<r}, we have
 $
\sum\limits_{k=0}^\infty \|\tilde f(Y_1,\ldots, Y_n)^k\|<\infty$
and, therefore,
\begin{equation}
\label{Cf-tilde}
 C_{f,T}(\Lambda_1,\ldots,
\Lambda_n)=\sum_{k=0}^\infty \tilde f(Y_1,\ldots, Y_n)^k,
\end{equation}
where the convergence is in the operator norm topology.

We remark that   $C_{f,T}(\Lambda_1,\ldots, \Lambda_n)$ is in
   $R_n^\infty(\cD_f)\bar\otimes B(\cH)$, the $WOT$-closed operator
    algebra generated by the spatial tensor product.
    Moreover, due to the results of Section \ref{Noncommutative} and
    Section \ref{Hardy},
    its  Fourier representation
    is
    \begin{equation}
    \label{Fourier-Cauc}
    C_{f,T}(\Lambda_1,\ldots, \Lambda_n)=\sum_{\beta\in \FF_n^+} (\Lambda_\beta\otimes b_{\tilde
    \beta}T_{\tilde \beta}^*),
    \end{equation}
where the coefficients $b_\alpha$, $\alpha\in \FF_n^+$ are given by
relation   \eqref{b_alpha}. In the particular case when
    $f$ is a polynomial,   the Cauchy kernel is
    in $\cR_n(\cD_f)\bar\otimes B(\cH)$.

In what follows we also use the notation
$C_{f,T}:=C_{f,T}(\Lambda_1,\ldots, \Lambda_n)$. Here are a few
properties of the noncommutative Cauchy kernel.

\begin{proposition}\label{Prop-Cauc}
Let $T:=(T_1,\ldots, T_n)\in B(\cH)^n$ be an   $n$-tuple of
operators with joint spectral radius $r_f(T_1,\ldots, T_n)<1$. Then
the following statements hold:
\begin{enumerate}
\item[(i)]
$\|C_{f,T}\|\leq \sum\limits_{k=0}^\infty \left\|
\Phi_{f,T}^k(I)\right\|^{1/2},$
\item[(ii)]
$C_{f,T}-C_{f,X}=C_{f,T}\left[\sum\limits_{|\alpha|\geq 1} a_{\tilde
\alpha} \Lambda_\alpha\otimes (T_{\tilde \alpha}^*-X_{\tilde
\alpha}^*)\right] C_{f,X},$ and
$$
\|C_{f,T}-C_{f,X}\|\leq \|C_{f,T}\| \|C_{f,X}\|\left\|
\sum_{|\alpha|\geq 1} a_\alpha
(T_\alpha-X_\alpha)(T_\alpha-X_\alpha)^*\right\|^{1/2}
$$
for any $n$-tuple  of operators $X:=(X_1,\ldots, X_n)\in B(\cH)^n$
with joint spectral radius $r_f(X_1,\ldots, X_n)<1$.
\end{enumerate}
 \end{proposition}

\begin{proof} Since $\|\Phi_{\tilde f,\Lambda}(I)\|\leq 1$, relation
\eqref{RRR} implies
$$
\|C_{f,T}\|\leq \sum_{k=0}^\infty\left\|  \tilde f(Z_1,\ldots,
Z_n)^k\right\| = \sum_{k=0}^\infty\left\|
\Phi_{f,T}^k(I)\right\|^{1/2}.
$$
To prove (ii), notice that
\begin{equation*}
\begin{split}
&C_{f,T}-C_{f,X}\\
&\quad= \left(I- \sum\limits_{|\alpha|\geq 1} a_{\tilde \alpha}
\Lambda_\alpha\otimes T_{\tilde \alpha}^*\right)^{-1} \left[ I-
\sum\limits_{|\alpha|\geq 1} a_{\tilde \alpha} \Lambda_\alpha\otimes
X_{\tilde \alpha}^*-\left(I- \sum\limits_{|\alpha|\geq 1} a_{\tilde
\alpha}
\Lambda_\alpha\otimes T_{\tilde \alpha}^*\right)\right]\\
&\qquad \qquad\qquad\left(I- \sum\limits_{|\alpha|\geq 1} a_{\tilde
\alpha}
\Lambda_\alpha\otimes X_{\tilde \alpha}^*\right)^{-1}\\
&\quad=C_{f,T}\left[ \sum\limits_{|\alpha|\geq 1} a_{\tilde \alpha}
\Lambda_\alpha\otimes (T_{\tilde \alpha}^*-X_{\tilde
\alpha}^*)\right] C_{f,X},
\end{split}
\end{equation*}
and

\begin{equation*}
\begin{split}
\|C_{f,T}-C_{f,X}\|&\leq \|C_{f,T}\| \|C_{f,X}\|\left\|
 \sum\limits_{|\alpha|\geq 1} a_{\tilde \alpha}
\Lambda_\alpha\otimes (T_{\tilde \alpha}^*-X_{\tilde
\alpha}^*)\right\|\\
&\leq\|C_{f,T}\| \|C_{f,X}\|\left\| \sum\limits_{|\alpha|\geq 1} a_{
\alpha} (T_{\alpha}-X_{\alpha}) (T_{ \alpha}^*-X_{ \alpha}^*)
\right\|^{1/2},
\end{split}
\end{equation*}
which completes the proof.
\end{proof}

We remark that, in particular, if the $n$-tuple $T:=(T_1,\ldots,
T_n)\in \cD_f^\circ(\cH)$, then $\|C_{f,T}\|\leq
\frac{1}{1-\|\Phi_{f,T}(I)\|^{1/2}}$. Indeed, in this case,  we have
$\|\Phi_{f,T}(I)\|<1$ and
$$\sum_{k=0}^\infty\left\|  \Phi_{f,T}^k(I)\right\|^{1/2} \leq \sum_{k=0}^\infty\left\|
 \Phi_{f,T}(I)\right\|^{k/2} =\frac{1}{1-\|\Phi_{f,T}(I)\|^{1/2}}.
$$
Now, the assertion follows from part (i) of Proposition
\ref{Prop-Cauc}.

 %\begin{proof} Using the classical Cauchy inequality, we %obtain
%\begin{equation*}
%\begin{split}
%\left<\left(\sum_{i=1}^k A_i\right)\left(\sum_{i=1}^k %A_i^*\right)
%x,x\right>&\leq \left(\sum_{i=1}^k \|A_i^*x\|\right)^2
%\leq k \sum_{i=1}^k \|A_i^*x\|^2\\
%&=\left<k\sum_{i=1}^k A_iA_i^*x,x\right>
%\end{split}
%\end{equation*}
%for any $x\in \cH$.
%\end{proof}

 Given an $n$-tuple of operators  $T:=(T_1,\ldots, T_n)\in B(\cH)^n$
  with joint spectral radius $r_f(T_1,\ldots, T_n)<1$,
  we define the {\it Cauchy transform} at
  $T$ to be the mapping
$$
\cC_{f,T}:B(F^2(H_n))\to B(\cH)
$$
defined by
$$
\left< \cC_{f,T}(A)x,y\right>:= \left<(A\otimes I_\cH)(1\otimes x),
C_{f,T}(\Lambda_1,\ldots, \Lambda_n)(1\otimes y)\right>
$$
for any $x,y\in \cH$, where $\Lambda_1,\ldots, \Lambda_n$ are the
weighted right creation operators associated with the noncommutative
domain $\cD_f$. The operator $\cC_{f,T}(A)$ is called the Cauchy
transform of $A$ at $T$.

Now, we provide a {\it free analytic functional calculus} for
$n$-tuples of operators with joint spectral radius $r_p(X_1,\ldots,
X_n)<1$, which now turns out to be continuous and unique.

\begin{theorem}\label{abel} Let $p$ be a positive regular
noncommutative polynomial and  let $T:=(T_1,\ldots, T_n)\in
B(\cH)^n$ be an $n$-tuple of operators with joint spectral radius
$r_p(T_1,\ldots,T_n)<1$.
  If \ $g:=\sum\limits_{\alpha\in \FF_n^+}
c_\alpha Z_\alpha$ is  a free holomorphic function on the
noncommutative  domain $\cD_p(\cH)$, then the series
$$
g(T_1,\ldots, T_n):=\sum\limits_{k=0}^\infty
\sum\limits_{|\alpha|=k} c_\alpha T_\alpha
$$ is convergent in the operator norm  of $B(\cH)$
and the mapping
 $$\Psi_{p,T}: Hol(\cD_p) \to B(\cH)\quad \text{ defined by } \quad
\Psi_{p,T}(g):=g(T_1,\ldots, T_n)
$$
is a continuous unital algebra homomorphism. Moreover, the free
analytic functional calculus  is uniquely determined by the mapping
$$
Z_i\mapsto T_i,\qquad i=1,\ldots,n.
$$
\end{theorem}

\begin{proof}  Let $T:=(T_1,\ldots, T_n)\in B(\cH)^n$ be such that
$r_p(T_1, \ldots, T_n)<1$.
 Using  relations \eqref{Fourier-Cauc}, \eqref{WbWb}, \eqref{WbWb-r},
    and the
definition of the reconstruction operator, we obtain
\begin{equation*}
\begin{split}
\left<\cC_{p,T}(W_\alpha)x,y\right>&= \left< (W_\alpha\otimes
I_\cH)(1\otimes x), C_{p,T}(\Lambda_1,\ldots, \Lambda_n)
(1\otimes y)\right>\\
 &=\left< \frac{1}{\sqrt b_\alpha}e_\alpha\otimes x, \sum_{\beta\in
\FF_n^+}
b_{\tilde \beta}(\Lambda_\beta\otimes T_{\tilde \beta}^*)(1\otimes y)\right>\\
&= \left< \frac{1}{\sqrt b_\alpha}e_\alpha\otimes x, \sum_{\beta\in
\FF_n^+} \sqrt {b_{\tilde \beta}}e_{\tilde \beta}\otimes
T_{\tilde \beta}^*y\right>\\
&=\left<T_\alpha x,y\right>
\end{split}
\end{equation*}
for any $ x,y\in \cH$. Hence we deduce that, for any polynomial
$q(W_1,\ldots, W_n)$,

$$
\left< q(T_1,\ldots, T_n)x,y\right>=\left< (q(W_1,\ldots,
W_n)\otimes I_\cH)(1\otimes x), C_{p,T}(\Lambda_1,\ldots,
\Lambda_n)(1\otimes y)\right>
$$
and, therefore,
\begin{equation}\label{ine-bound}
\|q(T_1,\ldots, T_n)\|\leq \|q(W_1,\ldots,
W_n)\|\|C_{p,T}(\Lambda_1,\ldots, \Lambda_n)\|
\end{equation}
for any polynomial $q(W_1,\ldots, W_n)$. Since $g
:=\sum\limits_{\alpha\in \FF_n^+} c_\alpha   Z_\alpha $ is a free
holomorphic function on $\cD_p$,  the series $g_r(W_1,\ldots,
W_n):=\sum\limits_{k=0}^\infty \sum\limits_{|\alpha|=k}c_\alpha
r^{|\alpha|} W_\alpha$ converges in the operator norm topology. Now,
using \eqref{ine-bound}, we deduce that $g_r(T_1,\ldots,
T_n):=\sum\limits_{k=0}^\infty \sum\limits_{|\alpha|=k}c_\alpha
r^{|\alpha|} T_\alpha$ converges in the operator norm topology of
$B(\cH)$,
\begin{equation}
\label{gr-ine}
 \|g_r(T_1,\ldots, T_n)\|\leq \|g_r(W_1,\ldots,
W_n)\|\|C_{p,T}(\Lambda_1,\ldots, \Lambda_n)\|,
\end{equation}
and
\begin{equation}
\label{gr-cauc} \begin{split}&\left< g_r(T_1,\ldots,
T_n)x,y\right>\\
&\quad\quad=\left< (g_r(W_1,\ldots, W_n)\otimes I_\cH)(1\otimes x),
C_{p,T}(\Lambda_1,\ldots, \Lambda_n)(1\otimes y)\right>
\end{split}
\end{equation}
for any $x,y\in \cH$.
 Now, we need to prove that if $r_p(T_1,\ldots,
T_n)<1$, then there exists a constant $t>1$ such that
$r_p(tT_1,\ldots, tT_n)<1$. Indeed, notice first that if
$p(T_1,\ldots, T_n)=\sum\limits_{1\leq |\alpha|\leq m} a_\alpha
T_\alpha$, then
 $$r_p(T_1,\ldots, T_n)=r\left(\sum_{1\leq |\alpha|\leq m}V_{(\alpha)}\otimes
\sqrt{a_\alpha} T_\alpha^*\right),
$$
where $\{V_{(\alpha)}\}_{1\leq |\alpha|\leq m}$ is any sequence of
isometries with orthogonal ranges. Since the spectrum of an operator
is upper semi-continuous, so is the spectral  radius. Consequently,
for any $\delta>0$ there exists $\epsilon>0$ such that if
$$\|[T_1-Y_1,\ldots, T_n-Y_n]\|<\epsilon,\quad  \text{ then } \quad
r_p(Y_1,\ldots, Y_n)<r_p(T_1,\ldots, T_n)+\delta.
$$
Since $r_p(T_1,\ldots, T_n)<1$, we deduce that there exists a
constant $t>1$ such that $r_p(tT_1,\ldots, tT_n)<1$. Applying
relation  \eqref{gr-ine} and \eqref{gr-cauc} when $r:=\frac{1}{t}$
and $T_i$ is replaced by $tT_i$, $i=1,\ldots, n$,  we deduce that
$$
g(T_1,\ldots, T_n):=\sum\limits_{k=0}^\infty
\sum\limits_{|\alpha|=k} c_\alpha T_\alpha
$$ is convergent in the operator norm  of $B( \cH)$
and
\begin{equation}
\label{gr2-cauc}\begin{split} &\left< g(T_1,\ldots,
T_n)x,y\right>\\
&\qquad =\left< (g_r(W_1,\ldots, W_n)\otimes I_\cH)(1\otimes x),
C_{p,\frac{1}{r}X}(\Lambda_1,\ldots, \Lambda_n)(1\otimes y)\right>
\end{split}
\end{equation}
for any $x,y\in \cH$.  Hence , we have \begin{equation} \label{ggc}
\|g(T_1,\ldots, T_n)\|\leq \|g_r(W_1,\ldots,
W_n)\|\|C_{p,\frac{1}{r}T}(\Lambda_1,\ldots, \Lambda_n)\|.
\end{equation}

 Now, to
prove the continuity of $\Psi_{p,T}$, let $g_m$ and $g$ be in
$Hol(\cD_p)$ such that $g_m\to g$, as $m\to \infty$,  in the metric
$\rho$ of $Hol(\cD_p)$.  This is equivalent to the fact that, for
each $r\in [0,1)$,
\begin{equation}\label{conv-S}
g_m(rW_1,\ldots, rW_n)\to g(rW_1,\ldots, rW_n),\quad \text{ as }\
m\to\infty,
\end{equation}
where the convergence is  in the operator norm of $B(F^2(H_n))$.
According to \eqref{ggc}, we have
 \begin{equation*}
 \begin{split}
&\|g_m(T_1,\ldots, T_n)-g(T_1,\ldots, T_n)\|\\
&\qquad \leq \|g_m(rW_1,\ldots, rW_n)-g(rW_1,\ldots,
rW_n)\|\|C_{p,\frac{1}{r}T}(\Lambda_1,\ldots, \Lambda_n)\|.
\end{split}
\end{equation*}
Now, using \eqref{conv-S}, we deduce that
\begin{equation}
\label{conv-f_m} \|g_m(T_1,\ldots, T_n)-g(T_1,\ldots, T_n)\|\to 0,
\quad \text{ as }\ m\to\infty.
\end{equation}
which proves the continuity of $\Psi_{p,T}$.

To prove the uniqueness of the free analytic functional calculus,
assume that $\Phi:Hol(\cD_p)\to B(\cH)$ is a continuous unital
algebra homomorphism such that $\Phi(Z_i)=T_i$, \ $i=1,\ldots, n$.
Hence, we deduce that
\begin{equation}\label{pol2}
\Psi_{p,T}(q(Z_1,\ldots, Z_n))=\Phi(q(Z_1,\ldots, Z_n))
\end{equation}
for any polynomial $q(Z_1,\ldots, Z_n)$ in $Hol(\cD_p)$.   Let
$g=\sum_{k=0}^\infty \sum_{|\alpha|=k} a_\alpha Z_\alpha$ be an
element in $Hol(\cD_p)$ and let $q_m:=\sum_{k=0}^m\sum_{|\alpha|=k}
a_\alpha Z_\alpha$, \ $m=1,2,\ldots$. Since
$$
g(rW_1,\ldots, rW_n)=\sum_{k=0}^\infty \sum_{|\alpha|=k} r^k
a_\alpha W_\alpha$$ and the series $\sum_{k=0}^\infty
r^k\left\|\sum_{|\alpha|=k} a_\alpha W_\alpha\right\|$ converges, we
deduce that
$$
q_m(rW_1,\ldots, rW_n)\to g(rW_1,\ldots, rW_n)
$$
in the operator norm, as $m\to\infty$. Therefore, $q_m\to g$ in the
metric $\rho$ of $Hol(\cD_p)$. Hence, using \eqref{pol2} and the
continuity of $\Phi$ and $\Psi_{p,T}$, we deduce that
$\Phi=\Psi_{p,T}$.
 This completes the proof.
 \end{proof}

\begin{corollary}\label{F-infty-cauc} If $(T_1,\ldots,
T_n)\in B(\cH)^n$  is an $n$-tuple of operators with
$r_p(T_1,\ldots,T_n)<1$ and $g(W_1,\ldots, W_n)\in
F_n^\infty(\cD_p)$, then
$$
g(T_1,\ldots, T_n)=\cC_{p,T}(g(W_1,\ldots, W_n)),
$$
where $g(T_1,\ldots, T_n)$ is defined by the free analytic
functional calculus.
\end{corollary}

\begin{proof}
Assume that $g:=\sum\limits_{k=0}^\infty\sum\limits_{|\alpha|=k}
a_\alpha W_\alpha$ is in $F_n^\infty(\cD_p)$ and   $0<r<1$. Then,
 we have
$$
\lim_{m\to\infty}\sum_{k=0}^m r^k \sum_{|\alpha|=k} a_\alpha
W_\alpha=g_r(W_1,\ldots, W_n)\in \cA_n(\cD_p)
$$
in the operator norm of $B(F^2(H_n))$, and
$$
\lim_{m\to\infty}\sum_{k=0}^m r^k \sum_{|\alpha|=k} a_\alpha
T_\alpha=g_r(T_1,\ldots, T_n)
$$
in the operator norm of $B(\cH)$. Now,  due to the continuity of the
noncommutative Cauchy transform in the operator norm, we deduce that
\begin{equation}\label{f_r-C}
g_r(T_1,\ldots, T_n)=\cC_{p,T}(g_r(W_1,\ldots, W_n)).
\end{equation}
Since $g(W_1,\ldots, W_n)\in F_n^\infty(\cD_p)$, we know that
$\lim\limits_{r\to 1} g_r(W_1,\ldots, W_n)=g(W_1,\ldots, W_n)$ in
the strong operator topology. Since $\|g_r(W_1,\ldots, W_n)\|\leq
 \|g(W_1,\ldots, W_n)\|$, we deduce that
$$\text{\rm SOT}-\lim\limits_{r\to 1}[g_r(W_1,\ldots, W_n)\otimes
I_\cH]=g(W_1,\ldots, W_n)\otimes I_\cH.
$$
On the other hand,   notice that  $\lim\limits_{t\to 1}
g_t(T_1,\ldots, T_n)=g(T_1,\ldots, T_n)$ in the operator norm.
Indeed, using relation \eqref{gr2-cauc}, we deduce that
\begin{equation*}
\begin{split}
 &\|g(T_1,\ldots,
T_n)-g_t(T_1,\ldots, T_n)\|\\
&\qquad\leq \|g_r(W_1,\ldots, W_n)-g_r(tW_1,\ldots,
tW_n)\|\|C_{p,\frac{1}{r}T}(\Lambda_1,\ldots, \Lambda_n)\|.
\end{split}
\end{equation*}
Since $\|g_r(W_1,\ldots, W_n)-g_r(tW_1,\ldots, tW_n)\|\to 0$, as
$t\to 1$, we prove our assertion.

 Now, passing to the limit, as $r\to 1$, in the equality
\begin{equation*}
\begin{split}
&\left<g_r(T_1,\ldots, T_n)x,y\right>\\
&\qquad =\left<(g_r(W_1,\ldots, W_n)\otimes I_\cH)(1\otimes x),
C_{p,T}(\Lambda_1,\ldots, \Lambda_n)(1\otimes y)\right>,
\end{split}
\end{equation*}
where $x,y\in \cH$, we obtain $g(T_1,\ldots,
T_n)=\cC_{p,X}(g(W_1,\ldots, W_n))$, which
  completes the proof.
   \end{proof}

Using the Cauchy representation provided   by Corollary
\ref{F-infty-cauc}, one can deduce the following result.

\begin{corollary}\label{conv-u-w*}
Let $(T_1,\ldots, T_n)\in B(\cH)^n$ be an $n$-tuple of operators
with $r_p(T_1,\ldots, T_n)<1$.
\begin{enumerate}
\item[(i)] If $\{g_k\}_{k=1}^\infty$ and $g$ are free holomorphic
functions  in $Hol(\cD_p)$ such that $\|g_k-g\|_\infty\to 0$, as
$k\to \infty$, then $g_k(T_1,\ldots, T_n)\to g(T_1,\ldots, T_n)$ in
the operator norm of $B(\cH)$.
\item[(ii)] If $\{g_k\}_{k=1}^\infty$ and $g$ are   in  the noncommutative Hardy  algebra
$F_n^\infty (\cD_f)$ such that
$$g_k(W_1,\ldots, W_n)\to g(W_1,\ldots,
W_n)$$ in the $w^*$-topology (or strong operator topology)  and
$\|g_k(W_1,\ldots, W_n)\|\leq M$ for any $k=1,2,\ldots$, then
$g_k(T_1,\ldots, T_n)\to g(T_1,\ldots, T_n)$ in the weak operator
topology.
\end{enumerate}
\end{corollary}

Using Theorem \ref{f-infty}, Theorem \ref{abel}, and the results
from   Section \ref{Functional}, one can deduce   the following.

\begin{proposition}
For  any  $n$-tuple of operators  $(T_1,\ldots, T_n)\in
\cD_p^\circ(\cH)$, the free analytic functional calculus  coincides
with the $F_n^\infty(\cD_p)$-functional calculus.
\end{proposition}

Now, we present a connection between free holomorphic functions on
noncommutative domains and holomorphic functions on Reinhardt
domains.

\begin{proposition}\label{Reinhardt} Let $p$ be a positive regular
noncommutative polynomial. If  \ $G(Z_1,\ldots,
Z_m):=\sum_{\alpha\in \FF_n^+} c_\alpha Z_\alpha$ is a free
holomorphic function on the noncommutative  domain $\cD_p$, then
~$g(\lambda_1,\ldots, \lambda_n):=\sum_{\alpha\in \FF_n^+} c_\alpha
\lambda_\alpha$ is a holomorphic function on the Reinhardt domain
$$
\cD_p^\circ(\CC):=\left\{ (\lambda_1,\ldots, \lambda_n)\in \CC^n: \
\sum_{1\leq |\alpha|\leq m} a_\alpha|\lambda_\alpha|^2<1\right\}.
$$
 Moreover,   the map $\Phi:
Hol(\cD_p)\to Hol(\cD_p^\circ(\CC))$ defined by $\Phi(G):=g$ is
continuous.
\end{proposition}
\begin{proof}
For each $\mu:=(\mu_1,\ldots, \mu_n)\in \cD_p^\circ(\CC)$, define
$z_\mu:=\sum_{\alpha\in \FF_n^+} \sqrt{b_\alpha}
\overline{\mu}_\alpha e_\alpha$. Since $(\mu_1,\ldots, \mu_n)$ is of
class $C_{\cdot 0}$ with respect to $\cD_p$, relation \eqref{I-Q}
implies
$$
\left(1-\sum_{|\alpha|\geq 1} a_\alpha
|\mu_\alpha|^2\right)\left(\sum_{\beta\in \FF_n^+} b_\beta
|\mu_\beta|^2\right)= 1,
$$
which shows that $z_\mu\in F^2(H_n)$. If
$\lambda:=(\lambda_1,\ldots, \lambda_n)\in \cD_p^\circ$, then there
exists $r\in (0,1)$ such that $(\frac{1}{r}\lambda_1,\ldots,
\frac{1}{r} \lambda_n)\in \cD_p^\circ(\CC)$. Applying the above
relation to the latter $n$-tuple of complex numbers, we deduce that
$\sum\limits_{\beta\in \FF_n^+}\frac{1}{r^{2|\beta|}} b_\beta
|\lambda_\beta|^2<\infty$. On the other hand, since $G$ is a free
holomorphic function on $\cD_p$, the operator $G(rW_1,\ldots, rW_n)$
is in $\cA_n(\cD_p)$. Hence, $G(rW_1,\ldots, rW_n)(1)\in F^2(H_n)$
and consequently $\sum\limits_{\beta\in \FF_n^+} r^{2|\beta|}
|c_\beta|^2 \frac{1}{b_\beta}<\infty$. Now, we have
\begin{equation*}
\begin{split}
\sum_{\beta\in \FF_n^+} |c_\beta|
|\lambda_\beta|\leq\left(\sum_{\beta\in \FF_n^+} r^{2|\beta|}
|c_\beta|^2
\frac{1}{b_\beta}\right)^{1/2}\left(\sum\limits_{\beta\in
\FF_n^+}\frac{1}{r^{2|\beta|}} b_\beta
|\lambda_\beta|^2\right)^{1/2}<\infty,
\end{split}
\end{equation*}
which clearly implies that $g(\lambda_1,\ldots, \lambda_n)$ is  a
holomorphic function on  the Reinhardt domain $\cD^\circ_p(\CC)$.

 To prove the last part of the proposition, let
$\{G_m\}_{m=1}^\infty$ and $G$ be in $Hol(\cD_p)$ and let
$\{g_m\}_{m=1}^\infty$ and $g$  be the corresponding representations
on $\CC$, respectively. Due to  to Proposition \ref{pro}, we have
\begin{equation*}
\begin{split}
\sup_{ (\lambda_1,\ldots, \lambda_n)\in \cD_{p,r}(\CC)}
|g_m(\lambda_1,\ldots, \lambda_n)&-g(\lambda_1,\ldots,
\lambda_n)|\\
&\leq \|G_m(rW_1,\ldots, rW_n)-G(rW_1,\ldots, rW_n)\|
\end{split}
\end{equation*}
for any $r\in [0,1)$. Hence, we deduce that if $G_m\to G$ in the
metric $\rho$ of $Hol(\cD_p)$, then $g_m\to g$ uniformly on compact
subsets of $\cD^\circ_p(\CC)$, which completes the proof.
\end{proof}

Using Proposition \ref{Prop-Cauc} and  Proposition
\ref{F-infty-cauc}, we can obtain the following.

\begin{corollary}\label{an=cauch}
Let $T:=(T_1,\ldots, T_n)\in B(\cH)^n$ be an   $n$-tuple of
operators with  joint spectral radius $r_f(T_1,\ldots, T_n)<1$.
Then,  the map $\Psi_{f,T}:F_n^\infty(\cD_f)\to B(\cH)$ defined by
$$
\Psi_{f,T}(g)=g(T_1,\ldots, T_n):=\cC_T(g(W_1,\ldots, W_n)),
$$
for any $g(W_1,\ldots, W_n)\in F_n^\infty(\cD_f)$,  is a unital  WOT
continuous  homomorphism such that $\Psi_{f,T}(W_\alpha)=T_\alpha$
for any $\alpha\in \FF_n^+$. Moreover,
$$
\|g(T_1,\ldots, T_n)\|\leq \left(\sum\limits_{k=0}^\infty\left\|
\Phi_{f,T}^k(I)\right\|^{1/2}\right) \|g(W_1,\ldots, W_n)\|.
$$
\end{corollary}

In the last part of this section,  we introduce  a  Banach space of
analytic functions on $\cD_p^\circ(\CC)\subset\CC^n$ and obtain  von
Neumann  and Bohr type inequalities in the commutative  setting.

  Let
  $p:=\sum_{1\leq|\alpha|\leq m} a_\alpha X_\alpha$ be a
  positive regular noncommutative polynomial
  and
  let ${\bf p}:=(p_1,\ldots, p_n)$ be a multi-index in $\ZZ_+^n$.
  If
$\lambda:=(\lambda_1,\ldots,\lambda_n)\in\CC^n$, then we set
$\lambda^{\bf p}:=\lambda_1^{p_1}\cdots \lambda_n^{p_n}$ and define
the symmetrized functional calculus
$$
(\lambda^{\bf p})_{\text{\rm sym}} (W_1,\ldots, W_n):=\frac {1}
{\gamma_{\bf p} }\sum_{\alpha\in \Lambda_{\bf p}} b_\alpha W_\alpha,
$$
where
$$
\Lambda_{\bf p}:=\{\alpha\in \FF_n^+: \lambda_\alpha= \lambda^{\bf
p} \text{ for any } \lambda\in \BB_n\},
$$
the coefficients $b_\alpha$ are defined by \eqref{b_alpha},  and
$W_1,\ldots, W_n$ are the weighted left creation operators
associated with $\cD_p$.
  Denote by $H_{\text{\rm sym}}(\cD_p^\circ(\CC))$ the set of all
analytic functions on $\cD_p^\circ(\CC)$ with scalar coefficients
 $$
g(\lambda_1,\ldots,\lambda_n):=\sum\limits_{\bf p\in \ZZ_+^n}
\lambda^{\bf p} a_{\bf p}, \quad  a_{\bf p}\in \CC,
$$
such that
 \begin{equation*}
g_{\text{\rm sym}}(rW_1,\ldots, rW_n) :=\sum_{k=0}^\infty
\sum\limits_{{\bf p}\in \ZZ_+^n,|{\bf p}|=k}
r^k a_{{\bf p}}[(\lambda^{\bf p})_{\text{\rm sym}} (W_1,\ldots,
W_n)]
\end{equation*}
is convergent in the operator norm for any $r\in [0,1)$.

Setting $g_{\text{\rm sym}}(rW_1,\ldots, rW_n)=\sum_{k=0}^\infty
\sum_{|\alpha|=k} r^{|\alpha|} c_{\alpha}W_\alpha$, we have
 $c_{0}:=a_{0}$ and $c_{\alpha}:= b_\alpha
 \frac {a_{{\bf p}}}{\gamma_{{\bf p}}}$,   where  ${\bf p}\in \ZZ_+^n$,
  ${\bf p}\neq (0,\ldots, 0)$,  is uniquely determined such that
  $\lambda_\alpha=\lambda^{\bf p}$ for any $\lambda\in \BB_n$.

 Notice that  $g_{\text{\rm sym}}(Z_1,\ldots, Z_n)$ is a free
holomorphic function on $\cD_p$. We define  $H_{\text{\rm
sym}}^\infty(\cD_p^\circ(\CC)) $ as the set of all functions $g\in
H_{\text{\rm sym}}(\cD_p^\circ(\CC))$ such that
$$\|g\|_{\text{\rm sym}}:=\sup_{0\leq r<1}
\left\| g_{\text{\rm sym}}(rW_1,\ldots, rW_n)\right\|<\infty.
$$
\begin{theorem}\label{sym}
$\left(H_{\text{\rm sym}}^\infty(\cD_p^\circ(\CC)),
\,\|\cdot\|_{\text{\rm sym}}\right)$ is a Banach space.
\end{theorem}

\begin{proof}
First notice that if $g\in H_{\text{\rm
sym}}^\infty(\cD_p^\circ(\CC))$ then $g_{\text{\rm
sym}}(rW_1,\ldots, rW_n)$ is norm convergent and $g_{\text{\rm
sym}}(Z_1,\ldots, Z_n)$ is a free holomorphic function on the
noncommutative domain $\cD_p$ . Using Theorem \ref{free-ho}, it is
easy to see that $H_{\text{\rm sym}}^\infty(\cD_p^\circ(\CC))$ is a
vector space and $\|\cdot \|_{\text{\rm sym}}$ is a norm. Let
$\{g_m\}_{m=1}^\infty$ be a Cauchy sequence of functions in
$H_{\text{\rm sym}}^\infty(\cD_p^\circ(\CC))$. According to Theorem
\ref{f-infty}, $(g_m)_{\text{\rm sym}}\in F_n^\infty(\cD_p)$ and
$\{(g_m)_{\text{\rm sym}}\}_{m=1}^\infty$  is a Cauchy sequence in
$\|\cdot\|_\infty$, the norm of the Banach algebra
$F_n^\infty(\cD_p)$. Therefore, there exists
 $\varphi(W_1,\ldots, W_n)\in F_n^\infty(\cD_p)$
such that $$\|(g_m)_{\text{\rm sym}}-\varphi(W_1,\ldots,
W_n)\|_\infty\to 0,\quad \text{ as }\ m\to\infty. $$
 If $
g_m(\lambda_1,\ldots,\lambda_n):=\sum\limits_{\bf p\in \ZZ_+^n}
a_{\bf p}^{(m)}\lambda^{\bf p}, \quad  a_{\bf p}\in \CC, $ then
$(g_m)_{\text{\rm sym}} =\sum_{k=0}^\infty\sum_{|\alpha|=k}
c_\alpha^{(m)} W_\alpha$, where $c_\alpha^{(m)}:=  b_\alpha
 \frac {a_{{\bf p}}^{(m)}}{\gamma_{{\bf p}}}$.
 If $\sum_{\alpha\in \FF_n^+} b_\alpha W_\alpha$ is  the Fourier
representation of the operator $\varphi(W_1,\ldots, W_n)$, then we
have
\begin{equation*}
\begin{split}
\frac{1}{b_\alpha}|c_\alpha^{(m)}-d_\alpha|&=
\left|\left<W_\alpha^*[(g_m)_{\text{\rm sym}}
(W_1,\ldots, W_n)-\varphi(W_1,\ldots, W_n)]1,1\right>\right|\\
&\leq\|(g_m )_{\text{\rm sym}}(W_1,\ldots, W_n)-\varphi(W_1,\ldots,
W_n)]\|_\infty.
\end{split}
\end{equation*}
Taking $m\to \infty$,  we deduce that $c_\alpha^{(m)}\to d_\alpha$
for each $\alpha\in \FF_n^+$.
Since $c_\alpha^{(m)}:=  b_\alpha
 \frac {a_{{\bf p}}^{(m)}}{\gamma_{{\bf p}}}$, we deduce that
 $d_\alpha=b_\alpha \frac {a_{{\bf p}}'}{\gamma_{{\bf p}}}$, where
 $a_{{\bf p}}'=\frac{d_\alpha \gamma_{{\bf p}}}{b_\alpha}$.
  Setting
$h(\lambda_1,\ldots, \lambda_n):=  \sum\limits_{\bf p\in \ZZ_+^n}
a_{\bf p}'\lambda^{\bf p}$,  we have $h_{\text{\rm
sym}}=\varphi(W_1,\ldots, W_n)$. Moreover, $\|h\|_{\text{\rm
sym}}=\|\varphi(W_1,\ldots, W_n)\|<\infty$. This shows that
$H_{\text{\rm sym}}^\infty(\cD_p^\circ(\CC))$ is a Banach space.
\end{proof}

Now, using Theorem \ref{abel}, we can deduce the following.

\begin{proposition}
  If  $T:=(T_1,\ldots, T_n)\in B(\cH)^n$  is  a commuting $n$-tuple
   of operators with  the joint spectral radius $r_p(T_1,\ldots, T_n)<1$  and
$g(\lambda_1,\ldots,\lambda_n):=\sum\limits_{\bf p\in \ZZ_+^n}
a_{\bf p} \lambda^{\bf p} $ is in $ H_{\text{\rm
sym}}(\cD_p^\circ(\CC))$, then
$$g(T_1,\ldots, T_n):=
\sum_{k=0}^\infty \sum\limits_{{\bf p}\in \ZZ_+^n,|{\bf p}|=k}a_{\bf
p} T^{\bf p}
$$
is a well-defined operator in $B(\cH)$, where the series is
convergent in the operator norm topology. Moreover, the map
$$
\Psi_T: H_{\text{\rm sym}}(\cD_p^\circ(\CC))\to B(\cH)\qquad
\Psi_T(g):=g(T_1,\ldots, T_n)
$$
is continuous and
$$
\|g(T_1,\ldots, T_n) \|\leq M\|g\|_{\text{\rm sym}},
$$
where $M=\sum_{k=0}^\infty \left\| \Phi_{p,T}^k(I)\right\|^{1/2}$.
\end{proposition}

Using now the  Bohr type inequality of Section \ref{Schwartz}, we
can obtain the following commutative version for $H_{\text{\rm
sym}}(\cD_p^\circ(\CC))$.

\begin{corollary}\label{commu}
Let $g(\lambda_1,\ldots,\lambda_n):=\sum\limits_{{\bf p}\in \ZZ_+^n}
\lambda^{{\bf p}} a_{{\bf p}}$ be  a  holomorphic function in
$\cD_p^\circ(\CC)$ such that
    $$
  \|
g_{\text{\rm sym}}(rW_1,\ldots, rW_n)\|\leq 1\ \  \text{ for } 0\leq
r<1.
$$
 Then
$$\sum\limits_{k=0}^\infty
 \sum\limits_{{\bf p}\in \ZZ_+^n, |{\bf p}|=k} |\lambda^{\bf p}|
  |a_{{\bf p}}| \leq  1
$$
for any $\lambda:=(\lambda_1,\ldots, \lambda_n)\in
\cD_{p,1/3}(\CC)$.
\end{corollary}

  In particular, we obtain the following Bohr type inequality for
  analytic polynomials.

\begin{corollary}
Let $N=2,3,\ldots$, and let
$$q(\lambda_1,\ldots,\lambda_n):=\sum\limits_{{\bf p}\in \ZZ_+^n,
|{\bf p}|\leq N-1} \lambda^{{\bf p}} a_{{\bf p}}, \quad a_{{\bf
p}}\in \CC,
$$
 be  an  analytic polynomial   such that
 $\|
q_{\text{\rm sym}}(W_1,\ldots, W_n)\|\leq 1 $.
 Then
$$\sum\limits_{k=0}^{N-1}
 \sum\limits_{{\bf p}\in \ZZ_+^n, |{\bf p}|=k}
  |\lambda^{\bf p}| |a_{{\bf p}}| \leq  1
$$
for any $\lambda:=(\lambda_1,\ldots, \lambda_n)\in
\cD_{p,1/3}(\CC)$.
\end{corollary}

\chapter{Model theory and unitary invariants    on
 noncommutative domains}
 \label{Unitary}

\section{Weighted shifts and  invariant subspaces}
\label{Weighted}

In this section we obtain a Beurling type characterization of the
invariant subspaces under the weighted left creation operators
$(W_1,\ldots, W_n)$ associated with the noncommutative domain
$\cD_f$. We deduce  a similar result for the model shifts
$(B_1,\ldots,B_n)$ associated with a noncommutative variety
$\cV_{f,J}\subset\cD_f$, generated by a $w^*$-closed two sided ideal
of the Hardy algebra $F_n^\infty(\cD_f)$.

 Define the positive linear mapping $\Phi_{f,W\otimes
I}:B(F^2(H_n)\otimes\cH)\to B(F^2(H_n)\otimes \cH)$ by setting
$$
\Phi_{f,W\otimes I}(Y):=\sum\limits_{|\alpha|\geq 1} a_\alpha
(W_\alpha\otimes I_\cH) Y(W_\alpha\otimes I_\cH)^*,
$$
where the convergence of the series  is in the weak operator
topology. We recall that an operator $M:F^2(H_n)\otimes \cH\to
F^2(H_n)\otimes \cK$  is called multi-analytic  with respect to
$W_1,\ldots, W_n$ if $M(W_i\otimes I_\cH)=(W_i\otimes I_\cK)M$ \,
for any $i=1,\ldots,n$. In case $M$ is a partial isometry, we call
it inner.

\begin{theorem}\label{Beurling}
Let $f$ be a positive regular free holomorphic function on
$B(\cH)^n$.
 If $Y\in
B(F^2(H_n)\otimes \cH)$, then the following statements are
equivalent.
\begin{enumerate}
\item[(i)]
There is a Hilbert space $\cE$ and a multi-analytic operator
$\Psi:F^2(H_n)\otimes \cE\to F^2(H_n)\otimes \cH$ with respect to
the weighted left creation operators $W_1,\ldots, W_n$,  such that
$Y=\Psi \Psi^*$.
\item[(ii)]
$\Phi_{f,W\otimes I}(Y)\leq Y$.
\end{enumerate}
\end{theorem}

\begin{proof}
Assume that (ii) holds. Then $Y-\Phi_{f,W\otimes I}^m(Y)\geq 0$ for
any $m=1,2,\ldots$. Since $(W_1,\ldots, W_n)$ is of class $C_{\cdot
0}$ with respect to $\cD_f$, we have  \
SOT-$\lim_{m\to\infty}\Phi_{f,W\otimes I}^m(Y)=0$, which implies
$Y\geq0$. Let $\cM:=\overline{\text{\rm range} Y^{1/2}}$ and define
\begin{equation}\label{ai}
A_i(Y^{1/2} x):=Y^{1/2} (W_i^*\otimes I_\cH)x,\quad x\in
F^2(H_n)\otimes \cH,
\end{equation}
 for any $i=1,\ldots, n$. Notice that
\begin{equation*}
\begin{split}
\sum_{|\alpha|\geq 1} a_\alpha\|A_{\tilde \alpha}Y^{1/2}x\|^2&\leq
\sum_{|\alpha|\geq 1} \|Y^{1/2} (\sqrt{a_\alpha} W_\alpha^*\otimes
I)x\|^2\\
&=\left< \Phi_{f,W\otimes I}(Y)x,x\right>\leq \|Y^{1/2} x\|^2
\end{split}
\end{equation*}
for any $x\in F^2(H_n)\otimes \cH$.  Hence, we have
$a_{g_i}\|A_iY^{1/2} x\|^2\leq \|Y^{1/2} x\|^2$, for any $x\in
F^2(H_n)\otimes \cH$.  Since $a_{g_i}\neq 0$, $A_i$ can be uniquely
be extended to a bounded operator (also denoted by $A_i$) on the
subspace $\cM$. Denoting $T_i:=A_i^*$, $i=1,\ldots, n$,  an
approximation argument shows that
$$
\sum_{|\alpha|\geq 1} a_\alpha T_\alpha T_\alpha^*\leq I_\cM.
$$
On the other hand, due to \eqref{ai}, we have
\begin{equation*}
\begin{split}
\left< \Phi_{f, T}^m(I) Y^{1/2}x, Y^{1/2} x\right> &= \left<
\Phi_{f,
W\otimes I}^m(Y)x,  x\right>\\
&\leq \|Y\| \left< \Phi_{f, W\otimes I}^m(I)x,  x\right>
\end{split}
\end{equation*}
for any $x\in F^2(H_n)\otimes \cH$. Since \
SOT-$\lim\limits_{m\to\infty}\Phi_{f,W\otimes I}^m(I)=0$, we have \
SOT-$\lim\limits_{m\to\infty}\Phi_{f,T}^m(I)=0$, which shows that
$T:=(T_1,\ldots, T_n)$ is of class $C_{\cdot 0}$ with respect to
$\cD_f(\cM)$. According to Section \ref{domain algebra}, the Poisson
kernel $K_{f,T}:\cM\to F^2(H_n)\otimes \cE$ ($\cE$ is an appropriate
Hilbert space) is an isometry with the property that
\begin{equation}
\label{int-KT} T_iK_{f,T}^*=K_{f,T}^* (W_i\otimes I_\cE),\qquad
i=1,\ldots,n.
\end{equation}

Let $\Psi:=Y^{1/2} K_{f,T}^*:F^2(H_n)\otimes \cE\to F^2(H_n)\otimes
\cH$ and notice that
\begin{equation*}
\begin{split}
\Psi(W_i\otimes I_\cE)&=Y^{1/2} K_{f,T}^*(W_i\otimes I_\cE)=Y^{1/2}
T_i K_{f,T}^*\\
&=(W_i\otimes I_\cH) Y^{1/2} K_{f,T}^* =(W_i\otimes I_\cH) \Psi
\end{split}
\end{equation*}

for any $i=1,\ldots, n$. Notice also that $\Psi \Psi^*=Y^{1/2}
K_{f,T}^* K_{f,T} Y^{1/2} =Y$.
\end{proof}

Now, we can obtain a Beurling \cite{Be} type characterization  of
the invariant subspaces under  the weighted left creation operators
associated  with a noncommutative domain $\cD_f$.

\begin{theorem} \label{inv-subs} Let $f$ be a positive regular free
 holomorphic function on $B(\cH)^n$, and  let $(W_1,\ldots, W_n)$ be
the weighted left creation operators associated  with a
noncommutative domain $\cD_f$. A subspace $\cM\subseteq
F^2(H_n)\otimes \cH$ is invariant under each operator $W_1\otimes
I_\cH,\ldots, W_n\otimes I_\cH$  if and only if there exists an
inner multi-analytic operator $\Psi:F^2(H_n)\otimes \cE\to
F^2(H_n)\otimes \cH$ with respect to the weighted left creation
operators $W_1,\ldots, W_n$ such that
$$
\cM=\Psi[F^2(H_n)\otimes \cE].
$$
\end{theorem}
\begin{proof}
Assume that $\cM\subseteq F^2(H_n)\otimes \cH$ is invariant under
each operator $W_1\otimes I_\cH,\ldots, W_n\otimes I_\cH$. Since
$P_\cM(W_i\otimes I_\cH)P_\cM=(W_i\otimes I_\cH)P_\cM$ for any
$i=1,\ldots, n$, and  $(W_1,\ldots, W_n)\in \cD_f(F^2(H_n))$, we
have
\begin{equation*}
\begin{split}
\Phi_{f, W\otimes I_\cH}(P_\cM)&=P_\cM\left[ \sum_{|\alpha|\geq 1}
a_\alpha (W_\alpha \otimes I_\cH) P_\cM (W_\alpha ^*\otimes
I_\cH)\right] P_\cM\\
&\leq P_\cM\left[ \sum_{|\alpha|\geq 1} a_\alpha (W_\alpha \otimes
I_\cH)  (W_\alpha ^*\otimes
I_\cH)\right] P_\cM\\
&\leq P_\cM
\end{split}
\end{equation*}
Applying Theorem \ref{Beurling} to our particular case,  there is
multi-analytic operator $\Psi:F^2(H_n)\otimes \cE\to F^2(H_n)\otimes
\cH$ with respect to the weighted left creation operators
$W_1,\ldots, W_n$ such that $P_\cM=\Psi \Psi^*$. Since $P_\cM$ is an
orthogonal projection, we deduce that $\Psi$ is a partial isometry
and $\cM=\Psi[F^2(H_n)\otimes \cE]$. The converse is obvious. The
proof is complete.
\end{proof}

Now, we turn our attention to noncommutative varieties
$\cV_{f,J}\subseteq \cD_f$ and the  constrained shifts associated
with them.

Let $J$ be a $w^*$-closed two-sided ideal of $F_n^\infty(\cD_f)$
such that $J\neq F_n^\infty(\cD_f)$, and define the subspaces of
$F^2(H_n)$ by
$$
\cM_J:=\overline{JF^2(H_n)}\quad \text{and}\quad
\cN_J:=F^2(H_n)\ominus \cM_J.
$$
It is easy  to see that the subspace $\cN_J$  is  invariant under
the operators
 $W_1^*,\ldots, W_n^*$ and  $\Lambda_1^*,\ldots, \Lambda_n^*$.
On the other hand, notice that $\cN_J\neq 0$ if and only if $ J\neq
F_n^\infty(\cD_f)$.
 Indeed, due to
Theorem \ref{F/J}, for any $\varphi\in F_n^\infty(\cD_f)$,
$$
d(\varphi, J)=\|P_{\cN_J}\varphi (W_1,\ldots, W_n)|\cN_J\|.
$$
Define the {\it constrained  weighted left} (resp.~{\it right}) {\it
creation operators} associated with the noncommutative variety
$\cV_{f,J}$ by setting
$$B_i:=P_{\cN_J} W_i|\cN_J \quad \text{and}\quad C_i:=P_{\cN_J} \Lambda_i|\cN_J,\quad i=1,\ldots, n.
$$
Let $F_n^\infty(\cV_{f,J})$ be the $w^*$-closed algebra generated by
$B_1,\ldots, B_n$ and the identity. We showed, in Section
\ref{Functional II}, that
$$
F_n^\infty(\cV_{f,J})=P_{\cN_J}F_n^\infty(\cD_f)
|\cN_J=\{g(B_1,\ldots, B_n):\ g(W_1,\ldots, W_n)\in
F_n^\infty(\cD_f)\},
$$
where,  according to the $F_n^\infty(\cD_f)$-functional calculus,
$$g(B_1,\ldots, B_n)=\text{\rm SOT-}\lim\limits_{r\to 1}g(rB_1,\ldots, rB_n).
$$
Note that if $\varphi\in J$, then $\varphi(B_1,\ldots, B_n)=0$.
Therefore, $(B_1,\ldots, B_n)$ is in $\cV_{f,J}(\cN_J)$ and, due to
the results of Section \ref{Functional II}, it plays the role of
universal model for  the noncommutative variety $\cV_{f,J}$.

Similar results hold for $R_n^\infty(\cV_{f,J})$, the $w^*$-closed
algebra generated by $C_1,\ldots, C_n$ and the identity. Moreover,
one  can prove that
\begin{equation}\label{commutant0}
F_n^\infty(\cV_{f,J})^\prime=R_n^\infty(\cV_{f,J})\ \text{ and } \
R_n^\infty(\cV_{f,J})^\prime=F_n^\infty(\cV_{f,J}),
\end{equation}
where $^\prime$ stands for the commutant. An operator $M\in
B(\cN_J\otimes \cK,\cN_J\otimes \cK')$ is called multi-analytic with
respect to the constrained shifts $B_1,\ldots, B_n$ if
$$
M(B_i\otimes I_{\cK})=(B_i\otimes I_{\cK'})M,\quad i=1,\ldots, n.
$$
If  $M$ is partially isometric, then we call it an inner
multi-analytic operator. Using  \eqref{commutant0}, we can  show
that the set of all multi-analytic operators with respect to
 $B_1,\ldots, B_n$ coincides  with
$$
R_n^\infty(\cV_{f,J})\bar\otimes B(\cK,\cK')=P_{\cN_J\otimes
\cK'}[R_n^\infty(\cD_f)\bar\otimes B(\cK,\cK')]|\cN_J\otimes \cK,
$$
and a similar result holds for  the algebra $F_n(\cV_{f,J})$.

The next result provides  a Beurling type  characterization of the
invariant subspaces under the constrained weighted shifts
$B_1,\ldots, B_n$, associated with the noncommutative variety
$\cV_{f,J}$.

\begin{theorem}\label{Beur} Let $J\neq F_n^\infty(\cD_f)$ be a $w^*$-closed
two-sided ideal of $F_n^\infty(\cD_f)$ and let $B_1,\ldots, B_n$ be
the corresponding constrained left creation operators associated
with  the noncommutative variety $\cV_{f,J}(\cN_J)$. A subspace
$\cM\subseteq \cN_J\otimes \cK$ is invariant under each operator
$B_i\otimes I_\cK$, \ $i=1,\ldots, n$, if and only if there exists a
Hilbert space  $\cG$ and an inner operator
$$\Theta(C_1,\ldots, C_n)\in R_n^\infty(\cV_{f,J})
\bar\otimes B(\cG,\cK)$$ such that
\begin{equation}\label{Beu}
\cM=\Theta(C_1,\ldots, C_n)\left[\cN_J\otimes \cG\right].
\end{equation}
\end{theorem}

\begin{proof}
Due to the remarks preceding the theorem, the subspace $\cN_J\otimes
\cK$ is invariant under each operator $W_i^*\otimes I_\cK$,
$i=1,\ldots, n$, and
$$(W_i^*\otimes I_\cK)|\cN_J\otimes \cK=B_i^*\otimes I_\cK,
 \quad i=1,\ldots,n. $$
Since the subspace $[\cN_J\otimes \cK]\ominus \cM$  is invariant
under
 $B_i^*\otimes I_\cK$, \ $i=1,\ldots, n$, we deduce that it is
  also invariant under  each  operator $W_i^*\otimes I_\cK$,
   \ $i=1,\ldots, n$.
Therefore, the subspace
\begin{equation}\label{E}
\cE:=[F^2(H_n)\otimes \cK]\ominus\{[\cN_J\otimes \cK]\ominus \cM\}=
[\cM_J\otimes \cK]\oplus \cM
\end{equation}
is invariant under $W_i\otimes I_\cK$, \ $i=1,\ldots, n$, where
$\cM_J:=F^2(H_n)\ominus \cN_J$. Using Theorem \ref{inv-subs}, we
find a Hilbert space $\cG$ and  an inner multi-analytic operator
$$\Theta(\Lambda_1,\ldots,\Lambda_n)\in R_n^\infty(\cD_f)\bar
 \otimes B(\cG,\cK)$$
such that
$$
\cE=\Theta(\Lambda_1,\ldots,\Lambda_n)[F^2(H_n)\otimes \cG].
$$
Since $\Theta(\Lambda_1,\ldots,\Lambda_n)$ is a partial isometry, we
have
\begin{equation}\label{proj-alta}
P_\cE=\Theta(\Lambda_1,\ldots,\Lambda_n)\Theta(\Lambda_1,\ldots,\Lambda_n)^*,
\end{equation}
where $P_\cE$ is the orthogonal projection of $F^2(H_n)\otimes \cK$
onto $\cE$. Notice that  the subspace $\cN_J\otimes \cK$ is
invariant under the operators $\Lambda_i^*\otimes I_\cK$,
$i=1,\ldots, n$.
Moreover, using the remarks preceding the theorem we have
$$\Theta(C_1,\ldots,C_n)=
P_{\cN_J\otimes \cK}\Theta(\Lambda_1,\ldots, \Lambda_n)|\cN_J\otimes
\cK.$$
Now, it is clear that equation \eqref{proj-alta}
$$
P_{\cN_J\otimes \cK}P_\cE|\cN_J\otimes
\cK=\Theta(C_1,\ldots,C_n)\Theta(C_1,\ldots,C_n)^*.
$$
Due to relation \eqref{E}, the left hand side of this equality is
equal to $P_\cM$, the orthogonal projection  of $\cN_J\otimes \cK$
onto $\cM$. Hence, we deduce that
$$
P_\cM=\Theta(C_1,\ldots,C_n)\Theta(C_1,\ldots,C_n)^*,
$$
the operator  $\Theta(C_1,\ldots,C_n)$ is a partial isometry, and
$$\cM=\Theta(C_1,\ldots, C_n)\left[ \cN_J\otimes \cG\right].$$
 The proof
is complete.
\end{proof}

From now on, throughout this section,  we assume    that $J$ is a
$w^*$-closed two-sided ideal of $F_n^\infty(\cD_f)$ such that
$1\in\cN_J$. We recall that a   subspace $\cH\subseteq \cK$ is
called co-invariant under $\cS\subset B(\cK)$ if $X^*\cH\subseteq
\cH$ for any $X\in \cS$.

\begin{theorem}\label{cyclic}
Let  $J$ be a $w^*$-closed two-sided ideal of $F_n^\infty(\cD_f)$
such that $1\in \cN_J$ and let $\cK$  be a Hilbert space. If
$\cM\subseteq \cN_J\otimes \cK$ is a co-invariant subspace under
$B_i\otimes I_\cK$, \ $i=1,\ldots, n$, then  there exists a subspace
$\cE\subseteq \cK$ such that
\begin{equation}\label{span}
\overline{\text{\rm span}}\,\left\{\left(B_\alpha\otimes
I_\cK\right)\cM:\ \alpha\in \FF_n^+\right\}=\cN_J\otimes \cE.
\end{equation}
\end{theorem}

\begin{proof}
  Define the subspace $\cE\subset \cK$ by
 $\cE:=(P_\CC\otimes I_\cK)\cM$ and  let $\varphi$ be an arbitrary
nonzero element of $\cM$ with  Fourier representation
$$\varphi=\sum_{\alpha\in \FF_n^+}e_\alpha\otimes h_\alpha,
\quad h_\alpha\in \cK.
$$
Let $\beta\in \FF_n^+$ be such that $h_\beta\neq 0$ and notice that
\begin{equation}\label{PB*}
(P_\CC\otimes I_\cK)(B_\beta^*\otimes I_\cK)\varphi=1\otimes
\frac{1}{\sqrt{b_\beta}} h_\beta.
\end{equation}
Indeed, since $1\in \cN_J$ and using relation \eqref{projW*}, we
have
\begin{equation*}
\begin{split}
&\left< (P_\CC\otimes I_\cK) (B_\beta^*\otimes
I_\cK)\varphi,1\otimes
h\right>\\
&\qquad=\left<(P_{\cN_J}\otimes I_\cK)
(W_\beta^*\otimes I_\cK)\varphi,1\otimes h\right>\\
&\qquad=\left<(W_\beta^*\otimes I_\cK)\varphi, 1\otimes
h\right>=\left< (P_\CC\otimes I_\cK) (W_\beta^*\otimes
I_\cK)\varphi,1\otimes
h\right>\\
&\qquad=\frac{1}{\sqrt{b_\beta}} \left<h_\beta, h\right>
\end{split}
\end{equation*}
for any $h\in \cK$. Since $\left< (P_\CC\otimes I_\cK)
(B_\beta^*\otimes I_\cK)\varphi,e_\gamma\otimes h\right>=0$ for any
$\gamma\in \FF_n^+$ with $|\gamma|\geq 1$, and $h\in \cK$, the
result follows.
 Now, since $\cM$ is  a
co-invariant subspace under $B_i\otimes I_\cK$, \ $i=1,\ldots, n$,
relation \eqref{PB*}  implies  that $h_\beta\in \cE$.  Hence, and
taking into account that $1\in\cN_J$, we  have
$$
(B_\beta\otimes I_\cK)(1\otimes
h_\beta)=\frac{1}{\sqrt{b_\beta}}(P_{\cN_J} e_\beta)\otimes
h_\beta\in \cN_J\otimes \cE, \quad \text{ for }\ \beta\in \FF_n^+.
$$
Since $\varphi\in \cM\subseteq \cN_J\otimes \cK$, we deduce that
$$
\varphi=(P_{\cN_J}\otimes
I_\cK)\varphi=\lim_{m\to\infty}\sum_{k=0}^m \sum_{|\alpha|=k}
P_{\cN_J}e_\beta\otimes h_\beta
$$
is in $\cN_J\otimes \cE$. Therefore $\cM\subset\cN_J\otimes \cE$ and
$$
\cY:= \overline{\text{\rm span}}\,\left\{\left(B_\alpha\otimes
I_\cK\right)\cM:\ \alpha\in \FF_n^+\right\}\subset\cN_J\otimes \cE.
$$

To prove  the reverse inclusion, we show first that $\cE\subset
\cY$. Since $\cN_J$ is an  invariant subspace under each operator
$W_i^*$, \ $i=1,\ldots,n$, and contains the constants, we have
\begin{equation*}
\begin{split}
I_{\cN_J}-\sum_{|\alpha|\geq 1} a_\alpha B_\alpha B_\alpha^*
 &=P_{\cN_J}\left(I-\sum_{|\alpha|\geq 1} a_\alpha W_\alpha W_\alpha^*\right)|\cN_J\\
&=P_{\cN_J}P_\CC|\cN_J\\
&=P^{\cN_J}_\CC,
\end{split}
\end{equation*}
where $ P^{\cN_J}_\CC$ is the orthogonal projection of $\cN_J$ onto
$\CC$.
 If
$h_0\in \cE$, $h_0\neq 0$, then there exists $g\in \cM\subset
\cN_J\otimes \cE$ such that $g=1\otimes
h_0+\sum\limits_{|\alpha|\geq 1} e_\alpha\otimes h_\alpha.$
Consequently, we have
\begin{equation*}
\begin{split}
1\otimes h_0&=(P_\CC\otimes I_\cK) g=(P^{\cN_J}_\CC\otimes I_\cK)
g\\
&=\left(I_{\cN_J\otimes \cK}-\sum_{|\alpha|\geq 1}
a_\alpha(B_\alpha\otimes I_\cK) (B_\alpha\otimes I_\cK)^*\right)g.
\end{split}
\end{equation*}
 Hence and since $\cM$ is co-invariant under $B_i\otimes I_\cK$, \ $i=1,\ldots, n$,
we deduce that $h_0\in\cY$ for any $h_0\in \cE$, i.e., $\cE\subset
\cY$. The latter inclusion  shows that $(B_\alpha\otimes
I_\cK)(1\otimes \cE)\subset \cY$ for any $\alpha\in \FF_n^+$, which
implies
\begin{equation}\label{PN}
\frac{1}{\sqrt{b_\alpha}}P_{\cN_J}e_\alpha\otimes \cE\subset
\cY,\quad \alpha\in \FF_n^+.
\end{equation}
Assume that $\varphi\in \cN_J\otimes \cE\subset F^2(H_n)\otimes \cE$
has  Fourier representation $\varphi=\sum\limits_{\alpha\in
\FF_n^+}e_\alpha\otimes y_\alpha,\quad y_\alpha\in \cE. $   Due to
\eqref{PN}, we deduce that
$$
\varphi=(P_{\cN_J}\otimes I_\cE)\varphi=\lim_{m\to
\infty}\sum_{k=0}^m\sum_{|\alpha|=k}P_{\cN_J}e_\alpha\otimes
y_\alpha\in \cY.
$$
Therefore, $\cN_J\otimes \cE\subseteq \cY$.  The proof is complete.
\end{proof}

Now, we can easily deduce the following result.

\begin{corollary}
Let  $J$ be a $w^*$-closed two-sided ideal of $F_n^\infty(\cD_f)$
such that $1\in \cN_J$ and let $\cK$  be a Hilbert space.
 A subspace $\cM\subseteq \cN_J\otimes \cK$ is reducing under each operator
 $B_i\otimes I_\cK$, \ $i=1,\ldots, n$,
if and only if   there exists a subspace $\cE\subseteq \cK$ such
that
\begin{equation*}
 \cM=\cN_J\otimes \cE.
\end{equation*}
\end{corollary}

\bigskip

\section{$C^*$-algebras associated with noncommutative varieties
and Wold decompositions} \label{C*-algebras}

Given a positive regular noncommutative polynomial $p$ and a
$w^*$-closed two-sided ideal $J$ of the Hardy algebra
$F_n^\infty(\cD_p)$, we associate with each noncommutative variety
$\cV_{p,J}\subset \cD_p$ the $C^*$-algebra $C^*(B_1,\ldots, B_n)$
generated by the model operators $B_1,\ldots, B_n$ and the identity.
We obtain  Wold type decompositions for nondegenerate
$*$-representations of $C^*(B_1,\ldots, B_n)$.

We  recall that the universal model $(B_1,\ldots, B_n)$ associated
with the noncommutative variety $\cV_{p,J}$ is defined by
$B_i:=P_{\cN_J}W_i|\cN_J$, $i=1,\ldots,n$, where
$\cN_J:=F^2(H_n)\ominus \overline{JF^2(H_n)}$ and $(W_1,\ldots,W_n)$
is the universal model for the domain $\cD_p$.

 Assume that the noncommutative domain $\cD_p$ is
generated by a positive regular noncommutative  polynomial. We keep
the notation of Section \ref{Weighted}.

\begin{theorem}\label{compact}
Let $J$ be a $w^*$-closed two-sided ideal of $F_n^\infty(\cD_p)$
such that $1\in\cN_J$, and let $(B_1,\ldots, B_n)$ be the universal
model associated with the noncommutative variety $\cV_{p,J}$. Then
all the compact operators in $B(\cN_J)$ are contained in the
operator space
$$\overline{\text{\rm span}}\{B_\alpha B_\beta^*:\ \alpha,\beta\in \FF_n^+\}.
$$
Moreover,   the $C^*$-algebra  $C^*(B_1,\ldots, B_n)$ is
irreducible.
\end{theorem}

\begin{proof}
Since $1\in \cN_J$ and  $\cN_J$ is an  invariant subspace  $W_i^*$,
\ $i=1,\ldots,n$,   we deduce that
\begin{equation*}
I_{\cN_J}-\sum_{|\alpha|\geq 1} a_\alpha B_\alpha B_\alpha^*
 =P_{\cN_J}\left(I-\sum_{|\alpha|\geq 1} a_\alpha W_\alpha W_\alpha^*\right)|\cN_J
=P^{\cN_J}_\CC,
\end{equation*}
where $ P^{\cN_J}_\CC$ is the orthogonal projection of $\cN_J$ onto
$\CC$.
Fix
$$g(W_1,\ldots, W_n):=\sum\limits_{|\alpha|\leq m} d_\alpha W_\alpha\quad
\text{ and } \quad \xi:=\sum\limits_{\beta\in \FF_n^+} c_\beta
e_\beta\in \cN_J\subset F^2(H_n),
$$
and note  that
\begin{equation*}
\begin{split}
P^{\cN_J}_\CC g(B_1,\ldots, B_n)^*\xi&=\sum_{|\alpha|\leq m}P_\CC
\overline{d}_\alpha W_\alpha^*\xi=\sum_{|\alpha|\leq
m}\frac{1}{\sqrt{b_\alpha}}\overline{d}_\alpha  c_\alpha
\\
&=\left< \xi, \sum_{|\alpha|\leq m}
\frac{1}{\sqrt{b_\alpha}}d_\alpha e_\alpha\right>= \left<
\xi,g(B_1,\ldots, B_n)(1)\right>.
 \end{split}
\end{equation*}
Consequently,
\begin{equation}\label{rankone}
q(B_1,\ldots, B_n)P^{\cN_J}_\CC g(B_1,\ldots, B_n)^*\xi= \left<
\xi,g(B_1,\ldots, B_n)(1)\right>q(B_1,\ldots, B_n)(1)
\end{equation}
for any polynomial $q(B_1,\ldots, B_n)$ in $F_n^\infty(\cV_{p,J})$.
Hence, we deduce that the operator $q(B_1,\ldots, B_n)P^{\cN_J}_\CC
g(B_1,\ldots, B_n)^*$
 has rank one. Moreover, since
$P_\CC^{\cN_J}=I_{\cN_J}-\sum\limits_{|\alpha|\geq 1} a_\alpha B_
\alpha B_\alpha^*$, it is clear  that the above operator is also in
the operator space $\overline{\text{\rm span}}\{B_\alpha B_\beta^*:\
\alpha,\beta\in \FF_n^+\}$. Since  the   set
$$\cE:=\left\{\left(\sum\limits_{|\alpha|\leq m}d_\alpha
B_\alpha\right)(1):\ m\in \NN, d_\alpha\in \CC\right\}$$ is  dense
in $\cN_J$,   relation \eqref{rankone} implies  that all compact
operators  in $B(\cN_J)$ are included in the operator space
$\overline{\text{\rm span}}\{B_\alpha B_\beta^*:\ \alpha,\beta\in
\FF_n^+\}$.

To prove the last part of this theorem, let $\cM\neq\{0\}$ be a
subspace of $\cN_J\subseteq F^2(H_n)$, which is jointly reducing for
each operator $B_i$, $i=1,\ldots,n$. Let $\varphi\in \cM$,
$\varphi\neq 0$, and assume that
$\varphi=c_0+\sum\limits_{|\alpha|\geq 1} c_\alpha e_\alpha.$  If
 $c_\beta$  is a nonzero coefficient of $\varphi$,
then
 $
P_\CC B_\beta^*\varphi= \frac{1}{\sqrt{b_\beta}} c_\beta$. Indeed,
since $1\in \cN_J$, one can use relation \eqref{WbWb} to deduce that
\begin{equation*}
\begin{split}
\left< P_\CC
B_\beta^*\varphi,1\right>&=\left<P_{\cN_J}W_\beta^*\varphi,1\right>\\
&=\left<W_\beta^*\varphi, 1\right>=\left< P_\CC
W_\beta^*\varphi,1\right>=\frac{1}{\sqrt{b_\beta}} c_\beta.
\end{split}
\end{equation*}
Since $\left< P_\CC B_\beta^*\varphi,e_\gamma\right>=0$ for any
$\gamma\in \FF_n^+$ with $|\gamma|\geq 1$, our assertion  follows.
On the other hand, since
$P_\CC^{\cN_J}=I_{\cN_J}-\sum\limits_{|\alpha|\geq 1} a_\alpha
B_\alpha B_\alpha^* $ and $\cM$ is reducing for $B_1,\ldots, B_n$,
we deduce that $c_\beta\in \cM$, so $1\in \cM$. Using once again
that $\cM$ is invariant under the operators $B_1,\ldots, B_n$, we
have
 $\cE\subseteq \cM$. On the other hand, since $\cE$ is dense in $\cN_J$, we deduce that
 $\cN_J\subset \cM$. Therefore $\cN_J= \cM$.
This completes the proof.
\end{proof}

We say that two $n$-tuples of operators $(T_1,\ldots, T_n)$, \
$T_i\in B(\cH)$, and $(T_1',\ldots, T_n')$, \ $T_i'\in B(\cH')$, are
unitarily equivalent if there exists a unitary operator $U:\cH\to
\cH'$ such that
$$T_i=U^* T_i' U\ \text{  for  any } \ i=1,\ldots, n.
$$
 If $(B_1,\ldots, B_n)$ is the universal model
  associated with the noncommutative variety $\cV_{p,J}$,
then the $n$-tuple $(B_1\otimes I_\cH,\ldots, B_n\otimes I_\cH)$ is
called constrained weighted shift with multiplicity $\dim \cH$.
Using Theorem \ref{compact}, we can prove the following.

\begin{proposition}\label{eq-mult}
Two constrained weighted shifts associated with the noncommutative
variety $\cV_{p,J}$ are unitarily equivalent if and only if their
multiplicities are equal.
\end{proposition}
\begin{proof}
One implication is trivial. Let $(B_1\otimes I_\cH,\ldots,
B_n\otimes I_\cH)$ and $(B_1\otimes I_{\cH'},\ldots, B_n\otimes
I_{\cH'})$  be two constrained shifts and let $U:\cN_J\otimes \cH\to
\cN_J\otimes \cH'$ be a unitary operator such that
$$
U(B_i\otimes I_\cH)=(B_i\otimes I_{\cH'})U,\quad i=1,\ldots, n.
$$
Since $U$ is unitary , we have
$$
U(B_i^*\otimes I_\cH)=(B_i^*\otimes I_{\cH'})U,\quad i=1,\ldots, n.
$$
Since the $C^*$-algebra $C^*(B_1,\ldots, B_n)$ is irreducible, we
must have $U=I_{\cN_J}\otimes A$, where $A\in B(\cH,\cH')$ is a
unitary operator . Therefore, $\dim \cH=\dim \cH'$. The proof is
complete.
 \end{proof}

In what follows, we   prove  a Wold type decomposition for
non-degenerate
  $*$-representations of the $C^*$-algebra $C^*(B_1,\ldots,
B_n)$, generated by the the constrained weighted shifts associated
with $\cV_{p,J}$,  and the identity.

\begin{theorem}\label{wold} Let
$p=\sum\limits_{1\leq|\alpha|\leq m} a_\alpha X_\alpha$ be a
positive regular noncommutative polynomial and let  $J$ be a
$w^*$-closed two-sided ideal of the noncommutative Hardy algebra
$F_n^\infty(\cD_p)$  such that $1\in \cN_J$. If  \
$\pi:C^*(B_1,\ldots, B_n)\to B(\cK)$ is  a nondegenerate
$*$-representation  of $C^*(B_1,\ldots, B_n)$ on a separable Hilbert
space  $\cK$, then $\pi$ decomposes into a direct sum
$$
\pi=\pi_0\oplus \pi_1 \  \text{ on  } \ \cK=\cK_0\oplus \cK_1,
$$
where $\pi_0$ and  $\pi_1$  are disjoint representations of
$C^*(B_1,\ldots, B_n)$ on the Hilbert spaces
\begin{equation*}
\begin{split}
\cK_0:&=\overline{\text{\rm span}}\left\{\pi(B_\beta)
\left(I_\cK-\sum_{1\leq |\alpha|\leq m} a_\alpha
\pi(B_\alpha)\pi(B_\alpha)^*\right) \cK:\ \beta\in
\FF_n^+\right\}\quad \text{ and } \\
 \cK_1:&=\cK_0^\perp,
\end{split}
\end{equation*}
 respectively. Moreover, up to an isomorphism,
\begin{equation}\label{shi}
\cK_0\simeq\cN_J\otimes \cG, \quad  \pi_0(X)=X\otimes I_\cG \quad
\text{ for } \  X\in C^*(B_1,\ldots, B_n),
\end{equation}
 where $\cG$ is a Hilbert space with
$$
\dim \cG=\dim \left[\text{\rm range}\,\left(I_\cK-\sum_{1\leq
|\alpha|\leq m} a_\alpha \pi(B_\alpha)\pi(B_\alpha)^*\right)\right],
$$
 and $\pi_1$ is a $*$-representation  which annihilates the compact operators   and
$$
\sum_{1\leq |\alpha|\leq m} a_\alpha
\pi_1(B_\alpha)\pi_1(B_\alpha)^*=I_{\cK_1}.
$$
If  $\pi'$ is another nondegenerate  $*$-representation of
$C^*(B_1,\ldots, B_n)$ on a separable  Hilbert space  $\cK'$, then
$\pi$ is unitarily equivalent to $\pi'$ if and only if
$\dim\cG=\dim\cG'$ and $\pi_1$ is unitarily equivalent to $\pi_1'$.
\end{theorem}

\begin{proof}
Since  $1\in\cN_J$, Theorem \ref{compact} shows  that all the
compact operators $\cC(\cN_J))$ in $B(\cN_J)$ are contained in the
$C^*$-algebra  $C^*(B_1,\ldots, B_n)$. Due to  standard theory of
representations of the $C^*$-algebras \cite{Arv-book}, the
representation $\pi$ decomposes into a direct sum $\pi=\pi_0\oplus
\pi_1$ on  $ \cK=\tilde\cK_0\oplus \tilde\cK_1$, where
$$\tilde\cK_0:=\overline{\text{\rm span}}\{\pi(X)\cK:\ X\in \cC(\cN_J)\}
\quad \text{ and  }\quad  \tilde\cK_1:=\tilde\cK_0^\perp,
$$
and the  representations $\pi_j:C^*(B_1,\ldots, B_n)\to
\tilde\cK_j$, $j=1,2$, are defined by $\pi_j(X):=\pi(X)|\cK_j$, \
$j=0,1$. It is clear that $\pi_0$, $\pi_1$  are disjoint
representations of $C^*(B_1,\ldots, B_n)$
 such that
 $\pi_1$ annihilates  the compact operators in $B(\cN_J)$, and
  $\pi_0$ is uniquely determined by the action of $\pi$ on the
  ideal $\cC(\cN_J)$ of compact operators.
Since every representation of $\cC(\cN_J)$ is equivalent to a
multiple of the identity representation, we deduce \eqref{shi}.

 Now,
we show that the space $\tilde\cK_0$ coincides with the space
$\cK_0$. Indeed, using Theorem \ref{compact} and its proof, one can
easily see
that
\begin{equation*}\begin{split}
\tilde\cK_0&:=\overline{\text{\rm span}}\{\pi(X)\cK:\ X\in \cC(\cN_J)\}\\
&=\overline{\text{\rm span}}\{\pi(B_\beta P_\CC^{\cN_J} B_\alpha^*)\cK:\
 \alpha, \beta\in \FF_n^+\}\\
&= \overline{\text{\rm span}}\left\{\pi(B_\beta) \left(I_\cK-
\sum_{1\leq |\alpha|\leq m} a_\alpha
\pi(B_\alpha)\pi(B_\alpha)^*\right) \cK:\ \beta\in
\FF_n^+\right\}\\
&=\cK_0.
\end{split}
\end{equation*}
On the other hand, since $P_\CC^{\cN_J}=I- \sum\limits_{1\leq
|\alpha|\leq m} a_\alpha B_\alpha B_\alpha^*$ is a  projection of
rank one in $C^*(B_1,\ldots, B_n)$, we have
$$ \sum\limits_{1\leq
|\alpha|\leq m} a_\alpha \pi_1(B_\alpha)\pi_1(B_\alpha)^*=I_{\cK_1},
$$
and
$$
\dim \cG=\dim \left[\text{\rm range}\,\pi(P_\CC^{\cN_J})\right]
=\dim \left[\text{\rm range}\,\left(I_\cK- \sum_{1\leq |\alpha|\leq
m} a_\alpha \pi(B_\alpha)\pi(B_\alpha)^*\right)\right].
$$

To prove the uniqueness, one can use  standard theory of
representations of $C^*$-algebras. Consequently, $\pi$ and $\pi'$
are unitarily equivalent if and only if
 $\pi_0$ and $\pi_0'$ (resp.~$\pi_1$ and $\pi_1'$) are unitarily equivalent.
  By Proposition \ref{eq-mult}, we deduce that $\dim\cG=\dim\cG'$.
   The proof is
 complete.
\end{proof}

\begin{corollary}\label{Wold}
Assume the hypotheses and notations of Theorem $\ref{wold}$. Setting
$V_i:=\pi(B_i)$, \ $i=1,\ldots, n$,  and
$$\Phi_{p,
V}(Y):=\sum\limits_{1\leq |\alpha|\leq m} a_\alpha V_\alpha
YV_\alpha^*, \quad Y\in B(\cK),
$$
we have:
\begin{enumerate}
\item[(i)] $Q:=I_\cK-\Phi_{p,
V}(I_\cK)$ is an orthogonal projection and
$Q\cK=\bigcap\limits_{i=1}^n \ker V_i^*$;
\item[(ii)] $\cK_0=\left\{ x\in \cK:\ \lim\limits_{k\to \infty}
\left< \Phi^k_{p, V}(I_\cK)x,x\right>=0\right\}$;
\item[(iii)] $\cK_1=\left\{ x\in \cK:\ \left< \Phi^k_{p, V}(I_\cK)x,x\right>
=\|x\|^2 \ \text{ for any } k=1,2, \ldots\right\}$;
\item[(iv)] $\text{\rm SOT-}\lim\limits_{k\to\infty}  \Phi^k_{p, V}(I_\cK)=P_{\cK_1}$;
\item[(v)]
$\sum\limits_{k=0}^\infty  \Phi^k_{p, V}(Q) =P_{\cK_0}$.
\end{enumerate}
\end{corollary}

\begin{proof} Due to the proof of Theorem \ref{compact},
the operator  $I_{\cN_J}-  \Phi_{p, B}(I)=P_\CC^{\cN_J}$ is an
orthogonal projection. Therefore, so is $Q=\pi(P_\CC^{\cN_J})$.
Taking into account that $a_\alpha\geq 0$, $\alpha\in \FF_n^+$, and
$a_\alpha\neq 0$ if $|\alpha|=1$, we deduce that
\begin{equation*}
\begin{split}
Q\cK&=\left\{ x\in \cK:\ \left(I-\Phi_{p,
V}(I)\right)x=x\right\}\\
&=\left\{ x\in \cK:\ \sum_{1\leq |\alpha|\leq m} a_\alpha V_\alpha V_\alpha^*x=0\right\}\\
&=\bigcap_{i=1}^n \ker V_i^*.
\end{split}
\end{equation*}
Therefore (i) holds. By  Theorem \ref{wold}, we have
\begin{equation}\label{matr}
\Phi^k_{p, V}(I_\cK)= \left[
\begin{matrix}
\Phi_{p,
B}^k(I_{\cN_J})\otimes I_\cG &0\\
0&I_{\cK_1}
\end{matrix}\right],\quad k=1,2,\ldots.
\end{equation}
  Since
$\cN_J$ is co-invariant under each operator $W_i$, $i=1,\ldots,n$,
and $B_i^*=W_i^*|\cN_J$, $i=1,\ldots, n$, we  deduce that
$$
\text{\rm SOT-}\lim_{k\to \infty}\Phi^k_{p, B}(I_{\cN_J})\otimes
I_\cG= \text{\rm SOT-}\lim_{k\to \infty} \left[P_{\cN_J}\left(
\Phi^k_{p, W}(I_{F^2(H_n)})\right)|\cN_J\right]\otimes I_\cG=0.
$$
The latter equality holds due to the fact that $(W_1,\ldots, W_n)$
is   of class $C_{\cdot 0}$ with respect to the noncommutative
domain $\cD_f$. This proves part (iv). Hence, and taking into
account that
$$\sum_{k=1}^m  \Phi^k_{p,
V}(I_\cK-\Phi_{p,V}(I_\cK))=I- \Phi^{m+1}_{p, V}(I_\cK),
$$
we deduce item (v).
To prove (ii), let $x\in \cK=\cK_0\oplus \cK_1$,\ $x=x_0+ x_1$, with
$x_0\in \cK_0$ and $x_1\in \cK_1$. Using \eqref{matr}, we have
 \begin{equation}
\label{VVBB} \left<\Phi^m_{p, V}(I_\cK)x,x\right>=\left<\left(
 \Phi^m_{p,
B}(I_{\cN_J})\otimes I_\cG\right)x_0, x_0\right>+\|x_1\|^2, \quad
m=1,2,\ldots.
\end{equation}
Consequently, $\lim\limits_{m\to \infty} \left<\Phi^m_{p,
V}(I_\cK)x,x\right>=0 $ if and only if $x_1=0$, i.e., $x=x_0\in
\cK_0$.
 Now item (ii) follows.  It remains to prove (iii).
 Relation \eqref{VVBB} shows
that $\left<\Phi^m_{p, V}(I_\cK)x,x\right>=\|x\|^2$ for any
$m=1,2,\ldots$, if and only if
$$
\left<\left( \Phi^m_{p, B}(I_{\cN_J})\otimes I_\cG\right)x_0,
x_0\right>=\|x_0\|^2 \quad \text{
 for any } \ m=1,2,\ldots.
 $$
  Since
$(B_1,\ldots, B_n)$ is  of class $C_{\cdot 0}$, we have $\text{\rm
SOT-}\lim\limits_{m\to\infty}  \Phi^m_{p, B}(I_{\cN_J})=0$.
Therefore, the above  equality holds for any $m=1,2,\ldots$, if and
only if $x_0=0$, which is equivalent to $x=x_1\in \cK_1$. This
completes the proof.
\end{proof}

  Let $\cS\subseteq B(\cK)$ be a set of operators acting on the
Hilbert space $\cK$. We call  a subspace $\cH$ cyclic for $\cS$   if
$$\cK=\overline{\text{\rm span}}\{Xh:\ X\in \cS, \ h\in \cH\ \}.
$$

\begin{corollary}\label{co-shi}
Let $\pi$  be a nondegenerate $*$-representation of $C^*(B_1,\ldots,
B_n)$ on a separable Hilbert space $\cK$, and let $V_i:=\pi(B_i)$, \
$i=1,\ldots, n$. Then  the following statements are equivalent:
\begin{enumerate}
\item[(i)]
$V:=(V_1,\ldots, V_n)$ is a constrained  weighted shift;
 \item[(ii)]
$ \cK=\overline{\text{\rm span}}\left\{V_\beta \left(I- \Phi_{p,
V}(I)\right) \cK:\  \beta\in \FF_n^+\right\}; $
\item[(iii)] \text{\rm SOT-}$\lim\limits_{m\to\infty}
 \Phi^m_{p, V}(I)=0$.
\end{enumerate}
In this case,   the multiplicity of $V$
 satisfies the equality
\begin{equation}\label{multiplicity}
\text{\rm mult}\,(V)=\dim [(I- \Phi_{p, V}(I)) \cK],
\end{equation}
and  it is  also equal to the minimum dimension of a cyclic subspace
for $V_1,\ldots, V_n$.
\end{corollary}

\begin{proof}
It is clear that the statements above are equivalent due to
Corollary \ref{Wold}. On the other hand, relation
\eqref{multiplicity} follows from Theorem \ref{wold}.  Using (ii)
and Corollary \ref{Wold}, we deduce that the subspace
$$ \cL:=\bigcap_{i=1}^n \ker V_i^*=(I- \Phi_{p, V}(I)) \cK
$$
 is  cyclic  for $V_1,\ldots, V_n$.  Let $\cE$ be a
 cyclic subspace for $V_1,\ldots, V_n$, i.e.,
  $\cK=\bigvee\limits_{\alpha\in \FF_n^+} V_\alpha \cE$, and denote
   $A:=P_\cL|\cE\in B(\cE, \cL)$, where $P_\cL$ is the orthogonal projection
    of $\cK$ onto $\cL$.
Assume that $x\in \cL\ominus T\cE$ and let $y\in \cE$. Note that
$$
\left<x,y\right>=\left<x,P_\cL y\right>=\left<Ax,y\right>=0.
$$
On the other hand, since $V_\alpha^*x=0$ for all $x\in \cL$, we have
$\left<V_\alpha x, y\right>=0$ for any $\alpha\in \FF_n^+$ with
$|\alpha|\geq 1$. Therefore,  $y\perp V_\alpha\cE$ for
 $\alpha\in
\FF_n^+$. Since $\cE$ is a cyclic  subspace for $V_1,\ldots, V_n$,
we deduce that $y=0$. Therefore, $\overline {A\cE}=\cL$. Hence, the
operator  $A^*\in B(\cL, \cE)$ is one-to-one and, consequently, we
have $\dim\cL\leq \dim\cE$. This completes the proof.
\end{proof}

We can use now  Corollary \ref{co-shi} and Proposition \ref{eq-mult}
to obtain  the following result.
 \begin{proposition}
Two constrained weighted shifts  associated to the noncommutative
variety $\cV_{p,J}$ are similar if and only if they are unitarily
equivalent.
\end{proposition}
\begin{proof} One implication is clear.
Let $V:=(V_1,\ldots, V_n)$, $V_i\in B(\cK)$, and $V':=(V_1',\ldots,
V_n')$, $V_i'\in B(\cK')$, be two constrained weighted  shifts
associated with $\cV_{p,J}$,  and let $X:\cK\to \cK'$ be an
invertible operator such that
$$XV_i=V_i'X,\quad i=1,\ldots, n.
$$
If $\cM$ is a cyclic subspace for $V_1,\ldots, V_n$, then we obtain
\begin{equation*}
\begin{split}
\cK'=X\cK=X\left(\bigvee_{\alpha\in \FF_n^+} V_\alpha \cM\right)
\subseteq \bigvee_{\alpha\in \FF_n^+}X V_\alpha \cM=
\bigvee_{\alpha\in \FF_n^+} V_\alpha'X \cM\subseteq \cK'.
\end{split}
\end{equation*}
Hence, $\cK'=\bigvee_{\alpha\in \FF_n^+} V_\alpha'X \cM$.  This
shows that the subspace  $X\cM$ is cyclic for $V'$. Since $X$ is an
invertible operator, $\dim \cM=\dim X\cM$. Hence, and by Corollary
\ref{co-shi}, we deduce that  the two constrained weighted shifts
have the same multiplicity. Now, using  Proposition \ref{eq-mult},
the result follows.
\end{proof}

\bigskip

 \section{Dilations on noncommutative domains and varieties}
\label{Dilations}

In this section, we develop a dilation theory for $n$-tuples of
operators  in the noncommutative domain $\cD_f(\cH)$  or in the
noncommutative variety $\cV_{f,J}(\cH)$ defined by
$$
\cV_{f,J}(\cH):=\left\{ (X_1,\ldots, X_n)\in \cD_f(\cH): \
q(X_1,\ldots, X_n)=0, q\in \cP_J\right\},
$$
where $J$ is a $w^*$-closed two-sided ideal of $F_n^\infty(\cD_f)$
generated by a set of noncommutative polynomials $\cP_J$. Under
natural conditions on the ideal $J$ and  the $C^*$-algebra
$C^*(B_1,\ldots, B_n)$ associated with the variety $\cV_{f,J}$, we
have uniqueness for the minimal dilation.

Let $f:=\sum_{|\alpha|\geq1} a_\alpha X_\alpha$ be a positive
regular free holomorphic function on $B(\cH)^n$, and
 let $J\neq F_n^\infty(\cD_f)$ be a $w^*$-closed two-sided ideal of
$F_n^\infty(\cD_f)$ generated by a family of polynomials
$\cP_J\subset F_n^\infty(\cD_f)$. Fix an $n$-tuple of operators
$T:=(T_1,\ldots, T_n)$   in the noncommutative variety
$$
\cV_{f,J}(\cH):= \left\{(X_1,\ldots, X_n)\in \cD_f(\cH):\
q(X_1,\ldots, X_n)=0\quad \text{ for any }\quad q\in \cP_J\right\}.
$$
Let $\cD$ and $\cK$  be Hilbert spaces  and let $(U_1,\ldots,
U_n)\in \cD_f(\cK)$  be such that
$$\sum_{|\alpha|\geq 1} a_\alpha U_\alpha U_\alpha^*=I_\cK. $$
An $n$-tuple $V:=(V_1,\ldots, V_n)$ of operators  having the matrix
representation
\begin{equation}\label{Vi}
V_i:=\left[\begin{matrix}
B_i\otimes I_\cD&0\\
0&U_i
\end{matrix}
\right], \qquad i=1,\ldots, n,
\end{equation}
where the $n$-tuple $(B_1,\ldots, B_n)$ is the  universal model
associated with the noncommutative variety $\cV_{f,J}$,   is called
{\it constrained} (or $J$-{\it constrained}) {\it dilation} of $T\in
\cV_{f,J}$ if the following conditions hold:
\begin{enumerate}
\item[(i)] $(V_1,\ldots, V_n)\in \cV_{f,J}((\cN_J\otimes \cD)\oplus \cK)$;
\item[(ii)] $\cH$ can be identified with a co-invariant subspace
under each operator
  $V_i$, $i=1,\ldots,n$,  such that
$$
T_i=P_\cH V_i|\cH,\quad i=1,\ldots, n.
$$
\end{enumerate}
The dilation is called {\it minimal} if
$$
(\cN_J\otimes \cD)\oplus \cK=\overline{\text{\rm span}}\left\{
V_\alpha \cH: \  \alpha\in \FF_n^+\right\}.
$$
The {\it dilation index} of $T$,  denoted by $\text{\rm
dil-ind}\,(T)$, is the minimum dimension of $\cD$ such that $V$ is a
constrained   dilation of $T$.

Our first dilation result   on noncommutative varieties
$\cV_{f,J}\subset \cD_f$  is the following.

\begin{theorem}\label{dil1} Let $f$ be a positive regular free
holomorphic function on $B(\cH)^n$ and  let $J\neq
F_n^\infty(\cD_f)$ be a $w^*$-closed two-sided ideal of \
$F_n^\infty(\cD_f)$ generated by a family of polynomials
$\cP_J\subset F_n^\infty(\cD_f)$. If $T:=(T_1,\ldots, T_n)$  is an
$n$-tuple of operators in the noncommutative variety
$\cV_{f,J}(\cH)$, then there exists a Hilbert space $\cK$  and and
$n$-tuple $(U_1\ldots, U_n)\in \cV_{f,J}(\cK)$ with
$$ \sum_{|\alpha|\geq 1} a_\alpha U_\alpha U_\alpha^*=I_\cK $$
and such that
\begin{enumerate}
\item[(i)]
 $\cH$ can be identified with a  co-invariant subspace of
   $\tilde\cK:=(\cN_J\otimes \overline{\Delta_{f,T}\cH})\oplus \cK$ under the operators
$$
V_i:=\left[\begin{matrix} B_i\otimes
I_{\overline{\Delta_{f,T}\cH}}&0\\0&U_i
\end{matrix}\right],\quad i=1,\ldots,n,
$$
where $\Delta_{f,T}:=\left(I_\cH-\sum\limits_{|\alpha|\geq 1}
a_\alpha T_\alpha T_\alpha^*\right)^{1/2};$
\item[(ii)]
$T_i^*=V_i^*|\cH$,\ $i=1,\ldots,n$.
\end{enumerate}
Moreover, $\cK=\{0\}$ if and only if $(T_1,\ldots, T_n)$ is  of
class $C_{\cdot 0}$ with respect to $\cD_f(\cH)$.
\end{theorem}

\begin{proof}
Define  the subspace $\cM\subset F^2(H_n)$ be setting
$$
\cM:=\overline{\text{\rm span}}\left\{W_\alpha q(W_1,\ldots, W_n)
W_\beta(1):\ q\in \cP_J, \alpha,\beta\in \FF_n^+\right\}.
$$
We prove that $\cM=\cM_J$.  Since  $\cM\subseteq \cM_J$, it remains
to prove that $\cM_J\subseteq \cM$. To this end,  it is enough to
show that $\cM^\perp\subseteq \cM_J^\perp$. Let $g\in F^2(H_n)$  be
such that
$$\left< g, W_\alpha p(W_1,\ldots, W_n)
W_\beta(1)\right>=0\quad \text{ for any }
\ p\in \cP_J, \alpha,\beta\in \FF_n^+.
$$
 Due to  Corollary \ref{w*-density}, ,
  for any $\varphi(W_1,\ldots, W_n)\in F_n^\infty(\cD_f)$,
 there is a sequence of
polynomials  such that $\varphi(W_1,\ldots, W_n)=\text{\rm
SOT}-\lim_{m\to\infty}\{q_m(W_1,\ldots, W_n)\}_{m=1}^\infty$. Hence,
$$\left< g,\varphi(W_1,\ldots, W_n) q(W_1,\ldots, W_n) W_\beta(1)\right>=0
$$ for any  $\varphi(W_1,\ldots, W_n)\in F_n^\infty (\cD_f)$, $q\in \cP_J$,
 and $\alpha,
\beta\in \FF_n^+.$  Consequently, $g\in \cM_J^\perp$ and, therefore,
$\cM_J=\cM$.

As in the proof of Theorem \ref{funct-calc2}, using the properties
of the Poisson kernel $K_{f,T}$, we obtain
$$
\left<  K_{f,T}x, W_\alpha q(W_1,\ldots, W_n) W_\beta(1)\otimes
y\right>=\left<x,T_\alpha q(T_1,\ldots,T_n)T_\beta
\Delta_{f,T}y\right>=0
$$
for any $x\in \cH$, $y\in \overline{\Delta_{f,T}\cH}$, and $q\in
\cP_J$. Since $\cM_J=\cM$, we deduce  that
\begin{equation}\label{range}
\text{\rm range}\,K_{f,T}\subseteq \cN_J\otimes
\overline{\Delta_{f,T}\cH}.
\end{equation}
The {\it constrained Poisson kernel} $K_{f,T, J}:\cH\to \cN_J\otimes
\overline{\Delta_{f,T} \cH}$  is defined by
$$K_{f,T,J}:=(P_{\cN_J}\otimes I_{\overline{\Delta_{f,T} \cH}})K_{f,T},
$$
where $K_{f,T}$ is the Poisson kernel associated with $T\in
\cD_f(\cH)$ (see Section \ref{domain algebra}). Due to  relation
\eqref{range} and  \eqref{ker-inter}, we obtain
\begin{equation}\label{KJ}
K_{f,T,J}T_\alpha^*=(B_\alpha^*\otimes I_\cH)K_{f,T,J},\quad
\alpha\in \FF_n^+.
\end{equation}
We introduce the operators   $Q:=\text{\rm
SOT-}\lim\limits_{m\to\infty} \Phi_{f,T}^m(I)$ and
$$Y:\cH\to \cK:=\overline{Q^{1/2}\cH}\quad \text{ defined by }\quad
Yh:=Q^{1/2}h, \ h\in \cH.
$$
For each $i=1,\ldots, n$, define the operator $L_i:Q^{1/2}\cH\to
\cK$ by setting
\begin{equation}\label{Z_i}
L_i Yh:=YT_i^*h,\quad h\in \cH.
\end{equation}
Note that $L_i$, $i=1,\ldots, n$, are well-defined  due to the fact
that
\begin{equation*}
\begin{split}
\|L_i Yh\|^2&=\left<   T_i QT_i^*h,h \right>\\
&=\frac{1}{a_{g_i}} \left< a_{g_i} T_iQ T_i^*h,h\right>\leq
 \frac{1}{a_{g_i}}\left< \Phi_{f,T}(Q)h,h\right>\\
 &=\frac{1}{a_{g_i}} \|Q^{1/2}h\|^2=\frac{1}{a_{g_i}}\|Yh\|^2.
\end{split}
\end{equation*}
We recall that, since $f$ is positive regular free holomorphic
function,  $a_{g_i}\neq 0$ for any $i=1,\ldots,n$.
 Consequently,  $L_i$ can be extended to a bounded operator on $\cK$, which will
also be denoted by $L_i$. Now, setting $U_i:=L_i^*$, $i=1,\ldots,
n$, relation  \eqref{Z_i} implies
\begin{equation}
\label{YZT} Y^*U_i=T_iY^*,\quad i=1,\ldots, n.
\end{equation}
Due to relation \eqref{YZT}, we have
\begin{equation*}
\begin{split}
Y^* \Phi_{f,U}(I_\cK)Y&=Y^*\left(\sum_{|\alpha|\geq 1} a_\alpha
U_\alpha U_\alpha^*\right)Y
=\sum_{|\alpha|\geq 1} a_\alpha T_\alpha YY^* T_\alpha\\
&=\Phi_{f,T}(YY^*)=\Phi_{f,T}(Q)=Q=YY^*.
\end{split}
\end{equation*}
Hence,
$$\left< \Phi_{f,U}(I_\cK)Yh,Yh\right>=
\left< Yh,Yh\right>,\quad h\in \cH
$$
which implies $\Phi_{f,U}(I_\cK)=I_\cK$.
 Now, using relation \eqref{YZT}, we obtain
$$Y^*q(U_1,\ldots, U_n)=q(T_1,\ldots, T_n)Y^*=0, \quad q\in\cP_J.
$$
Since $Y^*$ is injective on $\cK=\overline{Y\cH}$, we have
$q(U_1,\ldots, U_n)=0$ for any $q\in\cP_J$.
Let $V:\cH\to [\cN_J\otimes \cH]\oplus \cK$ be defined by
$$V:=\left[\begin{matrix}
K_{f,T,J}\\ Y
\end{matrix}\right].
$$
Notice  that $V$ is an isometry. Indeed, we have
\begin{equation*}\begin{split}
\|Vh\|^2&=\|K_{f,T,J}h\|^2+\|Yh\|^2\\
&=\|h\|^2-\text{\rm SOT-}\lim_{k\to\infty}\left< \Phi_{f,T}^k(I)h,h\right>+\|Yh\|^2\\
&=\|h\|^2-\left<Qh,h\right> + \left<Qh,h\right>\\
&=\|h\|^2
\end{split}
\end{equation*}
for any $h\in \cH$. Now, using relations \eqref{KJ} and \eqref{Z_i},
we obtain
\begin{equation*}
\begin{split}
VT_i^*&=\left[\begin{matrix} K_{f,T,J}\\ Y
\end{matrix}\right]T_i^*h=K_{f,T,J}T_i^*h\oplus YT_i^*h\\
&=(B_i^*\otimes I_\cH)K_{f,T,J}h\oplus U_i^*Yh\\
&=\left[\begin{matrix} B_i^*\otimes
I_{\overline{\Delta_{f,T}\cH}}&0\\0&U_i^*
\end{matrix}\right]Vh.
\end{split}
\end{equation*}
 Identifying   $\cH$  with $V\cH$ we
complete the proof of (i) and (ii). The  last part  of the theorem
is obvious.
\end{proof}

\begin{example}
\begin{enumerate}
\item[(i)]
When  $f$ is a positive regular noncommutative polynomial and the
two-sided ideal $J$ is generated by the polynomials
$$
W_iW_j-W_jW_i,\qquad  i,j=1,\ldots, n, $$ we obtain the dilation
theorem of  Pott \cite{Pot}.
\item[(ii)]
When $f=X_1+\cdots +X_n$ and $\cP_J=0$ we obtain a version of the
noncommutative dilation theorem for  row contractions (see
\cite{Po-isometric}).
 \item[(ii)] In the particular case when $n=1$,
$f=X$, and $\cP_J=0$, we obtain the classical isometric dilation
theorem for contractions obtained by Sz.-Nagy $($see \cite{Sz1},
\cite{SzF-book}$)$.
\end{enumerate}
\end{example}

We can evaluate the dilation index of an    $n$-tuple of operators
in the noncommutative variety $\cV_{f,J}$ and show that it does not
depend on the constraints.

\begin{corollary}\label{dil-ind}
Let $J\neq F_n^\infty(\cD_f)$ be a $w^*$-closed two-sided ideal of
$F_n^\infty(\cD_f)$ generated by a family of polynomials
$\cP_J\subset F_n^\infty(\cD_f)$ and let $T:=(T_1,\ldots, T_n)\in
\cV_{f,J}(\cH)$. Then the dilation index of $T$   is equal to $\rank
\Delta_{f,T}$.
\end{corollary}

\begin{proof}
Let
  $(U_1,\ldots, U_n)\in \cD_f(\cH)$ be  such that
$$ \sum_{|\alpha|\geq 1} a_\alpha U_\alpha U_\alpha^*=I_\cK $$
 and assume that   the $n$-tuple $V:=(V_1,\ldots, V_n)$ given by
\begin{equation}\label{ViBi}
V_i:=\left[\begin{matrix}
B_i\otimes I_\cD&0\\
0&U_i
\end{matrix}
\right], \qquad i=1,\ldots, n,
\end{equation}
is a  constrained dilation of $T$. Since the subspace $\cH$ is
co-invariant under each operator  $V_i$, $i=1,\ldots, n$, and
$\cN_J$ is co-invariant under the weighted left creation operators
$W_i$, $i=1,\ldots, n$,  associated with the noncommutative domain
$\cD_f$,  we have
$$
I_\cH- \sum_{|\alpha|\geq 1} a_\alpha T_\alpha
T_\alpha^*=P_\cH\left[
\begin{matrix}
\left[P_{\cN_J}\left( I-\sum\limits_{|\alpha|\geq 1} a_\alpha
W_\alpha W_\alpha^*\right)
|\cN_J\right]\otimes I_\cD& 0\\
0&0
\end{matrix}\right]|\cH.
$$
Hence,  and taking into account that  $I- \sum\limits_{|\alpha|\geq
1} a_\alpha W_\alpha W_\alpha^*= P_\CC$  is an operator of  rank
one, we have
\begin{equation*}
\rank \Delta_{f,T}\leq \rank \left[P_{\cN_J}
 P_\CC|\cN_J\otimes I_\cD\right]
\leq \dim \cD.
\end{equation*}
Due to  Theorem \ref{dil1}, we deduce that the dilation index of $T$
coincides with $\rank \Delta_{f,T}.$ The proof is complete.
\end{proof}

  Let $C^*(\cY)$ be the $C^*$-algebra generated by a set of
  operators
$\cY\subset B(\cK)$ and the identity. A  subspace $\cH\subseteq \cK$
is called $*$-cyclic  for $\cY$    if
$$\cK=\overline{\text{\rm span}}\left\{Xh:\ X\in C^*(\cY), \ h\in
\cH\right\}.
$$

\begin{theorem}\label{dil2} Let
 $p=\sum\limits_{1\leq |\alpha|\leq
m} a_\alpha X_\alpha$ be  a positive regular noncommutative
polynomial and let $J\neq F_n^\infty(\cD_p)$   be a $w^*$-closed
two-sided ideal of $F_n^\infty(\cD_p)$ generated by a family
 $\cP_J$  of homogenous polynomials. Let $\cH$ be a separable Hilbert space, and
 $T:=(T_1,\ldots, T_n) $  be an $n$-tuple of operators  in the noncommutative variety
 $\cV_{p,J}(\cH)$.
Then there exists  a $*$-representation $\pi:C^*(B_1,\ldots, B_n)\to
B(\cK_\pi)$  on a separable Hilbert space $\cK_\pi$,  which
annihilates the compact operators and
$$
\sum_{1\leq |\alpha|\leq m} a_\alpha \pi(B_\alpha)
\pi(B_\alpha)^*=I_{ \cK_\pi},
$$
such that
\begin{enumerate}
\item[(i)]
$\cH$ can be identified with a $*$-cyclic co-invariant subspace of
$\tilde\cK:=(\cN_J\otimes \overline{\Delta_{p,T}\cH})\oplus \cK_\pi$
under  each operator
$$
V_i:=\left[\begin{matrix} B_i\otimes
I_{\overline{\Delta_{p,T}\cH}}&0\\0&\pi(B_i)
\end{matrix}\right],\quad i=1,\ldots,n,
$$
where $\Delta_{p,T}:=\left(I_\cH-\sum\limits_{1\leq |\alpha|\leq
m}a_\alpha T_\alpha T_\alpha^*\right)^{1/2};$
\item[(ii)]
$ T_i^*=V_i^*|\cH,\quad i=1,\ldots, n. $
\end{enumerate}
  \end{theorem}
\begin{proof}
%\bigskip
According to  Theorem \ref{poiss-poly},   there is a unique unital
completely contractive linear map
$$\Psi_{p,T,J}:\overline{\text{\rm span}}\{B_\alpha B_\beta^*:\
 \alpha,\beta\in \FF_n^+\}\to B(\cH)
$$
such that $\Psi_{p,T,J}(B_\alpha B_\beta^*)=T_\alpha T_\beta^*$, \
$\alpha,\beta\in \FF_n^+$. Applying Arveson extension theorem
\cite{Arv-acta} to the map $\Psi_{p,T,J}$, we find a unital
completely positive linear map $\Psi_{p,T,J}:C^*(B_1,\ldots, B_n)\to
B(\cH)$. Let $\tilde\pi:C^*(B_1,\ldots, B_n)\to B(\tilde\cK)$ be a
minimal Stinespring dilation \cite{St} of $\Psi_{p,T,J}$. Then
$$\Psi_{p,T,J}(X)=P_{\cH} \tilde\pi(X)|\cH,\quad X\in C^*(B_1,\ldots, B_n),
$$
and $\tilde\cK=\overline{\text{\rm span}}\{\tilde\pi(X)h:\ h\in
\cH\}.$  Since , for each $i=1,\ldots, n$,
\begin{equation*}
\begin{split}
\Psi_{p,T,J}(B_i B_i^*)&=T_iT_i^*=P_\cH\tilde\pi(B_i)\tilde\pi(B_i^*)|\cH\\
&=P_\cH \tilde\pi(B_i)(P_\cH+P_{\cH^\perp})\tilde\pi(B_i^*)|\cH\\
&=\Psi_{p,T,J}(B_i B_i^*)+ (P_\cH
\tilde\pi(B_i)|_{\cH^\perp})(P_{\cH^\perp} \tilde\pi(B_i^*)|\cH),
\end{split}
\end{equation*}
 we deduce that $P_\cH \tilde\pi(B_i)|_{\cH^\perp}=0$ and
\begin{equation}\label{morph}
\begin{split}
\Psi_{p,T,J}(B_\alpha X)&=P_\cH(\tilde\pi(B_\alpha) \tilde\pi(X))|\cH\\
&=(P_\cH\tilde\pi(B_\alpha)|\cH)(P_\cH\tilde\pi(X)
|\cH)\\
&=\Psi_{p,T,J}(B_\alpha) \Psi_{p,T,J}(X)\end{split}
\end{equation}
for any $X\in C^*(B_1,\ldots, B_n)$ and $\alpha\in \FF_n^+$. Since
 $P_\cH \tilde\pi(B_i)|_{\cH^\perp}=0$,
  $\cH$ is an
invariant subspace under each $\tilde\pi(B_i)^*$, \ $i=1,\ldots, n$.
 Consequently, we have
  \begin{equation}\label{coiso}
\tilde\pi(B_i)^*|\cH=\Psi_{p,T,J}(B_i^*)=T_i^*,\quad i=1,\ldots, n.
\end{equation}

On the other hand, since  $1\in \cN_J$,   Theorem \ref{compact}
implies  that all the compact operators  in $B(\cN_J)$ are contained
in the $C^*$-algebra $C^*(B_1,\ldots, B_n)$. By Theorem \ref{wold},
the representation $\tilde\pi$ decomposes into a direct sum
$\tilde\pi=\pi_0\oplus \pi$ on $\tilde \cK=\cK_0\oplus \cK_\pi$,
where $\pi_0$, $\pi$  are disjoint representations of
$C^*(B_1,\ldots, B_n)$ on the Hilbert spaces $\cK_0$ and $\cK_\pi$,
respectively, such that
\begin{equation}\label{sime}
\cK_0\simeq\cN_J\otimes \cG, \quad  \pi_0(X)=X\otimes I_\cG, \quad
X\in C^*(B_1,\ldots, B_n),
\end{equation}
 for some Hilbert space $\cG$, and $\pi$ is a representation such that
$\pi(\cC(\cN_J))=0$. Now, taking into account that
 $P_\CC^{\cN_J}=I-
\sum\limits_{|\alpha|\geq 1} a_\alpha  B_\alpha B_\alpha^*$ is a
 projection  of rank one in $C^*(B_1,\ldots, B_n)$, we  have that
$$
\sum\limits_{|\alpha|\geq 1} a_\alpha
\pi(B_\alpha)\pi(B_\alpha^*)=I_{\cK_\pi}\quad \text{ and }\quad
\dim \cG=\dim (\text{\rm range}\,\tilde\pi(P_\CC^{\cN_J})).
$$
Since the  Stinespring representation $\tilde\pi$  is minimal and
using the proof of Theorem \ref{compact}, we deduce that
\begin{equation*}\begin{split}
\text{\rm range}\,\tilde\pi(P_\CC^{\cN_J})&=
\overline{\text{\rm span}}\{\tilde\pi(P_\CC^{\cN_J})\tilde\pi(X)h:\ X\in C^*(B_1,\ldots, B_n), h\in \cH\}\\
&=
\overline{\text{\rm span}}\{\tilde\pi(P_\CC^{\cN_J})\tilde\pi(Y)h:\ Y\in \cC(\cN_J), h\in \cH\}\\
&=
\overline{\text{\rm span}}\{\tilde\pi(P_\CC^{\cN_J})\tilde\pi(B_\alpha P_\CC^{\cN_J} B_\beta^*)h:\ \alpha,\beta\in \FF_n^+, h\in \cH\}\\
&= \overline{\text{\rm
span}}\{\tilde\pi(P_\CC^{\cN_J})\tilde\pi(B_\beta^*)h:\ \beta\in
\FF_n^+, h\in \cH\}.
\end{split}
\end{equation*}
Now, due to  relation \eqref{morph}, we have
\begin{equation*}\begin{split}
\left<\tilde\pi(P_\CC^{\cN_J})\tilde\pi(B_\alpha^*)h,
\tilde\pi(P_\CC^{\cN_J})\tilde\pi(B_\beta^*)k\right>
&=
\left<h,\pi(B_\alpha)\pi(P_\CC^{\cN_J})\pi(B_\beta^*)h\right>\\
&= \left<h,T_\alpha\left(I_\cK-\sum_{ |\gamma|\geq 1} a_\gamma
T_\gamma T_\gamma^*
\right)T_\beta^*h\right>\\
&= \left<\Delta_{p,T}T_\alpha^*h,\Delta_{p,T}T_\beta^*k\right>
\end{split}
\end{equation*}
for any $h, k \in \cH$. Consequently, there is   a unitary operator
$\Lambda:\text{\rm range}\,\tilde\pi(P_\CC^{\cN_J})\to
\overline{\Delta_{p,T}\cH}$  defined by
$$
\Lambda(\tilde\pi(P_\CC^{\cN_J})\tilde\pi(B_\alpha^*)h):=\Delta_{p,T}
T_\alpha^*h,\quad h\in \cH, \,\alpha\in \FF_n^+.
$$
 This shows that
$$
\dim[\text{\rm range}\,\pi(P_\CC^{\cN_J})]= \dim
\overline{\Delta_{p,T}\cH}=\dim \cG.
$$
The required dilation is obtained using relations \eqref{coiso} and
\eqref{sime},  and identifying    $\cG$ with
$\overline{\Delta_{p,T}\cH}$.
 The proof is complete.
\end{proof}

\begin{corollary} \label{part-cases}
Let $V:=(V_1,\ldots, V_n)$ be the dilation
 of Theorem $\ref{dil2}$. Then,
\begin{enumerate}
\item[(i)]
 $V$ is a constrained  weighted shift if and only if
 $T=(T_1,\ldots,T_n)\in \cV_{p,J}(\cH)$
 is  of class $C_{\cdot 0}$ with respect to $\cD_p$;
\item[(ii)]
  $
 \sum_{1\leq |\alpha|\leq m}a_\alpha V_\alpha V_\alpha^*=I
$ if and only if \ $
 \sum_{1\leq |\alpha|\leq m}a_\alpha T_\alpha T_\alpha^*=I.
$
\end{enumerate}
\end{corollary}
\begin{proof}
Since
$$
\Phi_{p,T}^k(I_\cH)=P_\cH \left[\begin{matrix}
\Phi_{p,B}^k(I_{\cN_J})\otimes I_{\overline{\Delta_{p,T}\cH}}&0\\0&
I_{\cK_\pi}\end{matrix}\right]|\cH
$$
we have
$$
\text{\rm SOT-} \Phi_{p,T}^k(I_\cH)=P_\cH \left[\begin{matrix}
 0&0\\0& I_{\cK_\pi}\end{matrix}\right]|\cH.
$$
Consequently,  $T$ is  of class $C_{\cdot 0}$ with respect to
$\cD_p$ if and only if $P_\cH P_{\cK_\pi} |\cH=0$. The latter
condition is equivalent to $\cH\perp (0\oplus \cK_\pi)$, which
implies $\cH\subset \cN_J\otimes\overline{\Delta_{p,T}\cH}$. On the
other hand, since $\cN_J\otimes\overline{\Delta_{p,T}\cH}$ is
reducing for  $V_1,\ldots, V_n$, and $\tilde \cK$ is the smallest
reducing subspace for   $V_1,\ldots, V_n$, which contains $\cH$, we
must have $\tilde\cK=\cN_J\otimes\overline{\Delta_{p,T}\cH}$.
Therefore, (i) holds.

Assume  now that   $\sum\limits_{|\alpha|\geq 1}^n a_\alpha V_\alpha
V_\alpha^*=I_{\tilde\cK}$. Since
$$ \Phi_{p,V}^k(I_\cK)=
\left[\begin{matrix}  \Phi_{p,B}^k(I_{\cN_J})\otimes
I_{\overline{\Delta_{p,T}\cH}}&0\\0& I_{\cK_\pi}\end{matrix}\right],
$$
we must have $\Phi_{p,B}^k(I_{\cN_J})\otimes
I_{\overline{\Delta_{p,T}\cH}}=I_{\cK_0}$ for any $k=1,2,\ldots$. On
the other hand, since $\text{\rm SOT-}\lim\limits_{k\to\infty}
 \Phi_{p,B}^k(I_{\cN_J})=0$, we deduce that
$\cK_0=\{0\}$. Now,  the proof of Theorem \ref{dil2} implies
$\cG=\{0\}$. Therefore, $\Delta_{p,T}=0$ and the  proof is complete.
\end{proof}

The dilation is unique up to a unitary equivalence, under additional
hypotheses on the $C^*$-algebra  $C^*(B_1,\ldots, B_n)$.

\begin{corollary}\label{unique}
Assume  the hypotheses of Theorem $\ref{dil2}$ and that
\begin{equation}\label{BB^*}
\overline{\text{\rm span}}\,\{B_\alpha B_\beta^*:\ \alpha,\beta\in
\FF_n^+\}=C^*(B_1,\ldots, B_n).
\end{equation}
Then the  dilation of Theorem $\ref{dil2}$ is minimal, i.e.,
$\tilde\cK=\bigvee\limits_{\alpha\in \FF_n^+} V_\alpha \cH$,   and
it is unique up to a unitary equivalence.
\end{corollary}

\begin{proof} If condition \eqref{BB^*} holds, then the map $\Psi_{f,T,J}$
in the proof of  Theorem \ref{dil2} is unique. The uniqueness of the
minimal Stinespring representation  \cite{St} implies  the
uniqueness of the minimal dilation of Theorem \ref{dil1}.
\end{proof}

Using Theorem \ref{cyclic} and  standard arguments  concerning
representation theory of $C^*$-algebras \cite{Arv-book}, one can
deduce the following result.

\begin{remark} In addition to  the hypotheses of Corollary
\ref{unique}, let $T':=(T_1',\ldots, T_n')\in \cV_{p,J}(\cH')$   and
let $V':=(V_1',\ldots, V_n')$ be the corresponding constrained
dilation. Then the $n$-tuples $T$ and $T'$ are unitarily equivalent
if and only if
$$\dim \overline{\Delta_{p,T}\cH}=\dim \overline{\Delta_{p,T'}\cH'}
$$
and there are unitary operators
 $U:\overline{\Delta_{p,T}\cH}\to \overline{\Delta_{p,T'}\cH'}$
  and $\Gamma:\cK_\pi\to \cK_{\pi '}$
such that $ \Gamma\pi(B_i)=\pi'(B_i) \Gamma\quad \text{ for } \quad
i=1,\ldots, n$, and
$$\left[ \begin{matrix} I_{\cN_J}\otimes U&0\\
0&\Gamma\end{matrix} \right] \cH=\cH'.
$$
\end{remark}

Let  $J$ be a two-sided ideal of  the Hardy algebra
$F_n^\infty(\cD_f$ generated by a set  $\cP_J$ of polynomials and
consider the noncommutative variety $\cV_{f,J}$. For any $n$-tuple
of operators $T:=(T_1,\ldots, T_n)\in \cD_f(\cH)$,  let
 $\cG_J$
be the largest subspace of $\cH$  which is invariant  under each
operator  $T_i^*$, $i=1,\ldots,n$,  such that $$(P_{\cG_J}T_1
|\cG_J,\ldots, P_{\cG_J}T_n|\cG_J)\in \cV_{f,J}.
$$
 We call the  $n$-tuple
$$(P_{\cG_J}T_1 |\cG_J,\ldots, P_{\cG_J}T_n|\cG_J)
$$
the {\it maximal  piece} of $(T_1,\ldots, T_n)$ in the
noncommutative variety $\cV_{f,J}$. It is easy to see that
$$
\cG_J= \overline{\text{\rm span}}\left\{ T_\alpha q(T_1,\ldots,
T_n)\cH:\
 q\in \cP_J, \alpha\in \FF_n^+\right\}^\perp.
$$
If  $T:=(T_1,\ldots, T_n)\in \cD_f(\cH)$ is c.n.c., then the
definition above makes sense when $J$ is  any $w^*$-closed two-sided
ideal of $F_n^\infty(\cD_f)$.

Using  Theorem \ref{symm-Fock} and its proof,  one can easily deduce
the following result.

\begin{proposition}\label{max-piece}
Let $J\neq F_n^\infty(\cD_f)$ be an arbitrary $w^*$-closed two-sided
ideal of $F_n^\infty(\cD_f)$
 and let $W_1,\ldots, W_n$ be the weighted left creation operators
 associated with the noncommutative domain $\cD_f$.
Then the universal model
  $(B_1,\ldots, B_n)$, associated with $\cV_{f,J}$,  is the maximal piece of $(W_1,\ldots, W_n)$
  in the noncommutative  variety $\cV_{f,J}$.
\end{proposition}

Now, we turn our attention to  the particular  case  when
$T:=(T_1,\ldots, T_n)\in \cD_f(\cH)$ is
 of class $C_{\cdot 0}$ with respect to $\cD_f(\cH)$. In  this case, the
results of this section can be extended to a larger class of
 noncommutative varieties.

\begin{theorem}\label{dil3}
Let $f$ be a positive regular free holomorphic function on
$B(\cH)^n$ and  let $J\neq F_n^\infty(\cD_f)$ be a $w^*$-closed
two-sided ideal of $F_n^\infty(\cD_{f})$. If $(T_1,\ldots, T_n)\in
\cV_{f,J}(\cH)$ is of class $C_{\cdot 0}$, then the following
statements hold:
\begin{enumerate}
\item[(i)]
The constrained Poisson kernel $K_{f,T,J}:\cH\to \cN_J\otimes
\overline{\Delta_{f,T}\cH}$ defined by setting
$$K_{f,T,J}:=(P_{\cN_J}\otimes I)K_{f,T}$$
 is an isometry, $K_{f,T,J}\cH$ is co-invariant under each operator
$B_i\otimes  I_{\overline{\Delta_{f,T}\cH}}$, $i=1,\dots, n$,   and
\begin{equation}\label{Pois}
T_i=K_{f,T,J}^*(B_i\otimes
I_{\overline{\Delta_{f,T}\cH}})K_{f,T,J},\quad i=1,\ldots, n.
\end{equation}
\item[(ii)]  The dilation index of $T$ coincides with $\text{\rm rank}\,\Delta_{f,T}$.
\item[(iii)]
If $1\in \cN_J$, then the $J$-constrained dilation $(B_1\otimes
I_{\overline{\Delta_{f,T}\cH}},\ldots, B_n\otimes
I_{\overline{\Delta_{f,T}\cH}})$, provided by \eqref{Pois}, is
minimal.
\end{enumerate}
If,  in addition, $f=p$ is a polynomial,  $1\in \cN_J$, and
\begin{equation}\label{C*}
\overline{\text{\rm span}}\{B_\alpha B_\beta^*:\ \alpha,\beta\in
\FF_n^+\}= C^*(B_1,\ldots, B_n),
\end{equation}
then  the following statements hold:
\begin{enumerate}
\item[(iv)] The minimal $J$-constrained dilation $(B_1\otimes
I_{\overline{\Delta_{p,T}\cH}},\ldots, B_n\otimes
I_{\overline{\Delta_{p,T}\cH}})$
 is unique up to an isomorphism.
\item[(v)] The minimal $J$-constrained dilation $(B_1\otimes
I_{\overline{\Delta_{p,T}\cH}},\ldots, B_n\otimes
I_{\overline{\Delta_{p,T}\cH}})$    is the maximal piece of the
   dilation of $T$ (provided by  Theorem \ref{dil2} when $J:=\{0\}$), i.e.,
   $(W_1\otimes
I_{\overline{\Delta_{p,T}\cH}},\ldots, W_n\otimes
I_{\overline{\Delta_{p,T}\cH}})$, in the noncommutative variety
$\cV_{p,J}$.
\end{enumerate}
\end{theorem}
 \begin{proof}
Part (i)  was proved in Section \ref{Functional II}. To prove (ii),
let $\cD$ be a Hilbert space such that $\cH$ can be identified with
a co-invariant subspace of $\cN_J\otimes \cD$ under each operator
$B_i\otimes I_\cD$, \ $i=1,\ldots, n$, and such that
 $T_i=P_\cH(B_i\otimes I_\cD)|\cH$ for  $i=1,\ldots, n$.
Then, using relation \eqref{proj}, we have
\begin{equation*}
\begin{split}
I_\cH-\sum_{|\alpha|\geq 1} a_\alpha T_\alpha T_\alpha^*&=
P_\cH^{\cN_J\otimes \cD}\left[ \left(I_{\cN_J}-\sum_{|\alpha|\geq 1}
a_\alpha B_\alpha B_\alpha^*\right)\otimes I_\cD
\right]|\cH\\
&= P_\cH^{\cN_J\otimes \cD}\left[ P_{\cN_J\otimes
\cD}(I-\sum_{|\alpha|\geq 1} a_\alpha W_\alpha
W_\alpha^*)|\cN_J\otimes I_\cD
\right]|\cH\\
&= P_\cH^{\cN_J\otimes \cD}\left[ P_{\cN_J}P_\CC|\cN_J\otimes I_\cD
\right]|\cH.
\end{split}
\end{equation*}
Hence, $\rank \Delta_{f,T}\leq \dim \cD$. Now, part (i) implies
that the dilation index of $T$ is equal to $\rank \Delta_{f,T}$.

 To prove (iii), assume that $1\in \cN_J$.
 Due to relation \eqref{WbWb}, $e_\alpha=\sqrt{b_\alpha} W_\alpha
 (1)$ and therefore
  $P_\CC^{\cN_J} P_{\cN_J}e_\alpha=0$ for any
$\alpha\in \FF_n^+$, $|\alpha|\geq 1$. On the other hand, the
definition of the constrained Poisson kernel $K_{f,T,J}$  implies
$$
P_0K_{f,T,J}h=\lim_{m\to \infty}\sum_{k=0}^m \sum_{|\alpha|=k}
P_\CC^{\cN_J} P_{\cN_J} e_\alpha\otimes \Delta_{f,T}
T_\alpha^*h,\quad h\in \cH,
$$
where $P_0:=P^{\cN_J}_\CC \otimes I_{\overline{\Delta_{f,T} \cH}}$.
Therefore, $P_0K_{f,T,J}\cH=\overline{\Delta_{f,T} \cH}$.  Using now
Theorem \ref{cyclic} in the particular case when $\cM:=K_{J,T}\cH$
and $\cD:=\overline{\Delta_{f,T} \cH}$, we deduce that the subspace
$K_{f,T,J}\cH$ is cyclic for  $B_i\otimes I_\cD$,\ $i=1,\ldots, n$,
which proves the minimality of the $J$-dilation provided by
\eqref{Pois}, i.e.,
\begin{equation}\label{minimal1}
\cN_J\otimes \overline{\Delta_{f,T} \cH}=\bigvee_{\alpha\in \FF_n^+}
(B_\alpha\otimes I_{\overline{\Delta_{f,T} \cH}}) K_{J,T}\cH.
\end{equation}

To prove the last part of the theorem, assume that $1\in\cN_J$,
$f=p$ is a polynomial,  and  relation \eqref{C*} holds. Consider
another minimal $J$-constrained dilation of $T$, i.e.,
\begin{equation}
\label{another} T_i=V^* (B_i\otimes I_\cD)V, i=1,\ldots,n,
\end{equation}
where $V:\cH\to \cN_J\otimes \cD$ is an isometry, $V\cH$ is
co-invariant under each operator $B_i\otimes I_\cD$,\ $i=1,\ldots,
n$, and
\begin{equation}\label{minimal2}
\cN_J\otimes \cD=\bigvee_{\alpha\in \FF_n^+} (B_\alpha\otimes
I_{\cD}) V\cH.
\end{equation}
Due to \eqref{Pois}, there exists a  unital completely positive
linear map
$$\Psi:\text{\rm span}\{B_\alpha B_\beta^*:\ \alpha,\beta\in \FF_n^+\}\to B(\cH)$$
 such that
$\Psi(B_\alpha B_\beta^*)=T_\alpha T_\beta^*$, \ $\alpha, \beta\in
\FF_n^+$.  Moreover, $\Psi$ has a unique extension to
$C^*(B_1,\ldots, B_n)$. Consider the $*$-representations
\begin{equation*}
\begin{split}
\pi_1:C^*(B_1,\ldots, B_n)\to B(\cN_J&\otimes \overline{\Delta_{p,T}
\cH}),\quad \pi_1(X)
:= X\otimes I_{\overline{\Delta_{p,T} \cH}}, \text{ and }\ \\
\pi_2:C^*(B_1,\ldots, B_n)\to B(\cN_J&\otimes \cD),\quad \pi_2(X):=
X\otimes I_{\cD}.
\end{split}
\end{equation*}
Due to relations  \eqref{Pois}, \eqref{another}, \eqref{C*}, and the
co-invariance of the subspaces $K_{p,T,J}\cH$ and $V\cH$ under
$B_i\otimes I_\cD$, $i=1,\ldots,n$, we have
$$
\Psi(X)=K_{p,T,J}^*\pi_1(X)K_{p,T,J}=V^*\pi_2(X)V\quad \text{ for }
\  X\in C^*(B_1,\ldots, B_n).
$$
Now, due to the minimality conditions \eqref{minimal1} and
\eqref{minimal2}, and relation \eqref{C*}, we deduce that $\pi_1$
and $\pi_2$ are  minimal Stinespring dilations of $\Psi$.  Since
they are unique, there exists a unitary operator $U:\cN_J\otimes
\overline{\Delta_T \cH}\to \cN_J\otimes \cD$ such that
\begin{equation}\label{int-U}
U(B_i\otimes I_{\overline{\Delta_{p,T} \cH}})=(B_i\otimes
I_\cD)U,\quad i=1,\ldots,n,
\end{equation}
and $UK_{p,T,J}=V$.  Since $U$ is unitary, we also have
$$
U(B_i^*\otimes I_{\overline{\Delta_{p,T} \cH}})=(B_i^*\otimes
I_\cD)U,\quad i=1,\ldots,n.
$$
On the other hand, since  $C^*(B_1,\ldots, B_n)$ is irreducible (see
Theorem \ref{compact}),  we must have $U=I_{\cN_J}\otimes Z$, where
$Z\in B(\overline{\Delta_{p,T} \cH},\cD)$ is a unitary operator.
Therefore, $\dim \overline{\Delta_{p,T} \cH}=\dim\cD$ and
$UK_{f,T,J}\cH=V\cH$, which proves that the two dilations are
unitarily equivalent.

In the particular case when $J=\{0\}$, item  (iv) shows that
$$ W:=(W_1\otimes I_{\overline{\Delta_{p,T} \cH}},\ldots, W_n\otimes
I_{\overline{\Delta_{p,T} \cH}})$$ is a realization of the minimal
dilation of $(T_1,\ldots, T_n)\in \cD_p(\cH)$.
 Using   Proposition \ref{max-piece}, one can easily see
that the maximal   piece of $W$ in the noncommutative variety
$\cV_{p,J}$ coincides with $(B_1\otimes I_{\overline{\Delta_{p,T}
\cH}},\ldots, B_n\otimes I_{\overline{\Delta_{p,T} \cH}})$, which
proves (v). The proof is complete.
\end{proof}

\begin{proposition}\label{rank-m} Let $p$ be a positive regular polynomial
and let $J$ be a $w^*$-closed ideal of $F_n^\infty(\cD_p)$ such that
  $1\in \cN_J$ and condition \eqref{C*} holds.
   A pure
$n$-tuple of operators $T\in \cV_{p,J}(\cH)$ has $\rank
\Delta_{p,T}=m$,\ $m=1,2,\ldots, \infty$, if and only if it is
unitarily equivalent to one obtained by compressing $(B_1\otimes
I_{\CC^m},\ldots, B_n\otimes I_{\CC^m})$ to a co-invariant subspace
$\cM\subset \cN_J\otimes \CC^m$  under each operator  $B_i\otimes
I_{\CC^m}$, $i=1,\ldots,n$,  with the property that $\dim P_0\cM=m$,
where $P_0$ is the orthogonal projection of $\cN_J\otimes \CC^m$
onto the subspace $1\otimes \CC^m$.
\end{proposition}
\begin{proof}
  The implication ``$\implies$'' follows from item
(i) of Theorem \ref{dil3}. Conversely, assume that
$$
T_i=P_\cH(B_i\otimes I_{\CC^m})|\cH,\quad i=1,\ldots,n,
$$
where $\cH\subset \cN_J\otimes \CC^m$ is a co-invariant subspace
under each operator $B_i\otimes I_{\CC^m}$,\  $i=1,\ldots,n$, such
that $\dim P_0\cH=m$.  First, notice  that $T:=(T_1,\ldots, T_n)\in
\cV_{p,J}(\cH)$ is  of class $C_{\cdot 0}$. Consider   the case when
$m<\infty$. Since $P_0\cH\subseteq \CC^m$ and $\dim P_0\cH=m$, we
deduce that $P_0\cH=\CC^m$.  Since $1\in \cN_J$, the latter
condition is equivalent  to
\begin{equation}
\label{perpend} \cH^\perp\cap \CC^m=\{0\}.
\end{equation}
Since $I_{\cN_J}-\sum\limits_{|\alpha|\geq 1} a_\alpha B_\alpha
B_\alpha^*=P_\CC^{\cN_J}$, we deduce that
\begin{equation*}\begin{split}
\rank \Delta_{p,T}&=\rank P_\cH\left[\left(I_{\cN_J}-
\sum\limits_{|\alpha|\geq 1}
a_\alpha B_\alpha B_\alpha^*\right)\otimes I_{\CC^m}\right]|\cH\\
&=\rank P_\cH P_0|\cH=\dim P_\cH P_0\cH\\
&=\dim P_\cH\CC^m.
\end{split}
\end{equation*}
If we assume that $\rank \Delta_{p,T}<m$, then there exists  $h\in
\CC^m$, $h\neq 0$, with $P_\cH h=0$.  This contradicts relation
\eqref{perpend}. Consequently, we must have $\rank \Delta_{p,T}=m$.

Now, we consider the case when $m=\infty$. According to Theorem
\ref{cyclic} and its proof,   we have
$$
\bigvee_{\alpha\in \FF_n^+} (B_\alpha\otimes
I_{\CC^m})\cH=\cN_J\otimes \cE,
$$
where $\cE:=P_0\cH$.  Since $\cN_J\otimes \cE$  is reducing for
$B_i\otimes I_{\CC^m}$, \ $i=1,\ldots, n$, we deduce that
$$
T_i=P_\cH(B_i\otimes I_\cE)|\cH,\quad i=1,\ldots, n.
$$
Due to  the uniqueness of the  minimal $J$-constrained dilation of
$T$ (see part (iv) of Theorem \ref{dil3}), we deduce that
$$\dim \overline{\Delta_{p,T} \cH}=\dim\cE=\infty.
$$
This completes the proof.
\end{proof}

 We can
characterize now the pure  $n$-tuples of operators in the
noncommutative variety $\cV_{p,J}$, having  rank one, i.e.,
$\rank\Delta_{p,T}=1$.

\begin{corollary}\label{rank1} Let $J$ be  a $w^*$-closed two-sided ideal of
 $F_n^\infty(\cD_p)$  such that  $1\in\cN_J$ and condition \eqref{C*} is satisfied.
 \begin{enumerate}
 \item[(i)]
If $\cM\subset\cN_J$ is  a co-invariant  subspace under $B_i$,
$i=1,\ldots,n$,  then
$$
T:=(T_1,\ldots, T_n), \quad  T _i:=P_\cM B_i|\cM,\ i=1,\ldots,n,
$$
is a pure    $n$-tuple of operators in $\cV_{p,J}(\cM)$ such that
$\rank\Delta_{p,T}=1$.
\item[(ii)]
If $\cM'$ is another co-invariant subspace under  $B_i$,
$i=1,\ldots,n$,  which gives rise to an $n$-tuple   $T'$, then $T$
and $T'$ are unitarily equivalent if and only if $\cM=\cM'$.
\item[(iii)]
Every pure  $n$-tuple $T\in \cV_{p,J}$ with $\rank\Delta_{p,T}=1$ is
unitarily equivalent to one obtained by compressing $(B_1,\ldots,
B_n)$ to a co-invariant subspace  under  $B_i$, $i=1,\ldots,n$.
\end{enumerate}
\end{corollary}
\begin{proof}
Since $\cM\subset\cN_J$ is  a co-invariant  subspace under each
operator  $B_i$, $i=1,\ldots,n$,  we can use the
$F_n^\infty(\cD_f)$-functional calculus and deduce that
$$
\varphi(T_1,\ldots, T_n)=P_\cM \varphi(B_1,\ldots, B_n)|\cM=0,\quad
\varphi\in J.
$$
This shows that $(T_1,\ldots, T_n)\in \cV_{p,J}(\cM)$.
 On the other hand, we have
\begin{equation*}\begin{split}
I_\cM-\sum_{|\alpha|\geq 1} a_\alpha T_\alpha T_\alpha^*&=
P_\cM\left(I_{\cN_J}-\sum_{|\alpha|\geq 1} a_\alpha B_\alpha B_\alpha^*\right)|\cM\\
&=P_\cM P_\CC^{\cN_J}|\cM.
\end{split}
\end{equation*}
Hence, $\text{\rm rank}\, \Delta_{p,T}\leq 1$. Since
$$
\Phi_{p,T}^k(I)=P_\cM \Phi_{p,B}^k(I_{\cN_J})|\cM, \quad
k=1,2,\ldots,
$$
and $(B_1,\ldots, B_n)$ is a pure  $n$-tuple in $\cV_{p,J}(\cN_J)$,
we deduce that $(T_1,\ldots, T_n)\in \cV_{p,J}(\cM)$ is pure. This
also implies that $\Delta_{p,T}\neq 0$, so $\text{\rm rank}\,
\Delta_{p,T}\geq 1$. Consequently, we have $\rank \Delta_{p,T}=1$.

To prove  (ii), notice that,  as in the proof of  Proposition
\ref{rank-m}, one can show that $T$ and $T'$ are unitarily
equivalent if and only if there exists a unitary operator
$\Lambda:\cN_J\to \cN_J$ such that
$$
\Lambda B_i=B_i \Lambda, \ \text{ for } i=1, \ldots, n, \ \text{ and
}  \ \Lambda \cM=\cM'.
$$
Hence $\Lambda B_i^*=B_i^* \Lambda$, $i=1,\ldots,n$. Since
$C^*(B_1,\ldots, B_n)$ is irreducible (see Theorem \ref{compact})
$\Lambda$ must be a scalar multiple of the identity. Hence, we have
$\cM=\Lambda \cM=\cM'$.

  Part (iii)  of this corollary follows from Theorem
\ref{dil2},   Corollary \ref{part-cases}, and Corollary
\ref{unique}.
\end{proof}

\bigskip

  \section{Characteristic functions and model theory}
\label{Characteristic}

Given a positive regular free holomorphic function $f$ on $B(\cH)^n$
and an $n$-tuple of operators $T:=(T_1,\ldots, T_n)$ in the
noncommutative domain  $\cD_f(\cH)$, we associate a characteristic
function $\Theta_{f,T}$, which is a multi-analytic operator with
respect to the universal model $(W_1,\ldots, W_n)$ of $\cD_f$. We
prove that the characteristic function is a complete unitary
invariant and obtain a functional model for the class of completely
non-coisometric elements of $\cD_f(\cH)$. Similar results are
obtained for  constrained characteristic functions $\Theta_{f,T,J}$
associated with the noncommutative varieties $\cV_{f,J}\subset
\cD_f$. In particular, the commutative case is discussed.

 Let $f=\sum_{|\alpha|\geq 1} a_\alpha X_\alpha$ a positive
 regular free holomorphic function on $B(\cH)^n$.
Given the  noncommutative domain
$$
\cD_f(\cH):=\left\{ (X_1,\ldots, X_n)\in B(\cH)^n: \
\left\|\sum_{|\alpha|\geq 1} a_\alpha X_\alpha
X_\alpha^*\right\|\leq 1\right\},
$$
define  $\Gamma:=\{\alpha\in \FF_n^+: \ a_\alpha \neq 0\}$ and
$N:=\text{\rm card}\, \Gamma$.
 If  $T:=(T_1,\ldots, T_n)\in \cD_f(\cH)$, we define the row operator
 $C(T):=[\sqrt{a_{\tilde\alpha}}T_{\tilde\alpha}: \ \alpha\in \Gamma]$,
 where the entries are arranged in the lexicographic order of
 $\Gamma\subset \FF_n^+$. Note that  $C(T)$ is an operator acting from
 $\cH^{(N)}$ (the direct sum of $N$ copies of $\cH$)  to $\cH$.
Let $(W_1,\ldots, W_n)$ (resp. $(\Lambda_1,\ldots, \Lambda_n)$ ) be
the weighted left (resp. right) creation operators associated with
the noncommutative domain $\cD_f$.

  We
define the characteristic function  of $T$ to be the multi-analytic
operator (with respect to the  operators $W_1,\ldots, W_n$)

$$
\Theta_{f,T}(\Lambda_1,\ldots, \Lambda_n):F^2(H_n)\otimes
\cD_{C(T)^*}\to F^2(H_n)\otimes \cD_{C(T)}
$$
with the formal Fourier representation
\begin{equation*}
\begin{split}
  -I_{F^2(H_n)}\otimes {C(T)}+
\left(I_{F^2(H_n)}\otimes \Delta_{C(T)}\right)&\left
(I_{F^2(H_n)\otimes \cH}-\sum_{|\alpha|\geq 1} a_{\tilde \alpha}
 \Lambda_\alpha\otimes T_{\tilde \alpha}^*\right)^{-1}\\
& \left[\sqrt{a_{\tilde\alpha}}\Lambda_\alpha\otimes I_\cH:\ \alpha
\in \Gamma \right] \left(I_{F^2(H_n)}\otimes
\Delta_{{C(T)^*}}\right),
\end{split}
\end{equation*}
where the defect operators  associated with the  row contraction
$C(T)$ are
\begin{equation*}
\begin{split}
 \Delta_{C(T)}&:=\left( I_\cH-C(T)C(T)^*\right)^{1/2}\in B(\cH)
\quad \text{ and }\\
 \Delta_{C(T)^*}&:=(I-C(T)^* C(T))^{1/2}\in
B(\cH^{(N)}),
\end{split}
\end{equation*}
  while the
defect spaces are $\cD_{C(T)}:=\overline{\Delta_{C(T)}\cH}$ and
$\cD_{C(T)^*}:=\overline{\Delta_{C(T)^*}\cH^{(N)}}$.

The following factorization result will play an important role in
our investigation.

\begin{theorem}\label{factor}
If  $T:=(T_1,\ldots, T_n)\in \cD_f(\cH)$, then the characteristic
function  $\Theta_{f,T}$  associated with $T$ is a well-defined
contractive multi-analytic operator with respect to the weighted
left creation operators associated with the  noncommutative domain
$\cD_f$. Moreover,
\begin{equation}\label{fa}
I-\Theta_{f,T}\Theta_{f,T}^*=K_{f,T}K_{f,T}^*,
\end{equation}
where  $K_{f,T}$ is the corresponding Poisson kernel associated with
$T$ and $\cD_f(\cH)$.
\end{theorem}
\begin{proof}

Consider the following operators:
\begin{equation*}
\begin{split}
 \widehat{C(T)} &:=
\left[I_{F^2(H_n)}\otimes \sqrt{a_{\tilde\alpha}}T_{\tilde\alpha}: \
\alpha\in \Gamma \right],\\
  \widehat{C(\Lambda)} &:=\left[
\sqrt{a_{\tilde\alpha}}\Lambda_{\alpha}\otimes I_\cH : \ \alpha\in
\Gamma \right],\\
\widehat{C(\Lambda)}_r &:=\left[
\sqrt{a_{\tilde\alpha}}r^{|\alpha|}\Lambda_{\alpha}\otimes I_\cH : \
\alpha\in \Gamma \right], \quad 0<r<1,
\end{split}
\end{equation*}
and define
$$
A:=\widehat{C(T)}^*,\ B:=\Delta_{\widehat{C(T)}^*}, \
C:=\Delta_{\widehat{C(T)}}, \ D:=- \widehat{C(T)}, \ \text{ and } \
Z:=\widehat{C(\Lambda)}_r, \ 0<r<1.
$$
Since $C(T)$ is a contraction, the operator matrix
\begin{equation*}
\left(\begin{matrix} A&B\\C&D
\end{matrix}
\right) = \left(\begin{matrix} \widehat{C(T)}^*
&\Delta_{\widehat{C(T)}^*}\\
\Delta_{\widehat{C(T)}}& - \widehat{C(T)}
\end{matrix}
\right)
\end{equation*}
is  unitary. Therefore,
\begin{equation}\label{unitar}
AA^*+BB^*=I,\ CC^*+DD^*=I, \text{ and } \ AC^*+BD^*=0.
\end{equation}
According to Section \ref{Noncommutative}, we have
$I-\sum\limits_{|\alpha|\geq 1} a_{\tilde \alpha} \Lambda_\alpha
\Lambda_\alpha^*=P_\CC$, which implies
$\|Z\|=\|\widehat{C(\Lambda)}_r\|\leq r<1$ for $0<r<1$. Since
$(T_1,\ldots, T_n)\in \cD_f(\cH)$, we have  $\|A\|\leq 1$ and
consequently $\|ZA\|<1$.

Therefore, the operator
$$\Phi(Z):=D+C(I-ZA)^{-1} ZB
$$
is well-defined and, using relation \eqref{unitar}, we have
\begin{equation*}
\begin{split}
I-\Phi(Z)\Phi(Z)^* &=
I-DD^* -C(I-ZA)^{-1}ZBD^*-DB^*Z^*(I-A^*Z^*)^{-1} C^*\\
&
\qquad -C(I-ZA)^{-1}ZBB^*Z^*(I-A^*Z^*)^{-1} C^*\\
&=
CC^*+C(I-ZA)^{-1}ZAC^*+CA^*Z^*(I-A^*Z^*)^{-1} C^*\\
&
\qquad -C(I-ZA)^{-1}ZZ^*(I-A^*Z^*)^{-1} C^*\\
&\qquad
+C(I-ZA)^{-1}ZAA^*Z^*(I-A^*Z^*)^{-1} C^*\\
&=
C(I-ZA)^{-1}\left[(I-ZA)(I-A^*Z^*)+ZA(I-A^*Z^*)\right.\\
& \qquad \left.+(I-ZA)A^*Z^*-ZZ^*+ZAA^*Z^*\right]
(I-A^*Z^*)^{-1} C^*\\
&= C(I-ZA)^{-1}(I-ZZ^*) (I-A^*Z^*)^{-1} C^*.
 \end{split}
\end{equation*}
Hence,
\begin{equation}\label{PZ}
I-\Phi(Z)\Phi(Z)^*= C(I-ZA)^{-1}(I-ZZ^*) (I-A^*Z^*)^{-1} C^*.
\end{equation}
Therefore,  according to our notations and the definition of
$\Theta_{f,T}(r\Lambda_1,\ldots, r\Lambda_n)$, the defect operator
$$I_{F^2(H_n)\otimes \cH}-\Theta_{f,T}(r\Lambda_1,\ldots, r\Lambda_n)
\Theta_{f,T}(r\Lambda_1,\ldots, r\Lambda_n)^*$$
 is equal to
the product
\begin{equation*}
\begin{split}
\Delta_{\widehat{C(T)}}\left(I_{F^2(H_n)\otimes
\cH}-\widehat{C(\Lambda)}_r \widehat{C(T)}^*\right)^{-1}
&\left(I_{F^2(H_n)\otimes \cH}-
\widehat{C(\Lambda)}_r\widehat{C(\Lambda)}_r^*\right)\\
&\left(I_{F^2(H_n)\otimes \cH}-\widehat{C(T)}
\widehat{C(\Lambda)}_r^* \right)^{-1} \Delta_{\widehat{C(T)}},
\end{split}
\end{equation*}
which is a positive operator for any $r\in (0,1)$. Consequently, we
have
$$
\sup\limits_{0\leq r<1}\|\Theta_{f,T}(r\Lambda_1,\ldots,
r\Lambda_n)\|\leq1.
$$
An operator-valued extension of Theorem
\ref{f-infty} shows that
 $\Theta_{f,T}$ is a multi-analytic operator
 with respect to the weighted left creation
  operators associated  with the domain $\cD_f$, i.e., in the
  operator space
$R_n^\infty(\cD_f)\bar\otimes B(\cD_{C(T)^*}, \cD_{C(T)})$, and the
characteristic function  $\Theta_{f,T}(\Lambda_1,\ldots, \Lambda_n)$
is equal to
$$
\text{\rm SOT-}\lim_{r\to 1}
\left[-\widehat{C(T)}+\Delta_{\widehat{C(T)}}\left(I_{F^2(H_n)\otimes
\cH}-\widehat{C(\Lambda)}_r
\widehat{C(T)}^*\right)^{-1}\widehat{C(\Lambda)}_r\Delta_{\widehat{C(T)}^*}
\right].
$$
Now, the above calculations  and relations \eqref{r<r} and
\eqref{Fourier-Cauc} of Section \ref{Cauchy}, reveal that
\begin{equation*}
\begin{split}
&I_{F^2(H_n)\otimes \cD_{C(T)}}-\Theta_{f,T}(r\Lambda_1,\ldots,
r\Lambda_n)\Theta_{f,T}(r\Lambda_1,\ldots,
r\Lambda_n)^*\\
&\qquad=\left(I_{F^2(H_n)}\otimes \Delta_{f,T}\right)
\left(I_{F^2(H_n)\otimes \cH}-\sum_{|\alpha|\geq 1} a_{\tilde\alpha}
r^{|\alpha|}\Lambda_\alpha \otimes T_{\tilde\alpha}^*\right)^{-1}\\
&\qquad\qquad \phantom{X}\left[\left(I_{F^2(H_n)}-
\sum_{|\alpha|\geq 1} a_{\tilde \alpha} r^{2|\alpha|}\Lambda_\alpha
\Lambda_\alpha^*\right)\otimes I_\cH\right]\\
&\qquad\qquad \phantom{X}\left(I_{F^2(H_n)\otimes
\cH}-\sum_{|\alpha|\geq 1} a_{\tilde \alpha}
r^{|\alpha|}\Lambda_\alpha ^*\otimes T_{\tilde\alpha}\right)^{-1}
\left(I_{F^2(H_n)}\otimes \Delta_{f,T}\right)\\
&\qquad= \left(\sum_{ \gamma\in \FF_n^+}
b_{\tilde\gamma}r^{|\gamma|} \Lambda_\gamma \otimes \Delta_{f,T}
T_{\tilde \gamma}^*\right)
 \left[\left(I_{F^2(H_n)}- \sum_{|\alpha|\geq 1}
a_{\tilde \alpha} r^{2|\alpha|}\Lambda_\alpha
\Lambda_\alpha^*\right)\otimes I_\cH\right]\\
&\qquad\qquad \phantom{X}
 \left(\sum_{ \beta\in \FF_n^+} b_{\tilde\beta}r^{|\beta|}
\Lambda_\beta^* \otimes  T_{\tilde \beta}\Delta_{f,T}\right).
\end{split}
\end{equation*}
Now, we recall from Section \ref{Noncommutative} that
$$
\Lambda_\beta^* e_\alpha =\begin{cases} \frac
{\sqrt{b_\gamma}}{\sqrt{b_{\alpha}}}e_\gamma& \text{ if }
\alpha=\gamma \tilde \beta \\
0& \text{ otherwise }
\end{cases} \quad \text{ and } \quad
P_\CC \Lambda_\beta^* e_\alpha =\begin{cases}
\frac {1}{\sqrt{b_{\beta}}}   & \text{ if } \alpha=\tilde\beta\\
0& \text{ otherwise }
\end{cases}.
$$
Hence, and using again that $I-\sum\limits_{|\alpha|\geq 1}
a_{\tilde \alpha} \Lambda_\alpha \Lambda_\alpha^*=P_\CC$, we deduce
that

\begin{equation*}
\begin{split}
\lim_{r\to 1}&\left<[I_{F^2(H_n)\otimes
\cD_{C(T)}}-\Theta_{f,T}(r\Lambda_1,\ldots,
r\Lambda_n)\Theta_{f,T}(r\Lambda_1,\ldots, r\Lambda_n)^*]
(e_\alpha\otimes h), e_\omega \otimes k\right>\\
&= \lim_{r\to 1}\sum_{{\gamma\in \FF_n^+}\atop {|\gamma|\leq
|\omega|}} \sum_{{\beta\in \FF_n^+} \atop { |\beta|\leq |\alpha|}}
\left< b_{\tilde \gamma} b_{\tilde
\beta}r^{|\gamma|+|\beta|}\Lambda_\gamma
 \left(I_{F^2(H_n)}- \sum_{|\alpha|\geq 1}
a_{\tilde \alpha} r^{2|\alpha|}\Lambda_\alpha
\Lambda_\alpha^*\right) \Lambda_\beta^* e_\alpha, e_\omega\right>\\
&\qquad \qquad \left< \Delta_{f,T}
T_{\tilde\gamma}^*T_{\tilde\beta}\Delta_{f,T} h,k\right>\\
&= \sum_{{\gamma\in \FF_n^+}\atop { |\gamma|\leq |\omega|}}
\sum_{{\beta\in \FF_n^+}\atop { |\beta|\leq |\alpha|}} \left<
b_{\tilde \gamma} b_{\tilde \beta}\Lambda_\gamma P_\CC
\Lambda_\beta^* e_\alpha, e_\omega\right> \left< \Delta_{f,T}
T_{\tilde\gamma}^*T_{\tilde\beta}\Delta_{f,T} h,k\right>\\
&= \sum_{{\beta\in \FF_n^+}\atop { |\beta|\leq |\alpha|}} \left<
b_{\tilde \gamma} b_{\alpha}\Lambda_\gamma
\left(\frac{1}{\sqrt{b_\alpha}}\right), e_\omega\right> \left<
\Delta_{f,T}
T_{\tilde\gamma}^*T_{\alpha}\Delta_{f,T} h,k\right>\\
&= \sum_{{\beta\in \FF_n^+}\atop{ |\beta|\leq |\alpha|}}  b_{\tilde
\gamma} \sqrt{b_\alpha} \frac{1}{\sqrt{b_{\tilde\gamma}}}\left<
e_{\tilde \gamma},e_\omega\right> \left< \Delta_{f,T}
T_{\tilde\gamma}^*T_{\alpha}\Delta_{f,T} h,k\right>\\
&=\sqrt{b_\omega}\sqrt{b_\alpha}\left< \Delta_{f,T}
T_{\omega}^*T_{\alpha}\Delta_{f,T} h,k\right>
\end{split}
\end{equation*}
for every $\alpha,\omega\in \FF_n^+$, \ $h\in \cD_{C(T)^*}$,
$k\in\cD_{C(T)}$.
An operator-valued version   of Theorem \ref{funct-calc}, shows that
$$\text{\rm SOT-}\lim_{r\to 1} \Theta_{f,T}(r\Lambda_1,\ldots, r\Lambda_n)^*
=\Theta_{f,T}(\Lambda_1,\ldots, \Lambda_n)^*.
$$
Therefore, taking into account  the above computations, we deduce
that

\begin{equation*}\begin{split}
&\left<[I_{F^2(H_n)\otimes
\cD_{C(T)}}-\Theta_{f,T}(\Lambda_1,\ldots, \Lambda_n)
\Theta_{f,T}(\Lambda_1,\ldots, \Lambda_n)^*](e_\alpha\otimes h),
 e_\omega \otimes k\right>\\
&\qquad =\sqrt{b_\omega}\sqrt{b_\alpha}\left< \Delta_{f,T}
T_{\omega}^*T_{\alpha}\Delta_{f,T} h,k\right>
\end{split}
\end{equation*}
for every $\alpha,\omega\in \FF_n^+$, \ $h,k\in \cD_{C(T)}$.

   On the other hand, using the definition
of the Poisson kernel associated with $(T_1,\ldots, T_n)\in
\cD_f(\cH)$, we deduce that $K_{f,T}^*(e_\alpha\otimes
h)=\sqrt{b_\alpha}T_\alpha \Delta_{f,T}$  for any $\alpha \in
\FF_n^+$, and
\begin{equation*}
\begin{split}
\left<K_{f,T}K_{f,T}^*(e_\alpha\otimes h), e_\omega\otimes k\right>
&=
\left< K_{f,T}\left(\sqrt{b_\alpha}T_\alpha \Delta_{f,T}h\right),
e_\omega\otimes k\right>\\
&= \left<\sum_{\beta\in\FF_n^+}\left(
\sqrt{b_\beta}\sqrt{b_\alpha}e_\beta \otimes\
 \Delta_{f,T} T_\beta^* T_\alpha \Delta_{f,T} h\right),\,  e_\omega\otimes k\right>\\
&= \sqrt{b_\alpha}\sqrt{b_\omega}\left<  \Delta_{f,T} T_\omega^*
T_\alpha \Delta_{f,T} h,  k\right>
\end{split}
\end{equation*}
for any $h,k\in \cD_T$ and $\alpha,\omega\in \FF_n^+$. Taking into
account the above relations,  we deduce that
$$
I_{F^2(H_n)\otimes \cD_{C(T)}}-\Theta_{f,T}(\Lambda_1,\ldots,
\Lambda_n)\Theta_{f,T}(\Lambda_1,\ldots,
\Lambda_n)^*=K_{f,T}K_{f,T}^*,
$$
which completes the proof.
\end{proof}

We recall from Section \ref{Cauchy} that the joint spectral radius
of an $n$-tuple of operators $X:=(X_1,\ldots, X_n)\in \cD_f(\cH)$,
with respect to the noncommutative domain $\cD_f$,  is defined by
$$
r_f(X_1,\ldots, X_n):=\lim_{k\to\infty}\left\| \Phi_{f,X}^k(I)
\right\|^{1/2k},
$$
where $\Phi_{f,X}(Y):=\sum\limits_{|\alpha|\geq 1}  a_{\alpha}
X_\alpha Y X_\alpha^*$ \,for any  $Y\in B(\cH)$. We also recall that
if $f=\sum_{|\alpha|\geq 1} a_\alpha X_\alpha$ is a positive regular
free holomorphic function on $B(\cH)^n$, then so is  $\tilde
f:=\sum_{|\alpha|\geq 1} a_{\tilde \alpha} X_\alpha$.

\begin{corollary}\label{fact2}
Let $T:=(T_1,\ldots, T_n)\in \cD_f(\cH)$  and let $\Theta_{f,T}$ be
its characteristic function. If $X:=(X_1,\ldots, X_n)\in \cD_{\tilde
f}(\cK)$ is an $n$-tuple of operators  with joint  spectral radius
$r_{\tilde f}(X_1,\ldots, X_n)<1$, then
\begin{equation*}
\begin{split}
I_{\cK\otimes \cD_{C(T)}}&-\Theta_{f,T}(X_1,\ldots,
X_n)\Theta_{f,T}(X_1,\ldots, X_n)^*\\
&\quad = \Delta_{\widehat{C( T)}}\left(I_{\cK\otimes
\cH}-\widehat{C( X)}\widehat{C( T)}^*\right)^{-1}\left(I_{\cK\otimes
\cH}-\widehat{C( X)} \widehat{C( X)}^*\right)\\
&\qquad\qquad \left(I_{\cK\otimes \cH}-\widehat{C( T)} \widehat{C(
X)}^*\right)^{-1} \Delta_{\widehat{C( T)}},
\end{split}
\end{equation*}
where $\widehat{C( X)}:=[\sqrt{a_{\tilde\alpha}} X_\alpha\otimes
I_\cH:\ \alpha\in \Gamma]$  (see the  notations   from the proof of
Theorem $\ref{factor}$).
\end{corollary}
\begin{proof}
As in Section \ref{Cauchy}, we deduce that
$$
r\left( \sum_{|\alpha|\geq 1} a_{\tilde\alpha} X_\alpha \otimes
T_{\tilde \alpha}^*\right) \leq r_{\tilde f}(X_1,\ldots,
X_n)r_f(T_1,\ldots, T_n).
$$
Since $(T_1,\ldots, T_n)\in \cD_f(\cH)$, we have $r_f(T_1,\ldots,
T_n)\leq 1$. The above inequality shows that the spectral radius
$r\left( \sum_{|\alpha|\geq 1} a_{\tilde\alpha} X_\alpha \otimes
T_{\tilde \alpha}^*\right)<1$ and the operator $\left(I_{\cK\otimes
\cH}-\sum_{|\alpha|\geq 1} a_{\tilde\alpha} X_\alpha \otimes
T_{\tilde \alpha}^*\right)^{-1} $ is bounded. Now, replacing $Z$ by
$\widehat{C( X)}$ in  relation \eqref{PZ} (see the proof of Theorem
\ref{factor}), we obtain the required factorization.
\end{proof}

We introduce now the constrained characteristic function associated
with $n$-tuples of operators in the noncommutative variety
$\cV_{f,J}\subset \cD_f$.

Let $J$ be a $w^*$-closed two-sided ideal of the  Hardy algebra
algebra $F_n^\infty(\cD_f)$.
As in our previous sections, define the subspaces of $F^2(H_n)$,
$$
\cM_J:=\overline{JF^2(H_n)}\quad \text{and}\quad
\cN_J:=F^2(H_n)\ominus \cM_J,
$$
and  the {\it constrained  weighted left} (resp.~{\it right}) {\it
creation operators} associated with  $\cD_f$ and $J$ by setting
$$B_i:=P_{\cN_J} W_i|\cN_J \quad \text{and}\quad
C_i:=P_{\cN_J} \Lambda_i|\cN_J,\quad i=1,\ldots, n.
$$
Now, we assume that the ideal $J$ is generated  by a family $\cP_J$
of polynomials in $W_1,\ldots, W_n$ and the identity. An $n$-tuple
of operators $T:=(T_1,\ldots, T_n)\in \cD_f(\cH)$ is called $J$-{\it
constrained} if
 it belongs to the noncommutative variety $\cV_{f,J}(\cH)$ defined by
 $$
 \cV_{f,J}(\cH):=\left\{ (X_1,\ldots, X_n)\in \cD_f(\cH): \ q(X_1,\ldots,
 X_n)=0\  \text{ for } q\in \cP_J\right\}.
 $$
 We define the {\it constrained
characteristic function} associated with  $T\in \cV_{f,J}(\cH)$ to
be the multi-analytic operator (with respect to the constrained
shifts $B_1,\ldots, B_n$)
$$
\Theta_{f,T,J}=\Theta_{f,T,J}(C_1,\ldots, C_n):\cN_J\otimes
\cD_{C(T)^*}\to \cN_J\otimes \cD_{C(T)}
$$
defined by the formal Fourier representation
\begin{equation*}
\begin{split}
 -I_{\cN_J}\otimes
C(T)+ \left(I_{\cN_J}\otimes
\Delta_{C(T)}\right)\Bigl(I_{{\cN_J}\otimes \cH}
&-\sum_{|\alpha|\geq 1} a_{\tilde \alpha} C_\alpha\otimes T_{\tilde
\alpha}^*
\Bigr)^{-1}\\
&\left[\sqrt{a_{\tilde \alpha}} C_\alpha\otimes I_\cH:\ \alpha\in
\Gamma \right] \left(I_{\cN_J}\otimes \Delta_{C(T)^*}\right),
\end{split}
\end{equation*}
where $C_i:= P_{\cN_J} \Lambda_i|_{\cN_J}$, $i=1,\ldots, n$. Taking
into account that $\cN_J$ is  an invariant subspace under
$\Lambda^*_1,\ldots, \Lambda^*_n$, we can see that

\begin{equation}\label{rest}
\begin{split}
\Theta_{f,T}(\Lambda_1,\ldots, \Lambda_n)^*(\cN_J\otimes\cD_{C(T)})
&\subseteq \cN_J\otimes \cD_{C(T)^*}\ \text{  and  }\\
P_{\cN_J\otimes \cD_{C(T)}}\Theta_{f,T}(\Lambda_1,\ldots,
\Lambda_n)|\cN_J\otimes \cD_{C(T)^*}&=\Theta_{f,T,J}(C_1,\ldots,
C_n),
\end{split}
\end{equation}
where $\Theta_{f,T}$ is the characteristic function of $T$ as an
element of the noncommutative domain $ \cD_f(\cH)$. Since
$\Theta_{f,T}\in R_n^\infty(\cD_f)\bar\otimes
B(\cD_{C(T)^*},\cD_{C(T)})$, it is easy to see that
$\Theta_{f,T,J}\in R_n^\infty(\cV_{f,J})\bar\otimes
B(\cD_{C(T)^*},\cD_{C(T)})$.

Let us remark that   the above definition of the constrained
characteristic function  makes sense (and has  the same properties)
when $J$ is an arbitrary $w^*$-closed two-sided ideal of
$F_n^\infty(\cD_f)$ and $T:=(T_1,\ldots, T_n)$ is an arbitrary
c.n.c.  $n$-tuple of operators in $\cV_{f,J}(\cH)$.

\begin{theorem}\label{J-factor}
Let $J\neq F_n^\infty(\cD_f)$ be a $w^*$-closed two-sided ideal of
$F_n^\infty(\cD_f)$ generated by a family of polynomials
$\cP_J\subset F_n^\infty(\cD_f)$. If $T:=(T_1,\ldots, T_n)$ is in
the noncommutative variety $ \cV_{f,J}(\cH)$, then
\begin{equation*}
I_{\cN_J\otimes
\cD_{C(T)}}-\Theta_{f,T,J}\Theta_{f,T,J}^*=K_{f,T,J}K_{f,T,J}^*,
\end{equation*}
where $\Theta_{f,T,J}$ is  the  constrained characteristic function
of   $T\in \cV_{f,J}(\cH)$ and $K_{f,T,J}$ is  the corresponding
constrained Poisson kernel.

Moreover,  the factorization above  holds when  $J$ is an arbitrary
$w^*$-closed two-sided ideal of $F_n^\infty(\cD_f)$ and
$T:=(T_1,\ldots, T_n)\in \cV_{f,T}(\cH)$ is a c.n.c. $n$-tuple of
operators.
\end{theorem}

\begin{proof}
 The constrained Poisson kernel associated with $T\in \cV_{f,J}(\cH)$  is the operator
  $K_{f,T,J}:\cH \to \cN_J\otimes \overline{\Delta_{C(T)}\cH}$ defined by
\begin{equation}\label{J-K}
K_{f,T,J}:=(P_{\cN_J}\otimes
I_{\overline{\Delta_{C(T)}\cH}})K_{f,T},
\end{equation}
where $K_{f,T}$ is the  Poisson kernel of $T\in \cD_f(\cH)$.
According to the proof of Theorem \ref{dil1}, we have $\text{\rm
range}\, K_{f,T}\subseteq \cN_J\otimes \overline{\Delta_{C(T)}\cH}$.
Using Theorem \ref{factor} and taking the compression of relation
\eqref{fa} to the subspace $\cN_J\otimes \cD_{C(T)}\subset
F^2(H_n)\otimes \cD_{C(T)}$, we obtain
\begin{equation*}
\begin{split}
I_{\cN_J\otimes \cD_{C(T)}}-P_{\cN_J\otimes
\cD_{C(T)}}\Theta_{f,T}(\Lambda_1,\ldots, \Lambda_n)
&\Theta_{f,T}(\Lambda_1,\ldots, \Lambda_n)^*|\cN_J\otimes
\cD_{C(T)}\\
&= P_{\cN_J\otimes \cD_{C(T)}}K_{f,T}K_{f,T}^*|\cN_J\otimes
\cD_{C(T)}.
\end{split}
\end{equation*}
Hence, and taking into account  relations
 \eqref{rest},  \eqref{J-K},  and that $C_i^*=\Lambda_i^*|\cN_J$,
  \ $i=1,\ldots, n$,  we infer that
$$
I_{\cN_J\otimes
\cD_{C(T)}}-\Theta_{f,T,J}(C_1,\ldots,C_n)\Theta_{f,T,J}(C_1,\ldots,C_n)^*
=K_{f,T,J}K_{f,T,J}^*.
$$
The last part of the theorem can be proved in a similar manner,
using Theorem \ref{dil3} and Theorem \ref{factor}. The proof is
complete.
\end{proof}

Now we present a functional  model for  c.n.c. $n$-tuples of
operators $(T_1,\ldots, T_n)$ in the noncommutative variety
$\cV_{f,J}(\cH)$, in terms of the constrained characteristic
functions.

\begin{theorem}\label{model}
Let $J\neq F_n^\infty(\cD_f)$ be a $w^*$-closed two-sided ideal of
the Hardy algebra $F_n^\infty(\cD_f)$ and  let $T:=(T_1,\ldots,
T_n)$  be a c.n.c.  $n$-tuple of operators in the noncommutative
variety
$$
\cV_{f,J}(\cH):=\left\{ (X_1,\ldots, X_n)\in \cD_f(\cH): \
\varphi(X_1,\ldots, X_n)=0,\quad \varphi\in J\right\}.
$$
Then $T:=(T_1,\ldots, T_n)$ is unitarily equivalent to the $n$-tuple
$\TT:=(\TT_1,\ldots, \TT_n)\in \cV_{f,J}(\HH_{f,T,J})$ on the
Hilbert space
\begin{equation*}
\begin{split}
\HH_{f,T,J}&:=\left[\left(\cN_J\otimes \cD_{C(T)}\right)\oplus
\overline{\Delta_{\Theta_{f,T,J}}(\cN_J\otimes
\cD_{C(T)^*})}\right]\\
& \qquad \qquad\ominus\left\{\Theta_{f,T,J}\varphi\oplus
\Delta_{\Theta_{f,T,J}}\varphi:\ \varphi\in \cN_J\otimes
\cD_{C(T)^*}\right\},
\end{split}
\end{equation*}
where $\Delta_{\Theta_{f,T,J}}:= \left(I-\Theta_{f,T,J}^*
\Theta_{f,T,J}\right)^{1/2}$ and the operator $\TT_i$, \
$i=1,\ldots, n$,  is uniquely defined by the relation
$$
\left( P_{\cN_J\otimes \cD_{C(T)}}|_{\HH_{f,T,J}}\right) \TT_i^*x=
(B_i^*\otimes I_{\overline{\Delta_{C(T)}\cH}})\left( P_{\cN_J\otimes
\cD_{C(T)}}|_{\HH_{f,T,J}}\right)x, \quad x\in \HH_{f,T,J},
$$
where $P_{\cN_J\otimes \cD_{C(T)}}|_{\HH_{f,T,J}}$ is an injective
operator, $ P_{\cN_J\otimes \cD_{C(T)}}$ is the orthogonal
projection of the Hilbert space $\left(\cN_J\otimes
\cD_{C(T)}\right)\oplus
\overline{\Delta_{\Theta_{f,T,J}}(\cN_J\otimes \cD_{C(T)^*})}$ onto
the subspace $\cN_J\otimes \cD_{C(T)}$, and $B_1,\ldots, B_n$ are
the constrained weighted left creation operators associated with
$\cV_{f,J}$.

Moreover,  $T$ is a pure  $n$-tuple of operators in
$\cV_{f,J}(\cH)$, if and only if the constrained characteristic
function $\Theta_{f,T,J}$ is an inner multi-analytic operator.
  In this case,   $T$ is unitarily equivalent to the   $n$-tuple
\begin{equation}\label{HH}
\left(P_ {\HH_{f,T,J}} (B_1\otimes
I_{\cD_{C(T)}})|\HH_{f,T,J},\ldots, P_{\HH_{f,T,J}} (B_n\otimes
I_{\cD_{C(T)}})|\HH_{f,T,J}\right),
\end{equation}
where $P_{\HH_{f,T,J}}$ is the orthogonal projection of
$\cN_J\otimes \cD_{C(T)}$ onto the Hilbert space
$$\HH_{f,T,J}=\left(\cN_J\otimes \cD_{C(T)}\right)\ominus
\Theta_{f,T,J}(\cN_J\otimes \cD_{C(T)^*}).
$$
 \end{theorem}

\begin{proof}
Define the Hilbert space
$$
\KK_{f,T,J}:=\left(\cN_J\otimes \cD_{C(T)}\right)\oplus
\overline{\Delta_{\Theta_{f,T,J}}(\cN_J\otimes \cD_{C(T)^*})}
$$
and denote by $P_{\HH_{f,T,J}}$ the orthogonal projection of
$\KK_{f,T,J}$ onto the subspace $\HH_{f,T,J}$, defined in the
theorem. In what follows, we show that there is a unique  unitary
operator $\Gamma:\cH\to \HH_{f,T,J}$ such that
\begin{equation}\label{Ga}
\Gamma(K_{f,T,J}^* g)=P_{\HH_{f,T,J}}(g\oplus 0), \quad  g\in
\cN_J\otimes \cD_{C(T)}.
\end{equation}
First,   let us prove that
\begin{equation}\label{K*P}
\|K_{f,T,J}^* g\|=\|P_{\HH_{f,T,J}}(g\oplus 0)\|,  \quad g\in
\cN_J\otimes \cD_{C(T)}.
\end{equation}
 Due to Theorem \ref{J-factor}, we have
\begin{equation}
\label{Ke-Te}
 \|K_{f,T,J}^* g\|^2+ \|\Theta_{f,T,J}^*g\|^2=\|g\|^2,
\quad g\in \cN_J\otimes \cD_{C(T)}. \end{equation}
 On the other
hand, note that the operator $\Phi: \cN_J\otimes \cD_{C(T)^*}\to
\KK_{f,T,J}$ defined by
$$
\Phi \varphi:=\Theta_{f,T,J} \varphi\oplus \Delta_{\Theta_{f,T,J}}
\varphi,\quad \varphi\in \cN_J\otimes \cD_{C(T)^*},
$$
 is an isometry and
\begin{equation}\label{fi}
\Phi^*(g\otimes 0)=\Theta_{f,T,J}^*g, \quad g\in \cN_J\otimes
\cD_{C(T)}.
\end{equation}
Consequently,
 we have
\begin{equation*}
\begin{split}
\|g\|^2&= \|P_{\HH_{f,T,J}}(g\oplus 0)\|^2+\|\Phi \Phi^*(g\oplus 0)\|^2\\
&=\|P_{\HH_{f,T,J}}(g\oplus 0)\|^2+\|\Theta_{f,T,J}^*g\|^2
\end{split}
\end{equation*}
for any $g\in \cN_J\otimes \cD_{C(T)}$. Hence, and using
\eqref{Ke-Te}, we deduce \eqref{K*P}.

Since $(T_1,\ldots, T_n)$ is  c.n.c. $n$-tuple in $\cV_{f,J}(\cH)$,
Proposition \ref{cnc} shows that $K_{f,T}$ is  a one-to-one
operator. On the other hand, $K_{f,T,J}=(P_{\cN_J}\otimes \cH)
K_{f,T}$ and the range of $K_{f,T}$  is included in the subspace
$\cN_J\otimes \overline{\Delta_{f,T}(\cH)}$.
 Consequently, the constrained
Poisson kernel $K_{f,T,J}$ is also one-to-one and $\text{\rm
range}\, K_{f,T,J}^*$ is dense in $\cH$. On the other hand,  let
$x\in \HH_{f,T,J}$ and assume that $x\perp P_{\HH_{f,T,J}}(g\oplus
0)$ for any $g\in \cN_J\otimes \cD_{C(T)}$. Using the definition of
$\HH_{f,T,J}$ and the fact that $\KK_{f,T,J}$ coincides with the
span
$$
\left\{g\oplus 0:\ g\in \cN_J\otimes \cD_{C(T)}\right\}\bigvee
\left\{\Theta_{f,T,J} \varphi\oplus \Delta_{\Theta_{f,T,J}}
\varphi,\quad \varphi\in \cN_J\otimes \cD_{C(T)^*}\right\},
$$
we deduce that $x=0$. This shows that
$$
\HH_{f,T,J}=\left\{P_{\HH_{f,T,J}}(g\oplus 0):\ g\in \cN_J\otimes
\cD_{C(T)}\right\}
$$
Hence, and using  relation \eqref{K*P},  we infer that there is a
unique unitary operator $\Gamma$ satisfying relation \eqref{Ga}.

For each $i=1,\ldots, n$, let $\TT_i:\HH_{f,T,J}\to \HH_{f,T,J}$ be
the transform of $T_i$ under the unitary operator $\Gamma:\cH\to
\HH_{f,T,J}$ defined  by \eqref{Ga}, that is,
$$\TT_i:=\Gamma
T_i\Gamma^*,\quad  i=1,\ldots, n.
$$
We prove now that, for each $i=1,\ldots, n$,
\begin{equation}\label{PN-intert}
 \left( P_{\cN_J\otimes \cD_{C(T)}}|_{\HH_{f,T,J}}\right) \TT_i^*x=
(B_i^*\otimes I_{\overline{\Delta_{f,T}\cH}})\left( P_{\cN_J\otimes
\cD_{C(T)}}|_{\HH_{f,T,J}}\right)x
\end{equation}
for any  $x\in \HH_{f,T,J}$
First, notice that using  Theorem \ref{J-factor}, relation \ref{fi},
and the fact that $\Phi$ is an isometry, we obtain
\begin{equation*}
\begin{split}
P_{\cN_J\otimes \cD_{C(T)}} \Gamma K_{f,T,J}^* g&= P_{\cN_J\otimes
\cD_{C(T)}} P_{\HH_{f,T,J}}(g\oplus 0)
 =
g-P_{\cN_J\otimes \cD_{C(T)}} \Phi \Phi^*(g\oplus 0)\\
&=g-\Theta_{f,T,J} \Theta_{f,T,J}^* g =K_{f,T,J} K_{f,T,J}^*g
\end{split}
\end{equation*}
for any $g\in \cN_J\otimes \cD_{C(T)}$. Consequently, using the fact
that  the range of $K_{f,T,J}^*$ is dense in $\cH$, we deduce that
\begin{equation}
\label{PGK} P_{\cN_J\otimes \cD_{C(T)}} \Gamma=K_{f,T,J}.
\end{equation}
Hence,  and taking into account that the constrained Poisson kernel
$K_{f,T,J}$ is one-to-one, we  can conclude that
\begin{equation}
\label{PKGa} P_{\cN_J\otimes \cD_{C(T)}} |_{\HH_{f,T,J}}=K_{f,T,J}
\Gamma^*
\end{equation}
is a one-to-one operator acting from $\HH_{f,T,J}$ to $\cN_J\otimes
\cD_{C(T)}$. Now, using relation \eqref{PGK} and  the properties of
the constrained Poisson kernel, we have
\begin{equation*}
\begin{split}
\left(P_{\cN_J\otimes \cD_{C(T)}} |_{\HH_{f,T,J}}\right)
\TT_i^*\Gamma h&= \left(P_{\cN_J\otimes \cD_{C(T)}}
|_{\HH_{f,T,J}}\right) \Gamma T_i^* h =K_{f,T,J} T_i^*h\\&=
\left( B_i^*\otimes I_{\overline{\Delta_{f,T}\cH}}\right) K_{f,T,J}h\\
&= \left( B_i^*\otimes I_{\overline{\Delta_{f,T}\cH}}\right)
\left(P_{\cN_J\otimes \cD_{C(T)}} |_{\HH_{f,T,J}}\right)\Gamma h
\end{split}
\end{equation*}
for any $h\in \cH$. Hence, we deduce   relation \eqref{PN-intert}.
We remark that, since the operator $P_{\cN_J\otimes
\cD_{C(T)}}|_{\HH_{f,T,J}}$ is one-to-one (see \eqref{PKGa}), the
relation \eqref{PN-intert} uniquely determines each operator
$\TT_i^*$, \ $i=1,\ldots, n$.

To prove that last part of the theorem, assume that $T:=(T_1,\ldots,
T_n)$ is a pure $n$-tuple in $\cV_{f,J}(\cH)$.
 Due to the proof of Theorem \ref{funct-calc2}, the constrained Poisson kernel
  $K_{f,T,J}:\cH\to \cN_J\otimes \overline{\Delta_{f,T}\cH}$
 is an isometry.
Consequently, $K_{f,T,J}K_{f,T,J}^*$ is the orthogonal projection of
$\cN_J\otimes \overline{\Delta_{f,T}\cH}$ onto $K_{f,T,J}\cH$. Using
Theorem \ref{J-factor}, we have
$$
K_{f,T,J}K_{f,T,J}^*+\Theta_{f,T,J}\Theta_{f,T,J}^*= I_{\cN_J\otimes
\overline{\Delta_{f,T}\cH}}.
$$
  This shows that $K_{f,T,J}K_{f,T,J}^*$ and
$\Theta_{f,T,J}\Theta_{f,T,J}^*$ are mutually orthogonal
projections. Therefore, $\Theta_{f,T,J}$ is a partial isometry,
i.e., an inner multi-analytic operator, with respect to $B_1,\ldots,
B_n$. Consequently,   $\Theta_{f,T,J}^* \Theta_{f,T,J}$ is a
projection. This implies that $\Delta_{\Theta_{f,T,J}}$ is the
projection on the orthogonal complement of $\text{\rm range}\,
\Theta_{f,T,J}^*$.

Note that  $u\oplus v\in \KK_{f,T,J}$ is in the subspace
$\HH_{f,T,J}$ if and only if
$$
\left< u\oplus v,\Theta_{f,T,J}\varphi\oplus
\Delta_{\Theta_{f,T,J}}\varphi\right>=0
$$
for any $\varphi\in \cN_J\otimes \cD_{C(T)^*}$, which is equivalent
to
\begin{equation}
\label{TH*} \Theta_{f,T,J}^*u+\Delta_{\Theta_{f,T,J}}v=0.
\end{equation}
According to the above observations, we have $\Theta_{f,T,J}^*u\perp
\Delta_{\Theta_{f,T,J}}v$. Consequently,  relation \eqref{TH*} holds
if and only if $\Theta_{f,T,J}^*u=0$ and $v=0$. Therefore,
$$\HH_{f,T,J}=\left(\cN_J\otimes \cD_{C(T)}\right)\ominus
\Theta_{f,T,J}(\cN_J\otimes \cD_{C(T)^*}).
$$
Notice also that, $P_{\cN_J\otimes \cD_{C(T)}}|_{\HH_{f,T,J}}$ is
the restriction operator and relation \eqref{PN-intert} implies
$$
\TT_i= P_ {\HH_{f,T,J}} (B_i\otimes I_{\cD_{C(T)}})|\HH_{f,T,J},
\quad i=1,\ldots, n.
$$
Conversely,  if we assume that  $\Theta_{f,T,J}$ is inner, then it
is a partial isometry. Once again, Theorem \ref{J-factor} implies
that $K_{f,T,J}$ is a partial isometry. On the other hand,
 since $T$ is
c.n.c., Proposition \ref{cnc}  shows that $K_{f,T}$ and,
consequently, $K_{f,T,J}$ are one-to-one partial isometries, and
therefore isometries. Due to  relation \eqref{K*K}, we have
$$
K_{f,T,J}^*K_{f,T,J}=K_{f,T}^*K_{f,T}=I_\cH-\text{\rm
SOT-}\,\lim_{k\to\infty} \Phi_{f,T}^k(I).
$$
Consequently, we have $\text{\rm SOT-}\,\lim_{k\to\infty}
\Phi_{f,T}^k(I)=0$, which proves that $T$ is a pure $n$-tuple of
operators with respect to $\cD_f(\cH)$.
    This completes the proof.
\end{proof}

Let  $\Phi\in R_n^\infty(\cV_{f,J})\bar\otimes B(\cK_1, \cK_2)$ and
$\Phi'\in R_n^\infty(\cV_{f,J})\bar\otimes B(\cK_1', \cK_2')$ be two
multi-analytic operators with respect to the constrained shifts
$B_1,\ldots, B_n$ associated with $\cV_{f,J}$. We say that $\Phi$
and $\Phi'$ coincide if there are two unitary operators $\tau_j\in
B(\cK_j, \cK_j')$ such that
$$
\Phi'(I_{\cN_J}\otimes \tau_1)=(I_{\cN_J}\otimes \tau_2) \Phi.
$$

Now we can show  that the constrained characteristic function
$\Theta_{f,T,J}$  is a complete unitary invariant for c.n.c.
$n$-tuples of operators in the noncommutative variety $\cV_{f,J}$.

\begin{theorem}\label{u-inv}
Let $J\neq F_n^\infty(\cD_f)$ be a $w^*$-closed two-sided ideal of
the Hardy algebra $F_n^\infty(\cD_f)$ and let $T:=(T_1,\ldots,
T_n)\in \cV_{f,J}(\cH)$ and $T':=(T_1',\ldots, T_n')\in
\cV_{f,J}(\cH')$ be two c.n.c. $n$-tuples of operators. Then $T$ and
$T'$ are unitarily equivalent if and only if their constrained
characteristic functions $\Theta_{f,T,J}$ and $\Theta_{f,T',J}$
coincide.
\end{theorem}

\begin{proof}
First assume that  the $n$-tuples $T$ and $T'$ are unitarily
equivalent and let $U:\cH\to \cH'$ be a unitary operator such that
$T_i=U^*T_i'U$ for any $i=1,\ldots, n$. We recall that
$\Gamma:=\{\alpha\in \FF_n^+: \ a_\alpha \neq 0\}$, \  $N:=\text{\rm
card}\, \Gamma$,  and
 $$
 C(T):=[\sqrt{a_{\tilde\alpha}}T_{\tilde\alpha}: \ \alpha\in
 \Gamma]$$
 is the the row operator
 with the entries are arranged in the lexicographic order of
 $\Gamma\subset \FF_n^+$. Note that  $C(T)$ is an operator acting from
 $\cH^{(N)}$ (the direct sum of $N$ copies of $\cH$)  to $\cH$.
Since $\Delta_{C(T)}^2=I-\sum_{|\alpha|\geq 1} a_\alpha T_\alpha
T_\alpha^*$,  we deduce  that $ U\Delta_{C(T)}=\Delta_{C(T')}U $.
Similarly, we obtain  $(\oplus_{i=1}^N
U)\Delta_{C(T)^*}=\Delta_{C(T)'^*}(\oplus_{i=1}^N U)$.
 Define the
unitary operators $\tau$ and $\tau'$ by setting
$$\tau:=U|\cD_{C(T)}:\cD_{C(T)}\to \cD_{C(T')} \ \text{ and }\
\tau':=(\oplus_{i=1}^N U)|\cD_{C(T)^*}:\cD_{C(T)*}\to \cD_{C(T')^*}.
$$
Using the definition of the constrained characteristic function, it
is easy to show that
$$
(I_{\cN_J}\otimes
\tau)\Theta_{f,T,J}=\Theta_{f,T',J}(I_{\cN_J}\otimes \tau').
$$

Conversely, assume that the constrained characteristic functions  of
$T$ and $T'$ coincide.  Then  there exist unitary operators
$\tau:\cD_{C(T)}\to \cD_{C(T')}$ and $\tau_*:\cD_{C(T)^*}\to
\cD_{{C(T')}^*}$ such that
\begin{equation}\label{com}
(I_{\cN_J}\otimes \tau)\Phi_{f,T,J}=\Phi_{f,T',J}(I_{\cN_J}\otimes
\tau_*).
\end{equation}
It is clear that relation \eqref{com} implies
$$
\Delta_{\Phi_{f,T,J}}=\left(I_{\cN_J}\otimes \tau_*\right)^*
\Delta_{\Phi_{f,T',J}}\left(I_{\cN_J}\otimes \tau_*\right)
$$
and
$$
\left(I_{\cN_J}\otimes
\tau_*\right)\overline{\Delta_{\Phi_{f,T,J}}(\cN_J\otimes
\cD_{C(T)^*})}= \overline{\Delta_{\Phi_{f,T',J}}(\cN_J\otimes
\cD_{{C(T')}^*})}.
$$
 Define now the unitary operator $U:\KK_{f,T,J}\to \KK_{f,T',J}$
by setting
$$U:=(I_{\cN_J}\otimes \tau)\oplus (I_{\cN_J}\otimes \tau_*).
$$
Simple computations reveal that the operator $\Phi:\cN_J\otimes
\cD_{C(T)^*}\to \KK_{f,T,J}$, defined  by $$ \Phi\varphi:=
\Theta_{f,T,J}\varphi \oplus \Delta_{\Theta_{f,T,J}}\varphi,\quad
\varphi\in \cN_J\otimes \cD_{C(T)^*}, $$
 and the corresponding $\Phi'$ satisfy the following relations:
\begin{equation}
\label{Uni1} U \Phi\left(I_{\cN_J}\otimes \tau_*\right)^*=\Phi'
\end{equation}
and
\begin{equation}
\label{Uni2} \left(I_{\cN_J}\otimes \tau\right) P_{\cN_J\otimes
\cD_{C(T)}}^{\KK_{f,T,J}} U^*=P_{\cN_J\otimes
\cD_{C(T')}}^{\KK_{f,T',J}},
\end{equation}
where $P_{\cN_J\otimes \cD_{C(T)}}^{\KK_{f,T,J}}$ is the orthogonal
projection of $\KK_{f,T,J}$ onto $\cN_J\otimes \cD_{C(T)}$. Note
also that relation \eqref{Uni1} implies
\begin{equation*}
\begin{split}
U\HH_{f,T,J}&=U\KK_{f,T,J}\ominus U\Phi(\cN_J\otimes \cD_{C(T)^*})\\
&=\KK_{f,T',J}\ominus \Phi'(I_{\cN_J}\otimes \tau_*)(\cN_J\otimes
 \cD_{C(T)^*})\\
&=\KK_{f,T',J}\ominus \Phi' (\cN_J\otimes \cD_{{C(T')}^*}).
\end{split}
\end{equation*}
Consequently, the operator $U|_{\HH_{f,T,J}}:\HH_{J,T}\to
\HH_{f,T',J}$ is unitary.

On the other hand, we have
\begin{equation}
\label{intertw} (B_i^*\otimes I_{\cD_{C(T')}})(I_{\cN_J}\otimes
\tau)= (I_{\cN_J}\otimes \tau)(B_i^*\otimes I_{\cD_{C(T)}}).
\end{equation}
Let $\TT:=(\TT_1,\ldots \TT_n)$ and $\TT':=(\TT_1',\ldots \TT_n')$
be the model operators provided by Theorem \ref{model}  for $T$ and
$T'$, respectively. Using the relation \eqref{PN-intert}  for $T'$
and $T$, as well as \eqref{Uni2} and \eqref{intertw}, we have
\begin{equation*}
\begin{split}
P_{\cN_J\otimes \cD_{C(T')}}^{\KK_{f,T',J}}{\TT_i'}^*Ux
&=  (B_i^*\otimes I_{\cD_{C(T')}}) P_{\cN_J\otimes
\cD_{C(T)}}^{\KK_{f,T,J}}Ux\\
&=(B_i^*\otimes I_{\cD_{C(T')}})(I_{\cN_J}\otimes \tau)
P_{\cN_J\otimes \cD_{C(T)}}^{\KK_{f,T,J}}x\\
&=(I_{\cN_J}\otimes \tau)(B_i^*\otimes I_{\cD_{C(T)}})
P_{\cN_J\otimes \cD_{C(T)}}^{\KK_{f,T,J}}x\\
&=(I_{\cN_J}\otimes \tau) P_{\cN_J\otimes \cD_{C(T)}}^{\KK_{f,T,J}}
\TT_i^*x\\
&= P_{\cN_J\otimes \cD_{C(T')}}^{\KK_{f,T',J}}U \TT_i^*x
\end{split}
\end{equation*}
for any $x\in \HH_{f,T,J}$ and $i=1,\ldots, n$. Using the fact that
$P_{\cN_J\otimes \cD_{C(T')}}^{\KK_{f,T',J}}|_{\HH_{f,T',J}}$ is an
one-to-one operator (see Theorem \ref{model}), we deduce that
$$
\left(U|_{\HH_{f,T,J}}\right)
\TT_i^*={\TT_i'}^*\left(U|_{\HH_{f,T,J}}\right),\quad i=1,\ldots,n.
$$
 Finally, using again Theorem \ref{model}, we conclude that the $n$-tuples
  $T$ and $T'$ are
 unitarily equivalent. The proof is complete.
\end{proof}

We remark that in the particular case when $J=\{0\}$, we obtain  the
characteristic function $\Theta_{f,T}$ and model theory for
$n$-tuples of operators in the noncommutative domain $\cD_f(\cH)$.

\smallskip

Now, we discuss the commutative case when $(T_1,\ldots, T_n)\in
\cD_f(\cH)$ and $$T_i T_j=T_jT_i\quad \text{ for } \quad
i,j=1,\ldots,n.
$$
Assume that $f=\sum_{|\alpha|\geq 1} a_\alpha X_\alpha$ is a
positive regular free holomorphic function on $B(\cH)^n$.  If $J_c$
is the $w^*$-closed two-sided ideal of $F_n^\infty(\cD_f)$ generated
by the polynomials
$$
 \{W_iW_j-W_jW_i:\ i,j=1,\ldots, n\},
$$
then $\cN_{J_c}=F^2_s(\cD_f)$, the symmetric weighted Fock space,
and
$$L_i:=P_{F_s^2(\cD_f)} W_i|_{F_s^2(\cD_f)}, \ i=1,\dots, n,
$$
 are the weighted creation operators on the symmetric  weighted
Fock space. We showed in  Section \ref{Symmetric},  that the
symmetric Fock space $F_s^2(\cD_f)$ can be identified with the
Hilbert space $H^2(\cD^\circ_f(\CC))$, which is the reproducing
kernel Hilbert space with reproducing kernel $K_f:
\cD_f^\circ(\CC)\times \cD_f^\circ(\CC)\to \CC$ defined by
 $$
 K_f(z,w):= \frac {1}
{1- \sum_{ |\alpha|\geq 1} a_\alpha z_\alpha\bar{w}_\alpha}, \qquad
z,w\in \cD_f^\circ(\CC),
$$
where
$$\cD_f^\circ(\CC):=\{z:=(z_1,\ldots, z_n)\in \CC^n:\
 \sum_{|\alpha|\geq 1}
a_\alpha |z_\alpha|^2 <1\}.
$$

 We proved in Section \ref{Symmetric} that
 $F_n^\infty(\cV_{f,J_c}):=P_{F_s^2(\cD_f)}
F_n^\infty(\cD_f)|_{F_s^2(\cD_f)},$
     is  the $w^*$-closed algebra
    generated by  the operators
   $L_i$, \ $i=1,\dots, n$, and the identity.
     Moreover,  we showed
   that $F_n^\infty(\cV_{f,J_c})$ can be identified with  the algebra of all
   multipliers  of $H^2(\cD_f^\circ(\CC))$.
   Under this identification, the weighted creation operators
    $L_1,\ldots, L_n$, generating $F_n^\infty(\cV_{f,J_c})$,  become the multiplication operators $M_{z_1},\ldots, M_{z_n}$
    by the coordinate functions $z_1,\ldots, z_n$ of $\CC^n$.
If the $n$-tuple of operators $T:=(T_1,\ldots, T_n)$ is in
$\cD_f(\cH)$, then $T\in \cV_{f,J_c}(\cH)$ if and only if
$$
T_iT_j=T_jT_i, \quad i,j=1,\ldots, n.
$$
Under the above-mentioned identifications, the constrained
characteristic function of $T\in \cV_{f,J_c}(\cH)$  is the
multiplication operator
$$M_{\Theta_{f,T,J_c}}:H^2(\cD_f^\circ(\CC))\otimes \cD_{C(T)^*}\to
H^2(\cD_f^\circ(\CC))\otimes \cD_{C(T)}$$ defined by
$$
\Theta_{f,T,J_c}(z):= -C(T)+\Delta_{C(T)}\left(I-\sum_{|\alpha|\geq
1} a_\alpha z_\alpha T_\alpha^*\right)^{-1}
\left[\sqrt{a_\alpha}z_\alpha I_\cH:\ \alpha \in \Gamma
\right]\Delta_{C(T)^*}
$$
for $z\in \cD_f^\circ(\CC)$.

\begin{remark} All the  results of this section can be written  in the particular
 case when $(T_1,\ldots,T_n)$ is a  completely non-coisometric $n$-tuple of commuting
  operators in $\cD_f(\cH)$.
\end{remark}

A few  comments about the results of this section are necessary.
Note that, in particular, if $p=\sum_{1\leq|\alpha|\leq m} a_\alpha
X_\alpha$
 is a positive
regular polynomial, then $\cD_p^\circ(\CC)$ is a Reinhardt domain in
$\CC^n$. We should remark that we can recover, in this  particular
setting, the results from \cite{BS}   when $(T_1,\ldots, T_n)\in
\cD_p(\cH)$ is a pure $n$-tuple of commuting operators.

When $p=X_1+\cdots + X_n$ and $(T_1,\ldots, T_n)\in \cD_p(\cH)$ is a
commuting $n$-tuple of operators, we obtain the characteristic
function introduced in \cite{Po-varieties} and \cite{BES1}, and
further  studied in \cite{Po-varieties2}, \cite{BES2} and
\cite{BT1}.

\bigskip

\section{Curvature invariant for $n$-tuples of operators in $\cD_p$}
\label{Curvature}

We introduce the curvature and $*$-curvature associated with
$n$-tuples of operators $T:=(T_1,\ldots, T_n)$ in the noncommutative
domain $\cD_p(\cH)$,  where $p$ is a positive regular noncommutative
polynomial. We prove the existence of these numerical invariants and
present basic properties. We show that both the curvature and
$*$-curvature can be express in terms of the characteristic function
$\Theta_{p,T}$ of $T\in \cD_p(\cH)$. The particular case when
$\cD_{p_e}(\cH)$ is the noncommutative ellipsoid
$$
\left\{(X_1,\ldots, X_n)\in B(\cH)^n: \ a_1 X_1X_1^*+\cdots+a_n
X_nX_n^*\leq I\right\}
$$
is discussed in greater details.

Throughout this section, we assume that  $p:=\sum_{1\leq|\alpha|\leq
m} a_\alpha X_\alpha$
 is a positive
regular polynomial, i.e., $a_\alpha\geq 0$ and  $a_\alpha>0$ if
$|\alpha|=1$.
As in the previous sections, we associate with the
 $n$-tuple of  operators $T:=(T_1,\dots, T_n)\in \cD_p(\cH) $
 the completely positive linear map
 $$
 \Phi_{p,T}(X):= \sum_{1\leq|\alpha|\leq
m} a_\alpha T_\alpha X T_\alpha^*,
 \quad X\in B(\cH).
 $$
 The {\it adjoint} of $\Phi_{p,T}$ is  defined by
 $
 \Phi_{p,T}^*(X):= \sum_{1\leq|\alpha|\leq
m} a_\alpha T^*_\alpha X T_\alpha,
 \quad X\in B(\cH).
 $

\begin{lemma}\label{ine-trac}
If ~$T:=(T_1,\dots, T_n) $ is
  in the noncommutative domain $\cD_p(\cH)$,
  then
\begin{equation}\label{a/b}
\text{\rm trace}\,\Phi_{p,T}(X)\leq \|\Phi_{p,T}^*(I)\| \text{\rm
trace}\,(X)\leq \left(\sum_{1\leq |\alpha|\leq m}
\frac{a_\alpha}{b_\alpha}\right)\text{\rm trace}\,(X)
 \end{equation}
 for any  positive trace class operator $X$,
  where the coefficients $b_\alpha$, $\alpha\in \FF_n^+$, are
  defined by \eqref{b_alpha}.
\end{lemma}

\begin{proof}
If $X\geq 0$ is a trace class operator, then
\begin{equation*}
\begin{split}
\text{\rm trace}\,\Phi_{p,T}(X)&= \text{\rm trace}\, \left(\sum_{
1\leq |\alpha|\leq m}X^{1/2} a_\alpha T_\alpha^* T_\alpha
X^{1/2}\right)=\text{\rm trace}\,
[X^{1/2} \Phi_{p,T}^*(I) X^{1/2}]\\
&\leq \|\Phi_{p,T}^*(I)\| \text{\rm trace}\, X.
\end{split}
\end{equation*}
Using  the von Neumann inequality of Theorem \ref{Poisson-C*},  we
have
\begin{equation*}
\begin{split}
\|\Phi_{p,T}^*(I)\|&= \|[\sqrt{a_\alpha} T_\alpha^*:\ 1\leq
|\alpha|\leq m]\|^2\\
&\leq \|[\sqrt{a_\alpha} W_\alpha^*:\ 1\leq
|\alpha|\leq m]\|^2\\
&=\left\| \sum_{ 1\leq |\alpha|\leq m} a_\alpha W_\alpha^*
W_\alpha\right\|\\
&\leq \sum_{1\leq |\alpha|\leq m} \frac{a_\alpha}{b_\alpha}
\end{split}
\end{equation*}
The latter inequality is due to the fact that
$\|W_\alpha\|=\frac{1}{\sqrt{b_\alpha}}$, $\alpha\in \FF_n^+$ (see
Section \ref{Noncommutative}). The proof is complete.
\end{proof}

 In what follows,
we define the curvature  $\text{\rm curv}_{p}:\cD_p(\cH)\to
[0,\infty)$   and show that it exists.

 \begin{theorem}\label{curva}Let $p:=\sum_{1\leq |\alpha|\leq m}
 a_\alpha X_\alpha $ be a positive regular noncommutative
 polynomial.
If ~$T:=(T_1,\dots, T_n) $ is
  in the noncommutative domain $\cD_p(\cH)$,
  then
\begin{equation}\label{curv}
 \text{\rm curv}_{p}(T):=
 \lim_{k\to\infty}
  {\frac {\text{\rm trace}\,\left\{K_{p,T}^* \left[ I-\Phi_{p,W}^{k+1}(I)
  \otimes I_{\overline{\Delta_{p,T}\cH}}\right]
   K_{p,T}\right\}} {\sum\limits_{j=0}^k\left(\sum\limits_{1\leq|\alpha|\leq m}
   \frac{a_\alpha}{b_\alpha}\right)^j}}
 \end{equation}
 exists, where   $K_{p,T}$ is the Poisson
 kernel associated with $T$ and the defect operator
 $\Delta_{p,T}:=\left(I-\Phi_{p,T}(I)\right)^{1/2}$.
 Moreover,
 $$\text{\rm curv}_{p}(T)<\infty \ \text{ if and only if }\
 \text{\rm trace}(I-\Phi_{p,T}(I))<\infty.
 $$
\end{theorem}
\begin{proof}
Due to the properties of the Poisson kernel $K_{p,T}$ (see Section
\ref{domain algebra}), we have $K_{p,T} T_i^*=(W_i^*\otimes I_\cH)
K_{p,T}$, $i=1,\ldots, n$, and $K_{p,T}^* K_{p,T} =I-Q_{p,T}$, where
$Q_{p,T}:=\lim_{m\to \infty} \Phi_{p,T}^m (I)$. Notice also that
$$
K_{p,T}^*(\Phi_{p,T}^{k+1}(I)\otimes I)K_{p,T}= \Phi_{p,T}^{k+1}
(I-Q_{p,T})= \Phi_{p,T}^{k+1}(I)-Q_{p,T}.
$$
Hence, we have
 \begin{equation*}\begin{split}
K_{p,T}^*  [ I-\Phi_{p,W}^{k+1}(I)\otimes
I_{\overline{\Delta_{p,T}\cH}}]
   K_{p,T}&=
   K_{p,T}^* K_{p,T}- K_{p,T}^*[\Phi_{p,W}^{k+1} (I)\otimes
    I_{\overline{\Delta_{p,T}\cH}}]
   K_{p,T}\\
   &= I-Q_{p,T}-\Phi_{p,T}^{k+1}(I) + Q_{p,T}\\
   &=  I-\Phi_{p,T}^{k+1}(I).
\end{split}
\end{equation*}
Since the sequence of positive operators $\{
I-\Phi_{p,T}^{k+1}(I)\}_{k=1}^\infty$ is increasing,
 it is clear that
 ~$\text{\rm curv}_p(T)=\infty$  whenever
 ~$\text{\rm trace}(I-\Phi_{p,T}(I))=\infty$.

 Assume now that
 ~$\text{\rm trace}(I-\Phi_{p,T}(I))<\infty$ and
denote $\gamma:=\sum_{1\leq |\alpha|\leq m}
\frac{a_\alpha}{b_\alpha}$. Due to relation \eqref{b_alpha}, we have
$b_{g_i}=a_{g_i}$, $i=1,\ldots,n$. Consequently, $\gamma\geq n$ and
the equality holds if and only if the positive regular  polynomial
$p$ has the form $a_1X_1+\cdots +a_nX_n$.

First we consider the case when $\gamma >1$.
Using the definition  \eqref{curv} and the calculations above, we
have
\begin{equation}\label{cur-new}
 \text{\rm curv}_p(T)=( \gamma-1)
 \lim_{k\to\infty}
  {\frac {\text{\rm trace}\,[ I-\Phi_{p,T}^k(I)]} {  \gamma^k}}.
 \end{equation}
We show  now that this limit exists.
  Since
\begin{equation}\label{dphi}
I-\Phi_{p,T}^{k+1}(I)=
I-\Phi_{p,T}(I)+\Phi_{p,T}(I-\Phi_{p,T}^k(I)),\ k=1,2,\ldots,
\end{equation}
we infer that $I-\Phi_{p,T}^{k+1}(I)$  is a trace class operator.
 From  Lemma \ref{ine-trac} and relation
and \eqref{dphi}, we obtain
\begin{equation}\label{tr-ine}
\text{\rm trace}\,[ I-\Phi_{p,T}^{k+1}(I)]\leq  \gamma
   \text{\rm trace}\, [ I-\Phi_{p,T}^{k}(I)]+  \text{\rm trace}\,
   [I-\Phi_{p,T}(I)].
\end{equation}
Therefore, setting
$$
x_k:= {\frac {\text{\rm trace}\, [  I-\Phi_{p,T}^{k}(I)]} {
\gamma^k}}- {\frac {\text{\rm trace}\, [ I-\Phi_{p,T}^{k-1}(I)]} {
\gamma^{k-1}}},
$$
relation \eqref{tr-ine} implies
$$
x_k\leq {\frac {\text{\rm trace}\, [ I-\Phi_{p,T}(I)]} { \gamma^k}},
\quad k=1,2,\ldots.
$$
Since $ \gamma>1$, it is clear  that $\sum\limits_{k:\, x_k\geq 0}
x_k <\infty$. Furthermore,
 every partial sum of the negative $x_k$'s is greater
 then or equal to ~$-\text{\rm trace}\, [I-\Phi_{p,T}(I)]$, so
 ~$\sum_{k:\, x_k<0} x_k$ converges.
Therefore, $\sum_{k} x_k$  is convergent and consequently
$$
\lim_{k\to \infty} {\frac {\text{\rm trace}\, [
I-\Phi_{p,T}^{k}(I)]} { \gamma^k}}
$$
exists, which, due to relation \eqref{cur-new}, implies the
existence of the limit defining the curvature.

Now, we consider the case when $\gamma=1$. Using  again Lemma
\ref{ine-trac},
   we get
 \begin{equation}\label{ine-ine}
 0\leq  \text{\rm trace}\, [\Phi_{p,T}^{k+1}( I-\Phi_{p,T}(I))]\leq
    \text{\rm trace}\, [\Phi_{p,T}^{k}( I-\Phi_{p,T}(I))]
 \end{equation}
 for any $k=0,1, \ldots$.
 Hence, the sequence
  $\{\text{\rm trace}\, [\Phi_{p,T}^{k}( I-\Phi_{p,T}(I))]\}_{k=0}^\infty$
  is decreasing and
   $$
   \lim_{k\to\infty} \text{\rm trace}\, [\Phi_{p,T}^{k}( I-\Phi_{p,T}(I))]
   $$
    exists.
      Consequently, using
      an elementary classical result and   the calculations presented in
      the beginning of the proof,  we obtain
    \begin{equation*}\begin{split}
 \text{\rm curv}_p(T):&=
 \lim_{k\to\infty}
  \frac {\text{\rm trace}\,\{K_{p,T}^*  [ I-\Phi_{p,W}^{k}(I)\otimes
  I_{\overline{\Delta_{p,T}\cH}}]
   K_{p,T}\}} { k}\\
  &=\lim_{k\to\infty}
   \text{\rm trace}\,\{K_{\varphi, T}^* [ (\Phi_{p,W}^{k-1}(I)
   -\Phi_{p,W}^{k}(I))\otimes I_{\overline{\Delta_{p,T}\cH}}]
   K_{p,T}\} \\
  &= \lim_{k\to \infty} \text{\rm trace}\, [\Phi_{p,T}^{k-1}(
  I-\Phi_{p,T}(I))].
  \end{split}
      \end{equation*}
      The proof is complete.
 \end{proof}

\begin{corollary}\label{curvac}
 Let ~$T:=(T_1,\dots, T_n) $ be
  in the noncommutative domain $\cD_p(\cH)$   and let
  $\gamma:=\sum_{1\leq |\alpha|\leq m}
\frac{a_\alpha}{b_\alpha}$.
  Then
\begin{equation*}
 \text{\rm curv}_p(T)=
( \gamma-1)
 \lim\limits_{k\to\infty}
  {\frac {\text{\rm trace}\,[ I-\Phi_{p,T}^k(I)]} {  \gamma^k}}
  \quad \text{ if }\   \gamma>1,
  \end{equation*}
  and
  \begin{equation*} \begin{split}
  \text{\rm curv}_p(T)&=\lim_{k\to \infty} {\frac {\text{\rm trace}\, [
I-\Phi_{p,T}^{k}(I)]} { k}}\\
&=
  \lim\limits_{k\to \infty} \text{\rm trace}\,
  [\Phi_{p,T}^{k}(I-\Phi_{p,T}(I))]
\end{split}
 \end{equation*}
  if    $\gamma=1$.
 Moreover,
 $$
\text{\rm curv}_p(T)\leq \text{\rm trace}\, [I-\Phi_{p,T}(I)] \leq
\text{\rm rank}\,[I-\Phi_{p,T}(I)].
$$
\end{corollary}
\begin{proof} The equalities above were   obtained in the proof of
Theorem \ref{curva}. To prove the last part,   assume first that
$\gamma>1$. Using the inequality \eqref{tr-ine}, we deduce
\begin{equation*}\begin{split}
{\frac {\text{\rm trace}\, [I-\Phi_{p,T}^{k+1}(I)]} { \gamma^{k+1}}}
&\leq \sum_{m=1}^{k+1} {\frac {\text{\rm trace}\, [I-\Phi_{p,T}(I)]}
{ \gamma^m}}\\
&= {\frac {\text{\rm trace}\, [I-\Phi_{p,T}(I)]} {
\gamma^{k+1}}}\cdot {\frac { \gamma^{k+1} -1} { \gamma-1}}.
\end{split}
\end{equation*}
Hence,   we  obtain the last part of the corollary. When $
\gamma=1$, one can use inequality \eqref{ine-ine} to complete the
proof.
\end{proof}

Due to  the Theorem \ref{curva} and Corollary \ref{curvac}, it is
easy to see that if we have $A:=(A_1,\ldots, A_n)\in \cD_p(\cH)$ and
$B:=(B_1,\ldots, B_n)\in \cD_p(\cK)$, then the direct sum $A\oplus
B:=(A_1\oplus B_1,\ldots, A_n\oplus B_n)\in \cD_p(\cH\oplus \cK)$
and
$$
\text{\rm curv}_{p}(A\oplus B)=\text{\rm curv}_{p}(A)+ \text{\rm
curv}_{p}(B).$$ In particular, if $\cK$ is finite dimensional, then
$\text{\rm curv}_{p}(A\oplus B)=\text{\rm curv}_{p}(A)$.

 We can  write now
the  curvature invariant in terms of the characteristic function
associated  with $T\in \cD_p(\cH)$.

\begin{theorem}\label{curva-charact}
If ~$T:=(T_1,\dots, T_n) $ is
  in the noncommutative domain $\cD_p(\cH)$ and rank $\Delta_{p,T}<\infty$,
  then
\begin{equation}\label{curv1}
 \text{\rm curv}_{p}(T)=
 \lim_{k\to\infty}
  {\frac {\text{\rm trace}\,\left[ \left(I-\Theta_{p,T} \Theta_{p,T}^*\right)
   \left( I-\Phi_{p,W}^{k+1}(I)\otimes I_{\overline{\Delta_{p,T}\cH}}
   \right)\right]
    } {\sum\limits_{j=0}^k\left(\sum\limits_{1\leq|\alpha|\leq m}
   \frac{a_\alpha}{b_\alpha}\right)^j}}
 \end{equation}
 exists, where   $\Theta_{p,T}$ is the   characteristic function
  associated with $T$.
\end{theorem}
\begin{proof}
Since $I-\Phi_{p,W}(I)=P_\CC$ (see Section \ref{Noncommutative}),
the operator
$$I-\Phi_{p,W}^{k+1}(I)= \sum_{m=0}^k \Phi_{p,W}^k(I-\Phi_{p,W}(I)),
\qquad k=0,1,\ldots$$ has finite rank. On the other hand, rank
$\Delta_{p,T}<\infty$ implies that $I-\Phi_{p,W}^{k+1}(I)\otimes
I_{\overline{\Delta_{p,T}\cH}}$ has also finite rank.
 Using the properties of the
trace and the factorization result of
 Theorem \ref{curva}, we deduce that
\begin{equation*}
\begin{split}
&\text{\rm trace}\,\left[K_{p,T}^*\left(
I\right.\right.\left.\left.-\Phi_{p,W}^{k+1}(I)\otimes
I_{\overline{\Delta_{p,T}\cH}}
   \right)K_{p,T}\right]\\
   &=\text{\rm trace}\,\left\{
\left( I-\Phi_{p,W}^{k+1}(I)\otimes I_{\overline{\Delta_{p,T}\cH}}
   \right)^{1/2}\left(I-\Theta_{p,T}
\Theta_{p,T}^*\right)\left( I-\Phi_{p,W}^{k+1}(I)\otimes
I_{\overline{\Delta_{p,T}\cH}}
   \right)^{1/2}\right\}\\
&= \text{\rm trace}\,\left[ \left(I-\Theta_{p,T}
\Theta_{p,T}^*\right)
   \left( I-\Phi_{p,W}^{k+1}(I)\otimes I_{\overline{\Delta_{p,T}\cH}}
   \right)\right].
\end{split}
\end{equation*}
According to  Theorem \ref{factor}, we have $I-\Theta_{p,T}
\Theta_{p,T}^*=K_{p,T} K_{p,T}^*$. Now, we can  complete the proof.
\end{proof}

  \begin{proposition}\label{ex1}
  Let \ $p_1:=a_1X_1+\cdots +a_nX_n$ be a positive regular  polynomial
  and let $W:=(W_1,\ldots, W_n)$ be the model operator associated
  with the noncommutative domain $\cD_{p_1}$. A pure $n$-tuple of
  operators $T:=(T_1,\ldots, T_n)\in \cD_{p_1}(\cH)$ is unitarily
  equaivalent to $(W_1\otimes I_\cK,\ldots, W_n\otimes I_\cK)$,
  where $\cK$ is
  a finite dimensional Hilbert space,  if and only if
  \begin{equation}\label{cur-dim}
\text{\rm curv}_{p_1}(T)=\text{\rm dim}
\left(I-\Phi_{p_1,T}(I)\right)^{1/2}\cH <\infty.
\end{equation}
  \end{proposition}
\begin{proof}
First, notice that due to the results of Section
\ref{Noncommutative}, $W_i=\frac{1}{\sqrt{a_i}}S_i$ for
$i=1,\ldots,n$, where $S_1,\ldots, S_n$ are the left creation
operators on the full Fock space $F^2(H_n)$. Consequently,
$\Phi_{p_1,W}^k(I)=\sum_{|\alpha|=k} S_\alpha S_\alpha^*$ for
$k=1,2,\ldots$, and $\left( I-\Phi_{p_1,W}(I)\right)=P_\CC$.

Assume that there is a unitary operator $U:\cH\to F^2(H_n)\otimes
\cK$ such that $T_i=U^*(W_i\otimes I_\cK)U$, $i=1,\ldots,n$. Since
$\Phi_{p_1,T}^k(I)=U^*\left[ \Phi_{p_1,W}(I)\otimes I_\cK\right]U$
for $k=1,2,\ldots$, we deduce that
 $
 \text{\rm dim}
\left(I-\Phi_{p_1,T}(I)\right)^{1/2}\cH=\text{\rm dim}\cK<\infty$
 and
$$
\text{\rm curv}_{p_1}(T)=
 \lim_{k\to\infty}
  {\frac {\text{\rm trace}\, \left[ \left(I- \sum_{|\alpha|=k}
  S_\alpha S_\alpha^*\right)\otimes
  I_\cK\right]
    } { 1+ n+ n^2+\cdots + n^k  }}=\text{\rm dim}\cK.
    $$

Conversely, assume that the $n$-tuple $T:=(T_1,\ldots, T_n)\in
\cD_{p_1}(\cH)$  is pure and satisfies the condition
\eqref{cur-dim}. If $Y_i:=\sqrt{a_i} T_i$, $i=1,\ldots, n$, then
$Y:=[Y_1,\ldots, Y_n]$ is a row contraction and $\text{\rm
curv}_{p_1}(T)=\text{\rm curv}(Y)$, where $\text{\rm curv}(Y)$ is
the curvature associated with a row contraction $Y$ (see
  \cite{Po-unitary}). Notice that due to relation \eqref{cur-dim}, we have
  $$
\text{\rm curv}(Y)=\text{\rm dim}(I-Y_1Y_1^*-\cdots
-Y_nY_n^*)^{1/2}\cH<\infty.
$$
Using  Theorem 3.4 from \cite{Po-unitary}, we deduce that
$Y:=[Y_1,\ldots, Y_n]$ is unitarily equivalent to
 $[S_1\otimes
I_\cK,\ldots, S_n\otimes I_\cK]$ for some finite dimensional Hilbert
space $\cK$. Therefore, $(T_1,\ldots, T_n)$ is unitarily equivalent
to
 $(W_1\otimes
I_\cK,\ldots, W_n\otimes I_\cK)$. This completes the proof.
\end{proof}

We introduce now the $*$-curvature associated with $T:=(T_1,\dots,
T_n)\in \cD_p(\cH)$.

 \begin{theorem}\label{curva*} If ~$T:=(T_1,\dots, T_n) $ is
  in the noncommutative domain $\cD_p(\cH)$,
  then
  \begin{equation}\label{curv**}
 \text{\rm curv}^*_{p}(T):=
 \lim_{k\to\infty}
  {\frac {\text{\rm trace}\,\{K_{p,T}^* [ I-\Phi_{p,W}^{k+1}(I)\otimes
  I]
   K_{p,T}\}} { 1+\|\Phi_{p,T}^*(I)\|+\cdots + \|\Phi_{p,T}^*(I)\|^k}}.
 \end{equation}
exists. Moreover, ~ $\text{\rm curv}^*_{p}(T)<\infty$ if and only if
 ~$\text{\rm trace}(I-\Phi_{p,T}(I))<\infty$.
\end{theorem}
\begin{proof}
 When $\|\Phi_{p,T}^*(I)\|\geq 1$, the proof is similar to the
 proof of Theorem \ref{curva}. All we have to do is to replace
 $\gamma$ by $\|\Phi_{p,T}^*(I)\|$ and use Lemma \ref{ine-trac}.

 Now,  assume $\|\Phi_{p,T}^*(I)\|<1$.
  Since
\begin{equation*}
I-\Phi_{p,T}^{k+1}(I)=
I-\Phi_{p,T}(I)+\Phi_{p,T}(I-\Phi_{p,T}^k(I)),\ k=1,2,\ldots,
\end{equation*}
Lemma \ref{ine-trac} implies
\begin{equation}\label{tr-ine2}
\text{\rm trace}\,[ I-\Phi_{p,T}^{k+1}(I)]\leq \|\Phi_{p,T}^*(I)\|
   \text{\rm trace}\, [ I-\Phi_{p,T}^{k}(I)]+  \text{\rm trace}\,
   [I-\Phi_{p,T}(I)].
\end{equation}
 Iterating relation \eqref{tr-ine2}, we deduce that the sequence
  $\{\text{\rm trace}\, [ I-\Phi_{p,T}^{k}(I)]\}_{k=1}^\infty$ is bounded.
  Since  the sequence  of operators$\{ I-\Phi_{p,T}^{k}(I)\}_{k=1}^\infty$ is increasing, we infer that
    $\lim_{k\to \infty} \text{\rm trace}\, [I-\Phi_{p,T}^{k}(I)]$
  exists.
According to the definition
 \eqref{curv**} and the calculations at the beginning of the proof
  of Theorem \ref{curva}, we get
 $$
  \text{\rm curv}_p^*(T)=(1-\|\Phi_{p,T}^*(I)\|)
  \lim_{k\to \infty} \text{\rm trace}\, [I-\Phi_{p,T}^{k}(I)].
  $$
 The proof is complete.
\end{proof}

The $*$-curvature can be expressed in terms of the characteristic
function of $T\in\cD_p(\cH)$. A result similar to Theorem
\ref{curva-charact} holds.  Since the proof is essentially the same,
we shall omit it.

\begin{corollary}\label{curvac2}
 If ~$T:=(T_1,\dots, T_n) $ is
  in the noncommutative domain $\cD_p(\cH)$,
  then
\begin{equation*}
 \text{\rm curv}_p^*(T)=
 \begin{cases}
 (\|\Phi_{p,T}^*(I)\|-1)
 \lim\limits_{k\to\infty}
  {\frac {\text{\rm trace}\,[ I-\Phi_{p,T}^k(I)]} { \|\Phi_{p,T}^*(I)\|^k}}&
  \quad \text{ if }\  \|\Phi_{p,T}^*(I)\|>1,\\
  \lim\limits_{k\to \infty} \text{\rm trace}\,
  [\Phi_{p,T}^{k}(I-\Phi_{p,T}(I))]& \quad
  \text{ if }\  \|\Phi_{p,T}^*(I)\|=1,\\
  (1-\|\Phi_{p,T}^*(I)\|)
  \lim\limits_{k\to \infty} \text{\rm trace}\, [I-\Phi_{p,T}^{k}(I)]&
  \quad
  \text{ if } \ \|\Phi_{p,T}^*(I)\|<1.
  \end{cases}
 \end{equation*}
 Moreover, if $\|\Phi_{p,T}^*(I)\|\geq 1$, then
 $$
\text{\rm curv}_p^*(T)\leq \text{\rm trace}\, [I-\Phi_{p,T}(I)] \leq
\text{\rm rank}\,[I-\Phi_{p,T}(I)].
$$
\end{corollary}
\begin{proof}
The proof is similar to that of Corollary \ref{curvac}.
\end{proof}

 Let $\Lambda$ be a nonempty  set of positive numbers $t>0$ such that
 $$
 \text{\rm trace}\,\Phi_{p,T}(X)\leq t \,\text{\rm trace}\,X
 $$
 for any positive trace class operator $X\in B(\cH)$. Notice that
    Theorem \ref{curva*} and  Corollary
  \ref{curvac2} remain true  (with exactly the same proofs)
   if we replace $\|\Phi_{p,T}^*(I)\|$ with  $t\in \Lambda$.
   The corresponding curvature is denoted by
   $\text{\rm curv}_{p}^t(T)$.

   When $d:=\inf \Lambda$,  the curvature
   $\text{\rm curv}_{p}^d(T)$ is called the distinguished
    curvature associated with $T\in \cD_p(\cH)$ and with respect to $\Lambda$.
  Now, using the analogues of Theorem \ref{curva*} and  Corollary
  \ref{curvac2}   for the curvatures  $\text{\rm curv}_{p}^t(T)$,
  $t\in \Lambda$, one can easily prove the following.
  If  $\text{\rm curv}_{p}^t(T)>0 $, then
  $$
   \text{\rm curv}_p^t(T)=
   \begin{cases}
    \text{\rm curv}_p^d(T)&\quad  \text{ if } ~t\geq 1,\\
    \frac{1-t} {1-d}
     \,\text{\rm curv}_p^d(T)&\quad  \text{ if } ~t< 1.
   \end{cases}
 $$
 As we will see latter,
 if $\text{\rm curv}_p^{t_0}(T)=0 $ for some $t_0\in \Lambda$, then,
 in general, $\text{\rm curv}_p^d(T)\neq 0 $. Therefore, the
 distinguished
    curvature $\text{\rm curv}_p^d(T)$ is a refinement of all the other
    curvatures $\text{\rm curv}_p^t(T)$, ~$t\in \Lambda$.

All the results of this section concerning the $*$-curvature  have
analogues (and similar proofs) for the distinguished
    curvature.

 \begin{corollary}\label{curv-curv}
   If  ~$\text{\rm trace}\,(I-\Phi_{p,T}(I))<\infty$ and
  $$
  0<\|\Phi_{p,T}^*(I)\|<\sum\limits_{1\leq |\alpha|\leq m}
  \frac{a_\alpha}{b_\alpha},
  $$
 then
  ~$\text{\rm curv}_p(T)=0$.
\end{corollary}
\begin{proof}
Assume $\text{\rm curv}_p(T)>0$ and  let $\gamma:=\sum\limits_{1\leq
|\alpha|\leq m}
  \frac{a_\alpha}{b_\alpha}$.   Due to  Theorem
\ref{curva*}, we  have
 $
  ~\text{\rm curv}_p^*(T)
 <\infty$, which implies $\frac{\text{\rm curv}_p^*(T)}{\text{\rm
 curv}_p(T)}<\infty$. On the other hand,
 due to the fact that $\|\Phi_{p,T}^*(I)\|<\gamma$,
 we can use  Theorem \ref{curva} and Theorem
\ref{curva*} to  deduce that
$$
\frac{\text{\rm curv}_p^*(T)}{\text{\rm
 curv}_p(T)}
 =\lim_{k\to\infty}\frac{1+\gamma+\cdots+\gamma^k}
 { 1+\|\Phi_{p,T}^*(I)\|+\cdots +
 \|\Phi_{p,T}^*(I)\|^k}=\infty,
 $$
 which is a contradiction. Therefore, we must have ~$\text{\rm
 curv}_p(T)=0$.
\end{proof}

Due to Corollary \ref{curv-curv}, the curvature invariant
${\text{\rm curv}_p(T)}$ does not
 distinguish  among the Hilbert modules
 over $\CC\FF_n^+$   with  $0<\|\Phi_{p,T}^*(I)\|<
 \sum\limits_{1\leq |\alpha|\leq m}
  \frac{a_\alpha}{b_\alpha}$.
   However, in this case,
   our $*$-curvature $\text{\rm curv}_p^*(T)$ in not zero in general.

\begin{proposition}\label{ex2}
  Let \ $p_1:=a_1X_1+\cdots +a_nX_n$ be a
  polynomial with $a_i>0$, $i=1,\ldots,n$.
  \begin{enumerate}
 \item[(i)]
  For any $t>0$ there exists $T:=(T_1,\ldots, T_n)\in
  \cD_{p_1}(\cH)$ such that $$\text{\rm curv}_{p_1}(T)=t.$$
  \item[(ii)]
For any $t>0$ there exists $A:=(A_1,\ldots, A_n)\in
  \cD_{p_1}(\cH)$ such that $$\text{\rm curv}_{p_1}(A)=0\quad \text{
  and }\quad
  \text{\rm curv}_{p_1}^*(A)=t.
  $$
  \end{enumerate}
\end{proposition}
\begin{proof}
According to Theorem 3.8 from \cite{Po-unitary}, there exists a row
contraction $Y:=[Y_1,\ldots, Y_n]\in B(\cH)^n$ such that $\text{\rm
dim}\,(I-Y_1Y_1^*-\cdots Y_nY_n^*)^{1/2}\cH<\infty$ and $\text{\rm
curv}(Y)=t$. Setting $T_i:=\frac{1}{\sqrt{a_i}}Y_i$, $i=1,\ldots,
n$, we have $T:=(T_1,\ldots, T_n)\in
  \cD_{p_1}(\cH)$. Due to  our Theorem \ref{curva} and Corollary 2.7
  from \cite{Po-unitary}, we deduce that
  $$
  \text{\rm curv}_{p_1}(T)=\text{\rm curv} (Y)=t.
  $$

To prove part (ii),  let $2\leq m\leq n-1$ and define the positive
regular polynomial $q:=c_1X_1+\cdots c_m X_m$, where $c_j:=a_j$ if
$j=1,\ldots, m-1$ and $c_m:=\sqrt{a_m+a_{m+1}+\cdots + a_n}$. Due to
part (i), we can find an $m$-tuple of operators $(T_1, \ldots,
T_m)\in \cD_q(\cH)$ such that $\text{\rm curv}_{q}(T)=t$. Using
Corollary \ref{curv-curv}, we deduce that $\|\Phi_{q,T}^*(I)\|=m$.
Now, define $A_i:=T_i$ if $i=1,\ldots, m-1$ and $A_j:=T_m$ if
$j=m,\ldots, n$. Notice that $\Phi_{q,T}(I)=\Phi_{p_1, A}(I)$ and
$\Phi_{q,T}^*(I)=\Phi_{p_1, A}^*(I)$. Since $\|\Phi_{q,T}^*(I)\|=m$,
we must have $\|\Phi_{p_1, A}^*(I)\|=m$. Hence, and using Theorem
\ref{curva*} and Theorem \ref{curva}, we have
\begin{equation*}
\begin{split}
\text{\rm curv}_{p_1}^*(A)
 &= \lim_{k\to\infty}
  {\frac {\text{\rm trace}\,  [ I-\Phi_{p_1,A}^{k+1}(I) ]
    } { 1+\|\Phi_{p,T}^*(I)\|+\cdots + \|\Phi_{p,T}^*(I)\|^k}}\\
    &=
    \lim_{k\to\infty}
  {\frac {\text{\rm trace}\,  [ I-\Phi_{q,T}^{k+1}(I) ]
    } { 1+ m+m^2\cdots +  m^k}}\\
    &=\text{\rm curv}_{q}(T)=t
\end{split}
\end{equation*}
On the other hand, since $m<n$, we have
\begin{equation*}
\begin{split}
\text{\rm curv}_{p_1}(A)
 &=
 \lim_{k\to\infty}
  {\frac {\text{\rm trace}\,  [ I-\Phi_{p_1,A}^{k+1}(I) ]
    } { 1+ n+n^2\cdots +  n^k}}\\
&= \lim_{k\to\infty}\left[
  {\frac {\text{\rm trace}\,  [ I-\Phi_{q,T}^{k+1}(I) ]
    } { 1+ m+m^2\cdots +  m^k}}\cdot
    \frac{1+ m+m^2\cdots +  m^k}{1+ n+n^2\cdots +  n^k}\right]\\
    &=\text{\rm curv}_{q}(T)\cdot0=0.
\end{split}
\end{equation*}
This completes the proof.
\end{proof}

An important problem that remains open is whether  Proposition
\ref{ex1} and Proposition \ref{ex2} remain true for noncommutative
domains $\cD_p$ generated by arbitrary positive regular polynomials.

\chapter{Commutant lifting and applications}

\section{Interpolation on
noncommutative domains} \label{Interpolation}

In this section  we provide a  Sarason type commutant lifting
theorem for pure  $n$-tuples of operators in noncommutative domains
$\cD_f$ (resp. varieties $\cV_{f,J}$) and we obtain Nevanlinna-Pick
and Schur-Carath\' eodory type interpolation
 results.

Let $f(X_1,\ldots, X_n):= \sum_{\alpha\in \FF_n^+} a_\alpha
X_\alpha$, \ $a_\alpha\in \CC$,  be a positive regular free
holomorphic function on $B(\cH)^n$. According to Section
\ref{Noncommutative},  $f$ satisfies the conditions
$$
\limsup_{k\to\infty} \left( \sum_{|\alpha|=k}
|a_\alpha|^2\right)^{1/2k}<\infty,
$$
   $a_\alpha\geq 0$ for any $\alpha\in \FF_n^+$, \ $a_{g_0}=0$,
 \ and  $a_{g_i}>0$, $i=1,\ldots, n$.

Using Parrott's lemma \cite{Par}, one can prove the following
commutant lifting theorem which extends Arias' result \cite{Ar1} to
our more general setting. Since the proof follows the same lines we
should omit it.

\begin{theorem}\label{CLT}
Let $f$ be a positive regular free holomorphic function on
$B(\cH)^n$ and let $(W_1,\ldots, W_n)$ and $(\Lambda_1,\ldots
\Lambda_n)$ be the weighted  left (resp.right) creation operators
associated with the noncommutative domain $\cD_f$. For each $j=1,2$,
let $\cE_j\subset F^2(H_n)\otimes \cK_j$ be a co-invariant subspace
under each operator $W_j\otimes I_{\cK_j}$, $i=1,\ldots,n$. If
$X:\cE_1\to \cE_2$ is a bounded operator such that
$$
X[P_{\cE_1}(W_i\otimes I_{\cK_1})|_{\cE_1}]=[P_{\cE_2}(W_i\otimes
I_{\cK_2})|_{\cE_2}]X,\quad i=1,\ldots,n,
$$
then there exists
$$\Phi(\Lambda_1,\ldots, \Lambda_n)\in R_n^\infty(\cD_f)
\bar\otimes B(\cK_1,\cK_2)$$ such that $\Phi(\Lambda_1,\ldots,
\Lambda_n)^*\cE_2\subseteq \cE_1$,
$$
\Phi(\Lambda_1,\ldots, \Lambda_n)^*|\cE_2=X^*, \quad \text{ and
}\quad \|\Phi(\Lambda_1,\ldots, \Lambda_n)\|=\|X\|.
$$
\end{theorem}

We remark that a similar result to Theorem \ref{CLT} holds if we
start with co-invariant subspaces under $\Lambda_i\otimes
I_{\cK_j}$, $i=1,\ldots,n$.

 Let
   $J$ be a $w^*$-closed two-sided ideal of the Hardy algebra
    $F_n^\infty(\cD_f)$ and assume
that   $T:=(T_1,\ldots, T_n) $    is a pure the $n$-tuple  of
operators in the noncommutative variety
\begin{equation}\label{J}
\cV_{f,J}(\cH):=\left\{ X_1,\ldots, X_n)\in \cD_f(\cH):
\varphi(T_1,\ldots, T_n)=0 \quad \text{ for any } \ \varphi\in
J\right\}.
\end{equation}
  Due to Theorem \ref{dil3}, the   $n$-tuple   $(T_1,\ldots, T_n)\in \cV_{f,J}(\cH)$ is unitarily
equivalent to the compression of $(B_1\otimes I_\cK,\ldots,
B_n\otimes I_\cK)$ to a co-invariant subspace $\cE$ under each
operator $B_i\otimes I_\cK$, $i=1,\ldots, n$. Therefore, we have
$$T_i=P_\cE(B_i\otimes I_\cK)|\cE,\quad i=1,\ldots,n.
$$

As in our previous sections, define the subspaces of $F^2(H_n)$ by
$$
\cM_J:=\overline{JF^2(H_n)}\quad \text{and}\quad
\cN_J:=F^2(H_n)\ominus \cM_J,
$$
and  the {\it constrained  weighted left} (resp.~{\it right}) {\it
creation operators} associated with  $\cD_f$ and $J$ by setting
$$B_i:=P_{\cN_J} W_i|\cN_J \quad \text{and}\quad
C_i:=P_{\cN_J} \Lambda_i|\cN_J,\quad i=1,\ldots, n.
$$

The following result is  a Sarason type \cite{S} commutant lifting
theorem for pure   $n$-tuples of operators in the noncommutative
variety $\cV_{f,J}$.

\begin{theorem}\label{CLT2}
Let $J$ be a $w^*$-closed two-sided ideal of  the Hardy algebra
$F_n^\infty(\cD_f)$  and let $(B_1,\ldots, B_n)$ and $(C_1,\ldots,
C_n)$ be the corresponding constrained weighted left (resp. right)
creation operators associated with the noncommutative variety
$\cV_{f,J}$. For each $j=1,2$, let $\cK_j$ be a Hilbert space and
$\cE_j\subseteq \cN_J\otimes \cK_j$ be a co-invariant subspace under
each operator $B_i\otimes I_{\cK_j}$, \ $i=1,\ldots, n$. If
$X:\cE_1\to \cE_2$ is a bounded operator such that
\begin{equation}\label{int}
X[P_{\cE_1}(B_i\otimes I_{\cK_1})|_{\cE_1}]=[P_{\cE_2}(B_i\otimes
I_{\cK_2})|_{\cE_2}]X,\quad i=1,\ldots,n,
\end{equation}
then there exists
$$G(C_1,\ldots, C_n)\in R_n^\infty(\cV_{f,J})\bar\otimes B(\cK_1,\cK_2)$$
such that $G(C_1,\ldots, C_n)^*\cE_2\subseteq \cE_1$,
$$
G(C_1,\ldots, C_n)^*|\cE_2=X^*, \quad \text{ and }\quad
\|G(C_1,\ldots, C_n)\|=\|X\|.
$$
\end{theorem}

\begin{proof}
 The subspace  $\cN_J\otimes \cK_j$
is  invariant under  each operator $W_i^*\otimes I_{\cK_j}$, \
$i=1,\ldots, n$,  and
$$
(W_i^*\otimes I_{\cK_j})|\cN_J\otimes \cK_j=B_i^*\otimes
I_{\cK_j},\quad  i=1,\ldots, n.
$$
Since $\cE_j\subseteq \cN_J\otimes \cK_j$ is invariant   under
$B_i^*\otimes I_{\cK_j}$ it is also invariant under $W_i^*\otimes
I_{\cK_j}$  and therefore
$$
(W_i^*\otimes I_{\cK_j})|\cE_j=  (B_i^*\otimes I_{\cK_j})|\cE_j,
\quad i=1,\ldots, n.
$$
Consequently, relation \eqref{int} implies
\begin{equation}\label{int3}
XP_{\cE_1}(W_i\otimes I_{\cK_1})|_{\cE_1}=P_{\cE_2}(W_i\otimes
I_{\cK_2})|_{\cE_2}X,\quad i=1,\ldots,n.
\end{equation}
It is clear that, for each $j=1,2$, the $n$-tuple $(W_1 \otimes
I_{\cK_j},\ldots, W_n \otimes I_{\cK_j})$ is  a  dilation of the
$n$-tuple
$$
[P_{\cE_j}(W_1\otimes I_{\cK_j})|\cE_j,\ldots, P_{\cE_j}(W_n\otimes
I_{\cK_j})|\cE_j].
$$
 Due to  Theorem \ref{CLT}, we can find $\Phi(\Lambda_1,\ldots, \Lambda_n)\in
R_n^\infty(\cD_f)\bar\otimes B(\cK_1,\cK_2)$,  a multi-analytic
operator
 with respect to $W_1,\ldots, W_n$,
   such that
$\Phi(\Lambda_1,\ldots, \Lambda_n)^*\cE_2\subseteq \cE_1$,
\begin{equation}\label{clt}
\Phi(\Lambda_1,\ldots, \Lambda_n)^*|\cE_2=X^*,\quad \text{ and }
\quad \|\Phi(\Lambda_1,\ldots, \Lambda_n)\|=\|X\|.
\end{equation}

Now, let $G(C_1,\ldots, C_n):=P_{\cN_J\otimes \cK_2}
\Phi(\Lambda_1,\ldots, \Lambda_n)|\cN_J\otimes \cK_1$. According to
Section \ref{Weighted}, we have
$$
G(C_1,\ldots, C_n)\in [P_{\cN_J}R_n^\infty(\cD_f) |\cN_J]\bar
\otimes B(\cK_1, \cK_2) =R_n^\infty(\cV_{f,J})\bar\otimes B(\cK_1,
\cK_2).
$$
On the other hand, since $\Phi(\Lambda_1,\ldots,
\Lambda_n)^*(\cN_J\otimes \cK_2)\subseteq \cN_J\otimes \cK_1$ and
$\cE_j\subseteq \cN_J\otimes \cK_j$, relation \eqref{clt} implies
$$
G(C_1,\ldots, C_n)^*\cE_2\subseteq \cE_1\quad \text{ and }\quad
G(C_1,\ldots, C_n)^*|\cE_2=X^*.
$$
Using again  relation \eqref{clt}, we have
$$
\|X\|\leq \|G(C_1,\ldots, C_n)\|\leq \|\Phi(\Lambda_1,\ldots,
\Lambda_n)\|=\|X\|.
$$
Therefore, $\|G(C_1,\ldots, C_n)\|=\|X\|$. The proof is complete.
\end{proof}

We remark that in the  particular case when   $\cE_j:=\cG\otimes
\cK_j$, where $\cG $ is a co-invariant subspace under each  operator
$B_i$ and $C_i$, \ $i=1,\ldots, n$,  the  implication of Theorem
\ref{CLT2} becomes an equivalence. Indeed, the implication
``$\implies$'' is clear from   the same theorem. For the converse,
let $X=P_{\cG \otimes \cK_2} \Psi(C_1,\ldots, C_n)|\cG \otimes
\cK_1$, where $\Psi(C_1,\ldots, C_n)\in
R_n^\infty(\cV_{f,J})\bar\otimes B(\cK_1, \cK_2)$. Since
$B_iC_j=C_jB_i$ for $i,j=1,\ldots, n$,  we have
$$
(B_i^*\otimes I_{\cK_1})\Psi(C_1,\ldots, C_n)^*=\Psi(C_1,\ldots,
C_n)^* (B_i^*\otimes I_{\cK_2}), \quad i=1,\ldots, n.
$$
Now, taking into account that $\cG$ is  an invariant subspace under
each of the operators $B_i^*$  and $C_i^*$, \ $i=1,\ldots, n$, we
deduce relation \eqref{int}. This proves our assertion.

\begin{corollary}\label{commutant}
Let $J$ be a $w^*$-closed two-sided ideal of  the  Hardy  algebra
$F_n^\infty(\cD_f)$  and let $(B_1,\ldots, B_n)$ and $(C_1,\ldots,
C_n)$ be the corresponding constrained weighted left (resp. right)
creation operators   associated with the noncommutative variety
$\cV_{f,J}$. If $\cK$ is a Hilbert space and $\cG\subseteq \cN_J$ is
an invariant subspace under each of the operators $B_i^*$ and
$C_i^*$, \ $i=1,\ldots, n$, then
$$
\left\{\left[P_\cG \cW(B_1,\ldots, B_n)|\cG\right]\otimes
I_\cK\right\}^\prime=
 [P_\cG \cW(C_1,\ldots, C_n)|\cG]\bar\otimes B(\cK).
$$

\end{corollary}

If $f(X_1,\ldots, X_n):= \sum_{\alpha\in \FF_n^+} a_\alpha
X_\alpha$, \ $a_\alpha\in \CC$,  is a positive regular  free
holomorphic function on $B(\cH)^n$, we denote
$$
\cD_f^\circ(\CC):=\left\{ (\lambda_1,\ldots, \lambda_n)\in \CC^n: \
\sum_{|\alpha|\geq 1} a_\alpha |\lambda_\alpha|^2<1\right\},
$$
and
$$
\cV_{f,J}^\circ(\CC):=\left\{ (\lambda_1,\ldots, \lambda_n)\in
\cD_f^\circ(\CC): \ g(\lambda_1,\ldots \lambda_n)=0 \text{ for }
g\in J \right\}.
$$

We recall, from Section \ref{Symmetric}, that the mapping
$K:\cD_f^\circ(\CC)\otimes \cD_f^\circ(\CC)\to \CC$ defined by
$$
K_f(\mu,\lambda):=\frac{1}{1-\sum_{|\alpha|\geq 1} a_\alpha
\mu_\alpha \overline{\lambda}_\alpha}\quad \text{ for all }\
\lambda,\mu\in \cD_f^\circ(\CC)
$$
is a positive definite kernel on $\cD_f^\circ(\CC)$.

 Now we can obtain  the following  Nevanlinna-Pick  interpolation
 result  in our setting.

\begin{theorem}\label{Nev}
Let $J$ be  a $w^*$-closed two-sided ideal of the Hardy algebra
$F_n^\infty(\cD_f)$ and let $(B_1,\ldots, B_n)$ and $(C_1,\ldots,
C_n)$ be the corresponding constrained weighted left (resp. right)
creation operators  associated with the noncommutative variety
$\cV_{f,J}$. Let $\lambda_1,\ldots, \lambda_k$ be $k$ distinct
points in the variety $\cV^\circ_{f,J}(\CC)$ and let $A_1,\ldots,
A_k\in B(\cK)$.

Then there exists $\Phi(C_1,\ldots, C_n)\in
R_n^\infty(\cV_{f,J})\bar\otimes B(\cK)$ such that
$$\|\Phi(C_1,\ldots, C_n)\|\leq 1\quad \text{and }\quad
\Phi(\lambda_j)=A_j,\ j=1,\ldots, k,
$$
if and only if the operator matrix
\begin{equation}\label{semipo}
\left[(I_\cK-A_i A_j^*) K_f(\lambda_i,\lambda_j)\right]_{k\times k}
\end{equation}
is positive semidefinite.
\end{theorem}

\begin{proof}

  Let $\lambda_j:=(\lambda_{j1},\dots,\lambda_{jn})\in \CC^n$, $j=1,\dots,k$,
   and denote $\lambda_{j\alpha}:=\lambda_{ji_1}\lambda_{ji_2}\dots
\lambda_{ji_m}$ if $\alpha=g_{i_1}g_{i_2}\dots g_{i_m}\in \FF_n^+$,
and $\lambda_{jg_0}:=1$.  As in Section \ref{Symmetric}, define the
vectors
$$
z_{\lambda_j}:=\sum_{\alpha\in \FF_n^+} \sqrt{b_\alpha}\overline
{\lambda}_{j\alpha} e_\alpha,\quad j=1,2,\dots,n.
$$
According to Theorem \ref{eigenvectors}, for any $g\in J$,
$\lambda\in \cV^\circ_{f,J}(\CC)$, and $\alpha,\beta\in \FF_n^+$, we
have
$$\left< z_\lambda, [W_\alpha g(W_1,\ldots, W_n) W_\beta](1)\right>
=\overline{\lambda}_\alpha \overline{g(\lambda)}
\overline{\lambda}_\beta=0,
$$
which implies  $z_\lambda\in \cN_J:=F^2(H_n)\ominus \overline{J(1)}$
for any $\lambda\in \cV^\circ_{f,J}(\CC)$. Note also that, since
$B_i^*=W_i^*|\cN_J$ and $C_i^*=\Lambda_i^*|_{\cN_J}$   for
$i=1,\ldots, n$, the same theorem implies
\begin{equation}\label{zj}
\begin{split}
B_i^*z_{\lambda_j}&=W_i^*z_{\lambda_j}=\overline{\lambda}_{ji}
z_{\lambda_j}\quad \text{for }\ i=1,\ldots, n; j=1,\ldots, k,\text{
and }\\
C_i^*
z_{\lambda_j}&=\Lambda^*_iz_{\lambda_j}=\overline{\lambda}_{ji}
z_{\lambda_j}\quad \text{for }\ i=1,\ldots, n; j=1,\ldots, k.
\end{split}
\end{equation}
Note that the subspace
$$
\cM:=\text{span} \{ z_{\lambda_j}:\ j=1,\dots,k\}
$$
 is invariant under $B_i^*$ and $\Lambda_i^*$ for any $i=1,\ldots,
 n$, and $\cM\subset \cN_J$.
  Define  the operators $X_i\in B(\cM\otimes\cK)$ by setting  $X_i=P_\cM B_i|_\cM\otimes I_\cK$, $i=1,\ldots, n$.
    Since $z_{\lambda_1},\dots, z_{\lambda_k}$ are linearly independent, we can define an operator
    $T\in  B(\cM\otimes\cK)$ by setting
    \begin{equation}\label{def-T*}
    T^*(z_{\lambda_j}\otimes h)=
z_{\lambda_j}\otimes A_j^* h
    \end{equation}
    for any $h\in \cK$ and   $j=1,\dots, k$.
    A simple calculation using  relations \eqref{zj} and \eqref{def-T*} shows that  $TX_i=X_iT$ for  $i=1,\dots, n $.

      Since $\cM$ is  a co-invariant subspace  under  each operator $B_i$, $ i=1,\dots,n$,  we can apply  Theorem
\ref{CLT2}
       and find
       $\Phi(\Lambda_1,\ldots, \Lambda_n)\in
       R_n^\infty(\cD_f)\bar \otimes B(\cK)$ such that  the operator
$$\Phi(C_1,\dots, C_n):=P_{\cN_J\otimes \cK}\Phi(\Lambda_1,\ldots,
\Lambda_n)|_{\cN_J\otimes \cK} \in R_n^\infty(\cV_{f,J})\bar\otimes
B(\cK)$$
 has the properties
\begin{equation}\label{PHI*}
\Phi(C_1,\dots, C_n)^* \cM\otimes \cK\subset \cM\otimes \cK,\quad
\Phi(C_1,\dots, C_n)^*|\cM\otimes \cK=T^*,
\end{equation}
 and  $ \|\Phi(C_1,\dots, C_n)\|=\|T\|$.

  Now, we  prove   that $\Phi(\lambda_j)=A_j$,\
$j=1,\ldots, k$, if and only if $$P_{\cM\otimes \cK}\Phi(C_1,\ldots,
C_n)|_{\cM\otimes \cK}=T.
$$
 Indeed,  notice that since $\Lambda^*_iz_{\lambda_j}=\overline{\lambda}_{ji}
z_{\lambda_j}$ we have
$$\varphi(\Lambda_1,\ldots,
\Lambda_n)^*z_\lambda=\overline{\varphi(\lambda)}z_\lambda$$
 for any
$\varphi(\Lambda_1,\ldots, \Lambda_n)\in R_n^\infty(\cD_f)$ and
$\lambda\in \cD_f^\circ(\CC)$. Consequently, we  obtain
\begin{equation}\label{Phi*-eq}
\Phi(\Lambda_1,\ldots, \Lambda)^*(z_\lambda\otimes
h)=z_\lambda\otimes \Phi(\lambda_1,\ldots, \lambda_n)^*h
\end{equation}
for any $\lambda\in \cD_f^\circ(\CC)$ and $h\in \cK$. Now,  due to
relations \eqref{zj} and \eqref{Phi*-eq}, we have
\begin{equation*}
\begin{split}
\left<\Phi(C_1,\dots, C_n)^*(z_{\lambda_j}\otimes x),
z_{\lambda_j}\otimes y\right>&= \left<\Phi(\Lambda_1,\ldots,
\Lambda)^*(z_{\lambda_j}\otimes x), z_{\lambda_j}\otimes y\right>\\
&=\left<z_{\lambda_j}\otimes \Phi(\lambda_j)^*x,
z_{\lambda_j}\otimes
y\right>\\
&= \left<z_{\lambda_j},
z_{\lambda_j}\right>\left<\Phi(\lambda_j)^*x,y\right>.
\end{split}
\end{equation*}
On the other hand, relation \eqref{def-T*} implies
$$\left<T^*(z_{\lambda_j}\otimes x), z_{\lambda_j}\otimes y\right>=
\left<z_{\lambda_j}, z_{\lambda_j}\right>\left<A_j^*x,y\right>.
$$
Due to Theorem \ref{eigenvectors}, we have $\left<z_{\lambda_j},
z_{\lambda_j}\right>= K_f(\lambda_j, \lambda_i)\neq 0$ for any
$j=1,\ldots, k$. Consequently,  the above relations imply  our
assertion.

 Now, since $ \|\Phi(C_1,\dots, C_n)\|=\|T\|$, it is clear that
 $\|\Phi(C_1,\dots, C_n)\|\leq 1$ if and only if $TT^*\leq I_\cM$.
Notice that, for any $h_1,\ldots, h_k\in \cK$, we have
 \begin{equation*}
 \begin{split}
 \left<\sum_{j=1}^k z_{\lambda_j}\otimes h_j, \sum_{j=1}^k z_{\lambda_j}\otimes
 h_j\right>&-\left<T^*\left(\sum_{j=1}^k z_{\lambda_j}\otimes h_j\right),
  T^*\left(\sum_{j=1}^k z_{\lambda_j}\otimes
 h_j\right)\right>\\
 &=\sum_{i,j=1}^k \left<z_{\lambda_i},
z_{\lambda_j}\right>\left< (I_\cK-A_jA_i^*)h_i, h_j\right>\\
&=\sum_{i,j=1}^k K_f(\lambda_j, \lambda_i)\left<
(I_\cK-A_jA_i^*)h_i, h_j\right>.
 \end{split}
 \end{equation*}
  Consequently, we
have $\|\Phi(C_1,\dots, C_n)\|\leq 1$ if and only if
  the matrix \eqref{semipo} is
positive semidefinite. This completes the proof.
\end{proof}

Let us make a few remarks. According to  Proposition \ref{tilde-f}
and  Proposition \ref{tilde-f2}, we have
$R_n^\infty(\cD_f)=U^*(F_n^\infty(\cD_{\tilde f}))U$. On the other
hand, an $n$-tuple $(\lambda_1,\ldots, \lambda_n)$ is in
$\cD_f^\circ(\CC)$ if and only if it is in $\cD_{\tilde
f}^\circ(\CC)$. Now, using Theorem \ref{symm-Fock} and  Theorem
\ref{Nev}, we can deduce that following result.

\begin{corollary} Let \, $f$ be a positive regular free holomorphic
function on $B(\cH)^n$  and  let $\lambda_1,\ldots, \lambda_k$ be
$k$ distinct points in  $ \cD_f^\circ(\CC)$. Given  $A_1,\ldots,
A_k\in B(\cK)$, the following statements are equivalent:
\begin{enumerate}\label{Nev-equi}
\item[(i)] there exists $\Psi \in  F_n^\infty(\cD_f)\bar\otimes B(\cK)$ such that
$$\|\Psi\|\leq 1\quad \text{and }\quad
\Psi(\lambda_j)=A_j,\ j=1,\ldots, k;
$$
\item[(ii)]
 there exists $\Phi \in  H^\infty(\cD_f^\circ(\CC))\bar\otimes B(\cK)$ such that
$$\|\Phi\|\leq 1\quad \text{and }\quad
\Phi(\lambda_j)=A_j,\ j=1,\ldots, k,
$$
where $H^\infty(\cD_f^\circ(\CC))$ is the algebra of multipliers of
$H^2(\cD_f^\circ(\CC))$;
 \item[(iii)] the operator matrix
\begin{equation*}
\left[ (I_\cK-A_i A_j^*) K_f(\lambda_i, \lambda_j)\right]_{k\times
k}
\end{equation*}
is positive semidefinite.
\end{enumerate}
\end{corollary}

Using this corollary, we can obtain the the following result.

\begin{theorem}\label{cara-hol}
Let \, $p$ be a  positive regular free polynomial on $B(\cH)^n$ and
let \,$\varphi $ be a complex-valued function defined on the
Reinhardt domain $\cD_p^\circ(\CC)\subset \CC^n$, such that
$$
|\varphi(z_1,\ldots, z_n)|<1\quad \text{ for all }\ (z_1,\ldots,
z_n)\in \cD_p^\circ(\CC).
$$
 Then there exists $F\in F_n^\infty(\cD_p)$ such that
$$\varphi(z_1,\ldots, z_n)=F(z_1,\ldots, z_n)\quad \text{ for all }\ (z_1,\ldots,
z_n)\in \cD_p^\circ(\CC),
$$
if and only if for each $k$-tuple of  distinct points
$\lambda_1,\ldots, \lambda_k\in \cD_p^\circ(\CC)$, the matrix
\begin{equation}
\label{varfi}
\left[(1-\varphi(\lambda_i)\overline{\varphi(\lambda_j)})
K_p(\lambda_i, \lambda_j) \right]_{k\times k}
\end{equation}
is positive semidefinite. In particular, if \eqref{varfi} holds,
then $\varphi$ is analytic on $\cD_p^\circ(\CC)$.
\end{theorem}

\begin{proof}
One implication follows from Corollary \ref{Nev-equi}. Conversely,
assume that $\varphi: \cD_p^\circ(\CC)\to \CC$ is such that the
matrix \eqref{varfi} is positive semidefinite for any  $k$-tuple of
distinct points $\lambda_1,\ldots, \lambda_k\in \cD_p^\circ(\CC)$.
Let $\{\lambda_j\}_{j=1}^\infty$ be a countable dense set in the
Reinhardt  domain $\cD_p^\circ(\CC)$. Applying  Theorem \ref{Nev},
for each $k\in\NN$, we find  $F_k\in F_n^\infty(\cD_p)$ such that
$\|F_k\|\leq 1$ and
\begin{equation}
\label{F-varfi} F_k(\lambda_j)=\varphi(\lambda_j)\quad \text{ for }
\ j=1,\ldots,k.
\end{equation}
Since the Hardy algebra $F_n^\infty(\cD_f)$ is $w^*$-closed
subalgebra in $B(F^2(H_n))$ and $\|F_k\|\leq 1$ for any $k\in \NN$,
we can use Alaoglu's theorem to find a subsequence
$\{F_{k_m}\}_{m=1}^\infty$ and $F\in F_n^\infty(\cD_p)$  such that
 $F_{k_m}\to F$, as $m\to\infty$, in the $w^*$-topology.
Since $\lambda_j:=(\lambda_{j1},\ldots, \lambda_{jn})\in
\cD_p^\circ(\CC)$, the $n$-tuple  is also of class $C_{\cdot 0}$.
Due to Theorem \ref{funct-calc}, the $F_n^\infty(\cD_p)$-functional
calculus  for c.n.c $n$-tuples of operators is $WOT$-continuous  on
bounded sets. Consequently, we deduce that $F_{k_m}(\lambda_j)\to
F(\lambda_j)$, as $m\to\infty$, for any $j\in \NN$. Hence and using
\eqref{F-varfi}, we obtain $\varphi(\lambda_j)=F(\lambda_j)$ for
$j\in \NN$. Given an arbitrary  element $z\in \cD_p^\circ(\CC)$, we
can apply again the above argument to find $G\in F_n^\infty(\cD_p)$,
$\|G\|\leq 1$ such that
$$
G(z)=\varphi(z) \ \text{ and }\ G(\lambda_j)=\varphi(\lambda_j),\
j\in \NN.
$$
Due to Proposition \ref{Reinhardt}, the maps $\lambda\mapsto
G(\lambda)$ and $\lambda\mapsto F(\lambda)$  are analytic on the
Reinhardt  domain  $\cD_p^\circ(\CC)$. Since they coincide on the
set $\{\lambda_j\}_{j=1}^\infty$, which is dense in
$\cD_p^\circ(\CC)$, we deduce that $G(\lambda)=F(\lambda)$ for any
$\lambda\in \cD_p^\circ(\CC)$. In particular, we have
$F(z)=\varphi(z)$. Since $z$ is an arbitrary element in
$\cD_p^\circ(\CC)$, the proof is complete.
\end{proof}

Let $\cP_m$ be the set of all polynomials in the Fock space
$F^2(H_n)$, of degree less then or equal to $m\in \NN$. Here is our
version of Schur-Carath\' eodory interpolation problem  on the Hardy
algebra $R_n^\infty(\cD_f)$.
\begin{theorem}
\label{Carate} Let $p(\Lambda_1,\ldots,
\Lambda_n):=\sum_{|\alpha|\leq m} \Lambda_\alpha \otimes
B_{(\alpha)}$ be a polynomial in $R_n^\infty(\cD_f)\bar\otimes
B(\cK_1, \cK_2)$. Then there is an element
$\varphi(\Lambda_1,\ldots, \Lambda_n)=\sum_{\alpha\in \FF_n^+}
\Lambda_\alpha \otimes D_{(\alpha)}$ in
$R_n^\infty(\cD_f)\bar\otimes B(\cK_1, \cK_2)$ with
$$\|\varphi(\Lambda_1,\ldots, \Lambda_n)\|\leq
1 \quad \text{ and } \quad B_{(\alpha)}=D_{(\alpha)}\ \text{ for any
} \ \alpha\in \FF_n^+,  \ |\alpha|\leq m,
$$
if and
only if
$$
\left\|P_{\cP_m\otimes \cK_2} p(\Lambda_1,\ldots,
\Lambda_n)|_{\cP_m\otimes \cK_1}\right\|\leq 1.
$$
\end{theorem}

\begin{proof}
Since $\cP_m$ is invariant under each operator $W_i^*$ and
$\Lambda_i^*$, $i=1,\ldots, n$, and $W_i\Lambda_j=\Lambda_j W_i$ for
any $i,j=1,\ldots,n$, we have
$$
(W_i^*\otimes I_{\cK_1})p(\Lambda_1,\ldots,
\Lambda_n)^*|_{\cP_m\otimes \cK_2}=p(\Lambda_1,\ldots, \Lambda_n)^*
(W_i^*\otimes I_{\cK_2})|_{\cP_m\otimes \cK_2}.
$$
Hence, we deduce that
$$
X[P_{\cP_m\otimes \cK_1}(W_i\otimes I_{\cK_1})|_{\cP_m\otimes
\cK_1}]=[P_{\cP_m\otimes \cK_2}(W_i\otimes I_{\cK_2})|_{\cP_m\otimes
\cK_2}]X,\quad i=1,\ldots,n,
$$
where $X:=P_{\cP_m\otimes \cK_2} p(\Lambda_1,\ldots,
\Lambda_n)|_{\cP_m\otimes \cK_1}$.
Applying Theorem \ref{CLT},
 we find
$$\Phi(\Lambda_1,\ldots, \Lambda_n)=\sum_{\alpha\in \FF_n^+}
\Lambda_\alpha \otimes D_{(\alpha)}\in R_n^\infty(\cD_f) \bar\otimes
B(\cK_1,\cK_2)$$ such that $\|\Phi(\Lambda_1,\ldots,
\Lambda_n)\|=\|X\|$, \ $\Phi(\Lambda_1,\ldots,
\Lambda_n)^*(\cP_m\otimes \cK_2)\subseteq \cP_m\otimes \cK_1$, and
$$
 \Phi(\Lambda_1,\ldots,
\Lambda_n)^*|\cP_m\otimes \cK_2=p(\Lambda_1,\ldots,
\Lambda_n)^*|\cP_m\otimes \cK_2.
$$
Using the definition of the weighted  right creation operators
associated with the noncommutative domain $\cD_f$, the latter
equality applied on vectors of the form $1\otimes h$, $h\in \cK_2$,
implies $D_{(g_0)}=B_{(g_0)}$. Now, using the same equality on
vectors of type $e_{g_i}\otimes h$, $i=1,\ldots, n$, we obtain
$\frac{1}{\sqrt{b_{g_i}}}\otimes
D_{(g_i)}^*h=\frac{1}{\sqrt{b_{g_i}}}\otimes B_{(g_i)}^*h$ and,
hence, $D_{(g_i)}=B_{(g_i)}$. Continuing this process and applying
the above equality on vector $e_\alpha\otimes h$ with
$|\alpha|=2,3,\ldots, m$, we deduce that $D_{(\alpha)}=B_{(\alpha)}$
for $|\alpha|\leq m$. This completes the proof.
\end{proof}

We remark that Theorem \ref{Nev}, Theorem \ref{cara-hol}, and
Theorem \ref{Carate} can be  recast   now  as statements for
multipliers of $H^2(\cD_f^\circ(\CC))$.

\bigskip

\section{Corona theorem for a class of Hardy algebras}
\label{Corona}

In this section, we prove a corona type theorem for  a class  of
Hardy  algebras  associated with the noncommutative domain $\cD_f$,
or the noncommutative variety $\cV_{f,J}$.

Let us recall that if $J$ be a $w^*$-closed two-sided ideal of  the
Hardy  algebra $F_n^\infty(\cD_f)$, then
  the {\it constrained  weighted left} (resp.~{\it right}) {\it
creation operators} associated with the noncommutative variety
$\cV_{f,J}$   are defined by setting
$$B_i:=P_{\cN_J} W_i|\cN_J \quad \text{and}\quad
C_i:=P_{\cN_J} \Lambda_i|\cN_J,\quad i=1,\ldots, n,
$$
where $ \cN_J:=F^2(H_n)\ominus \overline{JF^2(H_n)}$.

  \begin{theorem}\label{facto}
  Let $J$ be a $w^*$-closed two-sided ideal of  the  Hardy  algebra
$F_n^\infty(\cD_f)$  and let $(B_1,\ldots, B_n)$ and $(C_1,\ldots,
C_n)$ be the corresponding constrained weighted left (resp. right)
creation operators  associated with the noncommutative variety
$\cV_{f,J}$.
  If $A\in R_n^\infty(\cV_{f,J})\bar\otimes B(\cH, \cH')$
  and $B\in R_n^\infty(\cV_{f,J})\bar\otimes B(\cH'', \cH')$, then
  there exists a contraction  $G\in R_n^\infty(\cV_{f,J})
  \bar\otimes B(\cH, \cH'')$
  such that $A=BG$  if and only if
  $AA^*\leq BB^*$.
  \end{theorem}
  \begin{proof}
   Assume that $AA^*\leq BB^*$.
  Then, there is a contraction
  $Y:\cM:=\overline{B^*( \cN_J\otimes \cH')}\to  \cN_J\otimes\cH$
   satisfying $YB^*=A^*$.
  Since $B\in R_n^\infty(\cV_{f,J})\bar\otimes B(\cH'', \cH')$, Corollary
  \ref{commutant} implies
  $$
  B (B_i\otimes I_{\cH''})=(B_i\otimes I_{\cH'})B, \quad i=1,\ldots, n,
  $$
  where $B_i:=P_{  \cN_J} W_i |   \cN_J$, $i=1,\ldots, n$.
    Hence, the subspace $\cM\subseteq   \cN_J\otimes \cH''$ is invariant
     under each
  operator
  $B_i^*\otimes I_{{\cH''}}$, $i=1,\ldots, n$.
  Define the operators
   $$T_i:= P_\cM(B_i\otimes I_{\cH''})|\cM,\quad  i=1,\ldots, n,
   $$
   acting from $\cM$ to $\cM$.
  Since
  $A$ is in the operator space $R_n^\infty(\cV_{f,J})\bar\otimes B(\cH, \cH')$,
   while $B$ is in $ R_n^\infty(\cV_{f,J})\bar\otimes B(\cH'', \cH')$,
   and $YB^*=A^*$, we have
  \begin{equation*}
  \begin{split}
  Y(B_i^*\otimes I_{\cH''}) B^*x&=YB^* (B_i^*\otimes I_{\cH'})x=
  A^*(B_i^*\otimes I_{\cH'})x\\
  &=
  (B_i^*\otimes I_{\cH})A^*x=
  (B_i^*\otimes I_{\cH}) Y B^*x
  \end{split}
  \end{equation*}
  for any $x\in  \cN_J\otimes \cH'$.
  Hence,
  \begin{equation}
  \label{intert}
  Y^*(B_i\otimes I_{\cH})= T_iY^* ,\quad i=1,\ldots, n.
  \end{equation}
Applying the commutant lifting of Theorem \ref{CLT2} in our
particular case, we find
   an operator $G\in R_n^\infty(\cV_{f,J})\bar\otimes B(\cH,\cH'')$
   such that $G^*|_\cM=Y$ and $\|G\|=\|Y^*\|$.
Now, it is clear that
 $$A=BY^*=BP_\cM G=B(P_\cM+P_\cM^\perp)G=BG.
$$
This completes the proof of the theorem.
 \end{proof}

  Applying Theorem \ref{facto} to  the particular case
   when $\cH=\cH'$ and  $A=\delta I_{\cN_J\otimes \cH'}$,
   $\delta>0$,   we  obtain the following consequence.

  \begin{corollary}\label{coron}
  Let  $B\in R_n^\infty(\cV_{f,J})\bar\otimes B(\cH'', \cH')$. The
  following statements are equivalent:
  \begin{enumerate}
  \item[(i)] There is $D\in B(\cN_J\otimes \cH', \cN_J\otimes
  \cH'')$ such that
  $BD=I;$
  \item[(ii)] There is $\delta>0$ such that $\|B^*h\|\geq \delta\|h\|$,
   for any
  $h\in \cN_J\otimes \cH';$
  \item[(iii)] There is $G\in R_n^\infty(\cV_{f,J})\bar\otimes B(\cH',
   \cH'')$
  such that
  $BG=I$.
  \end{enumerate}
  \end{corollary}

  Another consequence  of Theorem \ref{facto} is the following corona type result
   for the algebra $R_n^\infty(\cV_{f,J})\bar\otimes B(\cH)$.

  \begin{corollary}\label{corona2}
  If  $\varphi_1,\ldots, \varphi_k\in
  R_n^\infty(\cV_{f,J})\bar\otimes B(\cH)$,
  then there exist  operators
  $g_1,\ldots, g_k\in R_n^\infty(\cV_{f,J})\bar\otimes B(\cH)$ such that
  $$
  \varphi_1g_1+\cdots +\varphi_kg_k=I
  $$
  if and only if there exists $\delta>0$ such that
  $$
  \varphi_1\varphi_1^*+\cdots +\varphi_k \varphi_k^*\geq \delta^2 I.
  $$
   \end{corollary}

  \begin{proof}
  Take
  \begin{equation*}
  \begin{split}
  B&:=[\varphi_1,\ldots, \varphi_k]\in R_n^\infty(\cV_{f,J})\bar \otimes
   B(\cH^{(k)}, \cH), \text{ and }\\
  C&:=\left[\begin{matrix}g_1\\\vdots\\
  g_k\end{matrix}\right]\in R_n^\infty(\cV_{f,J})
  \bar\otimes B(\cH,\cH^{(k)})
  \end{split}
  \end{equation*}
  and apply
  Corollary \ref{coron}.
  \end{proof}

We remark  that a similar result can be obtained for
$F_n^\infty(\cV_{f,J})\bar\otimes B(\cH)$.

Now, let us discuss a few particular cases. First,  we consider the
case when $J=\{0\}$. If  $\varphi_i  \in F_n^\infty(\cD_f)$, \
$i=1,\ldots, k$, we denote by $\sigma_r(\varphi_1,\ldots,
\varphi_k)$ the right joint spectrum
  with respect to the
   noncommutative analytic Toeplitz algebra $F_n^\infty(\cD_f)$.
 \begin{corollary}\label{spectru}
Let $(\varphi_1,\ldots, \varphi_k)$ be a $k$-tuple of operators in
  the  Hardy algebra $F_n^\infty(\cD_f)$.
     Then
 the following properties hold:
 \begin{enumerate}
  \item[(i)] $(\lambda_1,\ldots, \lambda_n)\notin
  \sigma_r(\varphi_1,\ldots, \varphi_k)$ if and only
  if there is $\delta>0$ such that
  $$(\lambda_1 I-\varphi_1)(\overline{\lambda}_1 I-\varphi_1^*)+\cdots +
  (\lambda_k I-\varphi_k)(\overline{\lambda}_k I-\varphi_k^*)\geq \delta^2 I;$$
  \item[(ii)]
  $\sigma_r(\varphi_1,\ldots, \varphi_k)$ is a  compact subset of $\CC^k$ and
  $$
  \{(\varphi_1(\lambda ), \ldots,
  \varphi_k(\lambda )):\  \lambda \in \cD_f^\circ(\CC)\}^{-}\subseteq
  \sigma_r(\varphi_1,\ldots, \varphi_k).
   $$
   \end{enumerate}
  \end{corollary}
\begin{proof} The first part of the corollary follows from the remarks
preceding  Proposition \ref{right-spec}, and Corollary
\ref{corona2}. According to Theorem \ref{eigenvectors}, we have
 $$
  \varphi(W_1,\ldots, W_n)^* z_\lambda=
  \overline{\varphi(\lambda )}
   z_\lambda \quad \text{ for all } \ \lambda\in \cD_f^\circ(\CC).
  $$
 Denote $\mu_j:= \varphi_j(\lambda )$, \ $j=1,\ldots, k$, and notice that
 $$
 \sum_{j=1}^k \|(\overline{\mu}_jI-\varphi(W_1,\ldots, W_n)^*) z_\lambda\|^2=0.
 $$
 Since $\|z_\lambda\|\neq 0$ for any $\lambda \in \cD_f^\circ(\CC)$ ,
 we can use  again  Corollary
 \ref{corona2} and deduce that
  $(\mu_1,\ldots, \mu_k)\in \sigma_r(\varphi_1,\ldots, \varphi_k)$.
\end{proof}

 Consider now the particular case when $J_c$ is the $w^*$-closed two-sided ideal
 of the Hardy algebra $F_n^\infty(\cD_f)$ generated by
 $W_iW_j-W_jW_i$, $i,j=1,\ldots,n$. Using Theorem \ref{symm-Fock}
 and Theorem \ref{W(L)}, we can
 write all the results of this section for
 $F_n^\infty(\cV_{f,J_c})$, the algebra    of
  all multipliers of the  Hilbert space $H^2(\cD_f^\circ(\CC))$. We
  leave this task to the reader.

      \bigskip

      %\Refs
      %\widestnumber\key{BFPQR}
      %\def\n{\key}
       %


\begin{thebibliography}{99}


\bibitem{Ag1} {\sc J.~Agler}, The Arveson extension theorem and
coanalytic models, {\it Integral Equations Operator Theory} {\bf 5}
(1982), 608-631.

\bibitem{Ag2} {\sc J.~Agler},
Hypercontractions and subnormality, {\it J. Operator Theory} {\bf
13} (1985), 203--217.


\bibitem{AMc2} {\sc J.~Agler and J.~E.~McCarthy},
 Complete Nevanlina-Pick kernels,
%
{\it J. Funct. Anal.}
  {\bf 175} (2000), 111--124.


      \bibitem{Ai} {\sc L.~Aizenberg},
      {Multidimensional analogues of Bohr's theorem on power series,}
      {\it Proc. Amer. Math. Soc.} {\bf 128} (2000), 1147--1155.


\bibitem{Ar1} {\sc A.~Arias},
{Projective modules of Fock spaces}
{\it J. Operator Theory}
{\bf 52} (2004), no. 1, 139--172.



\bibitem{ArPo} {\sc  A.~Arias and G.~Popescu},
Factorization and reflexivity on Fock spaces, {\it  Integr. Equat.
Oper. Th.}
 %
 {\bf  23} (1995),  268--286.



      \bibitem{ArPo1}   {\sc A.~Arias and G.~Popescu},
      {Noncommutative interpolation and Poisson transforms II,}
      {\it Houston J. Math.} {\bf 25} (1999), No. 1,  79--97.



      \bibitem{ArPo2}   {\sc A.~Arias and G.~Popescu},
      {Noncommutative interpolation and Poisson transforms,}
      {\it Israel J. Math.} {\bf 115} (2000), 205--234.

  %
\bibitem{Arv-acta} {\sc W.B.~Arveson}, Subalgebras of $C^*$-algebras,
{\it Acta.Math.} {\bf 123} (1969), 141--224.



\bibitem{Arv-book}  {\sc  W.B.~Arveson},
{\it An invitation to $C^*$-algebras}, Graduate Texts in Math., {\bf
39}. Springer-Verlag, New-York-Heidelberg, 1976.





\bibitem{Arv}    {\sc W.B.~Arveson},
      {Subalgebras of $C^*$-algebras III: Multivariable operator theory,}
{\it Acta Math.} {\bf 181} (1998), 159--228.


\bibitem{Arv2} {\sc W.B.~Arveson},
  The curvature invariant of a Hilbert module  over $\CC [z_1,\ldots, z_n]$,
 {\it J. Reine Angew. Math.}
 {\bf 522} (2000),  173--236.





\bibitem{BaGR} {\sc J.~A.~Ball, I.~Gohberg, and L.~Rodman},
{\it Interpolation of Rational Matrix Functions}, {\bf OT 45},
Birkh\"auser-Verlag, Basel-Boston, 1990.

\bibitem{BTV} {\sc J.~A.~Ball, T.~T.~Trent, and V.~Vinnikov}, Interpolation
and commutant
lifting for multipliers on reproducing kernels Hilbert spaces,
{\it Operator Theory and Analysis: The M.A. Kaashoek Anniversary Volume},
pages 89--138,
{\bf OT 122}, Birkhauser-Verlag, Basel-Boston-Berlin, 2001.

\bibitem{BB} {\sc J.A.~Ball and V.~Bolotnikov}, On bitangential interpolation
 problem for contractive-valued
 functions on the unit ball,
 {\it  Linear Algebra Appl.}
  {\bf 353} (2002), 107--147.



\bibitem{BB2} {\sc J.~A.~Ball and V.~Bolotnikov},
 Realization and
interpolation for Schur-Agler-class functions on domains with matrix
polynomial defining function in $\Bbb C\sp n$, {\it J. Funct. Anal.}
{\bf 213} (2004),  45--87.


\bibitem{BV} {\sc J.~A.~Ball   and V.~Vinnikov},
Lax-Phillips Scattering and Conservative Linear Systems: A
Cuntz-Algebra Multidimensional Setting, {\it Mem. Amer. Math. Soc.}
{\bf 837}, (2005).

 %
 \bibitem{BV3} {\sc J.~A.~Ball   and V.~Vinnikov},
 Functional models for representations of the Cuntz algebra,
 {\it Operator theory, systems theory and scattering theory:
  multidimensional generalizations}, 1--60,
Oper. Theory Adv. Appl., {\bf 157}, Birkhäuser, Basel, 2005.

\bibitem{BGM} {\sc J.~A.~Ball, G.~Groenewald, and T.~Malakorn},
 Conservative structured noncommutative multidimensional linear systems,
 {\it  The state space method generalizations and applications}, 179--223,
Oper. Theory Adv. Appl., {\bf 161}, Birkhäuser, Basel, 2006.




 \bibitem{BT1} {\sc C.~Benhida, and D.~Timotin},
Characteristic functions for multicontractions and automorphisms
 of the unit ball, {\it  Integr. Equat. Oper.Th.} {\bf 57} (2007), no. 2, 153--166.



 \bibitem{BT2} {\sc C.~Benhida, and D.~Timotin},
Some automorphism invariance properties for multicontractions, {\it
 Indiana Univ. Math. J.} {\bf 56} (2007), no. 1, 481--499.


\bibitem{Be} {\sc A.~Beurling}
{ On two problems concerning linear transformations in Hilbert
space}, {\it  Acta Math.} {\bf 81}, (1948), 239--251.





\bibitem{BBD} {\sc B.V.R.~Bhat, T.~Bhattacharyya, and S.~Dey},
Standard noncommuting and commuting dilations of commuting tuples,
{\it Trans. Amer. Math. Soc.} {\bf 356} (2003), 1551--1568.



\bibitem{BES1} {\sc  T.~Bhattacharyya, J.~Eschmeier, and J.~Sarkar},
Characteristic function of a   pure commuting  contractive tuple,
{\it Integral Equation Operator Theory} {\bf 53} (2005), no.1,
23--32.


\bibitem{BES2} {\sc  T.~Bhattacharyya, J.~Eschmeier, and J.~Sarkar},
On CNC commuting contractive tuples, {\it Proc. Indian Acad. Sci.
Math. Sci.} {\bf 116} (2006), no. 3, 299--316.




\bibitem{BS} {\sc  T.~Bhattacharyya  and J.~Sarkar},
Characteristic function for polynomially contractive commuting
tuples, {\it  J. Math. Anal. Appl.}  {\bf 321} (2006), no. 1,
242--259.


      \bibitem{BK} {\sc H.P. Boas and D. Khavinson},
      Bohr's power series theorem in several variables,
      {\it Proc. Amer. Math. Soc.} {\bf 125} 1997, 2975--2979.

      \bibitem{Bo} {\sc H. Bohr},
      A theorem concerning power series,
      {\it Proc. London Math. Soc.} (2){\bf 13} 1914, 1--5.



\bibitem{B} {\sc J.~W.~Bunce},
 Models for n-tuples of noncommuting operators, {\it J. Funct. Anal.}
  {\bf 57}(1984), 21--30.


\bibitem{Ca} {\sc C.~Carath\' eodory},
 \" Uber den Variabilit\" atsbereich der Koeffzienten von
  Potenzreinen die gegebene Werte nicht annehmen,
{\it Math. Ann.}
 {\bf 64}  (1907), 95--115.



\bibitem{Co} {\sc J.B.~Conway},
{\em Functions of one complex variable. I.} Second Edition. Graduate
Texts in Mathematics  {\bf 159}. { Springer-Verlag, New York}, 1995.


 \bibitem{CV1} {\sc R.E.~Curto, F.H.~Vasilescu}, Automorphism invariance
 of the operator-valued Poisson transform,
{\it Acta Sci. Math. (Szeged)} {\bf 57} (1993), 65--78.

\bibitem{CV2} {\sc R.E.~Curto, F.H.~Vasilescu}, Standard operator
models in the polydisc, {\it Indiana Univ. Math. J.} {\bf 42}
(1993), 791--810.

\bibitem{CV3} {\sc R.E.~Curto, F.H.~Vasilescu}, Standard operator
models in the polydisc II, {\it Indiana Univ. Math. J.} {\bf 44}
(1995), 727--746.






\bibitem{DP} {\sc K.~R.~Davidson and D.~Pitts},
 Nevanlinna-Pick interpolation for noncommutative analytic Toeplitz algebras,
 {\it  Integr. Equat. Oper.Th.}
{\bf 31} (1998), 321--337.









\bibitem{DP2} {\sc K.~R.~Davidson and D.~Pitts},
%
 The algebraic structure of  non-commutative
  analytic Toeplitz algebras,
{\it  Math. Ann.}
   {\bf 311} (1998),  275--303.



            \bibitem{DP1} {\sc K.R.~Davidson and D. ~Pitts},
      {Invariant subspaces and hyper-reflexivity for free semigroup
      algebras,}
      {\it Proc. London Math. Soc.} {\bf 78} (1999), 401--430.




\bibitem{DKS} {\sc K.R.~Davidson, D.W.~ Kribs,  and  M.E.~Shpigel},
 Isometric
  dilations of non-commuting finite rank $n$-tuples,
  {\it Canad. J. Math.}
  {\bf 53} (2001),
  506--545.


\bibitem{DKP}  {\sc K.R.~Davidson, E.~Katsoulis, and D.~Pitts},
  The structure of free semigroup algebras,
 {\it J. Reine Angew. Math.}
  {\bf 533} (2001), 99--125.




             \bibitem{DT} {\sc S.~Dineen and R.M.~Timoney},
      {On a problem of H. Bohr,}
      {\it Bull. Soc. Roy. Sci. Liege} {\bf 60} (1991), 401--404.


 \bibitem{Dr} {\sc S.~Drury}, A generalization of von Neumann's inequality
 to the complex
 ball,
 {\it Proc. Amer. Math. Soc.} {\bf 68} (1978), 300-304.


\bibitem{EP} {\sc J.~Eschmeier and M.~Putinar}, Spherical contractions and interpolation problems
on the unit ball, {\it J. Reine Angew. Math.} {\bf 542} (2002),
219--236.

\bibitem{FF-book}  {\sc C.~Foias, A.~E.~Frazho},
 {\it The commutant lifting approach to interpolation problems},
 Operator Theory: Advances and Applications,
Birh\"auser Verlag,  Bassel,
 1990.







\bibitem{FFGK-book}  {\sc C.~Foias, A.~E.~Frazho, I.~Gohberg, and M.~Kaashoek},
{\em  Metric constrained interpolation, commutant lifting and
systems},
  Operator Theory: Advances and Applications vol. 100,  Birh\"auser Verlag,
 Bassel, 1998.


 \bibitem{F} {\sc  A.~E.~Frazho},
 Models for noncommuting operators, {\it J. Funct. Anal.}
  {\bf 48} (1982), 1--11.
   %
\bibitem{G} {\sc J.~Garnett},
 {\it  Bounded analytic functions},
 Academic Press, New York, 1981.



\bibitem{H} {\sc K.~Hoffman},  {\em Banach Spaces of Analytic Functions},
Englewood Cliffs: Prentice-Hall, 1962.




\bibitem{He} {\sc W.J.~Helton},
``Positive'' noncommutative polynomials are sums of squares, {\it
Ann. of Math.} (2) {\bf 156} (2002), no.2, 675--694.

\bibitem{He-C}{\sc J.W.~ Helton and  S.A.~McCullough},
  A Positivstellensatz for non-commutative polynomials,
  {\it Trans.
Amer. Math. Soc.}  {\bf 356} (2004), no. 9, 3721--3737.


\bibitem{He-C-P2}{\sc J.W.~ Helton, S.A.~McCullough, and
M.~Putinar}, A non-commutative Positivstellensatz on isometries,
{\it  J. Reine Angew. Math.} {\bf  568} (2004), 71--80.



\bibitem{He-C-P}{\sc J.W.~ Helton, S.A.~McCullough, and M.~Putinar},
 Matrix representations for positive noncommutative
polynomials, {\it  Positivity} {\bf 10} (2006), no. 1, 145--163.



 \bibitem{He-C-V}{\sc J.W.~ Helton, S.A.~McCullough, and V.~
 Vinnikov},
  Noncommutative convexity arises from linear matrix inequalities,
  {\it J. Funct. Anal.} {\bf 240} (2006), no. 1, 105--191.







\bibitem{Kr} {\sc D.~Kribs},
  The curvature invariant of a  non-commuting $N$-tuple,
 {\it Integral Equations Operator Theory}  {\bf 41} (2001), no.\,4, 426--454.


\bibitem{MuSo1}  {\sc P.S.~Muhly and  B.~Solel},
 Tensor algebras over $C^*$-correspondences: representations,
dilations, and $C^*$-envelopes, {\it J. Funct. Anal.}
 {\bf 158} (1998),  389--457.


  \bibitem{MuSo5}  {\sc P.S.~Muhly and  B.~Solel},
Hardy algebras, $W^*$-correspondences and interpolation theory, {\it
Math. Ann.} {\bf 330} (2004),  353--415.


\bibitem{MuSo6}  {\sc P.S.~Muhly and  B.~Solel},
Canonical models for representations of Hardy algebras,
 {\it Integral Equations Operator Theory} {\bf 53} (2005), 411--452.



\bibitem{MuSo7}  {\sc P.S.~Muhly and  B.~Solel},
Schur class  operator functions and automorphisms of Hardy algebras,
preprint.

\bibitem{M} {\sc V.~M\"uller}, Models for operators using weighted
shifts, {\it J. Operator Theory} {\bf 20} (1988), 3--20.

\bibitem{MV} {\sc V.~M\"uller and F.H.~Vasilescu}, Standard models for
some commuting multioperators, {\it Proc. Amer. Math. Soc.} {\bf
117} (1993), 979--989.




\bibitem{N}  {\sc R.~Nevanlinna},
 \"Uber  beschr\"ankte Functionen, die in gegebenen Punkten
vorgeschribene Werte annehmen,
%
{\it Ann. Acad. Sci. Fenn. Ser A}
  {\bf 13} (1919),  7--23.




\bibitem{Par} {\sc S.~Parrott},
{On a quotient norm and the Sz.-Nagy--Foias lifting theorem}, {\it
J. Funct. Anal.}, {\bf 30} (1978), 311--328.


\bibitem{Pa1} {\sc V.I.~Paulsen}, Every completely polynomially bounded
operator is similar to a contraction, {\it J. Funct. Anal.}, {\bf55}
(1984), 1--17.


\bibitem{Pa-book} {\sc V.I.~Paulsen},
 {\it Completely Bounded Maps and Dilations},
Pitman Research Notes in Mathematics, Vol.146, New York, 1986.


\bibitem{PPoS}  {\sc V.I.~Paulsen, G.~Popescu,  and D.~Singh},
        {On Bohr's inequality,}
        {\it Proc. London Math. Soc.}
        {\bf 85} (2002), 493--512.

\bibitem{Pi} {\sc G.~Pisier}, \emph{ Similarity Problems and Completely Bounded Maps},
Springer Lect. Notes Math., Vol.1618, Springer-Verlag, New York,
1995.

             %%\bigskip

\bibitem{Po-models} {\sc G.~Popescu}, Models for infinite sequences of
noncommuting operators, {\it Acta. Sci. Math. (Szeged)} {\bf
53} (1989), 355--368.


\bibitem{Po-isometric} {\sc G.~Popescu}, Isometric dilations for infinite
sequences of noncommuting operators, {\it Trans. Amer. Math. Soc.}
{\bf 316} (1989), 523--536.

\bibitem{Po-charact} {\sc G.~Popescu}, Characteristic functions for infinite
sequences of noncommuting operators, {\it J. Operator Theory}
{\bf 22} (1989), 51--71.


\bibitem{Po-multi} {\sc G.~Popescu},
 Multi-analytic operators and some  factorization theorems,
 {\it Indiana Univ. Math.~J.}
 {\bf 38} (1989),   693--710.


      \bibitem{Po-von} {\sc G.~Popescu},
{Von Neumann inequality for $(B(H)^n)_1$,}
      {\it Math.  Scand.} {\bf 68} (1991), 292--304.
      %%\bigskip


\bibitem{Po-intert} {\sc G.~Popescu},  On intertwining dilations for
 sequences of noncommuting operators, {\it J. Math. Anal. Appl.}
{\bf 167} (1992), 382--402.

   %\leavevmode\vrule height 2pt depth -1.6pt width 44pt,



      \bibitem{Po-funct} {\sc G.~Popescu},
      {Functional calculus for noncommuting operators,}
       {\it Michigan Math. J.} {\bf 42} (1995), 345--356.



      \bibitem{Po-analytic} {\sc G.~Popescu},
      {Multi-analytic operators on Fock spaces,}
      {\it Math. Ann.} {\bf 303} (1995), 31--46.



      \bibitem{Po-disc}  {\sc G.~Popescu},
      {Noncommutative disc algebras and their representations,}
       {\it Proc. Amer. Math. Soc.} {\bf 124} (1996),  2137--2148.



\bibitem{Po-interpo} {\sc G.~Popescu},
  Interpolation problems in several variables, {\it  J. Math Anal. Appl.}
 {\bf 227}  (1998),    227--250.




      \bibitem{Po-poisson} {\sc G.~Popescu},
     {Poisson transforms on some $C^*$-algebras generated by isometries,}
       {\it J. Funct. Anal.} {\bf 161} (1999),  27--61.

\bibitem{Po-structure} {\sc G.~Popescu},
Structure and entropy for positive-definite Toeplitz kernels on free
semigroups, {\it J. Math. Anal. Appl.} {\bf 254} (2001), 191--218.

\bibitem{Po-curvature} {\sc  G.~Popescu},
  Curvature invariant for Hilbert modules over free semigroup algebras,
   {\it Adv. Math.}
 {\bf 158} (2001), 264--309.


\bibitem{Po-central} {\sc G.~Popescu},
  Central intertwining lifting, suboptimization, and interpolation
  in several variables, {\it  J. Funct. Anal.}  {\bf 189}  (2002), 132--154.


\bibitem{Po-nehari} {\sc G.~Popescu},
  Multivariable Nehari problem and interpolation,
   {\it  J. Funct. Anal.}
   {\bf 200} (2003), 536--581.


   \bibitem{Po-similarity} {\sc G.~Popescu},
    Similarity and ergodic theory of positive linear maps,
   {\it J. Reine Angew. Math.}
   {\bf 561} (2003), 87--129.



\bibitem{Po-entropy} {\sc G.~Popescu},
Entropy and Multivariable Interpolation,
 {\it Mem. Amer. Math. Soc.}
{\bf 184} (2006), No.868, vi+83 pp.


\bibitem{Po-charact-inv} {\sc G.~Popescu},
Characteristic functions and joint invariant suspaces, {\it J.
Funct. Anal.} {\bf 237} (2006), 277--320.

\bibitem{Po-varieties} {\sc  G.~Popescu},
Operator theory on noncommutative varieties, {\it Indiana Univ.
Math.~J.} {\bf 55}, no.2, (2006), 389--442.


\bibitem{Po-varieties2} {\sc  G.~Popescu},
Operator theory on noncommutative varieties II, {\it Proc. Amer.
Math. Soc.} {\bf 135} (2007), no. 7,  2151--2164.



      \bibitem{Po-unitary} {\sc G.~Popescu},
      {Unitary invariants in multivariable operator theory}, {\it Mem. Amer. Math. Soc.}, to appear.


      \bibitem{Po-holomorphic} {\sc G.~Popescu},
      {Free holomorphic functions on the unit ball of
      $B(\cH)^n$}, {\it J. Funct. Anal.}  {\bf 241} (2006),  268--333.

\bibitem{Po-Bohr} {\sc G.~Popescu},
{Multivariable Bohr inequalities}, {\it Trans. Amer. Math. Soc.}
{\bf 359} (2007), no. 11, 5283--5317.



\bibitem{Po-Berezin} {\sc G.~Popescu},
{  Noncommutative Berezin transforms and multivariable operator
model theory}, {\it J. Funct. Anal.},  {\bf 254} (2008), 1003--1057.



 \bibitem{Pot} {\sc S.~Pott},
 Standard models under polynomial positivity conditions,
 {\it J. Operator Theory} {\bf 41} (1999), no. 2, 365--389.


 %%\bigskip
\bibitem{S} {\sc D.~Sarason},
 Generalized interpolation in $H^\infty$,
{\it Trans. Amer. Math. Soc.}
  {\bf 127} (1967),   179--203.


\bibitem{Sh} {\sc A.L.~Shields},
{Weighted shift operators and analytic function theory}, {\it Topics
in operator theory}, pp. 49--128. Math. Surveys, No. 13, {\it Amer.
Math. Soc., Providence, R.I.,} 1974






\bibitem{Sc}  {\sc I.~Schur},
%
 \"Uber Potenzreihen die im innern des Einheitshreises beschr\"ankt sind,
{\it  J. Reine Angew. Math.}
 {\bf 148} (1918), 122--145.


  \bibitem{St} {\sc W.F.~Stinespring}, Positive functions on $C^*$-algebras,
  {\it Proc. Amer. Math. Soc.}
  {\bf 6} (1955) 211--216.


\bibitem{Sz1} {\sc B.~Sz.-Nagy}, Sur les contractions de l'espace de Hilbert,
{\it Acta. Sci.  Math. (Szeged)} {\bf 15} (1953), 87--92.


\bibitem{SzF-book} {\sc B.~Sz.-Nagy and C.~Foia\c{s}}, {\em
Harmonic Analysis of Operators on Hilbert Space}, North Holland, New
York 1970.


\bibitem{Va}{\sc F.H.~Vasilescu},
{An operator-valued Poisson kernel}, {\it J. Funct. Anal.}
  {\bf 110} (1992), 47--72.
        %

      \bibitem{vN}  {\sc J.~von Neumann},
      {Eine Spectraltheorie f\"ur allgemeine Operatoren eines unit\"aren
      Raumes,}
      {\it Math. Nachr.} {\bf 4} (1951), 258--281.
      %%\bigskip

       \end{thebibliography}
      \end{document}